\documentclass[a4paper,10pt]{report}
\usepackage[utf8x]{inputenc}
\usepackage{amsmath,amssymb,latexsym,mathrsfs}
\usepackage[all]{xy}
\usepackage{rotating}
\usepackage{changepage}

\newcommand{\Hom}{\mathbf{Hom}}
\newcommand{\Ext}{\mathbf{Ext}}
\newcommand{\rad}{\mathbf{rad}}
\newcommand{\soc}{\mathbf{soc}}
\newcommand{\Irr}{\mathbf{Irr}}

\newcommand{\End}{\mathbf{End}}
\newcommand{\dimk}{\mathbf{dim}}
\newcommand{\vdim}{\underline{\mathbf{dim}}}
\newcommand{\dimGr}{\mathbf{dimGr}}

\newcommand{\Ker}{\mathbf{Ker} \;}
\newcommand{\Img}{\mathbf{Im} \;}

\newcommand{\Rnk}{\mathbf{rk}}
\newcommand{\Ind}{\mathbf{ind}}
\newcommand{\ConR}{\mathbf{ConR}}
\newcommand{\ConL}{\mathbf{ConL}}
\newcommand{\dConR}{\mathbf{ConR}_0}
\newcommand{\dConL}{\mathbf{ConL}_0}

\newtheorem{lema}{Lemma}[chapter]
\newtheorem{proposicion}[lema]{Proposition}
\newtheorem{corolario}[lema]{Corollary}
\newtheorem{teorema}[lema]{Theorem}

\newtheorem{definicion}[lema]{Definition}

\def\bproof{\noindent{\usefont{OT1}{cmss}{m}{n}\selectfont\textbf{Proof. }}}
\def\eproof{\noindent{\hfill{$\square$}}\bigskip}

\newcommand{\CoefMax}{\xy 
  ( -1,0)="X" *{{}_{\bullet}};
  ( 1,0)="Y" *{{}_{\bullet}};
  ( 2,1)="X1" *{};
  ( 2,-1)="X2" *{};
  ( -2,-1)="X3" *{};
  ( -2,1)="X4" *{};
  "X1";"X2" **@{-};
  "X2";"X3" **@{-};
  "X3";"X4" **@{-};
  "X4";"X1" **@{-};
  \endxy}
\newcommand{\CoefAlaW}[1]{\xy 
  ( 0,0)="X" *{{}_{\bullet}};
  ( 0,-1.5)="Y" *{{}_{#1}};
  ( 0,3)="X1" *{};
  ( 3,-3)="X2" *{};
  ( -3,-3)="X3" *{};
  "X1";"X2" **@{-};
  "X2";"X3" **@{-};
  "X3";"X1" **@{-};
  \endxy}
\newcommand{\CoefAlaU}[1]{\xy 
  ( 0,0)="X" *{{}_{\bullet}};
  ( 1.5,0)="Y" *{{}_{#1}};
  ( 3,-3)="X1" *{};
  ( 3,3)="X2" *{};
  ( -3,0)="X3" *{};
  "X1";"X2" **@{-};
  "X2";"X3" **@{-};
  "X3";"X1" **@{-};
  \endxy}
\newcommand{\CoefAlaV}[1]{\xy 
  ( 0,0)="X" *{{}_{\bullet}};
  ( -1.5,0)="Y" *{{}_{#1}};
  (-3,-3)="X1" *{};
  (-3,3)="X2" *{};
  ( 3,0)="X3" *{};
  "X1";"X2" **@{-};
  "X2";"X3" **@{-};
  "X3";"X1" **@{-};
  \endxy}
\newcommand{\vctACi}[5]{\xy  
    (0,3)="B"       *{\scriptscriptstyle \mathbf{#1}};
    (0,1.5)="B1"      *{\scriptscriptstyle #2};
    (0,0)="B2"    *{\scriptscriptstyle #3};
    (0,-1.5)="B3"      *{\scriptscriptstyle #4};
    (0,-3)="B4"     *{\scriptscriptstyle \mathbf{#5}};
    \endxy}
\newcommand{\vctDn}[7]{\xy  
    (0,1)="B1"      *{\scriptscriptstyle #5};
    (0,2.5)="B2"    *{\scriptscriptstyle \mathbf{#6}};
    (1,3)="B3"      *{\scriptscriptstyle #7};
    (0,0)="p"       *{\scriptstyle {\cdots}};
    (0,-1)="A1"     *{\scriptscriptstyle #4};
    (0,-2.5)="A2"   *{\scriptscriptstyle #3};
    (-1,-3)="A3"    *{\scriptscriptstyle #1};
    (1,-3)="A4"     *{\scriptscriptstyle #2};
    \endxy}
\newcommand{\vctEse}[6]{\xy
    (0,0)="B1"          *{\scriptscriptstyle #1};
   (-1.3,0)="B2"        *{\scriptscriptstyle #2};
    (1.3,0)="B3"        *{\scriptscriptstyle #3};
  (-1.3,-1.5)="B4"      *{\scriptscriptstyle #4};
    (1.3,-1.5)="B5"     *{\scriptscriptstyle #5};
    (0,1.5)="B6"        *{\scriptscriptstyle \mathbf{#6}};
    \endxy}
\newcommand{\vctEsi}[7]{\xy
    (0,0)="B1"          *{\scriptscriptstyle #1};
   (-1.3,0)="B2"        *{\scriptscriptstyle #2};
    (1.3,0)="B3"        *{\scriptscriptstyle #3};
  (-1.3,-1.5)="B4"      *{\scriptscriptstyle \mathbf{#4}};
    (1.3,-1.5)="B5"     *{\scriptscriptstyle #5};
    (0,1.5)="B6"        *{\scriptscriptstyle #6};
    (1.3,-3)="B7"     *{\scriptscriptstyle #7};
    \endxy}
\newcommand{\vctEo}[8]{\xy
    (0,0)="B1"          *{\scriptscriptstyle #1};
   (-1.3,0)="B2"        *{\scriptscriptstyle #2};
    (1.3,0)="B3"        *{\scriptscriptstyle #3};
  (-1.3,-1.5)="B4"      *{\scriptscriptstyle #4};
    (1.3,-1.5)="B5"     *{\scriptscriptstyle #5};
    (0,1.5)="B6"        *{\scriptscriptstyle #6};
    (1.3,-3)="B7"       *{\scriptscriptstyle #7};
    (1.3,-4.5)="B7"     *{\scriptscriptstyle \mathbf{#8}};
    \endxy}
\newcommand{\CoefMaxDnO}{\xy  
    (  4,3.5)="m1" *{};
    (  4,-27)="m2" *{};
    ( -4,-27)="m3" *{};
    ( -4,3.5)="m4" *{};
    ( 0,  3.5)="A" *{};
    ( -1,  0)="B1" *{{}_{\bullet}};
    (  1,  0)="B2" *{{}_{\bullet}};
    ( -1, -1)="B'1" *{};
    (  1, -1)="C'1" *{};
    ( -1, -6)="B'2" *{};
    ( 1, -6)="C'2" *{};
    (0,-10.5)="dot" *{\cdots};
    ( -1,-15)="B'4" *{};
    (  1,-15)="C'4" *{};
    (-1,-20)="B'5" *{};
    ( 1,-20)="C'5" *{};
    (-1,-22)="B''5" *{};
    ( 1,-22)="C''5" *{};
    ( -4,-24.6)="B6" *{};
    ( 1,-27)="B7" *{};
    "A";"B1" **@{-};
    "A";"B2" **@{-};
    "B'1";"B'2" **@{-};
    "C'1";"C'2" **@{-};
    "B'4";"B'5" **@{-};
    "C'4";"C'5" **@{-};
    "B''5";"B6" **@{-};
    "C''5";"B7" **@{-};
   \endxy}
\newcommand{\CoefXUDn}{\xy 
    (  4,3.5)="m1" *{};
    (  4,-27)="m2" *{};
    ( -4,-27)="m3" *{};
    ( -4,3.5)="m4" *{}; 
    ( 0,  0)="B1" *{{}_{\bullet}};
    ( 0, -1)="B'1" *{};
    ( 0, -6)="B'2" *{};
    (0,-10.5)="dot" *{\cdots};
    ( 0,-15)="B'4" *{};
    ( 0,-20)="B'5" *{};
    ( 0,-22)="B''5" *{};
    ( 0,-27)="B6" *{};
    "B'1";"B'2" **@{-};
    "B'4";"B'5" **@{-};
    "B''5";"B6" **@{-};
   \endxy}
\newcommand{\CoefXDDn}{\xy
    (  4,3.5)="m1" *{};
    (  4,-27)="m2" *{};
    ( -4,-27)="m3" *{};
    ( -4,3.5)="m4" *{};  
    ( 0,  0)="B1" *{{}_{\bullet}};
    ( 0, -1)="B'1" *{};
    ( 0, -6)="B'2" *{};
    (0,-10.5)="dot" *{\cdots};
    ( 0,-15)="B'4" *{};
    ( 0,-20)="B'5" *{};
    ( 0,-22)="B''5" *{};
    ( -3,-24.6)="B6" *{};
    "B'1";"B'2" **@{-};
    "B'4";"B'5" **@{-};
    "B''5";"B6" **@{-};
   \endxy}
\newcommand{\CoefXTuDn}{\xy
    (  4,3.5)="m1" *{};
    (  4,-27)="m2" *{};
    ( -4,-27)="m3" *{};
    ( -4,3.5)="m4" *{};  
    ( 3,  3.5)="A" *{};
    ( 0,  0)="B1" *{{}_{\bullet}};
    ( 0, -1)="B'1" *{};
    ( 0, -6)="B'2" *{};
    (0,-10.5)="dot" *{\cdots};
    ( 0,-15)="B'4" *{};
    (0,-20)="B'5" *{};
    (0,-22)="B''5" *{};
    ( -3,-24.6)="B6" *{};
    ( 0,-27)="B7" *{};
    "A";"B1" **@{-};
    "B'1";"B'2" **@{-};
    "B'4";"B'5" **@{-};
    "B''5";"B6" **@{-};
    "B''5";"B7" **@{-};
   \endxy}
\newcommand{\CoefXTdDn}{\xy
    (  4,3.5)="m1" *{};
    (  4,-27)="m2" *{};
    ( -4,-27)="m3" *{};
    ( -4,3.5)="m4" *{};  
    ( 3,  3.5)="A" *{};
    ( 0,  0)="B1" *{{}_{\bullet}};
    ( 0, -1)="B'1" *{};
    ( 0, -6)="B'2" *{};
    (0,-10.5)="dot" *{\cdots};
    ( 0,-15)="B'4" *{};
    (-1,-20)="B'5" *{};
    ( 1,-20)="C'5" *{};
    (-1,-22)="B''5" *{};
    ( 1,-22)="C''5" *{};
    ( -4,-24.6)="B6" *{};
    ( 1,-27)="B7" *{};
    "A";"B1" **@{-};
    "B'1";"B'2" **@{-};
    "B'4";"B'5" **@{-};
    "B'4";"C'5" **@{-};
    "B''5";"B6" **@{-};
    "C''5";"B7" **@{-};
   \endxy}
\newcommand{\CoefXTnDn}{\xy  
    (  4,3.5)="m1" *{};
    (  4,-27)="m2" *{};
    ( -4,-27)="m3" *{};
    ( -4,3.5)="m4" *{};
    ( 3,  3.5)="A" *{};
    ( 0,  0)="B1" *{{}_{\bullet}};
    ( 0, -1)="B'1" *{};
    ( -1, -6)="B'2" *{};
    ( 1, -6)="C'2" *{};
    (0,-10.5)="dot" *{\cdots};
    ( -1,-15)="B'4" *{};
    (  1,-15)="C'4" *{};
    (-1,-20)="B'5" *{};
    ( 1,-20)="C'5" *{};
    (-1,-22)="B''5" *{};
    ( 1,-22)="C''5" *{};
    ( -4,-24.6)="B6" *{};
    ( 1,-27)="B7" *{};
    "A";"B1" **@{-};
    "B'1";"B'2" **@{-};
    "B'1";"C'2" **@{-};
    "B'4";"B'5" **@{-};
    "C'4";"C'5" **@{-};
    "B''5";"B6" **@{-};
    "C''5";"B7" **@{-};
   \endxy}
\newcommand{\CoefYUDn}{\xy
    (  4,3.5)="m1" *{};
    (  4,-27)="m2" *{};
    ( -4,-27)="m3" *{};
    ( -4,3.5)="m4" *{};  
    ( 3,  3.5)="A" *{};
    ( 0,  0)="B1" *{{}_{\bullet}};
    ( 0, -1)="B'1" *{};
    ( 0, -6)="B'2" *{};
    (0,-10.5)="dot" *{\cdots};
    ( 0,-15)="B'4" *{};
    ( 0,-20)="B'5" *{};
    ( 0,-22)="B''5" *{};
    ( -3,-24.6)="B6" *{};
    "A";"B1" **@{-};
    "B'1";"B'2" **@{-};
    "B'4";"B'5" **@{-};
    "B''5";"B6" **@{-};
   \endxy}
\newcommand{\CoefYDDn}{\xy
    (  4,3.5)="m1" *{};
    (  4,-27)="m2" *{};
    ( -4,-27)="m3" *{};
    ( -4,3.5)="m4" *{};  
    ( -3,3.5)="A" *{};
    ( 0,  0)="B1" *{{}_{\bullet}};
    ( 0, -1)="B'1" *{};
    ( 0, -6)="B'2" *{};
    (0,-10.5)="dot" *{\cdots};
    ( 0,-15)="B'4" *{};
    ( 0,-20)="B'5" *{};
    ( 0,-22)="B''5" *{};
    ( 0,-27)="B6" *{};
    "A";"B1" **@{-};
    "B'1";"B'2" **@{-};
    "B'4";"B'5" **@{-};
    "B''5";"B6" **@{-};
   \endxy}
\newcommand{\CoefYDDnPrime}{\xy
    (  4,3.5)="m1" *{};
    (  4,-27)="m2" *{};
    ( -4,-27)="m3" *{};
    ( -4,3.5)="m4" *{};  
    ( 0,3.5)="A" *{};
    ( 0,  0)="B1" *{{}_{\bullet}};
    ( 0, -1)="B'1" *{};
    ( 0, -6)="B'2" *{};
    (0,-10.5)="dot" *{\cdots};
    ( 0,-15)="B'4" *{};
    ( 0,-20)="B'5" *{};
    ( 0,-22)="B''5" *{};
    ( 0,-27)="B6" *{};
    "A";"B1" **@{-};
    "B'1";"B'2" **@{-};
    "B'4";"B'5" **@{-};
    "B''5";"B6" **@{-};
   \endxy}
\newcommand{\CoefYTuDn}{\xy
    (  4,3.5)="m1" *{};
    (  4,-27)="m2" *{};
    ( -4,-27)="m3" *{};
    ( -4,3.5)="m4" *{};  
    ( 0,  0)="B1" *{{}_{\bullet}};
    ( 0, -1)="B'1" *{};
    ( 0, -6)="B'2" *{};
    (0,-10.5)="dot" *{\cdots};
    ( 0,-15)="B'4" *{};
    ( 0,-20)="B'5" *{};
    "B'1";"B'2" **@{-};
    "B'4";"B'5" **@{-};
   \endxy}
\newcommand{\CoefYTdDn}{\xy
    (  4,3.5)="m1" *{};
    (  4,-27)="m2" *{};
    ( -4,-27)="m3" *{};
    ( -4,3.5)="m4" *{};  
    ( 0,  0)="B1" *{{}_{\bullet}};
    (0,-10.5)="dot" *{\cdots};
    ( 0, -1)="B'1" *{};
    ( 0, -6)="B'2" *{};
    "B'1";"B'2" **@{-};
   \endxy}
\newcommand{\CoefYTnDn}{\xy
    (  4,3.5)="m1" *{};
    (  4,-27)="m2" *{};
    ( -4,-27)="m3" *{};
    ( -4,3.5)="m4" *{};  
    ( 0,  0)="B1" *{{}_{\bullet}};
   \endxy}
\newcommand{\RaizM}{\xy
( -4,-5)="R" *{\mathcal{R}^{(-1)}_{1,1,1}[\ell]=};
(4,0)="Gr" *{\xymatrix@R=3pc@C=1pc{{\ell-1} \ar@<-.5ex>[d] \ar@<.5ex>[d] \\ {\ell}}};
\endxy}
\newcommand{\RaizC}{\xy
( -4,-5)="R" *{\mathcal{R}^{(0)}_{2,1,1}[\ell]=};
(4,0)="Gr" *{\xymatrix@!0@R=2pc@C=1pc{
{\ell} \ar@<-.5ex>[dd] \ar@<.5ex>[dd] 
\ar@/^10pt/[rrddd]
\\ \\ {\ell} \ar@{.>}[rrd]  \\ & &  1 }};
\endxy}
\newcommand{\RaizUc}{\xy
( -4,-5)="R" *{\mathcal{R}^{(1)}_{1,1,1}[\ell]=};
(4,0)="Gr" *{\xymatrix@R=3pc@C=1pc{{\ell+1} \ar@<-.5ex>[d] \ar@<.5ex>[d] \\ {\ell}}};
\endxy}
\newcommand{\RaizUu}{\xy
( -4,-5)="R" *{\mathcal{R}^{(1)}_{2,1,1}[\ell]=};
(4,0)="Gr" *{\xymatrix@!0@R=2pc@C=1pc{
{\ell+1} \ar@<-.5ex>[dd] \ar@<.5ex>[dd] 
\ar@/^10pt/[rrddd]
\\ \\ {\ell} \ar@{.>}[rrd]  \\ & &  1 }};
\endxy}
\newcommand{\RaizUd}{\xy
( -4,-5)="R" *{\mathcal{R}^{(1)}_{2,1,1}[\ell]=};
(4,0)="Gr" *{\xymatrix@!0@R=1pc@C=1pc{
{\ell+1} \ar@<-.5ex>[dd] \ar@<.5ex>[dd] 
\ar@/^10pt/[rrddd] \ar@/^20pt/[rrrrddddd]
\\ \\ {\ell} \ar@{.>}[rrd] \ar@/_6pt/@{.>}[rrrrddd]  \\ & &  1 \\ \\ & & & & 1 }};
\endxy}
\newcommand{\RaizUt}{\xy
( -4,-5)="R" *{\mathcal{R}^{(1)}_{2,2,2}[\ell]=};
(4,0)="Gr" *{\xymatrix@!0@R=1pc@C=1pc{
{\ell+1} \ar@<-.5ex>[dd] \ar@<.5ex>[dd]
\ar@/^10pt/[rrddd] \ar@/^15pt/[rrrdddd] \ar@/^20pt/[rrrrddddd]
\\ \\ {\ell} \ar@{.>}[rrd] \ar@/_3pt/@{.>}[rrrdd] \ar@/_6pt/@{.>}[rrrrddd]  \\ & &  1 \\ & & & 1 \\ & & & & 1 }};
\endxy}
\newcommand{\Renglon}[5]{\xy  
    (3,0)="B"        *{\scriptscriptstyle #5};
    (1.7,0)="B1"     *{\scriptscriptstyle #4};
    (0,0)="B2"       *{\scriptscriptstyle #3};
    (-1.7,0)="B3"    *{\scriptscriptstyle #2};
    (-3,0)="B4"      *{\scriptscriptstyle #1};
    \endxy}
\newcommand{\renglon}[3]{\xy  
    (1.7,0)="B1"     *{\scriptscriptstyle #3};
    (0,0)="B2"       *{\scriptscriptstyle #2};
    (-1.7,0)="B3"    *{\scriptscriptstyle #1};
    \endxy}

\begin{document}

\renewcommand{\bibname}{Bibliography}
\renewcommand{\contentsname}{Contents}
\renewcommand{\abstractname}{Abstract}
\renewcommand{\chaptername}{Chapter}
\renewcommand{\figurename}{Figure}
\renewcommand{\tablename}{Table}
\renewcommand{\appendixname}{Appendix}

% Title Page
%---------------------------------------------
%---------------------------------------------
%---------------------------------------------
\begin{titlepage}

\begin{adjustwidth}{-2cm}{-2cm}
%------------------------------
\begin{center}
\begin{tabular}{lcr}
\begin{tabular}{l}
%\\
\includegraphics[angle=0,width=0.1\linewidth]{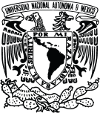}
\end{tabular}
& \textbf{\begin{tabular}{c}
{\large Posgrado Conjunto en Ciencias Matem\'aticas}\\[.3cm]
{\normalsize Universidad Nacional Aut\'onoma de M\'exico}\\[.1cm]
{\normalsize Universidad Michoacana de San Nicol\'as Hidalgo}\\
\end{tabular}} &
\begin{tabular}{l}
%\\
\includegraphics[angle=0,width=0.1\linewidth]{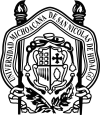}
\end{tabular}
\end{tabular}
\end{center}
%------------------------------
\end{adjustwidth}
%------------
\vspace{2.5cm}
\begin{center}
{\Large \textsc{Geometry in the indecomposable\\
modules of a hereditary algebra}}
\end{center}
%------------
\vspace{2cm}
\begin{center}
{Thesis submitted for the degree of\\
Doctor in Mathematical Sciences}
\end{center}
%------------
\vspace{1cm}
\begin{center}
{BY}
\end{center}
%------------
\begin{center}
{Jes\'us Arturo Jim\'enez Gonz\'alez}
\end{center}
%------------
\vspace{1cm}
\begin{center}
{THESIS ADVISOR}
\end{center}
%------------
\begin{center}
{Dr. Raymundo Bautista Ramos}
\end{center}
%------------
\vspace{3cm}
\begin{center}
{\small Morelia, Michoac\'an. \hspace{5cm} 2014}
\end{center}
\end{titlepage}
%---------------------------------------------
%---------------------------------------------
%---------------------------------------------

\tableofcontents

%\listoffigures
%\listoftables
%\chapter*{Agradecimientos}

%\begin{abstract}
%En esta tesis damos de manera pr\'actica y sistem\'atica formas matriciales para representaciones excepcionales de \'algebras hereditarias.
%Nuestra herramienta principal es la \emph{reducci\'on inducida por m\'odulos} (como se presenta en 
%Bautista-Salmer\'on-Zuazua \cite{BSZ09}), el correspondiente categ\'orico a las formas inducidas por ra\'ices 
%(ver por ejemplo Barot-de la Pe\~na \cite{BP06}). 
%
%En el cap\'itulo \ref{CH:meth} damos representaciones de m\'odulos posproyectivos
%de un diagrama de Dynkin extendido $\widetilde{\Delta}$. El m\'etodo a seguir se elemental.
%Consiste en, por un lado, especificar los $\Delta$-m\'odulos que conforman la restricci\'on $M|_{\Delta}$ de una
%$\widetilde{\Delta}$-representaci\'on $M$. Estos son representaciones de diagramas de Dynkin cuyas formas matriciales son bien 
%conocidas. Por otro lado, damos expl\'icitamente una transformaci\'on lineal
%$M(\alpha_0)$ (dos transformaciones en caso $\widetilde{\mathbf{A}}_n$) que pega estos $\Delta$-m\'odulos en una
%$\widetilde{\Delta}$-representaci\'on isomorfa a $M$. De forma an\'aloga se obtiene un carcaj de coeficientes para $M$ 
%pegando los correspondientes carcajes de la restricci\'on (\emph{bases distinguidas}, Ringel \cite{cmR98} y \cite{cmR12}).
%\end{abstract}

\pagenumbering{roman}
%----------------------------------------------------------------------
%----------------------------------------------------------------------
\chapter*{Introduction.}
\label{Cap(I)}
\addcontentsline{toc}{chapter}{Introduction}

Let $k$ be an arbitrary field and consider the path $k$-algebra $A$ of a finite quiver $Q$.
A fundamental problem in the representation theory of associative algebras is to explicitly give,
through collections of matrices, indecomposable finite dimensional $A$-representations.
In general this is a difficult problem, which has stimulated the theory and that is far from finding
a satisfactory answer. When $Q$ corresponds to a Dynkin diagram (with a fixed orientation of 
its arrows) Gabriel exhibited representations whose isomorphism classes are a complete list of 
indecomposable $A$-modules~\cite[1972]{pG72}. The problem was revisited by Ringel in 1998~\cite{cmR98}. 
For the extended Dynkin diagrams, main topic of this thesis, there are several partial results.
In 1890, when solving a problem posed by Weierstrass, Kronecker classified what nowadays are called
representations of the quiver with two arrows $\xymatrix{\cdot \ar@<.5ex>[r] \ar@<-.5ex>[r] & \cdot}$ 
(corresponding to the diagram $\widetilde{\mathbf{A_1}}$) whose path algebra is called Kronecker algebra.
The four subspace problem, which can be interpreted in terms of a quiver on the diagram
$\widetilde{\mathbf{D_4}}$, has been analized by Nazarova in 1967 \cite{aN67}, Gelfand and
Ponomarev in 1972 \cite{GP72}, and Medina and Zavadskij in 2004 \cite{MZ04}. The case $\widetilde{\mathbf{A_n}}$
was described by Gabriel and Roiter in 1997 \cite{GR97}, who also gave an algorithm to determine
regular representations of $\widetilde{\mathbf{E_6}}$ (for algebraically closed fields). 
This method can be easily adapted to any extended
Dynkin diagram $\widetilde{\Delta}$. The posprojective and preinjective representations of $\widetilde{\mathbf{D_n}}$ 
were given by Kussin and Meltzer in 2006 \cite{KM}, as well as the series of rank three representations for 
$\widetilde{\mathbf{E_6}}$. Their method makes use of tilting theory and the explicit knowledge of some
representations for the domestic canonical algebras (as defined by Ringel in \cite{cmR}) given by Komoda, 
Kussin and Meltzer in \cite{KM07} and \cite{KM08}. With the same method Kedzierski and Meltzer described in 2011
the series of maximal rank (six) representations of $\widetilde{\mathbf{E_8}}$ \cite{KM11}.

Exceptional modules (indecomposable representations without self-ex\-ten\-sions) 
have been studied from distinct perspectives and with different purposes. Exceptional sequences,
noncrossing partitions, tilting theory and all the reduction processes in this work are examples where the
exceptional modules play a main role. As is well known, the isomorphism classes of exceptional modules are
determined by their dimension vectors. Moreover, Ringel has shown~\cite{cmR98} that the exceptional representations 
can be given by matrices such that the number of nonzero coefficients is the minimal expected (tree modules).
In that case one can guarantee that the matrices contain only coefficients $0$ and $1$.

One way to deal with the problem of finding quiver representations is to enlarge the algebraic context from which
the representations are taken. With that in mind, and by the own importance of these topics, it has been considered 
representations of quivers with relations (Kussin and Meltzer and their representations of canonical algebras), 
of partially ordered sets (from the work initiated by Zavadskij), matrix problems (Nazarova) among others.
As proposed by Professor Bautista, advisor of this thesis, within the context of ditalgebras the
\emph{reduction by admissible modules} (as developed by Bautista, Salmer\'on and Zuazua in \cite{BSZ09})
proves to be a practical tool to systematically give quiver representations. Of particular interest,
considering the historical background above, is the case of extended Dynkin diagrams and their posprojective and preinjective
components. In response to this interest we establish the objectives of this work.
Let $A_0$ be the path algebra of a Dynkin quiver and $A$ the one-point extension of $A_0$.

\begin{quote}
\textbf{Objectives.} To describe the restriction to $A_0$-mod of indecomposable posprojective and preinjective 
$A$-modules in terms of their positions inside the Auslander-Reiten quiver. Furthermore, to determine certain reductions of the 
algebra $A$ with respect to admissible $A_0$-modules and to show how by means of the reduction functors one can
explicitly construct series of exceptional $A$-modules.
\end{quote}

In general terms, in chapter one we prove preliminar results in the context of ditalgebras
and discuss some general aspects of the Auslander-Reiten theory. In chapter two we analize
two key elements in the construction of $A$-modules: the module category of the Kronecker algebra
and the category of $A_0$-modules. In chapters three and four we study one-point extensions,
define and analize the rank of an $A$-representation and construct exceptional $A$-modules
starting from certain reduced representations of minor rank. In the appendix we give fundamental definitions
and results from the theory of ditalgebras and their module categories, as presented by Bautista, Salmer\'on 
and Zuazua \cite{BSZ09}. In what follows we give a detailed description by chapters.

%--------------------
\textbf{Chapter 1.} In section~\ref{(P)S:carAlg} we show some relations between integral bilinear forms
and the module category of a (bi)quiver with differential. Some notions and results necessary for the rest
of the work are established here.

In particular we define the Euler characteristic and extend the notion of Cartan matrix to the setting of
positively graded algebras. This notion is compatible with the Euler characteristic
(lemmas~\ref{(P)L:invers} and~\ref{(P)L:Euler}) and with the construction of the Coxeter matrix given by means of
simple reflections (section~\ref{(P)S:reflex}). In sections~\ref{(P)S:redElem} and~\ref{(P)S:redMod} we present
elementary tools for the construction of ditalgebras, such as regularization and reduction of an edge. Moreover, we generalize
some known properties of exceptional modules to the context of ditalgebras (for instance, every exceptional module is a tree module,
proposition~\ref{(P)P:excep}). In section~\ref{(DE)S:SCD} we present some classical results of the Auslander-Reiten theory
following Ringel \cite{cmR}. In~\ref{(DE)S:trasl} we show definitions and results about translation quivers,
sections and cosections. These are fundamental results for the proofs in section~\ref{(DE)S:rank}. 
Finally, in~\ref{(DE)S:comp} we describe the posprojective and preinjective components in the Auslander-Reiten
quiver of the category of quiver representations.

%--------------------
\textbf{Chapter 2.} We construct exceptional representations for the classical Kronecker algebra
(section~\ref{(P)S:KroDos}) and the first generalized case (section~\ref{(P)S:KroTres}). 
Although these representations are well known (cf. Ringel~\cite{cmR98} and~\cite{cmR10}), 
we provide systematic proofs that ilustrate the tools presented in the preliminars chapter. Moreover, the
explicit knowledge of the Auslander-Reiten quiver of the classical Kronecker algebra is essential for the construction
of the extended category $A$-mod. In section~\ref{(DE)S:ARDynkin} we study the Auslander-Reiten quiver
of $A_0$-mod following Ringel~\cite{cmR}.
%---------------------------------------------------------------------- 
\begin{displaymath}
\xy 0;/r.20pc/:
( 0,0)="A0" *{};
(20,0)="A2" *{};
(30,0)="A3" *{};
(40,0)="A4" *{};
(50,0)="A5" *{};
(60,0)="A6" *{};
(90,0)="A9" *{};
( 0,30)="B0" *{};
(30,30)="B3" *{};
(40,30)="B4" *{};
(50,30)="B5" *{};
(90,30)="B9" *{};
( 0,10)="C1" *{};
( 0,20)="C2" *{};
(90,10)="D1" *{};
(90,20)="D2" *{};
(25,5)="E1" *{};
(30,10)="E2" *{};
(35,15)="E3" *{};
(45,25)="E4" *{};
(35,25)="F1" *{};
(45,15)="F2" *{};
(50,10)="F3" *{};
(55,5)="F4" *{};
(40,20)="WB" *{{}_{\bullet}};
(46,20)="W" *{[W_0]};
(70,25)="CR" *{\ConR([W_0])};
(15,25)="CR" *{\ConL([W_0])};
%-----------------------
"A0";"A2" **@{-};
"A6";"A9" **@{-};
"A9";"B9" **@{-};
"A0";"B0" **@{-};
"B0";"B3" **@{-};
"B5";"B9" **@{-};
"B3";"A6" **@{-};
"B5";"A2" **@{-};
"B5";"A2" **@{-};
"E1";"A3" **@{-};
"E2";"A4" **@{-};
"E3";"A5" **@{-};
"E4";"B4" **@{-};
"F1";"B4" **@{-};
"F2";"A3" **@{-};
"F3";"A4" **@{-};
"F4";"A5" **@{-};
\endxy
\end{displaymath}
We show that the isomorphism class of an indecomposable module $W_0$ corresponding to the maximal positive root 
of the Dynkin diagram is a wing vertex in the (finite) Auslander-Reiten quiver of $A_0$-mod. 
The wings of $[W_0]$ (middle part in the figure above) separate the Auslander-Reiten quiver
of $A_0$-mod (cf. lemma~\ref{(DE)L:alaDyn}).

%--------------------
\textbf{Chapter 3.} The first step in the description of the category $A$-mod is to identify it
with the category of subspaces $\check{\mathcal{U}}(A_0\text{-mod},|-|)$ with respect to the functor
\[
|-|= \Hom_{A_0}(R,-):A_0\text{-mod} \longrightarrow k\text{-mod},
\]
where $R$ is the extension module of $A=A_0[R]$ (section~\ref{(DE)S:extPunt}).
In this way an $A$-module $M=(M_0,M_{\omega},\gamma_M)$ consists in an $A_0$-module $M_0$,
a $k$-vector space $M_{\omega}$ and a linear transformation $\gamma_M:M_{\omega} \to |M_0|$.
In proposition~\ref{(DE)P:extPuntRed} we analize reductions by admissible modules 
of one-point extensions. The lifting lemma (lemma~\ref{(DE)L:levAR}) relates almost split sequences in 
the categories $A_0$-mod and $A$-mod. One of our goals is to directly prove, without any knowledge of the regular
components, that the liftings of the almost split sequences in $\ConL([W_0])$ and $\ConR([W_0])$ (see figure above) lie in 
the posprojective and preinjective components of $A$-mod res\-pectively. This is shown in theorem~\ref{(DE)T:finito}.
In section~\ref{(DE)S:levantaFun} we give some functorial properties of the universal mappings and their liftings 
to $A$-mod, used by Ringel~\cite{cmR} in the point (17) of his main theorem~3.4.
In section~\ref{(DE)S:rank} we exhibit (through the reduction process) the Kronecker subcategory contained in $A$-mod
determined by the objects in $A$-mod whose restriction to $A_0$-mod is direct summ of copies of $W_0$. 
We also give answer to the first objective by means of the functors of section~\ref{(DE)S:levantaFun} (propositions~\ref{(DE)C:invar}
and~\ref{(DE)P:rank}, see also the table at the end of section~\ref{(DE)S:rank}).  We finally use 
modules determined by the left side $\dConL([W_0])$ and right side $\dConR([W_0])$ of the wings of $[W_0]$, to construct 
reduced ditalgebras $\mathcal{A}^X$ and $\mathcal{A}^Y$ of the extended Dynkin algebra $A$ (section~\ref{(DE)S:algRed}). 
Their interest lies in the fact that almost all isomorphism classes (up to a finite number) of posprojective and preinjective
modules belong in the image of the reduction functors $F^X$ and $F^Y$ respectively (theorem~\ref{(DE)T:finito}). 
Moreover, the reduced tensor algebras $A^X$ and $A^Y$ depend only on the tubular type of $A$, and their associated quadratic
forms coincide $q^{xy}$.  At the end of section~\ref{(DE)S:algRed} we can find a list of positive roots of the reduced 
quadratic form and some corresponding exceptional representations.

%--------------------
\textbf{Chapter 4.} This last chapter is dedicated to the presentation of exceptional modules of
extended Dynkin algebras through the reduction functors. Starting from the knowledge of the restriction to
$A_0$-mod of an $A$-module $M=(M_0,M_{\omega},\gamma_{M})$, we directly construct a coefficient quiver of $M$
from the coefficient quivers of the direct summands $Z_0$ of the restriction $M_0$ (depicted by $\CoefMax$, $\CoefAlaU{i}$ 
and $\CoefAlaV{j}$; the marked vertices correspond to a basis of the $k$-vector space $\Hom_{A_0}(R,Z_0)$, which has dimension
two when $Z_0\cong W_0$ and one in all other cases) and the matrices corresponding to the transformation $\gamma_M$. 
The matrices which conform $\gamma_M$ are stored in the reduced representation $\widetilde{M}$ (that is, $\widetilde{M}$ 
is a representation of the reduced ditalgebra $\mathcal{A}^Y$ such that $F^Y(\widetilde{M})\cong M$). For instance, 
in the following figure we show the coefficient quiver $\mathcal{C}(\widetilde{M})$ of a reduced $\mathcal{A}^Y$-representation 
$\widetilde{M}$ (left). Since the restriction of $M=F^Y(\widetilde{M})$ to $A_0$-mod has the form $Y \otimes \widetilde{M}$, 
the coefficient quiver of $M$ can be obtained by substituting in $\mathcal{C}(\widetilde{M})$ the coefficient quivers 
chosen for the direct summands of the reduction $A_0$-module $Y$ (right). 
\[
\xy
(-45,8)="L1" *{};
( 30,8)="L2" *{};
(-31,18)="Coef1" *{\scriptstyle \mathcal{C}(\widetilde{M})};
( 12,18)="Coef2" *{\scriptstyle \mathcal{C}(M)};
(-10,1)="FY" *{\scriptstyle F^Y};
(-11,-2)="mps" *{\mapsto};
(-30,-2)="Mred" *{\xy
    %-----
    ( 12,-15)="A1" *{};
    ( 16,-15)="A2" *{};
    ( 20,-15)="A3" *{};
    ( 24,-15)="A4" *{};
    ( 28,-15)="A5" *{};
    ( 10,-2)="B0" *{};
    ( 14,-2)="B1" *{};
    ( 18,-2)="B2" *{};
    ( 22,-2)="B3" *{};	
    ( 26,-2)="B4" *{};
    ( 10,-30)="C0" *{};
    ( 26,-35)="C4" *{};
    %-----------------------
    "A1";"B1" **@{-};
    "A2";"B2" **@{-};
    "A3";"B3" **@{-};
    "A4";"B4" **@{-};
    "B0";"A1" **@{-};
    "B1";"A2" **@{-};
    "B2";"A3" **@{-};
    "B3";"A4" **@{-};
    "B0";"C0" **@{-};
    "B4";"C4" **@{-};
    \endxy};
%-------------------------
( 0,-12)="E1"  *{\CoefAlaU{i}};
( 24,-18)="E5" *{\CoefAlaV{j}};
( 0,15)="A1" *{};
( 6,15)="A2" *{};
(12,15)="A3" *{};
(18,15)="A4" *{};
(24,15)="A5" *{};
( 2, 0)="B1" *{};
( 8, 0)="B2" *{};
(14, 0)="B3" *{};
(20, 0)="B4" *{};
(26, 0)="B5" *{};
( 4, 0)="C1" *{};
(10, 0)="C2" *{};
(16, 0)="C3" *{};
(22, 0)="C4" *{};
(28, 0)="C5" *{};
( 3, 0)="D1"  *{\CoefMax};
( 9, 0)="D2"  *{\CoefMax};
( 15, 0)="D3" *{\CoefMax};
( 21, 0)="D3" *{\CoefMax};
%-----------------------
"A1";"B1" **@{-};
"A2";"B2" **@{-};
"A3";"B3" **@{-};
"A4";"B4" **@{-};
"A2";"C1" **@{-};
"A3";"C2" **@{-};
"A4";"C3" **@{-};
"A5";"C4" **@{-};
"A1";"E1" **@{-};
"A5";"E5" **@{-};
"L1";"L2" **@{--};
\endxy
\]
The arrows that cross the dotted line correspond to the matrices of the transformation $\gamma_M$.
Deleting these arrows we obtain, on the left side, six isolated vertices (corresponding to the coefficient quiver
of a module over a semi-simple algebra $S$), and on the right side, the coefficient quivers of the restriction $M_0$.
Observe that the reduced representations $\widetilde{M}$ contain a preinjective Kronecker module 
(in the example shown, with dimension vector $(4,5)$). Corresponding coefficient quivers for dimension vectors
of the form $(\ell,\ell+1)$ produce a series of (rank one) indecomposable preinjective $A$-modules.
In section~\ref{(B)S:1} we show the case 
$\widetilde{\mathbf{A_n}}$ for an arbitrary (noncyclic) ordering of its arrows
(compare with theorem~11.1 of Gabriel and Roiter~\cite{GR97}). 
In section~\ref{(B)S:2} we analize the case $\widetilde{\mathbf{D_n}}$ for a particular ordering of its arrows
(the one used by Kussin and Meltzer in~\cite{KM}), while in section~\ref{(B)S:mayor} we exhibit the reduced ditalgebras
for the cases $\widetilde{\mathbf{E_m}}$ ($m=6,7,8$).

\pagenumbering{arabic}

%----------------------------------------------------------------------
%----------------------------------------------------------------------
%----------------------------------------------------------------------
\chapter{Preliminars.}
\label{Cap(P)}
%----------------------------------------------------------------------
%----------------------------------------------------------------------

The general setting of this work are differential tensor algebras, as developed by Bautista, Salmer\'on
and Zuazua in \cite{BSZ09}. In this chapter we establish some preliminar results. All topics in the following
sections are known at least in the classical case. Definitions and results of the theory of differential
tensor algebras and their representations can be found in the appendix.

\section{Quivers and algebras.} \label{(P)S:carAlg}
%----------------------------------------------------------------------
The basic combinatorial structure of this thesis is the oriented graph with two kinds of arrows, or biquiver,
which will be called simply quiver. 
A \textbf{quiver} $Q=(Q_0,Q_1,s,t,|\cdot|)$ consists of two sets, one of \textbf{vertices} $Q_0$ and other of \textbf{arrows} $Q_1$,
together with functions $s,t: Q_1 \to Q_0$ called \textbf{source} and \text{target} respectively, and a 
\textbf{degree} fuction $|\cdot|:Q_1 \to \{0,1\}$. A \textbf{morphism} $f$ between quivers $Q$ y $Q'$ consists of
functions $f_0:Q_0 \to Q'_0$ and $f_1:Q_1 \to Q'_1$ which commute with the source and target functions and which preserve
degrees. We will usually give a quiver in graphical form with a solid arrow $\xymatrix@C=1pc{s(\alpha) \ar[r] & t(\alpha)}$ 
for each $\alpha \in Q_1$ with $|\alpha|=0$ and a dotted arrow $\xymatrix@C=1pc{s(\alpha) \ar@{.>}[r] & t(\alpha)}$ if $|\alpha|=1$. 
A quiver with all its arrow having zero degree is called \textbf{solid quiver}. 
A quiver is \textbf{finite} if $Q_0$ and $Q_1$ are finite sets. The \textbf{direct successors} of a vertex $i$ is the set
of vertices $j$ for which there existes an arrow from $i$ to $j$. Its \textbf{direct predecessors} $j$ are those vertices for
which there exists an arrow from $j$ to $i$. A \textbf{tree} is a quiver with no cycles. We will say that a quiver is
\textbf{regular} if it does not contain subquivers of the following types
\begin{displaymath}
 \xymatrix{
\bullet \ar@<.5ex>[rr] \ar@{.>}@<-.5ex>[rr] & & \bullet & & \bullet \ar@(ul,dl)[] \ar@{.>}@(ur,dr)[]
}
\end{displaymath}
An ordering of the set of vertices $Q_0=\{1,\ldots,n\}$ will be \textbf{admissible} if for any arrow $\alpha$ the inequality
$s(\alpha) > t(\alpha)$ holds. Clearly every tree han an admissible ordering of its vertices.
We associate to each quiver $Q$ two algebraic structures, a $k$-algebra $kQ$ (for an arbitrary field $k$) and
an integral bilinear form $\langle \cdot,\cdot \rangle_Q$. A fundamental result in the representation theory of 
associative algebras states that the module category of $kQ$ is controled by the bilinear form of $Q$
(see for example~\cite{cmR}). The following is a short survey of this relation.

%-----------------------------------------------
\[
 \xymatrix@R=1pc@C=4pc{
& *+[F]\txt{$T$\\ \footnotesize Tensorial algebra} \ar[r]_-{\txt{\scriptsize $+ \delta$ \\ \scriptsize differential }} 
\ar@{-->}[dd]^-{\txt{\scriptsize Cartan \\ \scriptsize matrix}} 
&*+[F]{(T,\delta)\text{-mod}} \ar@{-->}@/_1pc/[l]_-{\varepsilon} \ar[dd]^-{\txt{\scriptsize Euler \\ \scriptsize charact.}} \\
*+[F]\txt{$Q$\\ \footnotesize Quiver} \ar@{<->}[ru]^-{\beta} \ar[rd]^-{\alpha} \\
&*+[F]\txt{$M$\\ \footnotesize Integral square \\ \footnotesize matrix} \ar@/^1pc/[lu]^-{\overline{\alpha}} \ar@{<->}[r]_-{\gamma} 
&*+[F]\txt{$\mathbb{Z}^n \times \mathbb{Z}^n \to \mathbb{Z}$ 
\\ \footnotesize Integral bilinear \\ \footnotesize form}  \\
}
\]

%-----------------------------------------------
Associated to any integral matrix $M=(m_{ij})_{i,j=1}^n$ there is a finite quiver as follows. The set of vertices
$Q(M)_0$ is the set $\{ 1,\ldots,n\}$. For $i \neq j$, if $m_{ij} > 0$ we add $m_{ij}$ dotted arrows from $i$ to $j$ and
if $m_{ij} < 0$ we add $-m_{ij}$ solid arrows from $i$ to $j$. Moreover, if $m_{ii}>1$ we add $m_{ii}-1$ dotted loops
in the vertex $i$, and if $m_{ii}<1$ we add $1-m_{ii}$ solid loops in $i$. 
Observe that $Q(M)$ is always a regular quiver.

On the other hand, the \textbf{incidence matrix} associated to a finite quiver $Q$ is the square matrix $M_Q=(m_{ij})$ 
whose entries are given by
\begin{equation*}
m_{ij} = \left\{
\begin{array}{l l}
-\sum_{\alpha:i \to j} (-1)^{|\alpha|}, & \text{if } i \neq j,\\
1-\sum_{\alpha:i \to i} (-1)^{|\alpha|}, & \text{if } i=j,
\end{array} \right.
\end{equation*}
(if the set of arrows from $i$ to $j$ is empty the sum $\sum_{\alpha:i \to j}$ takes the value zero). 
Clearly every integral square matrix $M$ is the incidence matrix of a quiver, specifically, $M=M_{Q(M)}$. 
\begin{lema} \label{(P)L:carReg}
 If $Q$ is a regular finite quiver then $Q \cong Q(M_Q)$.
\end{lema}
\bproof
By construction of $Q(M_Q)$ there is a bijection between vertices
\[
 f_0: Q(M_Q)_0 \to Q_0.
\]
We give by cases a bijection $f_1$ between the sets of arrows of $Q(M_Q)$ and $Q$ 
corresponding to the entry $m_{ij}$ of the matrix $M_Q$.\\
\underline{Case $i\neq j$.}
Assume that $m_{ij}>0$, so there are $m_{ij}$ dotted arrows in $Q(M_Q)$ form $i$ to $j$. 
Since $Q$ is regular and $m_{ij}=-\sum_{\alpha:i \to j} (-1)^{|\alpha|}$, there are exactly $m_{ij}$ arrows
from $i$ to $j$ en $Q$, all of them dotted. The case $m_{ij}<0$ is similar.
Hence we can take a bijection $f_1$ between the arrows in $Q(M_Q)$ from $i$ to $j$
and the arrws in  $Q$ from $i$ to $j$.\\
\underline{Case $i=j$.} Assume now that $i = j$ and that $m_{ij}>1$. 
Then there are $m_{ii}-1$ dotted loops in $Q(M_Q)$ over the vertex $i$. 
Since $Q$ is regular and $m_{ii}=1-\sum_{\alpha:i \to i} (-1)^{|\alpha|}$, there are exactly $m_{ii}$ loops
over $i$ in the quiver $Q$, all of them dotted. The case $m_{ii}<1$ is similar.
Hence we can take again a bijection $f_1$ between the loops in $Q(M_Q)$ in the vertex $i$
and the loops in $Q$ over $i$. 

This completes the proof since there are no more arrows in $Q(M_Q)$.
\eproof

%-----------------------------------------------
The \textbf{integral bilinear form} associated to an $n$ by $n$ matrix,
\[
\xymatrix{
\langle \cdot,\cdot \rangle_M: \mathbb{Z}^n \times \mathbb{Z}^n \longrightarrow \mathbb{Z},
}
\]
is given by $\langle x,y\rangle=x^tMy$. The matrix associated to a integral bilinear form has as coefficients
the values $m_{ij}=\langle \mathbf{e}_i,\mathbf{e}_j \rangle$ for $i,j=1,\ldots,n$, where $\mathbf{e}_1,\ldots,\mathbf{e}_n$ 
is the canonical basis in $\mathbb{Z}^n$. It is clear that all integral bilinear forms are determined by ther associated matrices
(relation $\gamma$ in the figure). Their \textbf{symmetrization} is $(x,y)=(1/2)(\langle x,y \rangle + \langle y,x \rangle$) and
the \textbf{associated quadratic form} is $q(x)=\langle x,x \rangle$. An integral vector $x$ such that $q(x)=1$ is called 
\textbf{root} of $q$. A vector is \textbf{positiv} whenever it is not zero and all its entries are nonnegative.

%----------------------------------------------
Let $k$ be an arbitrary field.
All graded algebras $T=\bigoplus T_i$ to be considered are freely generated by $A=T_0$ and $V=T_1$
(definition~\ref{(A)D:tens} in the appendix). Then $T$ is isomorphic to the tensorial algebra $T_A(V)$
(we will usually identify $T$ with $T_A(V)$). We will also require that $A \cong T_R(W_0)$ 
and that $V \cong A \otimes_R W_1 \otimes_R A$ (that is, $(R,W_0 \oplus W_1)$ is a  \textbf{layer} of $T$, 
see definition~\ref{(A)D:estrato}) for some trivial $k$-algebra $R$ (finite product of copies of the field) and some
finitely generated $R$-$R$-bimodules $W_0$ y $W_1$ with central action of the field. 
Such algebras will be called \textbf{elementary}. Let $\{e_1,\ldots,e_n\}$ be the set of orthogonal primitive central 
idempotents of $R$. An element $w$ of an $R$-$R$-bimodule $W$ is called \textbf{legible} if there exist $e_i$, $e_j$ 
such that $w=e_jwe_i$. An $R$-$R$-bimodule $W$ can be considered as a left $R^{op} \otimes R$-module, 
and since this is a semi-simple algebra, $W$ has a finite dual basis of legible elements.
The quiver $Q(T)$ associated to a basis of legible elements $\{ z \}_{z \in I}$ of the bimodules in the layer $(R,W_0\oplus W_1)$ 
has as vertices the set of idempotents $Q(T)_0=\{ e_i \}_{i=1}^n$. The set of arrows $Q(T)_1$ is the basis $\{ z \}_{z\in I}$.
The source and target of an arrow $z=e_jze_i$ are the vertices $e_i$ and $e_j$ respectively. The degree of an arrow corresponds
to the degree of $W_0$ or $W_1$. 

A \textbf{path} from $i$ to $j$ is a finite quiver $Q$ is an ordered set of arrows
$\gamma=\{ \alpha_1,\ldots,\alpha_r \}$ such that $s(\alpha_1)=i$, $t(\alpha_r)=j$ and $t(\alpha_{\ell})=s(\alpha_{\ell+1})$
for $\ell=1,\ldots, r-1$. The degree of a path $\gamma$ is given by $|\gamma|=\sum |\alpha_\ell|$.
We will use the notation $\gamma=\alpha_r \alpha_{r-1} \ldots \alpha_2 \alpha_1$. The idempotents $e_i$ will be considered
as trivial paths. The \textbf{path algebra} $kQ$ of a finite quiver $Q$ is the $k$-vector space with basis the paths of $Q$
and product given by concatenation. Observe that the subspace $R$ of $kQ$ generated by the trivial paths is a trivial
subalgebra of $kQ$. Denote by $W_i$ the subspace of $kQ$ generated by the arrows of degree $i$ ($i\in \{0,1\}$).
Then $W_i$ is a finitely generated $R$-$R$-bimodule with central action of the field and $kQ$ is freely generated
by $(R,W_0\oplus W_1)$, so that we identify $kQ$ with the tensorial algebra $T_R(W_0\oplus W_1)$. 
Crearly $T \cong k(Q(T))$ and $Q=Q(kQ)$ (relation $\beta$ in the figure) and the path algebra of a finite quiver
$Q$ in finite dimensional over the field $k$ if and only if $Q$ has no closed paths or \textbf{oriented cycles}.

%---------------------------------------------------------
A \textbf{differential tensor algebra} or \textbf{ditalgebra} $(T,\delta)$ consists of a graded algebra $T=\bigoplus T_i$ 
which is freely generated by $(T_0,T_1)$ and a differential $\delta:T \to T$, that is, a $k$-linear transformation with 
$\delta^2=0$ such that $\delta(T_i)\subset T_{i+1}$ and that satisfies the Leibniz rule:
\[
 \delta(ab)=\delta(a)b+(-1)^{|a|}a\delta(b),
\]
for all homogeneous elements $a$, $b$ in $T$. The \textbf{category of finite dimensional $(T,\delta)$-modules} $(T,\delta)$-mod 
is defined in the following way. Objects are given by finite dimensional $T_0$-modules. The group of morphisms $\Hom_{(T,\delta)}(M,N)$ 
between modules $M$ and $N$ is identified with the kernel of the following map (see definition~\ref{(A)D:hom} in the appendix)
\begin{equation} \label{(P)EQ:sigma}
\xymatrix@R=.3pc{
\txt{$\Hom_R(M,N)$ \\ $\oplus$ \\ $\Hom_{R\text{-}R}(W_1,\Hom_k(M,N))$} \ar[r]^-{\sigma} & \Hom_R(W_0 \otimes_R M,N) \\
(f^0,f^1) \ar@{|->}[r] & \left[ \txt{$w\otimes m \mapsto wf^0(m)$ \\ $-f^0(wm)-\widehat{f^1}(\delta(w))(m)$} \right],
}
\end{equation} 
where $\widehat{f^1}$ corresponds to $f^1$ under the isomorphism
\[
 \Hom_{R\text{-}R}(W_1,\Hom_k(M,N)) \cong \Hom_{T_0,T_0}(T_1,\Hom_k(M,N)),
\]
(cf. lemma~\ref{(A)L:extend}).
We say that a differential is \textbf{triangular on solid arrows} if there exists an ordering of these arrows in such a way
that the differential of a solid arrow $\alpha$ involves only solid arrows strictly smaller than $\alpha$.
In a similar way one can define \textbf{triangularity in dotted arrows} and we say that a differential is \textbf{triangular}
if it is triangular in both solid and dotted arrows (compare with definition~\ref{(A)D:triang}).
A $(T,\delta)$-module $M$ is \textbf{indecomposable} if whenever 
$M=M_1 \oplus M_2$ then one (and only one) of the representations $M_1$ or $M_2$ is zero. 
Denote by $\Ind (T,\delta)$ the set of isomorphism classes of indecomposable $(T,\delta)$-modules.
If the differential $\delta$ is triangular, then $(T,\delta)$-mod is a \textbf{Krull-Schmidt} category (that is, an additive 
$k$-category such that the ring of endomorphisms of every indecomposable module is local, cf. section~2.2 in Ringel \cite{cmR} 
and definition 5.10 and theorem 5.13 in Bautista, Salmer\'on and Zuazua \cite{BSZ09}). A full subcategory of a Krull-Schmidt
category $\mathcal{K}$ that is closed under isomorphisms, direct sums and direct summands is called an \textbf{object class} in 
$\mathcal{K}$.

The category $(T,\delta)$-mod has an \textbf{exact structure} $\mathcal{E}$ \cite[definition 6.3]{BSZ09}
as described as follows. For two $(T,\delta)$-modules $M$ and $N$ one defines the collection $\mathcal{E}(M,N)$ 
of pairs of morphisms $(f,g)$ that can be composed $\xymatrix@C=1.5pc{
N \ar[r]_-{f} & E \ar[r]_-{g} & M}$,
such that $gf=0$ and the exact sequence of $R$-modules
\[
 \xymatrix{
0 \ar[r] & N \ar[r]^-{f^0} & E \ar[r]^-{g^0} & M \ar[r] & 0,
}
\]
splits. Define the relation $(f,g) \sim (f',g')$ in $\mathcal{E}(N,M)$ whenever there exists 
an isomorphism $h:E \to E'$ such that the following diagram is commutative
\[
 \xymatrix{
N \ar[r]^-{f} & E \ar[r]^-{g} \ar[d]^-{h} & M \\
N \ar[r]_-{f'} \ar@{=}[u] & E' \ar[r]_-{g'} & M \ar@{=}[u].
}
\]
Let $\Ext^1_{(T,\delta)}(M,N)=\mathcal{E}(M,N)/\sim $ be the group of extensions from $M$ to $N$.
Then $\Ext^1_{(T,\delta)}(M,N)$ is isomorphic to the cokernel of the transformation $\sigma$
(lemma~\ref{(A)L:exacta} in the appendix).

We say that a $(T,\delta)$-module $M$ is \textbf{rigid} if $\Ext_{(T,\delta)}^1(M,M)=0$.
An indecomposable rigid module is called \textbf{exceptional}. In some case the tensorial algebra $T$ can
be recovered from the category $(T,\delta)$-mod (relation $\varepsilon$ in the figure) as shown in the following lemma. 
For each vertex $i$ denote by $S(i)$ the simple $T_0$-module of vertex $i$. Clearly $S(i)$ is still a simple module
when considered as a $(T,\delta)$-module.

\begin{lema} \label{(P)L:recupera}
 If $(T,\delta)$ is a ditalgebra whose differential satisfies
\[\delta(W_0) \subseteq \bigoplus_{m\geq 2}(W_0 \oplus W_1)^{\otimes m}\] then
\begin{eqnarray}
 \Ext^1_{(T,\delta)}(S(i),S(j)) & \cong & e_jW_0e_i, \nonumber \\
\Hom_{(T,\delta)}(S(i),S(j)) & \cong & \left\{
\begin{array}{l l}
e_jW_1e_i, & \text{if $i \neq j$},\\
k \oplus e_iW_1e_i, & \text{if $i=j$.}
\end{array} \right. \nonumber
\end{eqnarray}
\end{lema}

\bproof
Since the actions of $W_0$ in the simple modules $M=S(i)$ and $N=S(j)$ are null, the transformation $\sigma$ given 
in~(\ref{(P)EQ:sigma}) has the form
\[
\xymatrix@R=.5pc{
\txt{$\Hom_R(S(i),S(j))$ \\ $\oplus$ \\ $\Hom_{R\text{-}R}(W_1,\Hom_k(S(i),S(j)))$} \ar[r]^-{\sigma} & \Hom_R(W_0 \otimes_R S(i),S(j)) \\
(f^0,f^1) \ar@{|->}[r] & \left[ \txt{$w\otimes m \mapsto -\widehat{f^1}(\delta(w))(m)$} \right].
}
\]
By hypothesis and definition of $\widehat{f^1}$ (see lemma~\ref{(A)L:extend} in the appendix) we have that 
$\widehat{f^1}(\delta(w))(m)$ is zero for all $m \in M$, again because $M$ is simple. 
Then $\sigma$ is the zero transformation and hence its kernel and cokernel coincide with its domain and codomain respectively.
On the other hand it is clear that $\Hom_R(W_0 \otimes_R S(i),S(j)) \cong e_jW_0e_i$ and that 
$\Hom_{R\text{-}R}(W_1,\Hom_k(S(i),S(j))) \cong e_jW_1e_i$. Then $\Ext^1_{(T,\delta)}(S(i),S(j)) \cong e_jW_0e_i$ and 
\[
\Hom_{(T,\delta)}(S(i),S(j)) \cong \Hom_R(S(i),S(j)) \oplus e_jW_1e_i.
\]
The claim is now evident since
\begin{equation*}
\Hom_R(S(i),S(j)) \cong \left\{
\begin{array}{l l}
0, & \text{if $i \neq j$},\\
k, & \text{if $i=j$.}
\end{array} \right.
\end{equation*}
\eproof

Let $Q$ be a finite quiver. We are mainly interested in ditalgebras of the form $(kQ,\delta)$. In this case 
one asks that $\delta(R)=0$ where $R$ in the subalgebra of $kQ$ generated by trivial paths.
A $(kQ,\delta)$-module $M$ is determined by the vector spaces $M_i=e_iM$ for each $i \in Q_0$ and the linear transformations
\[
M_{\alpha}:\xymatrix@R=.5pc{M_i \ar[r] & M_j \\
m \ar@{|->}[r] & \alpha m=e_j\alpha e_im
}
\]
for each solid arrow $\alpha:i \to j$. In what follows we identify a $(kQ,\delta)$-module $M$ 
with the representation $(M_i;M_{\alpha})_{i \in Q_0}^{|\alpha|=0}$. 

%--------------------------------------------
For a finite dimensional graded $k$-algebra $T=\bigoplus_m T_m$ with orthogonal pri\-mi\-tive idempotents $\{ e_i \}_{i=1}^n$ 
define the \textbf{Cartan matrix} as the square matrix $C_T$ with coefficients the integers 
$m_{ji}=\dimGr_k (e_jTe_i)$, where
\[
 \dimGr_k (e_jTe_i)=\sum_m (-1)^m \dimk_k (e_jT_me_i).
\]
If $T$ is the path algebra of a finite quiver $Q$ with no oriented cycles, then the coefficients of $C_T=(m_{ji})$ are given by
\[
m_{ji}= \dimGr_k (e_jTe_i)=\sum (-1)^{|\gamma|},
\]
where the sum is taken over all paths $\gamma$ in $Q$ with source $i$ and target $j$.

\begin{lema} \label{(P)L:invers}
If $Q$ is a finite quiver with an admissible ordering of its vertices, then
\[
 Id=M_Q^tC_{kQ}.
\]
\end{lema}

\bproof
Let $i$ be an arbitrary vertex in $Q$ and $i_1,\ldots,i_r$ its direct predecessors 
(since the ordering is admissible, $i_j>i$ for each $j$).
Denote by $r_i$ the $i$-th row of $M_Q^t$.
Then $r_i$ has a $1$ in the position $i$ and the integer $a_j-b_j$ in the place $i_j$, where $a_j$ is the number of dotted arrows
from $i_j$ to $i$ and $b_j$ is the number of solid arrows from $i_j$ to $i$ ($j=1,\ldots,r$). 
All other entries in $r_i$ are zero. Let $c_{\ell}$ be the $\ell$-th 
column in $C_{kQ}$. Hence the product $r_ic_{\ell}$ has the form
\[
 r_ic_{\ell}=\dimGr_k[e_i(kQ)e_{\ell}]+\sum_{j=1}^r(a_j-b_j) \dimGr_k[e_{i_j}(kQ)e_{\ell}].
\]
If $i \geq \ell$ then $e_{i_j}(kQ)e_{\ell}=0$ for each $j$, and
\begin{equation*}
r_ic_{\ell} = \left\{
\begin{array}{l l}
0, & \text{if $i>\ell$},\\
1, & \text{if $i=\ell$}.
\end{array} \right.
\end{equation*}
Assume that $i<\ell$ and that there are $t_j$ paths from $\ell$ to $i_j$, say $\gamma^j_1,\ldots,\gamma^j_{t_j}$.
Observe that for an arrow $\alpha:i_j \to i$ one has either $|\alpha \gamma^j_u|=|\gamma^j_u|$ or
$|\alpha \gamma^j_u|=|\gamma^j_u|+1$, depending if $\alpha$ is solid or dotted respectively.
Hence
\begin{eqnarray}
 \dimGr_k [e_i(kQ)e_{\ell}] & = & \sum_{j=1}^r \left( -a_j\sum_{u=1}^{t_j}(-1)^{|\gamma_u^j|} 
                                  +b_j\sum_{u=1}^{t_j}(-1)^{|\gamma_u^j|} \right) =\nonumber \\
& = & -\sum_{j=1}^r(a_j-b_j) \sum_{u=1}^{t_j}(-1)^{|\gamma_u^j|} = \nonumber \\
& = & -\sum_{j=1}^r(a_j-b_j) \dimGr_k[e_{i_j}(kQ)e_{\ell}]. \nonumber
\end{eqnarray}
Then $r_ic_{\ell}=0$. This completes the proof for $r_ic_{\ell}=\delta_{i,\ell}$.
\eproof

The \textbf{dimension vector} of a representation $M$ is
the element of $\mathbb{Z}^{Q_0}$ given by $\vdim M=(\dimk_k M_i)_{i\in Q_0}$.
The \textbf{Euler characteristic} of $(kQ,\delta)$-mod is given, for two $(kQ,\delta)$-modules $M$ and $N$, by
the difference
\[
 \mathbb{E}(M,N)=\dimk_k \Hom_{(kQ,\delta)}(M,N)-\dimk_k \Ext_{(kQ,\delta)}^1(M,N).
\]
The following expression of the bilinear form associated to the quiver $Q$ in terms of its arrows
will be useful. For vectors $m=(m_i)$ and $n=(n_i)$ in $\mathbb{Z}^{Q_0}$ one has
\begin{eqnarray} 
\langle m,n \rangle_Q & = & m^tM_Qn = \sum_{i,j\in Q_0} m_in_jm_{ij}= \nonumber \\
& = & \sum_{i\in Q_0}m_in_im_{ii} + \sum_{\substack{i,j \in Q_0 \\ i \neq j}} m_in_jm_{ij} = \nonumber \\
& = & \sum_{i \in Q_0}m_in_i\left( 1-\sum_{\substack{\alpha \in Q_1\\ \alpha:i \to i}} (-1)^{|\alpha|} \right) 
+\sum_{i,j\in Q_0}m_in_j\left(-\sum_{\substack{\alpha \in Q_1\\ \alpha:i \to j}} (-1)^{|\alpha|} \right) \nonumber \\
& = & \sum_{i \in Q_0} m_in_i - \sum_{\substack{\alpha \in Q_1 \\ |\alpha|=0}} m_{s(\alpha)}n_{t(\alpha)}
+\sum_{\substack{\gamma \in Q_1 \\ |\gamma|=1}} m_{s(\gamma)}n_{t(\gamma)}. \nonumber
\end{eqnarray}

\begin{lema} \label{(P)L:Euler}
For any triangular differential $\delta$ of $kQ$, the Euler characteristic of $(kQ,\delta)$ coincides,
through the dimension vector, with the bilinear form 
$\langle \cdot,\cdot \rangle_Q$ associated to the incidence matrix $M_Q$, that is,
\[
\langle \vdim(M),\vdim(N) \rangle_Q =\dimk_k \Hom_{(kQ,\delta)}(M,N)-\dimk_k \Ext_{(kQ,\delta)}^1(M,N).
\]
\end{lema}

\bproof
The transformation $\sigma$ given in (\ref{(P)EQ:sigma}) generates an exact sequence of the form
\begin{equation*}
 \xymatrix@C=1pc{
0 \ar[r] &  \Hom_{(kQ,\delta)}(M,N) \ar[r] & 
\txt{$\Hom_R(M,N)$ \\ $\oplus$ \\ $\Hom_{R\text{-}R}(W_1,\Hom_k(M,N))$} \ar[r]^-{\sigma} &
}
\end{equation*}
\begin{equation*}
\xymatrix@C=1pc{
\ar[r]^-{\sigma} &  \Hom_R(W_0 \otimes_R M,N) \ar[r]^-{\eta} &  \Ext^1_{(kQ,\delta)}(M,N) \ar[r]  & 0,
}
\end{equation*}
(cf. lemma \ref{(A)L:exacta} in the appendix).
On the other hand, by the expression given above of the bilinear form in the vectors $m=\vdim M$ and $n=\vdim N$   
one has
\begin{displaymath}
\langle m,n \rangle_Q=\sum_{i \in Q_0} m_in_i - \sum_{\substack{\alpha \in Q_1 \\ |\alpha|=0}} m_{s(\alpha)}n_{t(\alpha)}
+\sum_{\substack{\gamma \in Q_1 \\ |\gamma|=1}} m_{s(\gamma)}n_{t(\gamma)}.
\end{displaymath}
Clearly $\dimk_k \Hom_R(M,N)=\sum_{i \in Q_0}m_in_i$. Since
\[
\Hom_R(W_0 \otimes_R M,N)\cong \bigoplus_{\substack{\alpha \in Q_1 \\ |\alpha|=0}} \Hom_k(M_{s(\alpha)},N_{t(\alpha)}),
\]
we have that $\dimk_k \Hom_R(W_0 \otimes_R M,N)=\sum_{|\alpha|=0}m_{s(\alpha)}n_{t(\alpha)}$.
Moreover we have an isomorphism
\[
\Hom_{R\text{-}R}(W_1,\Hom_k(M,N))\cong \bigoplus_{\substack{\gamma \in Q_1 \\ |\gamma|=1}} \Hom_k(M_{s(\gamma)},N_{t(\gamma)})
\]
which implies that $\dimk_k \Hom_{R\text{-}R}(W_1,\Hom_k(M,N))=\sum_{|\gamma|=1}m_{s(\gamma)}n_{t(\gamma)}$.
\eproof

\section{Reflections and Coxeter matrix.} \label{(P)S:reflex}
%----------------------------------------------------------------------
Assume $Q$ i a finite quiver with an admissible ordering of its vertices. Then the Cartan
matrix $C_{kQ}$ of $Q$ is an upper triangular matrix with only 1's in the diagonal.
In particular $Q$ has no loops and the canonical vectors $\mathbf{e}_i$ are roots
of the quadratic form $q_Q$ of the quiver $Q$.

The \textbf{Coxeter matrix} of $Q$ is given by 
\[
\Phi_{Q}=-C_{kQ}^tC_{kQ}^{-1}.
\]
In this section we give an alternative expression for $\Phi_Q$ in terms of reflections.
If $n=|Q_0|$ define the \textbf{simple reflection} $\sigma_i:\mathbb{Z}^n \to \mathbb{Z}^n$ as
\[
 \sigma_i(x)=x-2(x,\mathbf{e}_i)\mathbf{e}_i,
\]
where $(a,b)=(1/2)[\langle a,b \rangle + \langle b,a \rangle]$ is the symmetrization of the bilinear form
associated to $Q$. Since $q_Q(\mathbf{e}_i)=1$ we have that $\sigma_i(\mathbf{e}_i)=-\mathbf{e}_i$ for $1 \leq i\leq n$. 
Thus, for any $x \in \mathbb{Z}^n$
\begin{eqnarray}
 \sigma_i^2(x) & = & \sigma_i(x-2(x,\mathbf{e}_i)\mathbf{e}_i)=[\sigma_i(x)]-2(x,\mathbf{e}_i)\sigma_i(\mathbf{e}_i)=\nonumber \\
& = & [x-2(x,\mathbf{e}_i)\mathbf{e}_i]+2(x,\mathbf{e}_i)\mathbf{e}_i = \nonumber \\
& = & x,\nonumber
\end{eqnarray}
that is, $\sigma_i^2=Id_{\mathbb{Z}^n}$. If we denote also by $\sigma_i$ the matrix corresponding
to the transformation $\sigma_i$ respect to the canonical basis $\{ \mathbf{e}_i \}_{i=1}^n$, we 
want to prove that
\[
 \Phi_Q=\sigma_n \sigma_{n-1} \cdots \sigma_2 \sigma_1.
\]
For that purpose, define the vectors $p_{i,j}$ for $1 \leq i \leq j \leq n$ in the following way.
Take first $p_{i,i}=\mathbf{e}_i$ for $1 \leq i \leq n$ and for $1 \leq i<j \leq n$ take
\[
 p_{i,j}=\sigma_{i}\sigma_{i+1} \cdots \sigma_{j-1}\mathbf{e}_j.
\]
Observe that the $i$-th entry $x_i$ of the vector $x \in \mathbf{Z}^n$ can be obtained as the product $\mathbf{e}_i^tx$.
The next equality, which follows directly from the definition of $\sigma_{\ell}$, will be useful in the proof
of the following lemma,
\begin{equation}
\mathbf{e}^t_k\sigma_{\ell}(x) = \left\{
\begin{array}{l l}
\mathbf{e}_k^tx=x_k, & \text{if } k \neq \ell,\\
x_{\ell}-2(x,\mathbf{e}_{\ell}), & \text{if } k=\ell.
\end{array} \right.
\end{equation} \label{(P)E:reflexUno}

\begin{lema} \label{(P)L:lema}
 If $1\leq i \leq j \leq n$ then
\begin{equation*}
\mathbf{e}^t_kp_{i,j} = \left\{
\begin{array}{l l}
\sum_{\substack{\text{paths } \gamma \\ \text{from $j$ to $k$} }}(-1)^{|\gamma|}, & \text{if } i \leq k \leq j ,\\
0, & \text{otherwise}.
\end{array} \right.
\end{equation*}
\end{lema}

\bproof
From the equality (\ref{(P)E:reflexUno}) it is clear that if $k < i$ then
\begin{eqnarray} \label{(P)EQ:Pruebacox}
\mathbf{e}_k^tp_{i,j} & = & \mathbf{e}_k^t \sigma_{i}(\sigma_{i+1}\cdots \sigma_{j-1}(\mathbf{e}_j))= 
\mathbf{e}_k^t \sigma_{i+1}(\sigma_{i+2}\cdots \sigma_{j-1}(\mathbf{e}_j))=\nonumber \\
& = & \nonumber \mathbf{e}_k^t\sigma_{i+2}(\sigma_{i+3}\cdots \sigma_{j-1}(\mathbf{e}_j)) = \ldots = 
\mathbf{e}_k^t \sigma_{j-1}(\mathbf{e}_j)=\mathbf{e}_k^t\mathbf{e}_j =\nonumber \\
& = & 0,
\end{eqnarray}
and in a similar way if $j < k$ than $\mathbf{e}_k^tp_{i,j}=0$. 
We restrict ourselves to the case $i \leq k \leq j$ and prove the result inductively over the difference $j-i$.
When $i=j$ we have
\[
 \mathbf{e}_i^tp_{i,i}=\mathbf{e}_i^t\mathbf{e}_i=1
=\sum_{\substack{ \text{paths } \gamma \\ \text{from $i$ to $i$} }}(-1)^{|\gamma|},
\]
since the only path from $i$ to $i$ is the trivial path, which has zero degree. This is the base case.
Assume that the claim is valid for the vectors $p_{i',j'}$ with $j'-i' \leq \ell$ and let $1 \leq i < j \leq n$
be integers with $j-i=\ell+1$. Notice first that if $i < k \leq j$, using again the equality (\ref{(P)E:reflexUno}), then
\begin{eqnarray}
 \mathbf{e}_k^tp_{i,j} & = & \mathbf{e}_k^t \sigma_i(\sigma_{i+1}\cdots \sigma_{j-1}(\mathbf{e}_j)) = \ldots = \nonumber \\
& = & \mathbf{e}_k^t\sigma_k (\sigma_{k+1} \cdots \sigma_{j-1}(\mathbf{e}_j)) = \nonumber \\
& = & \mathbf{e}_k^t p_{k,j}, \nonumber
\end{eqnarray}
and since $j-k\leq \ell$, by induction hypothesis,
\[
 \mathbf{e}_k^tp_{i,j}=\mathbf{e}_k^tp_{k,j}
=\sum_{\substack{ \text{paths } \gamma \\ \text{from $j$ to $k$} }}(-1)^{|\gamma|}.
\]
Thus it is enough to analyze the case $k=i$. Since $\mathbf{e}_i^tp_{i+1,j}=0$ we have
\begin{eqnarray}
 \mathbf{e}_i^tp_{i,j} & = & \mathbf{e}_i^t \sigma_i(p_{i+1,j}) = \mathbf{e}_i^t[p_{i+1,j}-2(p_{i+1,j},\mathbf{e}_i)\mathbf{e}_i] = \nonumber \\
& = & -2(p_{i+1,j},\mathbf{e}_i)\mathbf{e}_i^t\mathbf{e}_i = \nonumber \\
& = & -( \langle p_{i+1,j},\mathbf{e}_i \rangle + \langle \mathbf{e}_i,p_{i+1,j} \rangle ). \nonumber 
\end{eqnarray}
Since the ordering of the vertex set is admissible, $\langle \mathbf{e}_i,\mathbf{e}_m \rangle=0$ if $i < m$. Then
\begin{eqnarray}
 \langle \mathbf{e}_i,p_{i+1,j} \rangle & = & \langle \mathbf{e}_i,\sum_{m=1}^n(\mathbf{e}_m^tp_{i+1,j})\mathbf{e}_m \rangle = \nonumber \\
& = &\sum_{m=1}^n\mathbf{e}_m^tp_{i+1,j}\langle \mathbf{e}_i,\mathbf{e}_m \rangle = \nonumber \\
& = & \sum_{m=1}^{i}\mathbf{e}_m^tp_{i+1,j}\langle \mathbf{e}_i,\mathbf{e}_m \rangle =0, \nonumber
\end{eqnarray}
for by equality~(\ref{(P)EQ:Pruebacox}), $\mathbf{e}_m^tp_{i+1,j}=0$ for $m=1,\ldots,i$. On the other hand,
we recall that for $m \neq i$, $\langle \mathbf{e}_m,\mathbf{e}_i \rangle=\mathbf{e}_m^tM_Q\mathbf{e}_i
=-\sum_{\substack{\text{arrows } \alpha \\ \text{from $m$ to $i$} }}(-1)^{|\alpha|}$. 
By induction hypothesis in the vector $p_{i+1,j}$ we have
\begin{eqnarray}
 \langle p_{i+1,j},\mathbf{e}_i \rangle & = &\sum_{m=1}^n\mathbf{e}_m^tp_{i+1,j}\langle \mathbf{e}_m,\mathbf{e}_i \rangle
= \sum_{m=i+1}^n\mathbf{e}_m^tp_{i+1,j}\langle \mathbf{e}_m,\mathbf{e}_i \rangle= \nonumber \\
& = & \sum_{m=i+1}^n \left( \sum_{\substack{\text{paths }\gamma \\ \text{from $j$ to $m$} }}(-1)^{|\gamma|} \right) 
\left( -\sum_{\substack{\text{arrows }\alpha \\ \text{from $m$ to $i$} }}(-1)^{|\alpha|} \right)= \nonumber \\
& = & -\sum_{\substack{\text{paths }\gamma \\ \text{from $j$ to $i$} }}(-1)^{|\gamma|}. \nonumber
\end{eqnarray}
Hence $\mathbf{e}_i^tp_{i,j} 
= -( \langle p_{i+1,j},\mathbf{e}_i \rangle + \langle \mathbf{e}_i,p_{i+1,j} \rangle )
=\sum_{\substack{\text{paths }\gamma \\ \text{froms $j$ to $i$} }}(-1)^{|\gamma|}$,
which completes the proof.
\eproof

Denote by $p_i$ the vector $p_{1,i}$, for $i=1,\ldots,n$. The following corollary is direct consequence
of the last result.
\begin{corolario} \label{(P)C:Cartan}
If $Q$ is a finite quiver with an admissible ordering of its vertices, then the Cartan matrix $C_{kQ}$
has as columns the vectors $p_1,\ldots,p_n$.
\end{corolario}
\bproof 
From the last lemma~\ref{(P)L:lema} the entries of the vector $p_{i}=p_{1,i}$ are given by
\[
 (p_i)_k=\mathbf{e}_k^tp_{1,i}=\sum_{\substack{\text{paths } \gamma \\ \text{from $i$ to $k$}}} (-1)^{|\gamma|},
\]
where the sum is considered zero when the set of paths from $i$ to $k$ is empty (in particular whan $k>i$).
This is precisely the definition of the entry $(k,i)$ in the matrix $C_{kQ}$.
\eproof

In particular, if $Q$ has only solid arrows, then $p_i=\vdim P(i)$ where $P(i)=kQ\mathbf{e}_i$ is the 
indecomposable projective $kQ$-module which is projective cover of the simple module of vertex $i$.
In this case, using the standar duality, the vector $q_i=C_{kQ}^t\mathbf{e}_i$ is dimension vector
of the injective envelope of the simple module of vertex $i$. In general, with dual arguments to the case $p_i$,
it can be proved that
\[
 q_i=\sigma_{n}\sigma_{n-1}\cdots \sigma_{i+1}(\mathbf{e}_i).
\]
Observe that the set $p_1,\ldots,p_n$ is basis of the free group $\mathbb{Z}^n$. We show by induction that 
$\langle p_1,\ldots,p_i \rangle=\langle \mathbf{e}_1,\ldots,\mathbf{e}_i \rangle$ for $1\leq i \leq n$.
The case $i=1$ is clear for $p_1=\mathbf{e}_1$. Assume that 
\[
G=\langle p_1,\ldots,p_{\ell-1} \rangle=\langle \mathbf{e}_1,\ldots,\mathbf{e}_{\ell-1} \rangle.
\]
Notice that if $i=1,\ldots,\ell-1$ then $\sigma_i(G) \subset G$ and $\sigma_i(\mathbf{e}_{\ell})=\mathbf{e}_{\ell}+a$ for 
some element $a \in G$. Hence
\begin{eqnarray}
 p_{\ell} & = & \sigma_1 \sigma_2\cdots \sigma_{\ell-1}(\mathbf{e}_{\ell}) = \nonumber \\
& = & \sigma_1 \sigma_2 \cdots \sigma_{\ell-2}(\mathbf{e}_{\ell}+a_2) = \nonumber \\
& = & \sigma_1\sigma_2 \cdots \sigma_{\ell-3}(\mathbf{e}_{\ell}+a_3) =\ldots = \nonumber \\
& = & \sigma_1(\mathbf{e}_{\ell}+a_{\ell-1}) = \nonumber \\
& = & \mathbf{e}_{\ell}+a_{\ell}, \nonumber 
\end{eqnarray}
for some elements $a_1,\ldots,a_{\ell}$ in $G$. Then 
$\langle p_1,\ldots,p_{\ell} \rangle=\langle \mathbf{e}_1,\ldots,\mathbf{e}_{\ell} \rangle$ whice completes the
inductive step.

\begin{corolario} \label{(P)C:Coxeter}
Let $Q$ be a finite quiver with an admissible ordering of its vertices. Denote by $\sigma_i$ the matrix corresponding
to the tranformation $\sigma_i$ respect to the canonical basis $\{ \mathbf{e}_i \}_{i=1}^n$. Then
\[
\Phi_Q= -C_{kQ}^tC_{kQ}^{-1}=\sigma_n \sigma_{n-1}\cdots \sigma_2\sigma_1.
\]
\end{corolario}

\bproof
As shown above the matrix $C_{kQ}$ satisfies $C_{kQ}\mathbf{e}_i=p_i$ and $C^t_{kQ}\mathbf{e}_i=q_i$, thus
\begin{eqnarray}
 \Phi_{kQ}p_i & = & -C_{kQ}^tC_{kQ}^{-1}(C_{kQ}\mathbf{e}_i)=-C_{kQ}^t\mathbf{e}_i= \nonumber \\
& = & -q_i. \nonumber
\end{eqnarray}
The same is valid for the composition $\sigma_n \cdots \sigma_1$, for since $\sigma_i^2=Id_{\mathbb{Z}^n}$ and 
$\sigma_i(\mathbf{e}_i)=-\mathbf{e}_i$ we have
\begin{eqnarray}
 \sigma_n \cdots \sigma_1p_i & = & \sigma_n \cdots \sigma_1(\sigma_1 \cdots \sigma_{i-1}(\mathbf{e}_i)) = \nonumber \\
& = & \sigma_n \cdots \sigma_{i+1}\sigma_i(\mathbf{e}_i)=-\sigma_n \cdots \sigma_{i+1}(\mathbf{e}_i) = \nonumber \\
& = & -q_i. \nonumber
\end{eqnarray}
Since the set $\{p_i\}_{i=1}^n$ is generator of $\mathbb{Z}^n$, we conclude that the transformations $\Phi_Q$ and
$\sigma_n \cdots \sigma_1$ coincide.
\eproof

\begin{lema} \label{(DE)L:radCero}
Let $Q$ be a finite quiver with an admissible ordering of its vertices. Let $(\cdot ,\cdot)_Q$ be the symmetric bilinear form
associated to $Q$, and let $\sigma_i$ be the simple reflection respecto to the vertex $i$ and $\Phi_Q$ the Coxeter matrix of $Q$.
For a vector $x \in \mathbb{Z}^n$ the following are equivalent
\begin{itemize}
 \item[a)] $\Phi_Qx=x$,
 \item[b)] $\sigma_ix=x$ for each $i \in Q_0$,
 \item[c)] $(x,y)_Q=0$ for each vector $y \in \mathbb{Z}^n$. 
\end{itemize}
\end{lema}

\bproof
Clearly $(b)$ is equivalent to $(x,\mathbf{e}_i)_Q=0$ for all canonical vectors $\mathbf{e}_i$, which is equivalent
to $(c)$. It is also clear, by corollary~\ref{(P)C:Coxeter}, that $(b)$ implies $(a)$. Hence it is enough 
to show that $(a)$ implies $(b)$. First, since $\Phi_Qx=x$, the first entry $x_1$ of the vector $x$ 
($x_1=\mathbf{e}_1^tx$) satisfies, because of the equality~(\ref{(P)E:reflexUno}) above,
\[
 \mathbf{e}_1^tx=\mathbf{e}_1^t\Phi_Qx=\mathbf{e}_1^t\sigma_n\cdots \sigma_1(x)=\mathbf{e}_1^t\sigma_1(x),
\]
thus $\sigma_1x=x$. Assume we have proved that $\sigma_ix=x$ for $i=1,\ldots,\ell-1$.
Then 
\[
\mathbf{e}_{\ell}^tx=\mathbf{e}_{\ell}^t\Phi_Qx=\mathbf{e}_{\ell}^t\sigma_{\ell}\sigma_{\ell-1}\cdots \sigma_1x
=\mathbf{e}_{\ell}^t\sigma_{\ell}x,
\]
and hence $\sigma_{\ell}x=x$. This completes the proof.
\eproof

\section{Elementary reductions.} \label{(P)S:redElem}
%----------------------------------------------------------------------
The following constructions over a ditalgebra $\mathcal{A}=(kQ,\delta)$, called elementary reductions, 
have been used in different contexts to give inductive proofs. Their associated reduction functors $F^z$ ($z \in \{c,r,d\}$) 
are always full and faithful.
\begin{itemize}
 \item[\textbf{(c)}] Change of basis, $F^c:\mathcal{A}^c\text{-mod} \longrightarrow \mathcal{A}$-mod.
 \item[\textbf{(r)}] Regularization, $F^r:\mathcal{A}^r\text{-mod} \longrightarrow \mathcal{A}$-mod.
 \item[\textbf{(d)}] Deletion of idempotents,  $F^d:\mathcal{A}^d\text{-mod} \longrightarrow \mathcal{A}$-mod.
\end{itemize}
The associated functors $F^z$ are induced by ditalgebra morphisms
$\varphi^z:\mathcal{A} \to \mathcal{A}^z$, that is, graded algebra morphisms $\varphi:T \to T'$ 
which commute with differentials
\[
\varphi \delta_{\mathcal{A}} = \delta_{\mathcal{A}'} \varphi.
\]
Then $F^z$ is given by restriction of scalars \cite[lemma 2.4]{BSZ09}. 
Recall that a functor $F$ between additive $k$-categories $\mathcal{C}$ and $\mathcal{C}'$ with exact
structures $\mathcal{E}$ and $\mathcal{E}'$ respectively is called \textbf{exact functor} if it determines
a mapping from $\mathcal{E}$ to $\mathcal{E}'$ \cite[definition 6.3]{BSZ09}. Observe that functors induced by ditalgebra
morphisms are exact.

\begin{lema} \label{(P)L:monoExt}
Let $F:(\mathcal{C},\mathcal{E}) \to (\mathcal{C}',\mathcal{E}')$ be an exact functor. Then for each pair
of objects $M$, $N$ in $\mathcal{C}$, the functor $F$ induces a transformation
\[
 F^*:\Ext^1_{\mathcal{C}}(M,N) \to \Ext^1_{\mathcal{C}'}(F(M),F(N)).
\]
If $F$ is full and faithful then $F^*$ is a monomorphism.
\end{lema}
\bproof
Consider exact pairs $e_1=(f_1,g_1)$ and $e_2=(f_2,g_2)$ in $\mathcal{E}$ which are equivalents
$e_1 \sim e_2$, and the image of this equivalence under the functor $F$,
\[
 \xymatrix{
N \ar[r]^-{f_1} \ar@{=}[d] & E \ar[r]^-{g_1} \ar[d]^{h} & M \ar@{=}[d] \ar@{}[rd]|-{\mapsto}
& F(N) \ar[r]^-{F(f_1)} \ar@{=}[d] & F(E) \ar[r]^-{F(g_1)} \ar[d]^{F(h)} & F(M) \ar@{=}[d] \\
N \ar[r]_-{f_2} & E' \ar[r]_-{g_2} & M & F(N) \ar[r]_-{F(f_2)} & F(E') \ar[r]_-{F(g_2)} & F(M).
}
\]
Since $F$ is exact and $F(h)$ is an isomorphism if $h$ is, we have that $F(e_1)$ is equivalent to $F(e_2)$ in the exact 
structure $\mathcal{E}'$. In this way $F$ induces a function in quotients
\[
F^*:\Ext^1_{\mathcal{C}}(M,N) \to \Ext^1_{\mathcal{C}'}(F(M),F(N)).
\]
Observe that $F^*$ is actually a vector space transformation, for $F$ commutes with the diagonal and codiagonal
morphisms $\Delta$ and $\nabla$ in the definition of Baer sums (cf. \cite[proposition 6.11]{BSZ09}).

Assume now that $F$ is full and faithful, and that $e=(f,g)$ is an exact pair such that $F^*([e])=0$,
\[
\xymatrix{
 N \ar[r]^-{f} & E \ar[r]^-{g} & M \\
 F(N) \ar[r]^-{F(f)} & F(E) \ar[r]^-{F(g)} & F(M). 
}
\]
Then $F(g)$ is a retraction in $\mathcal{C}'$, and hence there exists $h':F(M) \to F(E)$ such that $F(g)h'=Id_{F(M)}$.
Since $F$ is full there exists a morphism $h:M \to E$ such that $F(h)=h'$. Then
\[
 F(Id_{M})=Id_{F(M)}=F(g)h'=F(g)F(h)=F(gh),
\]
and since $F$ is faithful $gh=Id_M$. Thus $g$ is a retraction and $[e]=0$, that is, 
$F^*$ is a monomorphism.
\eproof
We will say that an exact functor $F:(\mathcal{C},\mathcal{E}) \to (\mathcal{C}',\mathcal{E}')$ is \textbf{rigid}
if for each pair of objects $M$ and $N$ in $\mathcal{C}$ the induced morphism
\[
  F^*:\Ext^1_{\mathcal{C}}(M,N) \to \Ext^1_{\mathcal{C}'}(FM,FN)
\]
is an isomorphism.

Define the \textbf{norm} $||M||$ of an $\mathcal{A}$-module $M$ as
\[
 ||M||=\sum_{\substack{\alpha \in Q_1 \\ |\alpha|=0}}(\dimk_k M_{s(\alpha)})(\dimk_k M_{t(\alpha)}).
\]
We are interested in how the reduction functors $F^z$ behave respect to dimension vectors, norms
and extension groups.

%--------------------
\subsubsection*{Change of basis.} \label{(E):bas}
Let $Q$ be a finite quiver and let $(kQ,\delta)$ be a ditalgebra with triangular differential. Consider a set 
of arrows of degree one $\gamma_1,\ldots,\gamma_r$ with same source vertex and same target vertex, and
such that $1\leq i<j\leq r$ implies that $\gamma_i$ is smaller than
$\gamma_j$ in the ordering of dotted arrows in which $\delta$ is triangular.
Let $c_1,\ldots,c_r$ be nonzero scalars and let $Q^c$ be a copy of $Q$. 
Let $W$ and $W^c$ be the arrow bimodules of $Q$ and $Q^c$, and consider the functions
$g_1:Q_1 \to W^c$ and $\widetilde{g}_1:Q_1^c \to W$ given by
\begin{equation*}
g_1(\beta) = \left\{
\begin{array}{l l}
\beta^c, & \text{if } \beta \neq \gamma_r,\\
\frac{1}{c_r} \left(\gamma^c_r - \sum_{i=1}^{r-1}c_i\gamma^c_i \right), & \text{if } \beta=\gamma_r;
\end{array} \right.
\end{equation*}
\begin{equation*}
\widetilde{g}_1(\beta^c) = \left\{
\begin{array}{l l}
\beta, & \text{if } \beta^c \neq \gamma^c_r,\\
\sum_ic_i\gamma_i, & \text{if } \beta^c=\gamma^c_r.
\end{array} \right.
\end{equation*}
Hence $g_1$ and $\widetilde{g}_1$ can be extended to bimodule morphisms $\xymatrix{W \ar@<.5ex>[r]^-{g_1} & W^c 
\ar@<.5ex>[l]^-{\widetilde{g}_1}}$ which are inverse from each other. Let $R$ and $R^c$ be the
trivial subalgebras of the path algebras $kQ$ and $kQ^c$. Denote by $g_0:R \to R^c$ the isomorphism of
algebras given by the correspondence of vertices in $Q$ and vertices in $Q^c$, and let $\widetilde{g}_0$ be its inverse. 
Then the isomorphisms $g_0$ and $g_1$ can be extended to a tensor algebra isomorphism $g:T_R(W) \to T_{R^c}(W^c)$,
with inverse given by the tensorial extension $\widetilde{g}$ of the morphisms $\widetilde{g}_0$ and $\widetilde{g}_1$.
Since we have identified the algebras $T_R(W)$ and $kQ$, we have an isomorphism
\[
 g:kQ \longrightarrow kQ^c,
\]
with inverse $\widetilde{g}$. Define $\delta^c=g\delta \widetilde{g}$, which is clearly a differential. 
Give the same ordering in the arrows of $Q^c$ as that in $Q$. Since the functions $\widetilde{g}_1$ and $g_1$
send an arrow $\beta$ to a linear combination of arrows smaller or equal to $\beta$, the differential $\delta^c$
is triangular in solid and dotted arrows. Then $(kQ^c,\delta^c)$ is a ditalgebra with triangular differential and 
$g:(kQ,\delta) \to (kQ^c,\delta^c)$ is an isomorphism. The functor it iduces is an equivalence of categories
\cite[lemma 2.4]{BSZ09}
\[
 F^c: \mathcal{A}^c\text{-mod} \longrightarrow \mathcal{A}\text{-mod}.
\]
Notice that if $\alpha$ is a solid arrow in $Q$ with $\delta(\alpha)=\sum_{i=1}^rc_i\gamma_i$, 
then the corresponding $\alpha^c$ satisfies
\[
 \delta^c(\alpha^c)=g\delta \widetilde{g}(\alpha^c)=g(\delta(\alpha))
=g(\sum_{i=1}^rc_i\gamma_i)=g(\widetilde{g}(\gamma^c_r))=\gamma^c_r.
\]
The functior $F^c$ crearly preserves dimension vectors $\vdim M=\vdim F^c(M)$ and 
norms $||M||=||F^c(M)||$. Since $F^c$ is an exact functor, by lemma~\ref{(P)L:monoExt} it induces injections in extension groups
$\Ext^1_{\mathcal{A}^c}(M,N) \to \Ext^1_{\mathcal{A}}(F^c(M),F^c(N))$. The same is valid for the induced functor of the inverse 
$g^{-1}$, hence we conclude that $F^c$ is a rigid functor.

%--------------------
\subsubsection*{Regularization.} \label{(E):reg}
Assume that $\alpha$ and $\gamma$ are arrows such that $\delta(\alpha)=\gamma$.
The ditalgebra $\mathcal{A}^r$ obtained by regularization of the arrows $\alpha$ and $\gamma$ 
is isomorphic to a ditalgebra of the form $(kQ^r,\delta^r)$ \cite[lemma 23.18]{BSZ09}.
The quiver $Q^r$ is obtained removing from $Q$ the arrows $\alpha$ and $\gamma$. The associated
ditalgebra morphism $\varphi:(kQ,\delta) \to (kQ^r,\delta^r)$ induces an equivalence of categories
\[
 F^r:\mathcal{A}^r\text{-mod} \longrightarrow \mathcal{A}\text{-mod},
\]
such that $||M|| \leq ||F^r(M)||$ \cite[lemma 25.3]{BSZ09}. Moreover, $\vdim M = \vdim F^r(M)$ and since 
$F^r$ is an exact functor, by lemma~\ref{(P)L:monoExt} it induces an inclusion in extension groups 
$\Ext_{\mathcal{A}^r}(M,N) \to \Ext_{\mathcal{A}}(F^r(M),F^r(N))$. We show that it is surjective.
If $\xymatrix{F^r(N) \ar[r]^-{u'} & E' \ar[r]^-{v'} & F^r(M)}$ is an exact pair, since $F^r$ is dense
there is an $\mathcal{A}^r$-module $E$ and an isomorphism $h:E' \to F^r(E)$. 
Then there exists an equivalence of exact pairs
\[
 \xymatrix{
F^r(N) \ar@{=}[d] \ar[r]^-{u'} & E' \ar[d]^-{h} \ar[r]^-{v'} & F^r(M) \ar@{=}[d] \\
F^r(N) \ar[r]_-{u} & F^r(E) \ar[r]_-{v} & F^r(M).
}
\]
Since $F^r$ is full there exist $f:N \to E$ and $g:E \to M$ such that $F^r(f)=u$ and $F^r(g)=v$.
Finally we verify that $\xymatrix{N \ar[r]^-{f} & E \ar[r]^-{g} & M}$ is an exact pair in $\mathcal{A}^r$.
On the one hand $gf=0$, since $F^r(gf)=F^r(g)F^r(f)=vu=0$ and $F^r$ is faithful. On the other hand, the sequence 
of $R^r$-modules given by restriction of scalars
\[
 \xymatrix{0 \ar[r] & N \ar[r]^-{f_0} & E \ar[r]^-{g_0} & M \ar[r] & 0}
\]
is exact, since $F^r$ is given by restriction of scalars. Then $F^r$ is rigid.

%--------------------
\subsubsection*{Deletion of idempotents.} \label{(E):del}
Consider a ditalgebra with triangular differential $(kQ,\delta)$ and let $Q^d$ be the 
full subquiver of $Q$ determined by a subset of vertices $Q'_0 \subset Q_0$.  Then there exists
a triangular differential $\delta^d$ in $kQ^d$ such that $(Q^d,\delta^d)$ is a
ditalgebra \cite[lemma 23.14]{BSZ09} and the canonical projection
\[
 \varphi^d:\mathcal{A}=(kQ,\delta) \longrightarrow \mathcal{A}^d=(kQ^d,\delta^d),
\]
induces a full and faithful functor $F^d:\mathcal{A}^d\text{-mod} \to \mathcal{A}\text{-mod}$. An $\mathcal{A}$-module
$N$ is isomorphic to an object in the image of $F^d$ if and only if the support of $N$ is contained in $Q'_0$. 
Again by lemma~\ref{(P)L:monoExt} there are inclusions
\[
 \Ext_{\mathcal{A}^d}(M,N) \to \Ext_{\mathcal{A}}(F^d(M),F^d(N)),
\]
for each pair of $\mathcal{A}^d$-modules $M$ and $N$. If there is an exact pair
\[
e=\xymatrix{F^d(N) \ar[r]^-{u'} & E' \ar[r]^-{v'} & F^d(M)},
\]
then the support of $E'$ is contained in $Q'_0$, and hence there is an isomorphism $h:E' \to F^d(E)$ 
for some $\mathcal{A}^d$-module $E$. 
Repeating the argument given in the regularization section, we have that the exact pair $e$ is equivalent
to the image of an exact pair in the category $\mathcal{A}^d$, and hence $F^d$ is rigid. Moreover, if $M$
is an $\mathcal{A}^d$-module and we denote by $I_{Q'_0}$ the matrix with $|Q'_0|$ columns, $|Q_0|$ rows 
and coefficients given by
\begin{equation*}
d_{ij} = \left\{
\begin{array}{ll}
1, & \text{if } i=j \in Q'_0,\\
0, & \text{otherwise, }\\
\end{array} \right.
\end{equation*}
then $\vdim F^d(M)=I_{Q'_0} \vdim M$ and $||M||=||F^d(M)||$ \cite[lemma 25.4]{BSZ09}.

\begin{lema} \label{(R)L:paso}
Let $Q$ be a finite quiver and $(kQ,\delta)$ a ditalgebra with triangular differential.
Then there exists a subquiver $Q'$ of $Q$ with $Q'_0=Q_0$, 
a triangular differential $\delta'$ in $kQ'$ and a ditalgebra morphism
\[
 \varphi:\mathcal{A}=(kQ,\delta) \longrightarrow \mathcal{A}'=(kQ',\delta'),
\]
such that the induced functor $F_{\varphi}:\mathcal{A}'\text{-mod} \to \mathcal{A}\text{-mod}$ is an equivalence
of categories and the quiver $Q'$ has either no solid arrows or has at least one solid arrow  $\alpha $ with $\delta'(\alpha)=0$.
Moreover, the functor $F_{\varphi}$ is rigid and for any $\mathcal{A}'$-module $M$ we have 
$\vdim M=\vdim F_{\varphi}(M)$ and $||M||\leq ||F_{\varphi}(M)||$.
\end{lema}
\bproof
We start with a ditalgebra $(kQ,\delta)$ which has solid arrows but non of them has zero differential.
By triangularity the smallest solid arrow $\alpha_0$ satisfies 
$\delta(\alpha_0)=\sum_{i=0}^rc_i\gamma_i$ for some dotted arrows $\gamma_i$ and some nonzero scalars $c_i$. 
Through a change of bases if necessary we can assume that $r=0$ and $c_0=1$.

Let $(kQ^1,\delta^1)$ be the regularization of the arrows $\alpha_0$ and $\gamma_0$ and let 
$\varphi^1:(kQ,\delta) \to (kQ^1,\delta^1)$ be the associated morphism, composed with the change of basis morphism
if performed. Then the induced functor $F_{\varphi^1}:(kQ^1,\delta^1)\text{-mod} \to (kQ,\delta)\text{-mod}$ 
is an equivalence of categories.

If $(kQ^1,\delta^1)$ has no solid arrows or has at least one solid arrow $\alpha$ with $\delta(\alpha)=0$, we take 
$(kQ',\delta')=(kQ^1,\delta^1)$. On the contrary we repeat the initial argument in the ditalgebra $(kQ^1,\delta^1)$ 
to obtain $(kQ^2,\delta^2)$ and a morphism of ditalgebras  $\varphi^2:(kQ^1,\delta^1) \to (kQ^2,\delta^2)$
such that the induced functor $F_{\varphi^2}$ is an equivalence. Since $Q$ is a finite quiver, after a finite number
of repetitions we end up with a ditalgebra $(kQ^{\ell},\delta^{\ell})$ which has no solid arrows or has at least one
solid arrow with zero differential and a morphism of ditalgebras $\varphi=\varphi^{\ell} \circ \ldots \circ \varphi^1$ 
such that the induced functor $F_{\varphi}$ is an equivalence of categories.
We take then $(Q',\delta')=(Q^{\ell},\delta^{\ell})$.

Now, the functors associated to the elementary reductions change of basis and regularization are rigid,
preserve dimension vectors and do not decrease norm. Hence the same is true for the functors 
$F_{\varphi^1},\ldots, F_{\varphi^{\ell}}$ and for their composition $F_{\varphi}$.
\eproof

\section{Reduction by admissible modules} \label{(P)S:redMod}
%----------------------------------------------------------------------
The main tool of this work is a particular case of the reduction by admissible modules as developed by
Bautista, Salmer\'on and Zuazua in \cite[chapter 12]{BSZ09}. Assume that $A$ is a finite dimensional $k$-algebra
and let $P$ be its radical. Assume also that there exists a subalgebra $S$ of $A$ with $A=S \oplus P$ as
$S$-$S$-bimodules. Then we say that \textbf{$A$ splits over its radical} (see definition~\ref{(A)D:escision}).

Let $Q$ be a finite quiver and let $\mathcal{A}=(kQ,\delta)$ be a ditalgebra with triangular differential. 
Assume moreover that $Q'_1$ is a subset of solid arrows in $Q$ with zero differential and denote by $B$
the subalgebra of $kQ$ generated by the trivial paths $e_1,\ldots,e_n$ and the arrows in $Q'_1$. 
Let $X$ be a finite dimensional left $B$-module and consider its opposed endomorphism algebra 
$\Gamma=\End_B(X)^{op}$. We will say that $X$ is an \textbf{admissible} $B$-module whenever $\Gamma$ splits
over its radical $\Gamma=S \oplus P$. In this work we will also assume that $S$ is a trivial $k$-algebra
(cf. definitions~\ref{(A)D:subditalg} and~\ref{(A)D:admisible}). Denote by $M^*=\Hom_S(M_S,S_S)$ the right $S$-dual 
of a right $S$-modulo $M_S$. The \textbf{reduction of $\mathcal{A}$ by the admissible $B$-module $X$}
is given by $\mathcal{A}^X=((kQ)^X,\delta^X)$, where $(kQ)^X$ is the tensor algebra
\[
 T_S(X^* \otimes_R W'' \otimes_R X \oplus P^*),
\]
and $W''$ is the arrow bimodule determined by the complement $Q'_1$ in $Q_1$ (definition~\ref{(A):tensRed}).
The differential $\delta^X$ is described in \cite[lemma 12.9]{BSZ09} (cf. lemmas~\ref{(A)L:sigma} and~\ref{(A)L:diferen}
in the appendix). The elements in $P^*$ have degree one, while those in $X^* \otimes W'' \otimes X$ get
their degree from $Q$. The associated reduction functor $F^X$,
\[
\xymatrix{ ((kQ)^X,\delta^X)\text{-mod} \ar[r]^-{F^X} & (kQ,\delta)\text{-mod},}
\]
is described in the proposition~\ref{(A)P:funRed}.
Lemma~13.3 and proposition~13.5 in \cite{BSZ09} state that, when $X$ is a finite dimensional admissible module, the
functor $F^X$ is full and faithful (\ref{(A)L:finita} and~\ref{(A)P:completo} in the appendix). Moreover, for any pair
of $\mathcal{A}^X$-modules $M$ and $N$, the functor $F^X$ induces an exact sequence in extension groups 
(see~\ref{(A)L:extRed}),
\[
 \xymatrix@C=.65pc{
0 \ar[r] & \Ext^1_{\mathcal{A}^X}(M,N) \ar[r] & \Ext^1_{\mathcal{A}}(F^X(M),F^X(N)) \ar[r] & \Ext^1_B(F^X(M),F^X(N)) \ar[r] & 0.
}
\]
\begin{lema} \label{(P)L:isoTens}
Let $Q$ be a finite quiver, $\mathcal{A}=(kQ,\delta)$ a ditalgebra with triangular differential and 
$B$ a subalgebra of $kQ$ as above. Let $X$ be a rigid admissible $B$-module with indecomposable decomposition of
nonisomorphic modules $X=X_1 \oplus \ldots \oplus X_r$ such that $\End_B(X_i)\cong k$ for any $i$.
Let $T$ be the matrix whose columns are given by the dimension vectors $\vdim X_i$ and let $Q^x$ be the quiver 
corresponding to the integral matrix $T^tM_QT$. Then $\vdim F^X(M)=T \vdim M$ and the Euler characteristic 
of $\mathcal{A}^X$ coincides, through the dimension vector, with the bilinear form associated to the quiver $Q^x$.
Moreover, if the quiver $Q((kQ)^X)$ associated to the graded tensor algebra
\[
(kQ)^X= T_S(X^* \otimes W'' \otimes X \oplus P^*)
\]
is regular, then it is isomorphic to $Q^x$, and hence there is an isomorphism of tensor algebras $(kQ)^X \cong kQ^x$.
\end{lema}
\bproof
We show first that $\vdim F^X(M)=T \vdim M$. Observe that $\vdim F^X(M)$ coincides with the dimension vector of the
$B$-module $X \otimes_S M$. So if $M$ is a simple module $S(j)$ then we have
\begin{eqnarray}
\vdim F^X(S(j)) & = & \vdim X \otimes_S S(j) = \vdim \left( \bigoplus_{i=1}^r X_i \right) \otimes_S S(j)= \nonumber \\
& = & \vdim X_j = T \vdim S(j), \nonumber
\end{eqnarray}
for $\vdim S(j)$ is the canonical vector $\mathbf{e}_j$.
In general, since any right $S$-module $M$ has the form $M\cong \bigoplus_{i=1}^n m_iS(i)$ for
some nonnegative integers $m_i$, we notice that
\begin{eqnarray}
 \vdim F^X(M) & = & \vdim X \otimes_S \left( \bigoplus_{i=1}^n m_iS(i) \right)
= \sum_{i=1}^n m_i \vdim X\otimes_S S(i) = \nonumber \\
& = &  \sum_{i=1}^n m_i T \vdim S(i) = T \vdim \left( \bigoplus_{i=1}^n m_iS(i) \right) = \nonumber \\
& = & T \vdim M. \nonumber
\end{eqnarray}

Since $F^X$ is full and faithful and induces an isomorphism in extension groups ($X$ is rigid and $\Ext^1_B(F^X(M),F^X(N))=0$),
using lemma~\ref{(P)L:Euler} we have the following equalities
 \begin{eqnarray}
 \dimk \Hom_{\mathcal{A}^X}(M,N) -\nonumber \\
- \dimk \Ext_{\mathcal{A}^X}(M,N) & = &  \dimk \Hom_{\mathcal{A}}(F^X(M),F^X(N))- \nonumber \\
&  & - \dimk\Ext_{\mathcal{A}}(F^X(M),F^X(N)) = \nonumber \\
& = & \langle \vdim F^X(M) , \vdim F^X(N) \rangle_Q =  \nonumber \\
& = & \langle T \vdim M, T \vdim N \rangle_Q = \nonumber \\
& = & (\vdim M)^tT^tM_QT(\vdim N)= \nonumber \\
& = & \langle \vdim M, \vdim N \rangle_{Q^x}. \nonumber
 \end{eqnarray}
This shows that the Euler characteristic of $\mathcal{A}^X$ correspondes to the bilinear form $\langle \cdot,\cdot \rangle_{Q^x}$,
that is, $M_{Q((kQ)^X)} =T^tM_QT$. If the quiver associated to $(kQ)^X$ is regular, by lemma~\ref{(P)L:carReg} we have that
$Q(M_{Q((kQ)^X)}) \cong Q((kQ)^X)$. Thus
\[
 Q^x=Q(T^tM_QT)=Q(M_{Q((kQ)^X)})\cong Q((kQ)^X).
\]
Using the isomorphism $kQ(A) \cong A$ for an elementary tensor algebra $A$ we have
\[
 kQ^x  \cong kQ((kQ)^X) \cong (kQ)^X,
\]
which completes the proof.
\eproof

In general, if the quiver $Q((kQ)^X)$ is not regular then the tensor algebras $kQ^x$ and $(kQ)^X$ are not
isomorphic. Examples of this can be found in the proofs of lemmas~\ref{(P)L:redDos} and \ref{(P)L:Kro3Dos} 
in the next chapter.

The first example of reduction consists in assuming that $B$ is the subalgebra of $kQ$ generated only by trivial paths
and that $X$ is direct sum of all simple modules $S(i)$ of vertex $i$. Then $\Gamma \cong R$, $(kQ)^X \cong kQ$ and 
$\delta^X \cong \delta$. One of the simplest nontrivial examples, called \textbf{reduction of an edge}, is described next.

%--------------------
\subsubsection*{Reduction of an edge.} \label{(E):eje}
Assume that $Q'_1$ consists of one solid arrow $ \alpha:i_0 \to j_0$ ($i_0 \neq j_0$) with zero differential. 
Then $B$ is the path algebra of the quiver 
\[
\mathbf{A_2}=\xymatrix{{}_{i_0} \ar[r]^-{\alpha} & {}_{j_0}}
\]
plus certain number of isolated vertices. Denote by $X_i$ the simple $B$-module of vertex $i$, 
$S(i)$, and fix a nonzero vector $x_{i}^i$ in $X_i$.  Let $X_z=E$ be the exceptional $\mathbf{A_2}$-module 
$E=\xymatrix{kx_{i_0}^z \ar[r]^-{1} & kx_{j_0}^z}$ considered as $B$-module, and notice that $\End_B(E) \cong k$. 
Take the reduction $B$-module $X$ as the direct sum
\[
 X=\left( \bigoplus_{i=1}^{n}X_i \right) \oplus X_z.
\]
Besides the endomorphisms of the indecomposable modules of $X$, one has the nonzero morphisms
$\xymatrix{
X_{j_0} \ar[r]^-{\mu} & X_z \ar[r]^-{\nu} & X_{i_0}}$ given by
\[
 \xymatrix{
0 \ar[r] \ar[d] & E_{i_0} \ar[r]^-{1} \ar[d]^-{1} & S(i_0)_{i_0} \ar[d] \\
S(j_0)_{j_0} \ar[r]_-{1} & E_{j_0} \ar[r] & 0.
}
\]
Consider the opposed endomorphism algebra $\Gamma$ of the $B$-module $X$ as a matrix algebra $[\Hom_{B}(X_i,X_j)]$
for $i,j \in \{1,\ldots,n,z\}$.
Then $\Gamma$ splits over its radical $\Gamma=S \oplus P$, where the diagonal $S$ 
consists of copies of the field (since $\End_B(X_i)\cong k$ for any $i$)
and $P$ is generated by the morphisms $\mu$ and $\nu$ considered as elements of $\Gamma$.
Thus the $B$-module $X$ is admissible. The reduction of $\mathcal{A}=(kQ,\delta)$ 
respect to $X$ is denoted by $\mathcal{A}^e=((kQ)^X,\delta^e)$. 

The reduced tensor algebra $(kQ)^X$ is isomorphic to the path algebra of a quiver $Q^e$ that can be described as follows.
Consider the following partition of the set of arrows $Q_1$,
\[
 Q_1=\{ \alpha \} \sqcup Q_1(\emptyset) \sqcup Q_1(s) \sqcup Q_1(t) \sqcup Q_1(s,t),
\]
where
\begin{eqnarray}
 Q_1(\emptyset) & = & \{ \xymatrix@C=1pc{\beta \in Q_1 \; | \;  s(\beta)\notin \{ i_0,j_0\},t(\beta)\notin \{ i_0,j_0\}} \}, \nonumber \\
 Q_1(s)   & = & \{ \xymatrix@C=1pc{\beta \in Q_1 \; | \;  s(\beta)\in \{ i_0,j_0\},t(\beta)\notin \{ i_0,j_0\}} \}, \nonumber \\
 Q_1(t)   & = & \{ \xymatrix@C=1pc{\beta \in Q_1 \; | \;  s(\beta)\notin \{ i_0,j_0\},t(\beta)\in \{ i_0,j_0\}} \}, \nonumber \\
 Q_1(s,t) & = & \{ \xymatrix@C=1pc{\beta \in Q_1 \; | \; \beta\neq \alpha,  s(\beta)\in \{ i_0,j_0\},t(\beta)\in \{ i_0,j_0\}} \}. \nonumber 
\end{eqnarray}
The vertex set in the reduced quiver $Q^e$ is $Q^e_0=Q_0 \sqcup \{z\}$.
The set of arrows is given by
\[
 Q^e_1= Q^e_1(\emptyset) \sqcup Q^e_1(s) \sqcup Q^e_1(t) \sqcup Q^e_1(s,t) \sqcup \{ \alpha_{\mu},\alpha_{\nu} \},
\]
where the arrows $\alpha_{\mu}:j_0 \to z$ and $\alpha_\nu:z \to i_0$ correspond to the elements 
$\mu^*$ and $\nu^*$ of $P^*$ respectively, and
\begin{eqnarray}
 Q_1^e(\emptyset) & = & Q_1(\emptyset), \nonumber \\
 Q_1^e(s) & = & Q_1(s) \sqcup Q_1(s)^z, \nonumber \\
 Q_1^e(t) & = & Q_1(t) \sqcup {}^zQ_1(t), \nonumber \\
 Q_1^e(s,t) & = & Q_1(s,t) \sqcup Q_1(s,t)^z \sqcup {}^zQ_1(s,t) \sqcup {}^zQ_1(s,t)^z. \nonumber 
\end{eqnarray}
In the last notation, if $D$ is a subset of $Q_1$ and $s^e, t^e$ are the source function and target function of $Q^e$, 
then $D^z$ denotes a copy of $D$ where each arrow $\beta^z \in D^z$, copy of $\beta \in D$, 
satisfies $s^e(\beta^z)=z$ and $t^e(\beta^z)=t(\beta)$. 
In a similar way ${}^zD$ is a copy of arrows in $D$ which consists of arrows ${}^z\beta \in {}^zD$, copy of $\beta \in D$, 
with $s^e({}^z\beta)=s(\beta)$ and $t^e({}^z\beta)=z$. Similarly can be defined ${}^zD^z$. 
The arrows $\alpha_{\mu}$ and $\alpha_{\nu}$ have degree one. All other arrows take their degree from the quiver $Q$
\cite[lemma 23.18]{BSZ09}. The triangular differential $\delta^e$, as well as the associated functor 
$F^e:\mathcal{A}^e\text{-mod} \to \mathcal{A}\text{-mod}$, are described in \cite[lemma 12.10]{BSZ09}. The functor $F^e$ is an 
equivalence of categories \cite[lemma 25.8]{BSZ09}.

In the following example $Q_1(\emptyset)$ is the empty set, $Q_1(s)=\{\beta_1\}$, $Q_1(t)=\{\beta_2\}$ and $Q_1(s,t)=\{\beta_0\}$,
\[
Q= \xymatrix{
{ \bullet } \ar[rd]_-{\beta_2} & & & & { \bullet } \\
& {\bullet }_{i_0} \ar@<.5ex>[rr]^-{\alpha} \ar@<-.5ex>@{<-}[rr]_-{\beta_0} & 
& {\bullet }_{j_0} \ar[ru]_-{\beta_1}
}
\]
Its reduced quiver $Q^e$ has a loop and is not regular,
\[
Q^e= \xymatrix{
{ \bullet } \ar[rd]_-{\beta_2} \ar[rr]^-{{}^z\beta_2} & & {\bullet}_z \ar@<-.5ex>@{.>}[dl]_-{\alpha_{\nu}} 
\ar@(lu,ru)[]^(.3){{}^z\beta^z_0} \ar@<.5ex>@{<.}[dr]^-{\alpha_{\mu}} \ar[rr]^-{\beta^z_1} & & { \bullet } \\
& {\bullet }_{i_0} \ar@{<-}[rr]_-{\beta_0} \ar@{<-}@<-.5ex>[ru]_-{\beta^z_0} & 
& {\bullet }_{j_0} \ar[ru]_-{\beta_1} \ar@<.5ex>[ul]^-{{}^z\beta_0}
}
\]

Let $M^e$ be a module of the reduced ditalgebra $\mathcal{A}^e$ and $M=F^e(M^e)$.
We give explicit transformations that constitude $M$ in terms of the transformations $M_{\beta^e}$
(for a solid arrow $\beta^e$ in $Q^e$). 
Let $\beta$ be a solid arrow in $Q$. Then the transformation $M_{\beta}$ is given by
\begin{equation} \label{(P)E:transformaRed}
M_{\beta} = \left\{
\begin{array}{l l}
M^e_{\beta}, & \text{si }\beta \in Q_1(\emptyset),\\
\left[ \begin{matrix} M^e_{\beta} & M^e_{\beta^z} \end{matrix} \right], & \text{si } \beta \in Q_1(s),\\
\left[ \begin{matrix} M^e_{\beta} \\ M^e_{{}^z\beta} \end{matrix} \right], & \text{si } \beta \in Q_1(t),\\
\left[ \begin{matrix} M^e_{\beta} & M^e_{\beta^z} \\ M^e_{{}^z\beta} & M^e_{{}^z\beta^z} \end{matrix} \right], 
& \text{si } \beta \in Q_1(s,t),\\
\left[ \begin{matrix} 0 & 0 \\ 0 & Id_{M^e_z} \end{matrix} \right], & \text{si } \beta =\alpha.\\
\end{array} \right.
\end{equation} 
This description is obtain by analysing the recipe of $F^e$.

Consider now the vectors $\vdim M^e=(m^e_i)\in \mathbb{Z}^{n+1}$ and $\vdim M=(m_i)\in \mathbb{Z}^{n}$. 
Observe that the left $S$-module $M^e$ has the form 
\[
 M^e=\left( \bigoplus_{i=1}^n m^e_iS(i) \right) \oplus m^e_zS(z).
\]
Then $\vdim E\otimes_S M^e=m^e_z \vdim X_{i_0}+m^e_z \vdim X_{j_0}$ (recall that 
$X_z=E$ is the nonsimple exceptional $\mathbf{A_2}$-module $\xymatrix{kx_{i_0}^z \ar[r]^-{1} & kx_{j_0}^z}$) and hence
\begin{eqnarray}
 \vdim F^e(M^e) & = & \vdim X \otimes_S M^e = \sum_{i=1}^n\vdim X_i \otimes_S M^e +\vdim X_z \otimes_S M^e = \nonumber \\
& = & \sum_{i=1}^nm^e_i\vdim X_i \otimes_S S(i) +\vdim E \otimes_S M^e = \nonumber \\
& = & \sum_{i=1}^nm^e_i\vdim X_i+ m^e_z \vdim X_{i_0}+m^e_z \vdim X_{j_0}, \nonumber
\end{eqnarray}
that is,
\begin{equation} \label{(P)E:redEje} 
m_i = \left\{
\begin{array}{ll}
m^e_i, & \text{if } i\neq i_0,j_0,\\
m^e_i+m^e_z, & \text{if } i = i_0,j_0.
\end{array} \right.
\end{equation}  
Using this equality and the partition of arrows in $Q_1$ it is possible to compare the norms of $M^e$ and $M$,
\begin{eqnarray}
 ||M|| & = & \sum_{\substack{\beta \in Q_1 \\ |\beta|=0}}m_{s(\beta)}m_{t(\beta)} = \nonumber \\
& = & m_{i_0}m_{j_0} + \sum_{\substack{\beta \in Q_1(\emptyset) \\ |\beta|=0}}m^e_{s(\beta)}m^e_{t(\beta)}+
\sum_{\substack{\beta \in Q_1(s) \\ |\beta|=0}}(m^e_{s(\beta)}+m^e_z)m^e_{t(\beta)}+ \nonumber \\ 
& & +\sum_{\substack{\beta \in Q_1(t) \\ |\beta|=0}}m^e_{s(\beta)}(m^e_{t(\beta)}+m^e_z)+  
\sum_{\substack{\beta \in Q_1(s,t) \\ |\beta|=0}}(m^e_{s(\beta)}+m^e_z)(m^e_{t(\beta)}+m^e_z). \nonumber 
\end{eqnarray}
Considering the partition of the set of arrows $Q^e_1$ one has
\begin{eqnarray}
||M||& = & m_{i_0}m_{j_0} + \left( \sum_{\substack{\beta \in Q_1(\emptyset) \\ |\beta|=0}}m^e_{s^e(\beta)}m^e_{t^e(\beta)} \right) + \nonumber \\
& & +\left( \sum_{\substack{\beta \in Q_1(s) \\ |\beta|=0}}m^e_{s^e(\beta)}m^e_{t^e(\beta)} +
 \sum_{\substack{\beta^z \in Q_1(s)^z \\ |\beta^z|=0}}m^e_{s^e(\beta^z)}m^e_{t^e(\beta^z)} \right) + \nonumber \\
&  & + \left( \sum_{\substack{\beta \in Q_1(t) \\ |\beta|=0}}m^e_{s^e(\beta)}m^e_{t^e(\beta)} +
 \sum_{\substack{{}^z\beta \in {}^zQ_1(t) \\ |{}^z\beta|=0}}m^e_{s^e({}^z\beta)}m^e_{t^e({}^z\beta)} \right) + \nonumber \\
& & + \left( \sum_{\substack{\beta \in Q_1(s,t) \\ |\beta|=0}}m^e_{s^e(\beta)}m^e_{t^e(\beta)} + 
\sum_{\substack{\beta^z \in Q_1(s,t)^z \\ |\beta^z|=0}}m^e_{s^e(\beta^z)}m^e_{t^e(\beta^z)} \right. + \nonumber \\
& & + \left. \sum_{\substack{{}^z\beta \in {}^zQ_1(s,t) \\ |{}^z\beta|=0}}m^e_{s^e({}^z\beta)}m^e_{t^e({}^z\beta)}+
\sum_{\substack{{}^z\beta^z \in {}^zQ_1(s,t)^z \\ |{}^z\beta^z|=0}}m^e_{s^e({}^z\beta^z)}m^e_{t^e({}^z\beta^z)} \right)  \nonumber \\
& = & m_{i_0}m_{j_0} + ||M^e||, \nonumber
\end{eqnarray}
that is, $||F^e(M^e)||=||M|| = m_{i_0}m_{j_0}+ ||M^e||$. In particular, if $m_{i_0}\neq 0$ and $m_{j_0}\neq 0$ 
then $||M^e||<||M||$.

\section{Exceptional representations.} \label{(P)S:repExc}
%----------------------------------------------------------------------
Let $Q$ be a finite quiver and $\delta$ a triangular differential of $kQ$.
Next we give some properties of the exceptional representations of $(kQ,\delta)$.

\begin{lema} \label{(P)L:sub}
Let $Q'$ be a solid subquiver of $Q$ with zero differential. If $M$ is a rigid $(kQ,\delta)$-module 
then its restriction $M|_{Q'}$ is rigid as $kQ'$-module. 
In parti\-cu\-lar, for any solid arrow $\alpha:i \to j$ of zero differential, the linear transformation
$M_{\alpha}$ is injective or surjective. And if $\alpha:i \to i $ is a loop with zero differential
then $M_{i}=0$.
\end{lema}

\bproof
Let $Q''$ be the quiver obtained by adding to $Q'$ all vertices in $Q$ which are not in $Q'_0$, in such a way that
there is in inclusion of algebras $kQ'' \rightarrow kQ$. Clearly a $Q''$-module $N$ is rigid if and only if
its restriction $N|_{Q'}$ is rigid, thus we can assume that $Q''=Q'$.
Let $W_0$ and $W_0'$ be the bimodules of solid arrows in $Q$ and $Q'$ and denote by $i:W'_0 \to W_0$ the bimodule inclusion.
Notice that the following diagram is commutative, where $\sigma$ and $\sigma'$ are the transformations given in~(\ref{(P)EQ:sigma}).
\[
\xymatrix{
\txt{$\Hom_R(M,N)$ \\ $\oplus$ \\ $\Hom_{R\text{-}R}(W_1,\Hom_k(M,N))$} \ar[d]_-{[I \; 0]} \ar[r]^-{\sigma} 
& \Hom_R(W_0\otimes_R M,N) \ar[d]^-{\Hom(i \otimes M,N)} \\
\Hom_R(M,N) \ar[r]_-{\sigma'} & \Hom_R(W'_0\otimes_R M,N).}
\]
Indeed, on the one hand $\sigma'[I \; 0](f^0,f^1)=\sigma'(f^0,0):w'\otimes m \mapsto w'f^0(m)-f^0(w'm)$.
On the other hand $\Hom(i\otimes M,N)\sigma(f^0,f^1)$ is the morphism that sends $w'\otimes m$ into 
$w'f^0(m)-f^0(w'm)-\widehat{f^1}(\delta(w'))(m)$. This is iqual to $w'f^0(m)-f^0(w'm)$ since the elements in $W'_0$ have 
zero differential. Then the diagram above is commutative, and obtaining cokernels one gets
\begin{equation*}
\xymatrix@C=1pc{
\Hom_R(W_0 \otimes_R M,N) \ar@{->>}[d]_-{\Hom(i \otimes M,N)} \ar[r]^-{\eta} &  \Ext^1_{(kQ,\delta)}(M,N) \ar@{->>}[d] \ar[r]  & 0 \\
\Hom_R(W'_0 \otimes_R M,N) \ar[r]^-{\eta} &  \Ext^1_{kQ'}(M|_{Q'},N|_{Q'}) \ar[r]  & 0. 
}
\end{equation*}
Moreover, the morphism $\Hom(i\otimes M,N)$ is surjective since $i$ is injective and $R$ is semi-simple. 
Hence the induced morphism in cokernels is also surjective and 
$\Ext^1_{(kQ,\delta)}(M,M)=0$ implies that $\Ext^1_{kQ'}(M|_{Q'},M|_{Q'})=0$.

In order to prove the second claim assume that $i \neq j$. The first claim implies that the $\mathbf{A}_2$-representation 
\[
M|_{\mathbf{A}_2}=(M_i,M_j;M_{\alpha})
\]
is rigid. If the transformation $M_{\alpha}$ is neither injective nor surjective, then $M|_{\mathbf{A}_2}$ contains as direct 
summand the module $S(i) \oplus S(j)$. But this is impossible since $S(i) \oplus S(j)$ is not rigid 
(for $\langle e_i,e_j \rangle_{\mathbf{A}_2}=-1$, see lemma~\ref{(P)L:Euler}).

Assume now that $i=j$ and let $L$ be the subquiver of $Q$ given by the vertex $i$ and the loop $\alpha$. 
Then $M|_{L}=(M_i;M_{\alpha})$ is a rigid $kL$-module, and since the bilinear form associated to $L$ is identically zero,
again by lemma~\ref{(P)L:Euler} we have
\[
 0=\langle \vdim M|_{L} , \vdim M|_{L} \rangle_L =\dimk_k\Hom_{kL}(M|_L,M|_L).
\]
Then $\End_{kL}(M|_L)=0$ and hence $M_i=0$.
\eproof

We will say that a module $M$ is \textbf{sincere} if all entries in the dimension vector $\vdim M$ are nonzero.

\begin{lema} \label{(P)L:determ}
Let $Q$ be a finite quiver and $\mathcal{A}=(kQ,\delta)$ a ditalgebra with triangular differential. 
If $M$ and $N$ are exceptional representations of $\mathcal{A}$ with $\vdim M=\vdim N$ then $M \cong N$.
\end{lema}

\bproof
Deleting all vertices in $Q$ where $M_i=0$ we can assume that $M$ (and hence $N$) is a sincere representation.
We proceed by induction over the norm $||M||$. Since $\vdim M=\vdim N$ we have $||M||=||N||$.
Moreover $||M||=0$ if and only if $M$ and $N$ are simple modules corresponding to a vertex $i$, and hence
the statment is valid for norm zero modules. Assume that it is valid for representations with norm igual
or smaller than $\ell$ and let $M$ be a module with $||M||=\ell+1$. 
Since $||M||>0$, by lemma~\ref{(R)L:paso} we can assume that there is a solid arrow $\alpha:i_0 \to j_0$ 
with zero differential (the equivalence $F_{\varphi}$ in lemma~\ref{(R)L:paso} preserves exceptional modules and
does not descrease norm). As a consequence of lemma~\ref{(P)L:sub}, the transformation $M_{\alpha}$ is not zero, 
for $M$ is a sincere module. In particular $\alpha$ is not a loop and the reduction of an edge $M^e$ respect to $\alpha$
of the module $M$  (that is, $M^e$ is an $\mathcal{A}^e$-module such that $F^e(M^e)\cong M$) satisfies $||M^e||\leq \ell$. 
Moreover $M^e$ and $N^e$ are exceptional modules because of the exact sequence before lemma~\ref{(P)L:isoTens}.
We show that $\vdim M^e=\vdim N^e$. The last equality~(\ref{(P)E:redEje}) implies that
\[
 (m_1,\ldots,m_{i_0},\ldots,m_{j_0},\ldots,m_n)=(m^e_1,\ldots,m^e_{i_0}+m^e_z,\ldots,m^e_{j_0}+m^e_z,\ldots,m^e_n),
\]
\[
 (n_1,\ldots,n_{i_0},\ldots,n_{j_0},\ldots,n_n)=(n^e_1,\ldots,n^e_{i_0}+n^e_z,\ldots,n^e_{j_0}+n^e_z,\ldots,n^e_n).
\]
As a consequence of lemma~\ref{(P)L:sub}, the transformation $M_{\alpha}$ is injective or surjective (depending on the value
of $m_{i_0}-m_{j_0}$, see description of $M_{\alpha}$ in~\ref{(P)E:transformaRed}) and the same is true for $N_{\alpha}$. 
Then, since $\vdim M=\vdim N$, we have $\vdim M^e=\vdim N^e$. So by inductive hypothesis $M^e \cong N^e$, 
and since $F^e$ is an equivalence we conclude that $M \cong N$.
\eproof

\begin{lema} \label{(P)L:endExc}
If $Q$ is a finite solid quiver with an admissible ordering of its vertices,
then the endomorphism algebra of any exceptional $kQ$-module is isomorphic to the field $k$.
\end{lema}
\bproof
Cf. corollary 1 in Ringel \cite{cmR94}.
\eproof

In this case, as a consequence of lemma~\ref{(P)L:determ}, the dimension vector $\vdim$ induces an injective function
between the set of isomorphism classes of exceptional modules and the set of positives roots of $Q$.

Given a $(kQ,\delta)$-representation $M$ we fix a basis $\mathscr{B}_i$ for each $k$-vector space
$M_i$, $i \in Q_0$. We say that the disjoint union $\mathscr{B}=\bigsqcup \mathscr{B}_i$ is a \textbf{basis of the module} $M$. 
We define the \textbf{coefficient quiver} $\mathcal{C}(M,\mathscr{B})$ of $M$ respect to the basis $\mathscr{B}$ as the 
solid quiver with vertex set $\mathscr{B}$ and arrows given in the following way. For each solid arrow $\alpha:i \to j$ in 
$Q$ we add a (s\'olida) arrow from $e^i \in \mathscr{B}_i$ to $e^j\in \mathscr{B}_j$ in $\mathcal{C}(M,\mathscr{B})$ 
if the coefficient of the matrix $M_{\alpha}$ corresponding to the position $(e^i,e^j)$ is different from zero.

\begin{lema} \label{(P)L:conexo}
A $(kQ,\delta)$-module $M$ is indecomposable if and only if for any basis $\mathscr{B}$ of $M$, 
the coefficient quiver $\mathcal{C}(M,\mathscr{B})$ is connected.  
\end{lema}
\bproof
We show first that if $\mathcal{C}(M,\mathscr{B})$ is disconnected for some basis $\mathscr{B}$ of $M$,
then $M$ splits. Let $\mathcal{C}_1$ be a connected component of $\mathcal{C}(M,\mathscr{B})$ (nonempty and different 
form the total) and let $\mathcal{C}_2=\mathcal{C}_1^c$ be its complement.
We will define an endomorphism $f$ of $M$ in $(kQ,\delta)$-mod, which will be a nontrivial idempotent.
Consider the mapping $f^0: \mathscr{B} \to M$ given by
\begin{equation*}
f^0(x) = \left\{
\begin{array}{ll}
x, & \text{if } x \in \mathcal{C}_1,\\
0, & \text{if } x \in \mathcal{C}_2.
\end{array} \right.
\end{equation*}
Let $A$ be the degree zero subalgebra of $kQ$. We will show that the linear extension to $M$ of $f^0$, which we will denote
with the same symbol $f^0$, is a morphism of $A$-modules. Let $\alpha$ be a solid arrow in $Q$. 
If $x \in \mathcal{C}_1$ then $M_{\alpha}x=\sum_{i=1}^r a^ix^i_1$ with $x^i_1 \in \mathcal{C}_1$ and $a^i \in k$, 
for $\mathcal{C}_1$ is a component of the coefficient quiver $\mathcal{C}(M,\mathscr{B})$. Then
\[
 \alpha f^0(x)=\alpha x=\sum_{i=1}^ra^ix_1^i = \sum_{i=1}^ra^if^0(x_1^i)= f^0(\sum_{i=1}^ra^ix_1^i)=f^0(\alpha x).
\]
On the other hand, if $x \in \mathcal{C}_2$ then $M_{\alpha}x=\sum_{j=1}^u b^jx^j_2$ with $x^j_2 \in \mathcal{C}_2$
and $b^j \in k$, and hence $\alpha f^0(x)=\alpha 0=0$ and
\[
 f^0(\alpha x)=f^0 \left( \sum_{j=1}^ub^jx_2^j \right) = \sum_{j=1}^ub^jf^0(x_2^j)=0.
\]
In this way $f^0$ is a morphism of $A$-modules, thus $f=(f^0,0)$ is an endomorphism of the $(kQ,\delta)$-module $M$. 
Clearly $f$ is a nontrivial idempotent. Since idempotents split in $(kQ,\delta)$-mod \cite[lemma 5.12]{BSZ09}, 
we conclude that $M$ splits.

Assume now that $M$ decomposes as $(kQ,\delta)$-module. We want to show that there exists a basis
$\mathscr{B}$ of $M$ such that $\mathcal{C}(M,\mathscr{B})$ is disconnected.  For that we verity that $M$ splits
as $A$-module. Indeed, if $M$ decomposes as $(kQ,\delta)$-module, then there is an endomorphism $g$ of $M$ which is
a nontrivial idempotent. By lemma~5.7 in \cite{BSZ09}, we can assume that $g=(g^0,0)$.  Then $g^0$ is an endomorphism
of $A$-modules, which is a nontrivial idempotent. Hence $M$ decomposes as $A$-module, that is, there exists a 
nontrivial decomposition $M = M_1 \oplus M_2$. Let $\mathscr{B}_1$ and $\mathscr{B}_2$ be bases of $M_1$ and $M_2$.
We will verify that the full subgraphs $\mathcal{C}_1$ and $\mathcal{C}_2$ of the coefficient quiver
$\mathcal{C}=\mathcal{C}(M,\mathscr{B}_1 \sqcup \mathscr{B}_2)$ determined by the partition 
$\mathscr{B}_1 \sqcup \mathscr{B}_2$ separate $\mathcal{C}$. For since $M_1$ and $M_2$ are submodules of $M$,
if there exists an arrow $\alpha:x \to y$ in $\mathcal{C}$ then $x \in \mathcal{C}_i$ implies $y \in \mathcal{C}_i$ 
for $i \in \{ 1,2 \}$.
\eproof

A $(kQ,\delta)$-module $M$ is called \textbf{tree module} if it has a basis $\mathscr{B}$ 
respect to which the coefficient quiver $M$ is a tree. As Ringel has shown \cite[property 2]{cmR98}, 
in this case bases can be adjust in such a way that the matrices $M_{\alpha}$ exhibit only $0$ and $1$ coefficients.
The following result is a generalization to the ditalgebras setting of a well known theorem of Ringel~\cite{cmR98}.
\begin{proposicion} \label{(P)P:excep}
Let $Q$ be a finite quiver and $\mathcal{A}=(kQ,\delta)$ a ditalgebra with triangular differential. If $M$
is an exceptional representation of $\mathcal{A}$ then $M$ is a tree module.
\end{proposicion}

\bproof
Deleting some vertices we can assume that $M$ is a sincere representation.
We proceed by induction on the norm of $M$.
The statement is evident for norm zero exceptional modules (the simple modules $S(i)$). 
Assume that the claim is valid for exceptional modules with norm less or equal to $\ell$.
Let $M$ be an exceptional module with $||M||=\ell+1$. 
Since $||M||>0$, by lemma~\ref{(R)L:paso} we can assume the existence of a solid arrow  $\alpha:i_0 \to j_0$ 
with zero differential. As a consequence of lemma~\ref{(P)L:sub}, the transformation $M_{\alpha}$ is not null, 
for $M$ is a sincere module.
In particular $\alpha$ is not a loop and the reduction $M^e$ of the module $M$, which is also exceptional, 
satisfies $||M^e||\leq \ell$.

Let $\mathscr{B}^e=\bigsqcup_{i=1}^{n,z} \mathscr{B}^e_i$ be a basis of the reduced representation $M^e$ such that the
coefficient quiver $T^e=\mathcal{C}(M^e,\mathscr{B}^e)$ is a tree.
Together with the selected basis $\{ x^i_{i},x^z_{i_0},x^z_{j_0} \}_{i \in Q_0}$ of $X$ (see description of the reduction of
an edge in page~\pageref{(E):eje}), the set $\mathscr{B}^e$ determines a basis of the representation $X \otimes_S M^e \cong M$. 
Indeed, for $i \notin \{i_0,j_0 \}$ we can take
\[
 \mathscr{B}_{i}=\{ x^i_{i} \otimes b \; | \; b \in \mathscr{B}^e_{i}\},
\]
and for $i\in \{i_0,j_0\}$ take
\[
\mathscr{B}_i=\{ x^i_i \otimes b,x^z_i \otimes b \; | \; b  \in \mathscr{B}^e_i \}.
\]
Then it is clear that $\mathscr{B}=\bigsqcup_{i=1}^n \mathscr{B}_{i}$ is a basis of $M$ (since $M$, considered as a
$B$-module, is given by $X\otimes M^e$, see proposition~\ref{(A)P:funRed}). 
Since the module $M^e$ is indecomposable, by lemma~\ref{(P)L:conexo} the quiver $T^e$ is a connected tree. 
Hence $\dimk_kM^e=|T^e_0|=1+|T^e_1|$. 
Take $\vdim(M^e)=(m_1^e,\ldots,m^e_n,m^e_z)$ and $\vdim M=(m_1,\ldots,m_n)$. By the equality~(\ref{(P)E:redEje}) we have that
$\dimk_k M=\dimk_k M^e + m^e_z$, that is, $|T_0|=|T^e_0|+m^e_z$ (where $m^e_z=\dimk_k M^e_z$).

On the other hand, in the description of $M$ given in equation~(\ref{(P)E:transformaRed}) can be observed that
the only transformation $M_{\beta}$ wich is not constructed with the trasformations $M^e_{\beta}$,
occurs when $\beta$ is the reduction edge $\alpha$, and in that case
\[
 M_{\alpha}=\left[ \begin{matrix} 0 & 0 \\ 0 & Id_{M^e_z} \end{matrix} \right].
\]
Then the coefficient quiver $T=\mathcal{C}(M,\mathscr{B})$ has $|T_0|=|T^e_0|+m^e_z$ vertices and
$|T_1|=|T^e_1|+m^e_z$ arrows. Thus $|T_0|=1+|T_1|$ ans since $M$ is indecomposable, $T$ is connected. 
Hence $T$ is a tree. 
\eproof

\section{Almost split sequences.} \label{(DE)S:SCD}
%----------------------------------------------------------------------
%------------------------------------------------------------------
%-----------------------------------------------------------------
Let $\mathcal{K}$ be a Krull-Schmidt category. For two indecomposable objects $M$ and $N$
the \textbf{radical} $\rad(M,N)$ of $\Hom(M,N)$ is define to be the set of noninvertible morphisms from $M$ to $N$. 
For direct sums $M=\bigoplus_iM_i$, $N=\bigoplus_jN_j$ consider a morphism $f:M \to N$ as a matrix
with entries $f_{i,j}:M_i \to N_j$ (to be precise $f_{i,j}=\mu_jf\sigma_i$ where $\mu_j$, $\sigma_i$ 
are the corresponding canonical projections and inclusions). By definition the radical $\rad(M,N)$ is given by matrices 
whose entries $f_{i,j}$ belong te the radicals $\rad(M_i,N_j)$. Define the square radical 
$\rad^2(M,N)$ as the set of morphisms of the form $gf$ where $f \in \rad(M,X)$ and $g \in \rad(X,N)$ for some
object $X$. Define the \textbf{$\End_{\mathcal{K}}(N)$-$\End_{\mathcal{K}}(M)$-bimodule of irreducible morphisms} 
from $M$ to $N$ as
\[
 \Irr(M,N)=\rad(M,N)/\rad^2(M,N).
\]
A morphism $f:M \to N$ is called \textbf{irreducible} if it is neither a section nor a retraction, and for any factorization
$f=f'f''$ then either $f'$ is a retraction or $f''$ is a section. When $M$ and $N$ are indecomposable, $f$ is irreducible 
if and only if $f \in \rad(M,N)-\rad^2(M,N)$ (cf. Ringel \cite[2.2]{cmR}). 

A morphism $f:M \to E$ is \textbf{left almost split} if it satisfies the following properties
\begin{itemize}
 \item[a)] $f$ is not a section;
 \item[b)] if $h:M \to Z$ is not a section, then there exists $h':E \to Z$ such that $h=h'f$.
\end{itemize}
The morphism $f$ is \textbf{left minimal} if
\begin{itemize}
 \item[c)] for any $\gamma \in \End(E)$ such that $\gamma f=f$ then $\gamma$ is an automorphism.
\end{itemize}
We notice that if $f$ is left almot split then $M$ is indecomposable.
Left minimal almost split morphisms are determined, up to isomorphism, by the inde\-com\-po\-sable $M$.
Dually, a morphism $g:E \to N$ is \textbf{right almost split} if it satisfies 
\begin{itemize}
 \item[a')] $g$ is not a retraction;
 \item[b')] if $h:Z \to N$ is not a retraction, then there exists $h':Z \to E$ such that $h=gh'$.
\end{itemize}
The morphism $g$ is \textbf{right minimal} if
\begin{itemize}
 \item[c')] for any $\gamma \in \End(E)$ such that $g \gamma=g$ then $\gamma$ is an automorphism.
\end{itemize}
Again, if $g$ is a right almost split morphism then $N$ is indecomposable. 
Right minimal almost split morphisms are determined, up to isomorphism, by the indecomposable $N$.
The following is a well known result. See for example Ringel \cite{cmR}2.2(lemma 3).

\begin{lema} \label{(DE)L:irred}
Assume there is a left minimal almost split morphism in $\mathcal{K}$ beginning in $M$. Let $E_1,\ldots, E_r$ be nonisomorphic 
indecomposable objets and assume that there are morphisms $f_{ij}:M \to E_i$ for $j=1,\ldots,d_i$ with
residual classes $\overline{f_{ij}}$ in $\Irr(M,E_i)$. Then the morphism
\[
 f=(f_{ij})_{ij}:M \to \bigoplus_{i=1}^r d_iE_i
\]
is a left minimal almost split morphism if the set $\{ \overline{f_{ij}}\}_{j=1}^{d_i}$ is a basis of
$\Irr(M,E_i)$ for each $i=1,\ldots,r$, and any indecomposable object $E'$ such that $\Irr(M,E') \neq 0$ is isomorphic to $E_i$ 
for some $i$.
\end{lema}
The following observation will be used in section~\ref{(DE)S:rank}.
\begin{lema} \label{(DE)L:minimal}
Let $Q$ be a finite solid quiver with an admissible ordering of its vertices. 
\begin{itemize}
 \item[a)]  If $g:E \to N$ is a right almost split morphism in $kQ$-mod, then there exists a decomposition of the form
\[
 \xymatrix{
E \ar[r]^-{g} & N \\
E' \oplus E'' \ar[u]^-{\cong} \ar[ru]_-{[g'\; 0]} 
}
\]
where $g':E' \to N$ is a right minimal almost split morphism.
 \item[b)] Dually, if $f:M \to E$ is a left almost split morphism in $kQ$-mod, the there exists a decomposition of the form
\[
 \xymatrix{
M \ar[r]^-{f} \ar[rd]_-{\left[ \begin{smallmatrix} f' \\ 0 \end{smallmatrix}\right] } & E \\
& E' \oplus E'' \ar[u]^-{\cong} 
}
\]
where $f':M \to E'$ is a left minimal almost split morphism.
\end{itemize}
\end{lema}

\bproof
We show $(a)$, the proof of $(b)$ is similar. 
Assume that $g:E \to N$ is a left almost split morphism. Let $\mathbb{E}$ be the class of subobjects $\sigma_{E'}:E' \to E$ 
of $E$ for which there exists an epimorphism $h_{E'}:E \to E'$ such that $g=g\sigma_{E'} h_{E'}$. 
The collection $\mathbb{E}$ is not empty for it contains $E$. 
Let $E'$ be an element of $\mathbb{E}$. Then $g\sigma_{E'}$ is right almost split. Indeed, $g\sigma_{E'}$ is not a retraction
for $g$ is not a retraction. Assume that $v:Z \to N$ is not a retraction.
Then there exists $t:Z \to E$ such that $v=gt$ and hence $v=gt=g\sigma_{E'}h_{E'}t$, that is, $g \sigma_{E'}$ is right almost split.

Assume now that $E'$ is an element in $\mathbb{E}$ with minimal dimension over $k$. Then $g'=g\sigma_{E'}$ is
minimal. Indeed, if $h':E' \to E''$ is an epimorphism such that $g'= g' \sigma'_{E''}h'$ (where $\sigma'_{E''}:E'' \to E'$
is the subobject inclusion) then $(h'h)$ satisfies
$g\sigma_{E''}(h'h)=g\sigma_{E'}\sigma'_{E''}(h'h)=(g'\sigma'_{E''}h')h=g'h=g$. In this way $E''$ is an element of $\mathbb{E}$,
and by minimality $\dimk_k E''\geq \dimk_k E'$. Hence $h'$ is an automorphism.

Let $E' \in \mathbb{E}$ be minimal. Then the morphism $g'=g\sigma_{E'}$ is irreducible (lemma~\ref{(DE)L:irred})
thus $\sigma_{E'}$ is a section (for $g$ is not a retraction). 
Let $E''$ be a direct summand of $E$ such that $E\cong E' \oplus E''$ and $g''=g\sigma_{E''}$. 
Then $g$ has the form $g=(g' \; g'')$ and $g''$ is not a retraction for $g$ is not a retraction. Hence there exists $h'':E'' \to E'$ 
such that $g''=g'h''$, thus the following diagram is commutative
\[
 \xymatrix{
E'\oplus E'' \ar[dd]_-{\left( \begin{smallmatrix}I_{E'} & h'' \\ 0 & I_{E''} \end{smallmatrix} \right)} 
\ar[rrd]^-{\left( \begin{smallmatrix} g' & g'h'' \end{smallmatrix} \right)} \\
& & N \\
E'\oplus E'' \ar[rru]_-{\left( \begin{smallmatrix} g' & 0 \end{smallmatrix} \right)}
}
\]
This completes the proof.
\eproof

A pair $(\mathcal{K},\mathcal{E})$ is called \textbf{Krull-Schmidt category with short exact sequences}
if $\mathcal{K}$ is a Krull-Schmidt category and $\mathcal{E}$ is a collection of pairs $(f,g)$ of composable morphisms
in $\mathcal{K}$ such that $f$ is kernel of $g$ and $g$ is cokernel of $f$. For instance, if $Q$ is a finite quiver,
$(kQ,\delta)$ is a ditalgebra with triangular differential and $\mathcal{E}$ is the exact structure associated to 
$(kQ,\delta)$-mod, then $((kQ,\delta)\text{-mod},\mathcal{E})$ is a Krull-Schmidt category with short exact sequences.
We say that an exact pair $(f,g)$ is \textbf{almost split} if $f$ is left minimal almost split and
$g$ is right minimal almost split.
\begin{lema} \label{(DE)L:ARequiv}
Let $Q$ be a finite solid quiver with admissible ordering of its vertices. 
For an exact sequence $\xymatrix{0 \ar[r] & M \ar[r]^-{f} & E \ar[r]^-{g} & N \ar[r] & 0 }$ in $kQ$-mod the following
are equivalent
\begin{itemize}
 \item[a)] $(f,g)$ is an almost split exact pair.
 \item[b)] $g$ is a right minimal almost split morphism. 
 \item[b')] $f$ is a left minimal almost split morphism.
 \item[c)] $g$ is a right almost split morphism and $M$ is indecomposable.
 \item[c')] $f$ is a left almost split morphism and $N$ is indecomposable.
\end{itemize}
\end{lema}

\bproof
See for example theorem~IV.1.13 in \cite{ASS06}.
\eproof

Let $(\mathcal{K},\mathcal{E})$ be a Krull-Schmidt category with short exact sequences.
Given an almost split exact pair $\xymatrix{M \ar[r]^-{f} & E \ar[r]^-{g} & N}$ one observes that the isomorphism
class $[M]$ is determined by $[N]$, and the isomorphism class $[N]$ is determined by $[M]$. We denote then
$[M]$ by $\tau[N]$ and call $\tau[N]$ the \textbf{Auslander-Reiten translation} of $[N]$. If we have a decomposition of $E$
in nonisomorphic indecomposables $E = \bigoplus_{i=1}^rd_iE_i$ then by lemma~\ref{(DE)L:irred}
\[
\dimk_k \Irr(E_i,N)= d_i=\dimk_k \Irr(M,E_i).
\]
Define the \textbf{Auslander-Reiten quiver} $\Gamma(\mathcal{K},\mathcal{E})$ in the following way. 
The vertex set is the set of isomorphism classes of indecomposable objects. For each right minimal almost split morphism 
$g:\bigoplus d_iE_i \to N$  we add $\dimk_k \Irr(E_i,N)= d_i$ solid arrows from $[E_i]$ to $[N]$. Similarly, for each left
minimal almost split morphism $f:M \to \bigoplus d_iE_i$ we add
$\dimk_k \Irr(M,E_i)= d_i$ solid arrows from $[M]$ to $[E_i]$. Clearly $(\Gamma(\mathcal{K},\mathcal{E}),\tau)$ is a 
translation quiver, as defined as follows.

\section{Translation quivers and sections.} \label{(DE)S:trasl}
%------------------------------------------------------------------
%------------------------------------------------------------------
Let $\Gamma$ be a solid quiver. For each vertex $x$ denote by $x^-$ the set of direct predecessors of $x$ and by 
$x^+$ the set of direct successors of $x$. A quiver is called \textbf{locally finite} if for every vertex $x$ the number of arrows
between $x$ and its \textbf{neighbors} (direct successors and predecessors of $x$) is finite. 
A \textbf{translation} $\tau$ of the quiver $\Gamma$ is a bijection between two subsets $\Gamma'_0$ and $\Gamma''_0$ of $\Gamma_0$, 
such that for every $x \in \Gamma'_0$ and every $y \in \Gamma_0$ the
number of arrows from $y$ to $x$ is equal to the number of arrows from $\tau(x)$ to $y$. In particular $x^-=\tau(x)^+$.
A \textbf{translation quiver} $(\Gamma,\tau)$ consists of a locally finite solid quiver $\Gamma$ 
together with a translation $\tau$ of $\Gamma$. The vertices in $\Gamma_0-\Gamma'_0$ are called 
\textbf{projective} and those in $\Gamma_0-\Gamma_0''$ \textbf{injective}. A translation quiver $(\Gamma,\tau)$ is called
\textbf{hereditary} if the direct predecessors of a projective vertex are projective and the direct successors of an 
injective vertex are injective. The main motivation for this definitions is the Auslander-Reiten $\Gamma(A)$ of the 
category of $A$-modules of a finite dimensional $k$-algebra and the Auslander-Reiten translation.

For every vertex $x$ in a translation quiver $(\Gamma, \tau)$ there exists an interval $I_x$ of
$\mathbb{Z}$ that contains zero and such that $\tau^i x$ is defined if and only if $i \in I_x$. The \textbf{orbit} of $x$
is defined as the subset of $\Gamma_0$ given by $\mathcal{O}_x=\{ \tau^i x \}_{i \in I_x}$. 
Observe that each orbit $\mathcal{O}_x$ contains at most one projective vertex and at most one injective vertex.
If $y$ is an element in $\mathcal{O}_x$, that is, if there exists $i_0 \in I_x$
such that $y=\tau^{i_0} x$, then $I_y=I_x-i_0$. Hence
\[
 \mathcal{O}_y=\{ \tau^jy \}_{j \in I_y}=\{ \tau^{i-i_0}(\tau^{i_0}x) \}_{i \in I_x} =\mathcal{O}_x,
\]
and the set of orbits in $(\Gamma,\tau)$ determines a partition in the vertex set of $\Gamma$. 
We distinguish four types of orbits.
\begin{itemize}
 \item An orbit $\mathcal{O}$ is \textbf{injective-projective} if it contains an injective vertex $q$ and a projective
vertex $p$. In that case $\mathcal{O}=\{ \tau^iq \}_{i=0}^n$ (with $\tau^nq=p$) is a finite set.
 \item An orbit $\mathcal{O}$ is \textbf{injective} if it contains an injective vertex $q$ but contains no projective vertex.
Then $\mathcal{O}=\{\tau^iq\}_{i=0}^{\infty}$ in an infinite set.
 \item An orbit $\mathcal{O}$ is \textbf{projective} if it contains a projective vertex $p$ but it contains no injective vertex. 
Then $\mathcal{O}=\{\tau^{-i}p\}_{i=0}^{\infty}$ is an infinite set.
 \item An orbit $\mathcal{O}$ is \textbf{stable} if it does not contain neither injective nor projective vertices.
Finite stable orbits are called \textbf{periodic}.
\end{itemize}
The \textbf{orbit graph} $G_{ob}(\Gamma,\tau)$ of a translation quiver $(\Gamma,\tau)$ has as vertices the set of 
orbits of $(\Gamma,\tau)$ and there is an edge between $\mathcal{O}$ and $\mathcal{O}'$ if there exist 
elements $x\in \mathcal{O}$ and $y\in \mathcal{O}'$ which are neighbors in $\Gamma$. Clearly $G_{ob}(\Gamma,\tau)$ is
a connected graph if $\Gamma$ is a connected quiver. A translation quiver $(\Gamma,\tau)$ is called \textbf{preinjective}, 
\textbf{posprojective} or \textbf{stable} if all its orbits are inyective, projective or stable respectively.
We say that a translation quiver $(\Gamma,\tau)$ is \textbf{proper} if for any nonprojective vertex $x$ the
set $x^-$ is nonempty, and is called \textbf{directed} if $\Gamma$ has no oriented cycles.

\begin{lema} \label{(DE)L:dirigido}
If $(\Gamma,\tau)$ is a proper directed translation quiver then it has no periodic orbits.
\end{lema}

\bproof
Assume that $(\Gamma,\tau)$ is a proper translation quiver with a periodic orbit
$x_0,\tau x_0,\ldots,\tau^nx_0=x_0$ ($n\geq 0$). Then no element in $\mathcal{O}_{x_0}$ is projective, 
hence for each $0 \leq i < n$ there exists a vertex $y_i$ in $(\tau^ix_0)^-$. In this way there are arrows
\[
 x_0=\tau^nx_0 \to y_{n-1} \to \tau^{n-1}x_0 \to y_{n-2} \to \cdots \to y_1 \to x_1 \to y_0 \to x_0,
\]
that is, $(\Gamma,\tau)$ is not a directed quiver.
\eproof

A \textbf{section}  $\mathcal{S}$ in a translation quiver $(\Gamma,\tau)$ is a nonempty subset of vertices
$\mathcal{S} \subset \Gamma_0$ such that
\begin{itemize}
 \item[i)]  if $x \in \mathcal{S}$ is not a projective vertex then $\tau(x) \notin \mathcal{S}$,
 \item[ii)] if $x \in \mathcal{S}$ and $x \to m$ is an arrow in $\Gamma$ then either $m$ or $\tau(m)$ belongs to $\mathcal{S}$.
\end{itemize}
Dually, a \textbf{cosection}  $\mathcal{S}$ in a translation quiver $(\Gamma,\tau)$ is a nonempty subset of vertices 
$\mathcal{S} \subset \Gamma_0$ such that
\begin{itemize}
 \item[i')]  if $y \in \mathcal{S}$ is not an injective vertex then $\tau^{-1}(y) \notin \mathcal{S}$,
 \item[ii')] if $y \in \mathcal{S}$ and $n \to y$ is an arrow in $\Gamma$ then either $n$ or $\tau^{-1}(n)$ belongs to $\mathcal{S}$.
\end{itemize}

\begin{lema} \label{(DE)L:seccion}
Let $(\Gamma,\tau)$ be a proper directed hereditary translation quiver.
\begin{itemize}
 \item[a)] Every section in $(\Gamma,\tau)$, which does not contain injective vertices, is a cosection.
 \item[b)] Every cosection in $(\Gamma,\tau)$, which does not contain projective vertices, is a section.
\end{itemize}
\end{lema}

\bproof
We show $(b)$, the proof of $(a)$ can be obtained with dual arguments. By the last lemma $(\Gamma,\tau)$ has no 
periodic orbits. Let $\mathcal{S}$ be a cosection in $(\Gamma,\tau)$ which does not contain projective vertices and let
$x$ be an element of $\mathcal{S}$. Since $x$ is not projective then $\tau (x)$ is not injective and
$\tau(x) \notin \mathcal{S}$ (for since $\mathcal{S}$ is a cosection, if $\tau(x) \in \mathcal{S}$ then 
$x=\tau^{-1}(\tau(x)) \notin \mathcal{S}$, a contradiction). Then $\mathcal{S}$ satisfies condition $(i)$ in the  
definition of section. On the other hand, let $x \to m$ be an arrow in $\Gamma$. The vertex $m$ cannot be projective,
for $(\Gamma,\tau)$ is hereditary and $x$ is not projective. Then there is an arrow
$\tau(m) \to x$, thus either $\tau(m)$ or $m=\tau^{-1}(\tau(m))$ is an element in $\mathcal{S}$. We conclude that $\mathcal{S}$
satisfies $(ii)$ and is hence a section.
\eproof

\begin{lema} \label{(DE)L:partSeccion}
Let $(\Gamma,\tau)$ be a proper directed hereditary translation quiver and let $\mathcal{S}$ be a section 
(cosection) that intersects every orbit in $(\Gamma,\tau)$ in exactly one vertex. Then there is a partition
$\Gamma_0=A \sqcup \mathcal{S} \sqcup B$ of the vertex set in $\Gamma$, where 
\[
 A=\{ x \; | \; \tau^{-i}(x) \in \mathcal{S} \text{ for some }i\geq 1 \},
\]
\[
 B=\{ y \; | \; \tau^{i}(y) \in \mathcal{S} \text{ for some }i\geq 1 \}.
\]
Moreover, every path $\gamma$ in $\Gamma$ which starts in $A$ and ends in $B$ crosses $\mathcal{S}$. 
\end{lema}

\bproof
We show the case of a section, the proof for a cosection is similar. The partition shown is clear, for
$(\Gamma,\tau)$ has no periodic orbits (lemma~\ref{(DE)L:dirigido}). 
To prove the second claim assume that $\gamma$ is a path that starts in $A$ and ends in $B$, and such that
non of the arrows that conform $\gamma$ starts in an element of $\mathcal{S}$.
Then one of the arrows in $\gamma$ has the form $\alpha:x \to y$ with $x \in A$ and $y \notin A$.
We show that $y \in \mathcal{S}$.

Assume to the contrary that $y \in B$, that is, that there exists $j \geq 1$ such that $\tau^j(y) \in \mathcal{S}$.
Notice first that $\tau^{j-2}(x)$ cannot be projective, for $\Gamma$ is hereditary and $\tau^{j-1}(y)$ is not
projective.
\[
 \xymatrix@!0{
\tau^j(y) \ar[rd] &  & \tau^{j-1}(y) \ar[rd] & & \cdots & & y \\
& \tau^{j-1}(x) \ar[ru] & & \tau^{j-2}(x) & \cdots & x \ar[ru]
}
\]
Then $\tau^{j-1}(x)$ is defined and we have an arrow $\tau^j(y) \to \tau^{j-1}(x)$. Since $\tau^j(y) \in \mathcal{S}$
and $\mathcal{S}$ is a section, either $\tau^{j-1}(x) \in \mathcal{S}$ or $\tau(\tau^{j-1}(x))=\tau^j(x) \in \mathcal{S}$.
This contradicts the assumption that $\mathcal{S}$ intersects each orbit in exactly one vertex (for $\tau^{-i}(x) \in \mathcal{S}$).
Hence $y \in \mathcal{S}$, which completes the proof.
\eproof

\begin{lema} \label{(DE)L:PosPre}
Let $(\Gamma,\tau)$ be a connected and hereditary translation quiver.
\begin{itemize}
 \item[a)] If $(\Gamma,\tau)$ contains a projective vertex then all orbits in $(\Gamma,\tau)$ contains a projective vertex.
 \item[b)] If $(\Gamma,\tau)$ contains an injective vertex then all orbits in $(\Gamma,\tau)$ contains an injective vertex.
\end{itemize}
\end{lema}
\bproof
We show $(a)$, the proof of $(b)$ is similar.
Let $B$ be the subset of $\Gamma_0$ given by vertices $x$ such that $\tau^i(x)$ is projective for some $i\geq 0$.
By hypothesis $B$ is nonempty. Assume that $B$ is not the total vertex set in $\Gamma$. By connectedness there exist
$y \notin B$ and $x \in B$ such that $y$ and $x$ are neighbors in $\Gamma$. Since $y \notin B$, we can actually assume
there is an arrow $y \to x$. Let $i\geq 0$ such that $\tau^{i}(x)$ is projective. Since $\tau^j(y)$ is defined for all 
$j\geq 0$, we have an arrow $\tau^i(y) \to \tau^i(x)$. Because $(\Gamma,\tau)$ is hereditary, $\tau^i(y)$ is projective, 
which contradicts the fact that $y \notin B$. Hence $B=\Gamma_0$ which completes the proof. 
\eproof

Let $\Gamma=(\Gamma_0,\Gamma_1,\tau)$ be a translation quiver. The \textbf{translation subquiver}
$\widehat{\Gamma}=(\widehat{\Gamma}_0,\widehat{\Gamma}_1,\widehat{\tau})$ determined by a subset of vertices
$\widehat{\Gamma}_0$ consists of the full subquiver $\widehat{\Gamma}=(\widehat{\Gamma}_0,\widehat{\Gamma}_1)$
defined over $\widehat{\Gamma}_0$ and the restriction $\widehat{\tau}$ of $\tau$ to the subset of nonprojective vertices
$\widehat{\Gamma}_0'$ given by
\[
\widehat{\Gamma}'_0=\{ x \in \widehat{\Gamma}_0 \cap \Gamma'_0 \; | \; \tau(x) \in \widehat{\Gamma}_0 \}.
\]
The noninjective vertices of $\widehat{\Gamma}$ are given by
\[
\widehat{\Gamma}_0''=\{ y \in \widehat{\Gamma}_0 \cap \Gamma_0'' \; | \; \tau^{-1}(y) \in \widehat{\Gamma}_0 \}.
\]

If $x$ is a nonprojective vertex in a translation quiver $\Gamma=(\Gamma_0,\Gamma_1,\tau)$ then 
the full subquiver determined by the set $\{x,\tau(x)\} \cup x^-$ is called \textbf{mesh} of $(\Gamma,\tau)$. 
A full translation subquiver $\widehat{\Gamma}$ of $\Gamma$ is \textbf{mesh clomplete} if for any nonprojective vertex
$x$ in $\widehat{\Gamma}$ the set $x_{\widehat{\Gamma}}^-$ (direct predecessors of $x$ in the quiver 
$\widehat{\Gamma}$) coincides with the set $x_{\Gamma}^-$ (direct predecessors of $x$ respect to $\Gamma$). 
\begin{lema} \label{(DE)L:mallas}
Let $(\Gamma,\tau)$ be a translation quiver and $(\widehat{\Gamma},\widehat{\tau})$ a mesh complete full translation subquiver
of $(\Gamma,\tau)$.
\begin{itemize}
 \item[a)] If $(\Gamma,\tau)$ is proper then $(\widehat{\Gamma},\widehat{\tau})$ is proper.
 \item[b)] If $(\Gamma,\tau)$ is hereditary then $(\widehat{\Gamma},\widehat{\tau})$ is hereditary.
 \item[c)] If $(\Gamma,\tau)$ has an admissible ordering of its vertices then 
              $(\widehat{\Gamma},\widehat{\tau})$ admits an admissible ordering of its vertices.
\end{itemize}
\end{lema}
\bproof
The point $(a)$ is clear for if $x$ is a nonprojective vertex in $\widehat{\Gamma}$ then
$x$ is nonprojective in $\Gamma$ and $x_{\widehat{\Gamma}}^-=x_{\Gamma}^- \neq \emptyset$. 
To prove $(b)$ let $p$ and $q$ be vertices in $\widehat{\Gamma}_0$ with an arrow $q \to p$.
Assume that $q$ is nonprojective in $(\widehat{\Gamma},\widehat{\tau})$. In particular $q$ in nonprojective
in $(\Gamma,\tau)$. Since $(\Gamma,\tau)$ is a hereditary quiver, $p$ is not projective
in $(\Gamma,\tau)$. Hence $\tau(p) \in q_{\Gamma}^-=q_{\widehat{\Gamma}}^-$, and thus $p$ is not projective in
$(\widehat{\Gamma},\widehat{\tau})$. Finally, claim $(c)$ is evident for any subquiver of $\Gamma$.
\eproof

We say that a vertex $x$ in $\Gamma$ \textbf{precedes} the vertex $y$, which we denote by $x \preceq y$, 
if there exists a path $\gamma:x \to y$ in $\Gamma$. Observe that if $\Gamma$ is a quiver with no oriented cycles, 
then the precedence relation $\preceq$ is a partial order. Indeed, $\preceq$ is reflexive because of the
existence of trivial paths, transitive through concatenation of paths and anti-symmetric for in $\Gamma$
there are no oriented cycles.

Given a finite solid quiver $Q$ we define the translation quiver $\mathbb{Z}Q$ in the following way. The vertex set
consists of pairs $(x,\ell)$ with $x \in Q_0$ and $\ell \in \mathbb{Z}$. To each arrow $a:x \to y$ in $Q$ corresponds
two series of arrows in $\mathbb{Z}Q$,
\[
 (a,\ell):(x,\ell) \to (y,\ell) \qquad \text{and} \qquad \sigma(a,\ell):(y,\ell) \to (x,\ell+1).
\]
The bijection $\tau:(x,\ell) \mapsto (x,\ell-1)$ is clearly a translation. Each interval $I \subset \mathbb{Z}$ defines
a full translation subquiver $IQ$ of $\mathbb{Z}Q$ through the vertices of the form $(x,\ell)$ with $\ell$ in $I$. In particular,
$\mathbb{N}_0Q$ and $(-1)\mathbb{N}_0Q$ are posprojective and preinjective translation subquivers respectively.
All translation quivers which we are interested in are mesh complete translation subquiver of some $\mathbb{Z}Q$.
Last lemma~\ref{(DE)L:mallas}, together with the following lemma, give an account of some of their properties.

\begin{lema} \label{(DE)L:orden}
Let $Q$ be a finite solid connected quiver, with more than one vertex and with an admissible ordering of its vertices.
Then $\mathbb{Z}Q$ is a connected proper hereditary translation quiver which admits an admissible ordering of its vertices. 
In particular, $\mathbb{Z}Q$ is a directed quiver.  
\end{lema}
\bproof
That $\mathbb{Z}Q$ is connected, proper and hereditary is clear.
Fix an admissible ordering $x_1,\ldots,x_n$ in the vertex set of $Q$ and define the following relation in
$Q_0 \times \mathbb{Z}$: two elements are related $(x_i,\ell) \leq (x_j,m)$ 
if they satisfy one of the following conditions,
\begin{itemize}
 \item[$i)$] $\ell > m$ in $\mathbb{Z}$, or
 \item[$ii)$] $\ell=m$ and $i \leq j$ in $Q$. 
\end{itemize}
Then $\leq$ is a linear order in the vertices of $\mathbb{Z}Q$. Indeed, it is transitive since if 
$(x_{i_1},\ell_1) \leq (x_{i_2},\ell_2) \leq (x_{i_3},\ell_3)$ and $\ell_1 > \ell_2$ or $\ell_2 > \ell_3$
then $\ell_1 > \ell_3$, and hence $(x_{i_1},\ell_1) \leq (x_{i_3},\ell_3)$. In case $\ell_1=\ell_2=\ell_3$
we have $i_1 \leq i_2 \leq i_3$, and again $(x_{i_1},\ell_1)\leq (x_{i_3},\ell_3)$.
It is anti-symmetric, for if $(x_{i},\ell)\leq (x_{j},m) \leq (x_{i},\ell)$ then $\ell=m$ and $i\leq j \leq i$.
Moreover, it is clear that two arbitrary vertices $(x_i,\ell)$ and $(x_j,m)$ are related, thus $\leq $ is a linear order.

With the definition of arrows in $\mathbb{Z}Q$ given above it is clear that the order $\leq$ defined in the vertex set
$Q_0 \times \mathbb{Z}$ of $\mathbb{Z}Q$ is admissible.
\eproof

By lemmas~\ref{(DE)L:seccion} and~\ref{(DE)L:orden}, in the context of the following lemma we can replace
the concept of section by that of cosection.

\begin{lema} \label{(DE)L:orbital}
Let $Q$ be a finite connected solid quiver with an admissible ordering of its vertices. 
For a subset $\mathcal{S}$ of vertices in $\mathbb{Z}Q$ the following conditions are equivalent
\begin{itemize}
 \item[a)] $\mathcal{S}$ is a connected section in $\mathbb{Z}Q$,
 \item[b)] $\mathcal{S}$ intersects each orbit in $\mathbb{Z}Q$ in exactly one vertex, 
and two elements $(x,\ell)$, $(y,m)$ in $\mathcal{S}$ are neighbors in $\mathbb{Z}Q$ if and only if 
$x$, $y$ are neighbors in $Q$.
\end{itemize}
In particular, the underlying graph of a connected section $\mathcal{S}$ of
$\mathbb{Z}Q$ is isomorphic to the underlying graph of $Q$.
\end{lema}
\bproof 
We can assume that $Q$ has more that one vertex, so that $\mathbb{Z}Q$ is a proper translation quiver. \\
\underline{Step 1.} \textit{We show that $(b)$ implies $(a)$}. The first condition $(i)$ in the definition of section
is satisfied by hypothesis. To show condition $(ii)$ let $(x_i,\ell_i)\to (y,m)$ be
an arrow in $\mathbb{Z}Q$ with $(x_i,\ell_i)$ in $\mathcal{S}$. By construction of the arrows in 
$\mathbb{Z}Q$ the vertices $y$ and $x_i$ are neighbors in $Q$. Observe that $(x_i,\ell_i)$ have only two neighbors
in the orbit of $(y,m)$,
\[
\xy
(-24,0) *{\text{when $x_i \to y$ in $Q$,}};
(24, 0) *{\text{when $y \to x_i$ in $Q$,}};
(-15,-3) *{\xymatrix@!0{& (x_i,\ell_i) \ar[rd] \\ (y,\ell_i-1) \ar[ru] & & (y,\ell_i)}};
(  7,-3) *{\xymatrix@!0{& (x_i,\ell_i) \ar[rd] \\ (y,\ell_i) \ar[ru] & & (y,\ell_i+1)}};
\endxy
\vspace{1.3cm}
\]
In the first case we have that $m=\ell_i$ and hence either $(y,\ell_i)=(y,m) \in \mathcal{S}$ or 
$(y,\ell_i-1)=(y,m-1)=\tau(y,m) \in \mathcal{S}$, for by hypothesis an element in the orbit of $(y,m)$ 
which is neighbor of $(x_i,\ell_i)$ belongs to $\mathcal{S}$. In the second case we have that $m=\ell_i+1$ and 
for the same reason either $(y,\ell_i+1)=(y,m) \in \mathcal{S}$ or $(y,\ell_i)=(y,m-1)=\tau(y,m) \in \mathcal{S}$.
Hence $\mathcal{S}$ is a section. Its connectedness is direct consequence of the hypothesis in $(b)$.\\
\underline{Step 2.} \textit{We show that if $\mathcal{S}$ is a connected section which contains no injective vertex and
$\tilde{\mathcal{S}}$ is a subset of $\mathcal{S}$ which is also a section,
then $\tilde{\mathcal{S}}=\mathcal{S}$}. Assume that $\tilde{\mathcal{S}}\neq \mathcal{S}$. Then by connectedness
there exist neighbor vertices $x \in \mathcal{S}-\tilde{\mathcal{S}}$ and $y \in \tilde{\mathcal{S}}$.
In case $y \to x$, since $\tilde{\mathcal{S}}$ is section and $x \notin \tilde{\mathcal{S}}$ we have
$\tau(x) \in \tilde{\mathcal{S}} \subseteq \mathcal{S}$, which is impossible for $x \in \mathcal{S}$ and $\mathcal{S}$ is 
section. In case $x \to y$ there is an arrow $y \to \tau^{-1}(x)$, and hence 
$\tau^{-1}(x)\in \tilde{\mathcal{S}} \subseteq \mathcal{S}$ again because $\tilde{\mathcal{S}}$ is section and 
$x \notin \tilde{\mathcal{S}}$. This is also impossible for in that case $x$ and $\tau^{-1}(x)$ both belong to
$\mathcal{S}$. Thus $\tilde{\mathcal{S}}=\mathcal{S}$.\\
\underline{Step 3.} \textit{We prove that $(a)$ implies $(b)$}.
We will show that if $\mathcal{S}$ is a connected section in $\mathbb{Z}Q$ then 
$\mathcal{S}$ has a subset $\mathcal{S}_n$ which satisfies the hypothesis in $(b)$.
To give $\mathcal{S}_n$ we successively construct sets $\mathcal{S}_i$ for $1 \leq i \leq n$.
Let $u_1$ be an arbitrary element in $\mathcal{S}$ and make $\mathcal{S}_1=\{u_1\}$. Assume that the subset 
$\mathcal{S}_i$ has been constructed ($i=|\mathcal{S}_i|$) intersecting each orbit in at most one vertex,
and such that two if its vertices $(x_{i_1},\ell_{1})$ and $(x_{i_2},\ell_2)$ are neighbors if and only if
$x_{i_1}$ and $x_{i_2}$ are neighbors in $Q$. If $\mathcal{S}_i$ does not intersect all orbits in $\mathbb{Z}Q$
then there exist a vertex $v$ in $\mathbb{Z}Q$ such that $\tau^{\ell}(v)\notin \mathcal{S}_i$ 
for all $\ell \in \mathbb{Z}$ and an arrow $u_s \to v$ with $u_s \in \mathcal{S}_i \subseteq \mathcal{S}$.
Since $\mathcal{S}$ is a section, either $v \in \mathcal{S}$ or $\tau(v) \in \mathcal{S}$.
Take $u_{i+1}=u$ or $u_{i+1}=\tau(v)$ depending on the case and make $\mathcal{S}_{i+1}=\mathcal{S}_i \cup \{u_{i+1}\}$. 
Observe that by construction of $\mathbb{Z}Q$, since $u_s=(x,\ell)$ and $u_{i+1}=(y,m)$ are neighbors in $\mathbb{Z}Q$
then $x$ and $y$ are neighbors in $Q$, and hence $\mathcal{S}_{i+1}$ has the same properties as $\mathcal{S}_i$.
Thus, if $\mathbb{Z}Q$ has $n$ orbits ($n=|Q_0|$), we can successively construct a subset $\mathcal{S}_n$ which satisfies
the hypothesis of $(b)$. Then by step 1, $\mathcal{S}_n$ is a connected section and by step 2, $\mathcal{S}=\mathcal{S}_n$,
that is, $\mathcal{S}$ satisfies $(b)$.
\eproof

Assume that $(\Gamma,\tau)$ is a translation quiver and that $f:\Gamma_0 \to \mathbb{Z}$ is a function. We say that
$f$ is an \textbf{additive function} if for any mesh $\{ x,\tau(x) \}\cup x^-$ in $\Gamma$ it satisfies
\[
 f(x)+f(\tau(x))=\sum_{y \in x^-}f(y)(-m_{y,x}),
\]
where $(-m_{y,x})$ is the number of arrows from $y$ to $x$ in $\Gamma$. 
For a finite section (cosection) $\mathcal{S}=\{x_i\}_{i=1}^n$ 
in $(\Gamma,\tau)$ denote by $f(S)$ the integral vector $(f(x_i))_{i=1}^n$. 

\begin{lema} \label{(DE)L:aditCox}
Let $(\Gamma,\tau)$ be a proper directed hereditary translation quiver.
\begin{itemize}
 \item[a)] Assume that $\mathcal{S}$ is a section in $(\Gamma,\tau)$ such that the full subquiver 
$Q=Q_{\mathcal{S}}$ in $\Gamma$ determined by $\mathcal{S}$
is finite and admits an admissible ordering of its vertices. Assume further that $\mathcal{S}$
contains no injective vertices and that it intersects each orbit in $(\Gamma,\tau)$ in at most one vertex. 
Then $\tau^{-1} \mathcal{S}=\{ \tau^{-1}(x) \}_{x \in \mathcal{S}}$ is also a section in $(\Gamma,\tau)$, 
$Q_{\tau^{-1}\mathcal{S}} \cong Q_{\mathcal{S}}$ and for any additive function $f$ in $(\Gamma,\tau)$ we have that
\[
 f(\tau^{-1} \mathcal{S})=\Phi_Q^{-1} f(\mathcal{S}),
\]
where $\Phi_{Q}^{-1}$ is the inverse of the Coxeter matrix of $Q$ respect to the admissible ordering of its vertices.
 \item[b)] Assume that $\mathcal{S}$ is a cosection in $(\Gamma,\tau)$ such that the full subquiver 
$Q=Q_{\mathcal{S}}$ in $\Gamma$ determined by $\mathcal{S}$
is finite and admits an admissible ordering of its vertices. Assume further that $\mathcal{S}$
contains no projective vertices and that it intersects each orbit in $(\Gamma,\tau)$ in at most one vertex. 
Then $\tau \mathcal{S}=\{ \tau(x) \}_{x \in \mathcal{S}}$ is also a cosection in $(\Gamma,\tau)$, 
$Q_{\tau\mathcal{S}} \cong Q_{\mathcal{S}}$ and for any additive function $f$ in $(\Gamma,\tau)$ we have that
\[
 f(\tau \mathcal{S})=\Phi_Q f(\mathcal{S}),
\]
where $\Phi_{Q}$ is the Coxeter matrix of $Q$ respect to the admissible ordering of its vertices.
\end{itemize}
\end{lema}

\bproof We will show $(b)$, the proof of $(a)$ is similar.
Fix an additive function $f$. We prove first the following claim.\\
\underline{Step 1.} \textit{Assume that $x$ is a nonprojective element of a cosection
$\mathcal{S}_0$ and that $x$ is a sink in the full subquiver $Q_{\mathcal{S}_0}$ determined by $\mathcal{S}_0$
(that is, if there is an arrow $x \to m$ in $\Gamma$ then $m \notin \mathcal{S}_0$).
Denote by $\sigma_x\mathcal{S}_0$ the subset of $\Gamma_0$ obtained from $\mathcal{S}_0$ by replacing
$x$ with $\tau(x)$,
\[
 \sigma_x\mathcal{S}_0=\{ \tau(x) \} \cup (\mathcal{S}_0-\{x\}),
\]
and maintaining the same order in vertices. Then $\sigma_x\mathcal{S}_0$ is a cosection.}
On the one hand, if $\mathcal{S}_0$ has at most one vertex in each orbit, then $\sigma_x\mathcal{S}_0$ satisfies the
same condition. Hence for each noninjective element $y$ in $\sigma_x\mathcal{S}_0$ we have $\tau^{-1}(y) \notin 
\sigma_x\mathcal{S}_0$. On the other hand, assume that $n \to y$ is an arrow in $\Gamma$ and that $y \in \sigma_x\mathcal{S}_0$.
Consider first the case $y \neq \tau(x)$. Then $y \in \mathcal{S}_0$, and since $\mathcal{S}_0$ is a cosection,
either $n \in \mathcal{S}_0$ or $\tau^{-1}(n) \in \mathcal{S}_0$. If $n \in \mathcal{S}_0$ then $n \neq x$ for
$x$ is a sink in $Q_{\mathcal{S}_0}$, and hence $n \in \sigma_x\mathcal{S}_0$. Assume now that $n \notin \mathcal{S}_0$,
thus $\tau^{-1}(n) \in \mathcal{S}_0$. If $\tau^{-1}(n) \neq x$ then $\tau^{-1}(n) \in \sigma_x\mathcal{S}_0$.
And if $\tau^{-1}(n)=x$ then $n=\tau(\tau^{-1}(n))=\tau(x) \in \sigma_x\mathcal{S}_0$. In any case one of the vertices
$n$ or $\tau^{-1}(n)$ belongs to $\sigma_x\mathcal{S}_0$. Consider finally the case $y=\tau(x)$,
that is, when there is an arrow $n \to \tau(x)$ in $\Gamma$. 
Notice first that $n$ is not injective, for $(\Gamma,\tau)$ is hereditary
and $\tau(x)$ is not injective. Then there is an arrow $\tau^{-1}(n) \to x$. Since $\mathcal{S}_0$ is a cosection and 
$\tau^{-2}(n)$ cannot be in $\mathcal{S}_0$ (again because $x$ is a sink in $Q_{\mathcal{S}_0}$) 
we conclude that $\tau^{-1}(n)$ belongs to $\mathcal{S}_0$. 
Since $\Gamma$ is directed, $\tau^{-1}(n)$ is different from $x$ and hence $\tau^{-1}(n)$ is an element in
$\sigma_x\mathcal{S}_0$. This shows that $\sigma_x\mathcal{S}_0$ is a cosection. Observe
that the quiver $Q_{\sigma_x\mathcal{S}_0}$ is obtained from $Q_{\mathcal{S}_0}$ by changing the orientation of all
arrows which start or end in $x$.\\
\underline{Step 2.} \textit{If $x$ is a nonprojective element of a cosection $\mathcal{S}_0$ such that $x$ is a sink
in the full subquiver $Q_{\mathcal{S}_0}$, and $Q_{\mathcal{S}_0}$ admits an admissible ordering of its vertices, then
\[
f(\sigma_x\mathcal{S}_0)=\sigma_x(f(\mathcal{S}_0)),
\]
where the $\sigma_x$ in the right hand of the equation is the simple reflection of vertex $x$ associated to 
$Q_{\mathcal{S}_0}$.} Since the change $\mathcal{S}_0 \mapsto \sigma_x\mathcal{S}_0$ only modifies the vertex $x$, 
in the same way as the reflection $\sigma_x$ only changes the entry $x$, we can focus in this component. 
Since $(\mathbf{e}_x,\mathbf{e}_x)=1$ because $Q_{\mathcal{S}_0}$ has no oriented cycles, notice that 
\begin{eqnarray}
 \sigma_x(f(\mathcal{S}_0))_x & = & f(\mathcal{S}_0)_x -2(f(\mathcal{S}_0),\mathbf{e}_x) = \nonumber \\
& = & f(x)-\sum_{y \in \mathcal{S}_0}f(y)2(\mathbf{e}_y,\mathbf{e}_x)= \nonumber \\
& = & f(x)-2f(x)-\sum_{y \in \mathcal{S}_0-\{x\}}f(y)(\langle \mathbf{e}_y,\mathbf{e}_x 
\rangle+\langle \mathbf{e}_x,\mathbf{e}_y \rangle). \nonumber 
\end{eqnarray}
Now, since $x$ is a sink in $Q_{\mathcal{S}_0}$ we have that $\langle \mathbf{e}_x,\mathbf{e}_y \rangle=m_{x,y}=0$ for all
$y \in \mathcal{S}_0-\{x\}$ and since $\mathcal{S}_0$ is a cosection we have that if $y \in \mathcal{S}_0-\{x\}$ then
$\langle \mathbf{e}_y,\mathbf{e}_x \rangle=m_{y,x} \neq 0$ if and only if $y \in x^-$. Hence
\begin{eqnarray}
 \sigma_x(f(\mathcal{S}_0))_x & = & -f(x)-\sum_{y \in \mathcal{S}_0-\{x\}}f(y)\langle \mathbf{e}_y,\mathbf{e}_x \rangle= \nonumber \\
& = & -f(x)+\sum_{y \in x^-}f(y)(-m_{y,x})= \nonumber \\
& = & f(\tau(x)). \nonumber
\end{eqnarray}
In this way we have that $f(\sigma_x\mathcal{S}_0)=\sigma_xf(\mathcal{S}_0)$.\\
\underline{Step 3.} \textit{We prove $(b)$}.
Let $\mathcal{S}$ be a cosection as in the statement of the lemma and fix an admissible ordering $\{x_1,\ldots,x_n\}$ of
the vertices in $Q_{\mathcal{S}}$. Then $x_1$ is a sink and by the steps above $\mathcal{S}_1=\sigma_{x_1}\mathcal{S}$ 
is a cosection. The full subquiver of $\Gamma$ which $\mathcal{S}_1$ determines admits an admissible ordering 
(making $\tau(x_1)$ greater than the rest of the elements in $\mathcal{S}$) and has the vertex $x_2$ as nonprojective sink.
We can then construct $\mathcal{S}_2=\sigma_{x_2}\mathcal{S}_1$. Assume we have constructed in this way the cosection
$\mathcal{S}_i$ where $x_{i+1}$ is a nonprojective sink. Then $\mathcal{S}_{i+1}=\sigma_{x_{i+1}}\mathcal{S}_i$ is a cosection. 
Successively we obtain a cosection $\mathcal{S}_n=\sigma_{x_n}\mathcal{S}_{n-1}$, and by step 2 we have
\[
f(\mathcal{S}_n)=\sigma_{x_n}f(\mathcal{S}_{n-1})=\ldots=\sigma_{x_n}\sigma_{x_{n-1}}\cdots \sigma_{x_2} \sigma_{x_1}(f(\mathcal{S}))
=\Phi_Qf(\mathcal{S}).
\]
Since all simple reflections $\sigma_x$ have been used exactly once, it is clear that $\mathcal{S}_n=\tau \mathcal{S}$
and that $Q_{\tau \mathcal{S}} \cong Q_{\mathcal{S}}$. This completes the proof.
\eproof

\section{Posprojective and preinjective components.} \label{(DE)S:comp}
%------------------------------------------------------------------
%------------------------------------------------------------------
We end this chapter studying those components of the Auslander-Reiten quiver $\Gamma(kQ)$ which contain
the projective (injective) $kQ$-modules, where $Q$ is a finite solid quiver with an admissible ordering
of its vertices which is not a Dynkin diagram. The following result is well known, 
see for instance~\cite[theorem V.7.8]{ARS95}.

\begin{lema} \label{(P)L:KSfinita}
Let $A$ be a finite dimensional $k$-algebra of finite representation type.
If $M$ and $N$ are nonisomorphic indecomposable $A$-modules then every nonzero morphism $f:M \to N$
is sum of compositions of irreducible morphisms between indecomposable $A$-modules.
In particular, there exists a path from $[M]$ to $[N]$ in the Auslander-Reiten quiver $\Gamma(A)$ of $A$-mod.
\end{lema}

\begin{lema} \label{(DE)L:radMax}
Let $A=kQ$ be the path algebra of a finite solid quiver with an admissible ordering of its vertices
and denote by $P(i)$ the indecomposable projective $A$-module $Ae_i$ (with $e_i$ the trivial path 
over the vertex $i\in Q_0$). For each path $\gamma:i \to j$ in $Q$ the mapping $g_{\gamma}:P(j) \to P(i)$ given by 
$ae_j \mapsto ae_j\gamma e_i$ is a monomorphism of $A$-modules.
\begin{itemize}
 \item[i)] The set $\{ g_{\gamma} \; | \; \gamma \text{ is path from $i$ to $j$} \}$ is basis of $\Hom_{A}(P(j),P(i))$.
\item[ii)] Take $P'(i)=\bigoplus_{\substack{\alpha \in Q_1 \\ s(\alpha)=i}}P(t(\alpha))$ and consider the morphism
$g:P'(i) \to P(i)$ given by $g=[g_{\alpha}]_{\substack{\alpha \in Q_1 \\ s(\alpha)=i}}$ where $g_{\alpha}:P(t(\alpha)) \to P(i)$.
If $j\neq i$ then every nonzero morphism $f:P(j) \to P(i)$ factors through $g$.
 \item[iii)] $\Img g$ is the unique maximal submodule of $P(i)$. 
\end{itemize}
\end{lema}

\bproof
Recall that if $M$ is an $A$-module and $e$ is an idempotent in $A$, then there is an isomorphism of $k$-vector spaces
\[
 \Phi:\Hom_A(Ae,M) \to eM,
\]
defined by $f \mapsto f(e)$. Its inverse is given by $\Psi(m):a \mapsto am$. Then 
\[
\Hom_A(P(j),P(i))=\Hom_A(Ae_j,Ae_i) \cong e_jAe_i,
\]
and each morphism $g_{\gamma}$ defined above corresponds to $\Psi(\gamma)$. If we order the paths from $i$ to $j$ in $Q$,
$\gamma_1,\ldots,\gamma_r$, then the morphisms $g_{\gamma_1},\ldots,g_{\gamma_r}$ form a basis of $\Hom_A(P(j),P(i))$.
This proves $(i)$. Moreover, if $\gamma\tilde{\gamma}$ is a path in $Q$, then 
$g_{\tilde{\gamma}}(g_{\gamma}(a))=g_{\tilde{\gamma}}(a\gamma)=a\gamma\tilde{\gamma}=g_{\gamma\tilde{\gamma}}(a)$.
Assume that $g:P'(i)\to P(i)$ is as in the point $(ii)$ and that $f:P(j) \to P(i)$ is a nonzero morphism with 
$j\neq i$. By $(i)$ there is a nonzero linear combination $\sum_{\ell=1}^rc_{\ell}\gamma_{\ell}$ such that
$f=\Psi(\sum_{\ell=1}^rc_{\ell}\gamma_{\ell})=\sum_{\ell=1}^rc_{\ell}g_{\gamma_{\ell}}$. 
For each arrow $\alpha$ with source $i$ take 
\[I_{\alpha}=\{ \ell \in \{1,\ldots, r\} \; | \; 
\gamma_{\ell}=\widetilde{\gamma}_{\ell}\alpha \text{ for some path } \widetilde{\gamma}_{\ell} \}, \]
and define $h_{\alpha}:P(j) \to P(t(\alpha))$ as $h_{\alpha}=\sum_{\ell \in I_{\alpha}}c_{\ell}g_{\tilde{\gamma}_{\ell}}$.
Take the morphism $h=[h_{\alpha}]^t_{\substack{\alpha \in Q_1 \\ s(\alpha)=i}}:P(j) \to P'(i)$ and notice that
\[
 gh= \sum_{\substack{\alpha \in Q_1 \\ s(\alpha)=i}}g_{\alpha}h_{\alpha}=\sum_{\substack{\alpha \in Q_1 \\ s(\alpha)=i}}
\sum_{\ell \in I_{\alpha}} c_{\ell}g_{\alpha}g_{\tilde{\gamma}_{\ell}}=
\sum_{\substack{\alpha \in Q_1 \\ s(\alpha)=i}} \sum_{\ell \in I_{\alpha}} c_{\ell}g_{\tilde{\gamma}_{\ell} \alpha} =
\sum_{\ell=1}^r c_{\ell}g_{\gamma_{\ell}}=f,
\]
since the set $\{I_{\alpha} \; | \; \alpha \in Q_1, s(\alpha)=i \}$ is a partition of the interval $\{1,\ldots,r\}$.
This shows $(ii)$. Assume finally that $M$ is a proper submodule of $P(i)$. Since the algebra $A$ is hereditary,
every indecomposable direct summand of $M$ is isomorphic to $P(j)$ for some $j \neq i$. By point $(ii)$ the composition
$f:P(j) \to M \to P(i)$ factors through $g$. Since this is true for any direct summand of $M$, we conclude that $M$ 
is contained in $\Img g$. This completes the proof.
\eproof

\begin{lema} \label{(DE)L:rad}
Let $Q$ be a finite solid quiver with an admissible ordering of its vertices and let $A=kQ$ be the path $k$-algebra of $Q$.
\begin{itemize}
 \item[a)] Let $P$ be an indecomposable projective $A$-module. A morphism of $A$-modules $g:M \to P$ is
right minimal almost split if and only if $g$ is a monomorphism and its image is $\rad P$.
 \item[b)] Let $I$ be an indecomposable injective $A$-module. A morphism of $A$-modules $f:I \to M$ is
left minimal almost split if and only if $f$ is an epimorphism and its kernel is $\soc I$.
 \item[c)] The Auslander-Reiten quiver $\Gamma(A)$ of $A$-mod is a hereditary translation quiver.
\end{itemize}
\end{lema}
\bproof
We prove $(a)$, the proof of $(b)$ is similar (cf. proposition IV.3.5 in \cite{ASS06}).
By unicity it is enough to prove that $g: \rad P \to P$ is a right minimal almost split morphism.
If $h:\rad P \to \rad P$ is such that $gh=g$, then $g(h-Id_{\rad P})=0$ and since $g$ is injective, $h=Id_{\rad P}$. 
Hence $g$ is right minimal. Assume that $v:V \to P$ is not a retraction. Then $v$ is not an epimorphism
(for $P$ is projective), and thus $\Img v$ is a proper submodule of $P$. 
Considering point $(iii)$ in lemma~\ref{(DE)L:radMax}, the projective $P$ has a unique maximal submodule $\rad P$, 
hence $v$ factors through $g$. This shows that $g$ is a right almost split morphism.

We prove now $(c)$. Recall that $\Gamma(A)_0=\Ind A$. Moreover, for each right minimal almost split morphism
$g: \bigoplus d_iE_i \to N$ (with all $E_i$ nonisomorphic indecomposables) there are 
$d_i=\dimk_k\Irr(E_i,N)$ solid arrows from $[E_i]$ to $[N]$, and for each left minimal almost split morphism
$f:M \to \bigoplus d_iE_i$ there are $d_i=\dimk_k \Irr(M,E_i)$ solid arrows from $[M]$ to $[E_i]$.

By construction $\Gamma(A)$ is a solid quiver which is clearly locally finite. Let $\Gamma(A)_0'$ be the set of 
vertices for which there exists an almost split exact pair $(f,g)$ of the form
\[
\xymatrix{
0 \ar[r] & M \ar[r]^-{f} & E \ar[r]^-{g} & N \ar[r] & 0,
}
\]
and $\Gamma(A)_0''$ be the set of elements $[M]$ determined by such pairs. By definition, the Auslander-Reiten 
is a bijection between $\Gamma(A)_0'$ and $\Gamma(A)_0''$. By equivalences $(a)$, $(b')$ and $(c')$
in lemma~\ref{(DE)L:ARequiv} and the definition of arrows in $\Gamma(A)$, it follows that the Auslander-Reiten 
translation is indeed a quiver translation. Moreover, because of the existence theorem of almost split sequences
($A$ is an artinian $k$-algebra) the set $\Gamma(A)_0-\Gamma(A)_0'$ consists in the isomorphism classes of
indecomposable projective $A$-modules, and $\Gamma(A)_0-\Gamma(A)_0''$ consists in the isomorphism classes of
indecomposable injective $A$-modules. 

To verify that $\Gamma(A)$ is a hereditary translation quiver assume we have an arrow
$[E_i] \to [P]$ where $[P]$ is a projective vertex. Then there is an $A$-module $E$ from which
$E_i$ is direct summand and a right minimal almost split morphism $g:E \to P$. By the point 
$(a)$ the module $E$ is isomorphic to the radical $\rad P$. By lemma~\ref{(DE)L:radMax}$(ii)$ all direct summands of
$E$ are projective and in particular $[E_i]$ is a projective vertex. Using point $(b)$ and a dual version of
lemma~\ref{(DE)L:radMax}, one shows that if $[I] \to [E_i]$ is an arrow with $[I]$ an injective vertex, 
then $[E_i]$ is also injective.
\eproof

\begin{lema} \label{(DE)L:finito}
Let $Q$ be a finite solid connected quiver with an admissible ordering of its vertices and let $A=kQ$ be the path $k$-algebra
of $Q$. The following properties are equivalent.
\begin{itemize}
 \item[a)] $A$ is of finite representation type.
 \item[b)] There exist an indecomposable projective $A$-module $P$ and an integer $n\geq 0$ such that $\tau^{-n}P$ is injective. 
 \item[c)] There exist an indecomposable injective  $A$-module $I$ and an integer $n\geq 0$ such that $\tau^{n}I$ is projective.
\end{itemize}
\end{lema}
\bproof
The equivalence of $(b)$ and $(c)$ is evident. By lemma~\ref{(DE)L:rad}$(c)$ the Auslander-Reiten quiver $\Gamma(A)$ of
$A$-mod is a hereditary translation quiver. 
Observe that, because of lemmas~\ref{(DE)L:radMax} and \ref{(DE)L:rad}, the full subquiver
of $\Gamma(A)$ determined by the projective vertices is connected. In a similar way one notice that the full subquiver
of $\Gamma(A)$ determined by the injective vertices is also connected.
If $A$ is of finite representation type, since all indecomposable $A$-modules $M$ have projective covers $P \to M$, by 
lemma~\ref{(P)L:KSfinita} it follows that $\Gamma(A)$ is a connected quiver. By the desciption of orbits given at the beginning of
section~\ref{(DE)S:trasl}, since every orbit in $\Gamma(A)$ is finite, then every orbit is either periodic or contains a projective
and an injective vertex. This, lemma~\ref{(DE)L:PosPre} and the existence of projective and injective modules in $A$-mod show
that $(a)$ implies $(b)$ and $(c)$. 

We show now that $(b)$ implies $(a)$.  
Let $\mathcal{C}$ be the connected component in $\Gamma(A)$ containing all projective vertices.
By hypothesis in $(b)$, some orbit in $\mathcal{C}$ has an injective vertex, and by connectedness $\mathcal{C}$ contains all
injective vertices. Then by lemma~\ref{(DE)L:PosPre} all orbits in $\mathcal{C}$ are projective-injective, and hence 
$\mathcal{C}$ is a finite set. Observe finilly that for any indecomposable $A$-module $M$ there exists an element $[N]$ in 
$\mathcal{C}$ such that $\Hom_A(M,N) \neq 0$ (since $\mathcal{C}$ contains the isomorphism classes of all indecomposable injective
modules). From lemma 6 in section 2.2 in \cite{cmR} it follows that $[M]$ belongs to $\mathcal{C}$,
and hence $\Gamma(A)_0=\mathcal{C}$ is a finite set.
\eproof

As a reference for the following result consider proposition 10.2 in Gabriel and Roiter~\cite{GR97}.

\begin{lema} \label{(DE)L:componentes}
Let $Q$ be a finite solid connected quiver with an admissible ordering of its vertices and let $A=kQ$ be the path algebra
of $Q$. Assume that $kQ$-mod is of infinite representation type.
\begin{itemize}
 \item[a)] There is a connected component $\mathcal{P}$ of the Auslander-Reiten quiver of $kQ$-mod which contains the isomorphism 
classes of all indecomposable projective modules. This component is isomorphic to $\mathbb{N}_0Q^{op}$ and hence $\mathcal{P}$ 
is a directed posprojective quiver.
 \item[b)] There is a connected component $\mathcal{I}$ of the Auslander-Reiten quiver of $kQ$-mod which contains the isomorphism 
classes of all indecomposable injective modules. This component is isomorphic to $(-\mathbb{N}_0)Q^{op}$ and hence $\mathcal{I}$ 
is a directed preprojective quiver.
\end{itemize}
\end{lema}
\bproof
We show $(a)$. Let $P(i)=Ae_i$ be the projective $A$-module where $e_i$ is the trivial path (idempotent) of 
vertex $i\in \{1,\ldots,n\}$. By the definition of the Auslander-Reiten quiver of $A$-mod
and as consequence of lemmas~\ref{(DE)L:radMax}$(iii)$ and \ref{(DE)L:rad}(a), the full subquiver of $\Gamma(kQ)$ 
determined by the classes $[P(1)],\ldots,[P(n)]$ is isomorphic to $Q^{op}$. 

Let $\mathcal{P}$ be the connected component of $\Gamma(kQ)$ which contains all projective vertices.
Since $\Gamma(A)$ is hereditary (\ref{(DE)L:rad}$(c)$), by lemma~\ref{(DE)L:mallas}$(b)$ the subquiver $\mathcal{P}$
is also hereditary. By lemma~\ref{(DE)L:PosPre}, for any vertex $[M]$ in $\mathcal{P}$ there exists
$\ell \geq 0$ such that $\tau^{\ell}[M]$ is projective. Moreover, by lemma~\ref{(DE)L:finito}, for each projective 
$P(i)$ the translations $\tau^{-\ell}[P(i)]$ are defined for all $\ell\geq 0$. Then there is a bijection between the vertices
in $\mathcal{P}$ and those in $\mathbb{N}_0Q^{op}$ given by $\tau^{-\ell}[P(i)] \mapsto (i,\ell)$ for $i$ in $Q_0$ and $\ell \geq 0$.
Since the quiver determined by the indecomposable projective modules is isomorphic to $Q^{op}$  and by definition of traslation, 
the bijection above is a translation quiver isomorphism. In particular $\mathcal{P}$ is a posprojective directed quiver.
\eproof

The following lemma is a particular case of corollary~IV.2.9 in \cite{ASS06}.
\begin{lema} \label{(P)L:coxPP}
Let $Q$ be a finite solid connected quiver with an admissible ordering of its vertices and let $kQ$ be the path algebra
of $Q$. Assume that $kQ$-mod is of infinite representation type. Let $\Phi_Q$ be the Coxeter matrix associated to $Q$.
\begin{itemize}
 \item[a)] If $M$ is a preinjective $kQ$-module then $\vdim \tau M=\Phi_Q \vdim M$.
 \item[b)] If $N$ is a posprojective $kQ$-module then $\vdim \tau^{-1} N=\Phi_Q^{-1} \vdim N$.
\end{itemize}
\end{lema}
%----------------------------------------------------------------------
%----------------------------------------------------------------------
\chapter{Kronecker algebras and Dynkin diagrams.}
\label{Cap(D)}
%----------------------------------------------------------------------
For $n \geq 2$ consider the quiver $K_n$ with two vertices and $n$ arrows in the same direction and its associated matrix 
$M_{K_n}$,
\[
K_n=\xymatrix{{}_{\bullet_2} \ar@<-1ex>[d]_-{a_1} \ar@{}[d]|-{\cdots} \ar@<1ex>[d]^-{a_n} \\ {}_{\bullet_1} }
\qquad 
M_{K_n}=\left( \begin{matrix} 1 & 0 \\ -n & 1 \end{matrix} \right).
\]
The \textbf{$\mathbf{n}$-generalised Kronecker algebra} is the path algebra $A_n=kK_n$ of the quiver $K_n$.
The purpose of sections~\ref{(P)S:KroDos} and \ref{(P)S:KroTres} is to exhibit algorithms to construct 
exceptional modules for the classical Kronecker algebra $A_2$ and the first generalised case $A_3$ respectively.
All other cases can be treated in a similar way.  Although this representations are well known (cf. Ringel~\cite{cmR98} 
and \cite{cmR10}), the method used illustrates the notions presented in the preliminars chapter. Moreover,
the module category of the classical Kronecker algebra is fundamental for the construction of representations for
other extended Dynkin quivers. In this chapter we denote by $I$ the identity morphisms. 

%----------------------------------------------------------------------
\section{Classical Kronecker algebra $n=2$.} \label{(P)S:KroDos}
%----------------------------------------------------------------------
Consider first the classical case $A_2$. 
\[
K_2=\xymatrix{{}_{\bullet_2} \ar@<-.5ex>[d]_-{a_1} \ar@<.5ex>[d]^-{a_2} \\ {}_{\bullet_1} }
\qquad 
M_{K_2}=\left( \begin{matrix} 1 & 0 \\ -2 & 1 \end{matrix} \right).
\]
The corresponding quadratic form is given by $q_{K_2}(d_1,d_2)=(d_1-d_2)^2$. 
Hence the set of positive roots of $K_2$ consists in the vectors
\[
 p^t_2=(t+1,t) \qquad \text{and} \qquad q^t_2=(t,t+1),
\]
for $t$ greater or equal to zero. We will perform two reductions in the Kronecker algebra $A_2$,
each of them will produce exceptional representation of dimension vectors $p_2^t$ and $q_2^t$. 
We need the following technical lemma. Recall that a functor is call rigid if it is exact and induces isomorphisms
in extension groups.

\begin{lema} \label{(P)L:semiEje}
Let $B$ be the subalgebra of $A_2$ generated by the arrow $a_1$ and consider the admissible $B$-modules
$X=X^1 \oplus X^{\omega}$ and $Y=Y^{\omega} \oplus Y^2$, where $X^1,Y^2$ are the simple modules of vertex $1$, $2$ and $X^{\omega}$,
$Y^{\omega}$ are copies of the exceptional $B$-representation $\xymatrix{k \ar[r]^-{1} & k}$.
Fix bases $\{ x^1_1\}$ of $X^1$, $\{x_1^{\omega}=y_1^{\omega},x_2^{\omega}=y_2^{\omega} \}$ of $X^{\omega}=Y^{\omega}$ and
$\{ y^2_2\}$ of $Y^2$ such that $a_1x^{\omega}_2=x^{\omega}_1$ and $a_1y^{\omega}_2=y^{\omega}_1$. There are irreducible
morphisms $\xymatrix{X^1 \ar[r]^-{\sigma} & X^{\omega}=Y^{\omega} \ar[r]^-{\pi} & Y^2}$ given by 
\[
 \xymatrix{
0 \ar[d] \ar[r] & kx^{\omega}_2 \ar[d]^-{X^{\omega}_{a_1}}_-{1} & 
ky^{\omega}_2 \ar[d]_-{Y^{\omega}_{a_1}}^-{1} \ar[r]^-{\pi_2}_-{1} & ky^2_2 \ar[d] \\
kx^1_1 \ar[r]_-{\sigma_1}^-{1} & kx^{\omega}_1        & ky^{\omega}_1 \ar[r] & 0 \\ 
}
\]
The opposed endomorphism algebras of $X$ and $Y$ split over their radical,
\[
 \End_B(X)^{op}\cong S \oplus P \qquad \text{and} \qquad \End_B(Y)^{op} \cong S' \oplus P',
\]
where $P=k\sigma$, $P'=k\pi$, $S=\End(X^1) \oplus \End(X^{\omega}) \cong kf_1 \times kf_{\omega}$ and 
$S'=\End(Y^{\omega}) \oplus \End(Y^2) \cong kf'_{\omega} \times kf'_{2}$.
The matrices of dimension vectors of $X$ and $Y$ are given by
\[
 T^X=\left( \begin{matrix} 1&1\\0&1 \end{matrix} \right) \quad \text{and} \quad 
 T^Y=\left( \begin{matrix} 1&0\\1&1 \end{matrix} \right).
\]
In the dual basis $\{ x^1_1,x^{\omega}_1,x^{\omega}_2;(x^1_1)^*,(x^{\omega}_1)^*,(x^{\omega}_2)^* \}$ of $X$ the 
left coaction $\lambda:X^* \to P^* \otimes_S X^*$ and right coaction $\rho:X \to X \otimes_S P^*$ have the following form,
\begin{equation*}
\lambda(u) = \left\{
\begin{array}{l l}
\sigma^* \otimes (x^1_1)^*, & \text{if $u=(x_1^{\omega})^*$},\\
0, & \text{otherwise.}
\end{array} \right.
\end{equation*}
\begin{equation*}
\rho(z) = \left\{
\begin{array}{l l}
x_1^{\omega}\otimes \sigma^*, & \text{if $z=x^1_{1}$},\\
0, & \text{otherwise.}
\end{array} \right.
\end{equation*}
On the other hand, in the dual basis $\{ y^1_1,y^{\omega}_1,y^{\omega}_2;(y^1_1)^*,(y^{\omega}_1)^*,(y^{\omega}_2)^* \}$ 
of $Y$ the left coaction $\lambda':Y^* \to (P')^* \otimes_{S'} Y^*$ and right coaction $\rho':Y \to Y \otimes_{S'} (P')^*$ 
have the following form
\begin{equation*}
\lambda'(u) = \left\{
\begin{array}{l l}
\pi^*\otimes (y_2^{\omega})^*, & \text{if $u=(y^2_{2})^*$}, \\
0, & \text{otherwise.}
\end{array} \right.
\end{equation*}
\begin{equation*}
\rho'(z) = \left\{
\begin{array}{l l}
y^2_{2} \otimes \pi^*, & \text{if $z=y_2^{\omega}$}, \\
0, & \text{otherwise.}
\end{array} \right.
\end{equation*}
Moreover, the $B$-modules $X$ and $Y$ are rigid, hence the reduction functors
$F^X:A_2^X\text{-mod} \to A_2\text{-mod}$ and $F^Y:A_2^Y\text{-mod} \to A_2\text{-mod}$ are rigid functors.
\end{lema}
\bproof
For the expression of the left and right coactions confer formulas in lemma~\ref{(A)D:coacciones}.
To verify that the functor $F^X$ is rigid, recall the exact sequence on extension groups induced by the reduction
by the module $X$,
\[
 \xymatrix@C=.65pc{
0 \ar[r] & \Ext^1_{\mathcal{A}^X}(M,N) \ar[r] & \Ext^1_{\mathcal{A}}(F^X(M),F^X(N)) \ar[r] & \Ext^1_B(F^X(M),F^X(N)) \ar[r] & 0.
}
\]
Since the modules $F^X(M)$ and $F^X(N)$, considered as $B$-modules, have the form $X\otimes_S M$ and $X \otimes_{S}N$, 
and these are by hypothesis rigid modules, we conclude that $F^X$ is a rigid functor. The case $F^Y$ is similar.
All other claims in the statement are clear.
\eproof

Keeping in mind the notation in last lemma, consider the changes of basis
\[
(T^X)^tM_{K_2}T^X= \left( \begin{matrix} 1&1\\-1&0 \end{matrix} \right) \quad \text{and} \quad 
(T^Y)^tM_{K_2}T^Y= \left( \begin{matrix} 0&1\\-1&1 \end{matrix} \right).
\]
Their corresponding associated quivers, which will be denoted by $K^x_2$ and $K^y_2$, have the form
\[
\xymatrix{{}_{\bullet_2} \ar@(lu,ru)[]^(.2){a} \ar@<.5ex>@{<.}[d]^-{b} \ar@<-.5ex>[d]_-{a'} \\ {}_{\bullet_1} } \qquad \text{and} \qquad
\xymatrix{{}_{\bullet_2} \ar@<.5ex>@{<.}[d]^-{b} \ar@<-.5ex>[d]_-{a'} \\ {}_{\bullet_1} \ar@(ld,rd)[]_(.2){a} }
\]
\begin{lema} \label{(P)L:redUno}
The reduced ditalgebras $A^X_2$ and $A^Y_2$ of the Kronecker algebra $A_2$ are isomorphic to 
$(kK^x_2,\delta^x)$ and $(kK^y_2,\delta^y)$, where the differentials are given by
\[
 \delta^x(a)=ba', \qquad \qquad \delta^y(a)=-a'b,
\]
and zero in all other arrows. The reduction functors $F^X:A_2^X\text{-mod} \to A_2\text{-mod}$ and 
$F^Y:A_2^Y\text{-mod} \to A_2\text{-mod}$ have the following explicit form in objects respectively,
\[
\xymatrix{M_2 \ar@(lu,ru)[]^(.2){M_a} \ar@<.5ex>@{<.}[d] \ar@<-.5ex>[d]_-{M_{a'}} \\ M_1 } \quad \mapsto \quad
\xymatrix{M_2 \ar@<-.5ex>[d]_-{\left[ \begin{smallmatrix} 0 \\ I \end{smallmatrix} \right]} 
\ar@<.5ex>[d]^-{\left[ \begin{smallmatrix} M_{a'}\\ M_a \end{smallmatrix} \right]} \\ M_1 \oplus M_2 }
\quad \text{and} \quad
\xymatrix{N_2 \ar@<.5ex>@{<.}[d] \ar@<-.5ex>[d]_-{N_{a'}} \\ N_1 \ar@(ld,rd)[]_(.2){N_a} } \quad \mapsto \quad
\xymatrix{N_1 \oplus N_2 \ar@<-.5ex>[d]_-{\left[ \begin{smallmatrix} I & 0 \end{smallmatrix} \right]} 
\ar@<.5ex>[d]^-{\left[ \begin{smallmatrix} N_a & N_{a'} \end{smallmatrix} \right]} \\ N_1. }
\]
\end{lema}
\bproof
It is easy to see that the reduced tensor algebras
\[
 T_S(X^*\otimes_R ka_2 \otimes_R X \oplus k\sigma^*) \qquad \text{and} \qquad 
 T_{S'}(Y^*\otimes_R ka_2 \otimes_R Y \oplus k\pi^*),
\]
\[
\xymatrix{{}_{\bullet_{\omega}} \ar@(lu,ru)[]^(.5){(x_1^{\omega})^* \otimes a_2 \otimes x_2^{\omega}} \ar@<.5ex>@{<.}[d]^-{\sigma^*} 
\ar@<-.5ex>[d]_-{(x^1_{1})^* \otimes a_2 \otimes x_2^{\omega}} \\ {}_{\bullet_1} } \qquad \text{} \qquad
\xymatrix{{}_{\bullet_2} \ar@<.5ex>@{<.}[d]^-{\pi^*} \ar@<-.5ex>[d]_-{(y_1^{\omega})^* \otimes a_2 \otimes y^2_{2}} \\ 
{}_{\bullet_{\omega}} \ar@(ld,rd)[]_(.5){(y_1^{\omega})^* \otimes a_2 \otimes y_2^{\omega}} }
\]
have regular associated quivers. By lemma~\ref{(P)L:isoTens} these algebras are isomorphic to $kK_2^x$ and $kK_2^y$ respectively. 
We compute now their differentials. Since the differential in the Kronecker algebra $A_2$ is null, the reduced differentials
have the following form
\[
 \delta^x(x^* \otimes a_2 \otimes x)=\lambda(x^*) \otimes a_2 \otimes x - x^*\otimes a_2 \otimes \rho(x),
\]
\[
 \delta^y(y^* \otimes a_2 \otimes y)=\lambda'(y^*) \otimes a_2 \otimes y - y^*\otimes a_2 \otimes \rho'(y).
\]
Substituting the values of the coactions and considering that the coproducts $\mu$ and $\mu'$ in $P^*$ and $(P')^*$ 
are zero, the result follows. 

To justify the recovery algorithm corresponding to the reduction by the module $X$ described in the statement,
consider an $A_2^X$-representation $M=(M_1,M_2;M_a,M_{a'})$. Let $e_1$, $e_2$ be the trivial paths in $A_2$. Then $F^X(M)$,
considered as $B$-module, is given by $X\otimes_S M$, and is constituted by the vector spaces
\[
 e_1(X \otimes_S M)\cong e_1(x^1_1 \otimes M_1 \oplus e^{\omega}_1 \otimes_SM_2 \oplus x^{\omega}_2 \otimes_S M_2) \cong
 M_1 \oplus M_2,
\]
\[
 e_2(X \otimes_S M)\cong e_2(x^1_1 \otimes M_1 \oplus e^{\omega}_1 \otimes_SM_2 \oplus x^{\omega}_2 \otimes_S M_2) \cong
 M_2.
\]
The action of $a_1$ in $F^X(M)$ is given by the action of $a_1$ in $X^{\omega}$, that is, by the matrix 
$\left[ \begin{smallmatrix} 0\\I\end{smallmatrix} \right]$ (see~\ref{(A)P:funRed}). The action of $a_2$
is obtained through the actions of $a$ and $a'$ in $M$, hence they correspond to the action of the matrix
$\left[ \begin{smallmatrix} M_{a'}\\M_a\end{smallmatrix} \right]$.
The reduction by the module $Y$ can be described in a similar way.
\eproof

We give now algorithms for the construction of exceptional representations of $A_2^X$ and $A_2^Y$.

\begin{lema} \label{(P)L:redDos}
\begin{itemize}
 \item[a)] There is a full, faithful and rigid functor $G^X:A_2^X\text{-mod} \to A_2^X\text{-mod}$ which increases dimensions
and whose explicit form in objects is
\[
\xymatrix{M_2 \ar@(lu,ru)[]^(.2){M_a} \ar@<.5ex>@{<.}[d] \ar@<-.5ex>[d]_-{M_{a'}} \\ M_1 } \qquad \mapsto \qquad
\xymatrix{M_1 \oplus M_2 \ar@(lu,ru)[]^(.2){\left[ \begin{smallmatrix} 0 & M_{a'} \\ 0 & M_a \end{smallmatrix} \right]} 
\ar@<.5ex>@{<.}[d] \ar@<-.5ex>[d]_-{\left[ \begin{smallmatrix} I&0 \end{smallmatrix} \right]} \\ M_1. }
\]
 \item[b)] There is a full, faithful and rigid functor $G^Y:A_2^Y\text{-mod} \to A_2^Y\text{-mod}$ which increases dimensions
and whise explicit form in objects is
\[
\xymatrix{N_2 \ar@<.5ex>@{<.}[d] \ar@<-.5ex>[d]_-{N_{a'}} \\ N_1 \ar@(ld,rd)[]_(.2){N_a} } \qquad \mapsto \qquad
\xymatrix{N_2 \ar@<.5ex>@{<.}[d] \ar@<-.5ex>[d]_-{\left[ \begin{smallmatrix} 0\\I \end{smallmatrix} \right]} 
\\ N_1 \oplus N_2. \ar@(ld,rd)[]_(.2){\left[ \begin{smallmatrix} N_a & N_{a'} \\ 0 & 0 \end{smallmatrix} \right]} }
\]
\end{itemize} 
\end{lema}

\bproof
We show $(a)$, the proof of $(b)$ is similar. 
Take $B$ to be the subalgebra of $A^X$ generated by the arrow $a'$ and use lemma~\ref{(P)L:semiEje}.
Observe first that 
\[
 (T^Y)^tM_{K^x}T^Y=\left( \begin{matrix}1&1\\ 0&1 \end{matrix} \right)
\left( \begin{matrix} 1&1 \\-1&0 \end{matrix} \right)
\left( \begin{matrix} 1&0 \\ 1&1 \end{matrix} \right)=M_{K^x},
\]
so that we expect an equivalence of categories
\[
 (A_2^X)^Y\text{-mod}\longrightarrow A_2^X\text{-mod}.
\]
To be precise, we will give full, faithful and rigid functors
\[
 \xymatrix{
A_2^X\text{-mod} \ar[r]^-{F_3} & (A_2^X)^Y\text{-mod} \ar[r]^-{F_2} 
& (A_2^X)^Y\text{-mod} \ar[r]^-{F_1} & A_2^X\text{-mod},
}
\]
where $F_1=F^Y$ is the functor associated to the reduction of $A^X_2$ using $Y$, $F_2=F^c$ is induced by a change of basis 
and $F_3=F^r$ is a regularization. For the description of $F_1$ observe that the quiver associated to the reduced graded tensor 
algebra $(kK^X)^Y$ has the form
\[
 \xymatrix@C=4pc@R=5pc{
{}_{\bullet_2} \ar@(lu,ru)[]^(.2){\txt{\scriptsize $\underline{(y^2_2)^*ay^2_2}$\\ \scriptsize $\alpha$}} 
\ar@<.5ex>@{<.}[d]^-{\txt{\scriptsize $\underline{(y^2_2)^*by_1^{\omega}}$\\ \scriptsize $\beta'_1$}} 
\ar@{<-}@<-.5ex>[d]_-{\txt{\scriptsize $\underline{(y^2_2)^*ay_2^{\omega}}$ \\ \scriptsize $\alpha'_1$ }} \\ 
{}_{\bullet_{\omega}} \ar@{<-}@/^4pc/[u]^-{\txt{\scriptsize $\underline{(y_2^{\omega})^*ay^2_2}$ \\ \scriptsize $\alpha'$ }} 
\ar@{.>}@/_4pc/[u]_-{\txt{\scriptsize $\underline{\pi^*}$ \\ \scriptsize $\beta$ }} 
\ar@{<-}@(l,d)[]_(.4){\txt{\scriptsize $\underline{(y_2^{\omega})^*ay_2^{\omega}}$ \\ \scriptsize $\alpha_1$ }} 
\ar@{.>}@(d,r)[]_(.6){\txt{\scriptsize $\underline{(y_2^{\omega})^*by_1^{\omega}}$ \\ \scriptsize $\beta_1$ }}
}
\]
where, for the sake of simplicity, we suppress all tensor product signs $\otimes $ and give an alternative notation for each arrow
(depicted in the figure as denominators). Using the formulas for the coactions given in lemma~\ref{(P)L:semiEje} we compute the
reduced differential.  Its values in arrows are
\begin{center}
 \begin{tabular}{c | l  c |  l}
Arrow & Differential & Arrow & Differential \\
\hline
$\alpha $ & $\beta\otimes \alpha'$, & $\beta $ & $0$, \\
$\alpha'$ & $0$, & $\beta_1 $ & $0$, \\
$\alpha_1 $ & $\beta_1-\alpha' \otimes \beta$, & $\beta'_1 $ & $\beta\otimes \beta_1$. \\
$\alpha'_1 $ & $\beta'_1 + \beta\otimes \alpha_1-\alpha \otimes \beta$, \\
 \end{tabular}
\end{center}
The reduced functor $F_1$ has the following explicit form in representations,
using again the expression of $F^e$ given in equation~(\ref{(P)E:transformaRed})
\[
 \xymatrix{
{M_2} \ar@(lu,ru)[]^(.2){M_{\alpha}} \ar@<.5ex>@{<.}[d] 
\ar@{<-}@<-.5ex>[d]_-{M_{\alpha'_1}} \\ 
{M_{\omega}} \ar@{<-}@/^4pc/[u]^-{M_{\alpha'}} \ar@{.>}@/_4pc/[u] \ar@{<-}@(l,d)[]_(.4){M_{\alpha_1}} 
\ar@{.>}@(d,r)[]
}
\quad \qquad \qquad  \mapsto \quad
\xymatrix{
{M_2 \oplus M_{\omega}} \ar@(lu,ru)[]^(.3){\left[ \begin{smallmatrix} M_{\alpha} & M_{\alpha'_1} \\ M_{\alpha'} & M_{\alpha_1} \end{smallmatrix} \right] } \ar@<.5ex>@{<.}[d]
\ar@<-.5ex>[d]_-{[0 \; I]} \\ 
{M_{\omega}}
}
\]
Before regularization we perform the following change of basis. Define a copy $\widehat{Q}$ of the quiver $Q$ associated to
the reduced ditalgebra $(A_2^X)^Y$ and define functions $g:\widehat{Q} \to kQ$ and $h:Q \to k\widehat{Q}$ in the following way
\begin{equation*}
g(\widehat{\gamma}) = \left\{
\begin{array}{l l}
\gamma, & \text{if $\widehat{\gamma} \neq \widehat{\beta}_1,\widehat{\beta}'_1$},\\
\beta_1-\alpha' \otimes \beta, & \text{if $\widehat{\gamma}=\widehat{\beta}_1$}, \\
\beta'_1 + \beta\otimes \alpha_1-\alpha \otimes \beta, & \text{if $\widehat{\gamma}=\widehat{\beta}'_1$.}
\end{array} \right.
\end{equation*}
\begin{equation*}
h(\gamma) = \left\{
\begin{array}{l l}
\widehat{\gamma}, & \text{if $\gamma \neq \beta_1,\beta'_1$},\\
\widehat{\beta}_1+\widehat{\alpha}' \otimes \widehat{\beta}, & \text{if $\gamma=\beta_1$}, \\
\widehat{\beta}'_1 - \widehat{\beta} \otimes \widehat{\alpha}_1+\widehat{\alpha} \otimes \widehat{\beta}, & \text{if $\gamma=\beta'_1$.}
\end{array} \right.
\end{equation*}
Then the functions $g$ and $h$ extend to algebra morphisms which are inverse from each other. Observe that
$g$ and $h$ preserve degree, hence the transformation $\widehat{\delta}=h\delta g$ is a differential. The ditalgebra
isomorphism $h:(kQ,\delta) \to (k\widehat{Q},\widehat{\delta})$ induces an equivalence of categories 
$F_2$. The differential $\widehat{\delta}$ has the form
\begin{center}
 \begin{tabular}{c | l  c |  l}
  Arrow & Differential & Arrow & Differential \\
\hline
$\widehat{\alpha} $ & $\widehat{\beta} \otimes \widehat{\alpha'}$, &  $\widehat{\beta} $ & $0$, \\
$\alpha'$ & $0$, & $\widehat{\beta_1} $ & $0$, \\
$\widehat{\alpha_1} $ & $\widehat{\beta_1}$,  & $\widehat{\beta'_1} $ & $0$.\\
$\widehat{\alpha'_1} $ & $\widehat{\beta'_1}$, \\
 \end{tabular}
\end{center}
Equivalence $F_3$ results from regularization of the arrows $(\widehat{\alpha_1},\widehat{\beta_1})$ and 
$(\widehat{\alpha'_1},\widehat{\beta'_1})$. Using the expression of $F_1$ given above it is clear that the functor 
$G^X=F_3 \circ F_2 \circ F_1$ has the form given in the statement. Moreover, since the regularization $F_3=F^r$,
the change of basis $F_2=F^c$ and the semi-reduction of an edge $F_1=F^e$ are rigid, the functor $G^X$ is also rigid.
\eproof

For an $a \times b$ matrix $M$ denote by $M^{\rightarrow}$ (respectively $M^{\leftarrow}$)
the $a \times (b+1)$ matrix obtained by adding a column of zeros to the right 
(respectively left) of $M$. In a similar way define $M^{\uparrow}$ and $M^{\downarrow}$ 
adding a row of zeros to the top or bottom of $M$ respectively.
\begin{proposicion} \label{(P)P:Kro2}
For $t \in \mathbb{N}_0$ denote by $P_2^t$ and $Q_2^t$ the representations of $A_2$ given respectively by the matrices
\[
P_2^t= \xymatrix@!0@C=3pc@R=3pc{k^t \ar@<.5ex>[d]^-{I^{\uparrow}} \ar@<-.5ex>[d]_-{I^{\downarrow}} \\ 
k^{t+1} 
}
\qquad \text{ and } \qquad
Q_2^t= \xymatrix@!0@C=3pc@R=3pc{k^{t+1} \ar@<.5ex>[d]^-{I^{\rightarrow}} \ar@<-.5ex>[d]_-{I^{\leftarrow}} \\ 
k^t. 
}
\]
The representations $P^t_2$ and $Q^t_2$ ($t\geq 0$) are exceptional and form a complete list of posprojective and
preinjective $A_2$-modules respectively.
\end{proposicion}
\bproof
Successively applying the functors $G^X$ and $G^Y$ to the simple $A_2^X$-module $S(1)$ and to the simple $A_2^Y$-module $S(2)$ 
respectively, we obtain the following exceptional representations (see lemma~\ref{(P)L:redDos}), 
which we denote by $\widetilde{P_2^t}$ and $\widetilde{Q_2^t}$.
To the right of each representation is shown the coefficient quiver with respect to the canonical basis.
\[
\widetilde{P_2^t}=\xymatrix@C=1.3pc@R=1.5pc{k^t \ar@(lu,ru)[]^(.2){(I^{\leftarrow})^{\downarrow}} \ar@<.5ex>@{<.}[d] 
\ar@<-.5ex>[d]_-{\left[ \begin{smallmatrix} 1 & 0 \ldots 0 \end{smallmatrix} \right]} 
& & {}_{\bullet_1} \ar[dl] & \ar[l]  {}_{\bullet_2} & \ar[l] 
 {}_{\bullet_3} \ldots {}_{\bullet} & \ar[l] {}_{\bullet_{t-1}} & \ar[l] {}_{\bullet_t}\\
k &  {}_{\bullet_0}}
\]
\[
\widetilde{Q_2^t}= \xymatrix@C=1.3pc@R=1.5pc{k \ar@<.5ex>@{<.}[d] \ar@<-.5ex>[d]_-{\left[ \begin{smallmatrix} 0 \\ \vdots \\ 1 \end{smallmatrix} \right]} 
& & & & & & {}_{\bullet_0} \ar[dl] \\ k^t \ar@(ld,rd)[]_(.2){(I^{\leftarrow})^{\downarrow}} 
& {}_{\bullet_1} & \ar[l]  {}_{\bullet_2} & \ar[l] 
 {}_{\bullet_3} \ldots {}_{\bullet} & \ar[l] {}_{\bullet_{t-1}} & \ar[l] {}_{\bullet_t}
}
\]
Using now the reduction functors $F^X$ and $F^Y$ in lemma~\ref{(P)L:redUno}
we produce exceptional modules for the classical Kronecker algebra $A_2$, 
$P_2^t=F^X(\widetilde{P_2^t})$ and $Q_2^t=F^Y(\widetilde{Q_2^t})$
for $t \geq 0$. We show the corresponding coefficient quivers,
\[
P_2^t= \xymatrix@!0@C=3pc@R=3pc{k^t \ar@<.5ex>[d]^-{I^{\uparrow}} \ar@<-.5ex>[d]_-{I^{\downarrow}} 
& & {}_{\bullet_1} \ar[ld] \ar[d] & {}_{\bullet_2} \ar[ld] \ar[d] & {}_{\bullet_3} \ar[ld] \ar[d]  
& {}_{\bullet_{t-1}} \ar@{}[ld]|-{\cdots} \ar[d] & {}_{\bullet_t} \ar[ld] \ar[d]
\\ k^{t+1} & {}_{\bullet_0} & {}_{\bullet_1} & {}_{\bullet_2} & {}_{\bullet_3} & {}_{\bullet_{t-1}} & {}_{\bullet_t}
}
\]
\[
Q_2^t= \xymatrix@!0@C=3pc@R=3pc{k^{t+1} \ar@<.5ex>[d]^-{I^{\rightarrow}} \ar@<-.5ex>[d]_-{I^{\leftarrow}} 
& {}_{\bullet_1} \ar[d] & {}_{\bullet_2} \ar[d] \ar[dl] & {}_{\bullet_3} \ar[d] \ar[dl] & {}_{\bullet_{t-1}} \ar[d] \ar@{}[dl]|-{\cdots} 
& {}_{\bullet_t} \ar[d] \ar[dl] & {}_{\bullet_0} \ar[dl] 
\\ k^t & {}_{\bullet_1} & {}_{\bullet_2} & {}_{\bullet_3} & {}_{\bullet_{t-1}} & {}_{\bullet_t} 
}
\]
Observe that the dimension vectors are $\vdim P_2^t=p^t_2$ and $\vdim Q_2^t=q^t_2$.
By lemmas~\ref{(P)L:determ} and \ref{(P)L:endExc} there is an inclusion of isomorphism classes
of exceptional $A_2$-modules to the set of positive roots of the quadratic form $q_{K_2}$.
Hence the collection of modules $P_2^t$ and $Q_2^t$ ($t\geq 0$) determines a complete list of
exceptional representations for $A_2$. Moreover, computing the Cartan matrix $C_{A_2}$ and the Coxeter matrix
$\Phi_{K_2}$ corresponding to the quiver $K_2$,
\[
 C_{A_2}=\left( \begin{smallmatrix} 1&2\\0&1 \end{smallmatrix} \right),
\; \Phi_{K_2}=-C_{A_2}^tC_{A_2}^{-1}=\left( \begin{smallmatrix} -1&2\\-2&3 \end{smallmatrix} \right)
\; \text{and}
\; \Phi_{K_2}^{-1}=-C_{A_2}C_{A_2}^{-t}=\left( \begin{smallmatrix} 3&-2\\2&-1 \end{smallmatrix} \right),
\]
we notice that
\[
 \Phi_{K_2}q^{t-1}_2=\left( \begin{smallmatrix} -1&2\\-2&3 \end{smallmatrix} \right)
\left( \begin{smallmatrix} t-1\\t \end{smallmatrix} \right)=\left( \begin{smallmatrix} t+1\\t+2 \end{smallmatrix} \right)
=q^{t+1}_2,
\]
and in a similar way
\[
 \Phi_{K_2}^{-1}p^{t-1}_2=\left( \begin{smallmatrix} 3&-2\\2&-1 \end{smallmatrix} \right)
\left( \begin{smallmatrix} t\\t-1 \end{smallmatrix} \right)=\left( \begin{smallmatrix} t+2\\t+1 \end{smallmatrix} \right)
=p^{t+1}_2.
\]
Using lemma~\ref{(P)L:coxPP} and observing that $P_2^0$ and $P_2^1$ are indecomposable projective $A_2$-modules,
we conclude that the set $\{P^t_2\}_{t \geq 0}$ constitudes a complete list of posprojective $A_2$-modules. In a similar way,
since $Q_2^0$ and $Q_2^1$ are indecomposable injective $A_2$-modules, the set $\{Q^t_2\}_{t \geq 0}$ is a complete list 
of preinjective $A_2$-representations.
\eproof

By lemma~\ref{(DE)L:componentes} the posprojective component $\mathcal{P}$ of the Auslander-Reiten quiver $\Gamma(A_2)$
has the form
\[
\xymatrix@C=1pc@R=1pc{
[P_2^0] \ar@{--}[rr] \ar@<.5ex>[rd] \ar@<-.5ex>[rd] & & [P_2^2] \ar@{--}[rr] \ar@<.5ex>[rd] \ar@<-.5ex>[rd] 
& & [P_2^4] \ar@<.5ex>[rd] \ar@<-.5ex>[rd] & \cdots \\
& [P_2^1] \ar@{--}[rr] \ar@<.5ex>[ru] \ar@<-.5ex>[ru] & & [P_2^3] \ar@{--}[rr] \ar@<.5ex>[ru] \ar@<-.5ex>[ru] & & [P_2^5] & \cdots
}
\]
while the preinjective component $\mathcal{I}$ has the following form,
\[
\xymatrix@C=1pc@R=1pc{
\cdots & [Q_2^4] \ar@{--}[rr] \ar@<.5ex>[rd] \ar@<-.5ex>[rd] & & [Q_2^2] \ar@{--}[rr] \ar@<.5ex>[rd] \ar@<-.5ex>[rd] 
& & [Q_2^0] \ar@<.5ex>[rd] \ar@<-.5ex>[rd] \\
& \cdots & [Q_2^5] \ar@{--}[rr] \ar@<.5ex>[ru] \ar@<-.5ex>[ru] & & [Q_2^3] \ar@{--}[rr] \ar@<.5ex>[ru] \ar@<-.5ex>[ru] & & [Q_2^1].
}
\]

\section{First generalised case $n=3$.} \label{(P)S:KroTres}
%----------------------------------------------------------------------
%----------------------------------------------------------------------
we study now the first generalised case $A_3$. 
\[
K_3=\xymatrix{{}_{\bullet_2} \ar@<-2ex>[d]_-{a_1} \ar@<-.4ex>[d]^-{a_2} \ar@<2.3ex>[d]^-{a_3} \\ {}_{\bullet_1} }
\qquad 
M_{K_3}=\left( \begin{matrix} 1 & 0 \\ -3 & 1 \end{matrix} \right).
\]
\begin{lema} \label{(P)L:redZ}
Let $C$ be the subalgebra of $A_3$ generated by the arrows $a_2$ and $a_3$ ($C\cong A_2$). 
Let $Z={}_CC=Z^1 \oplus Z^2$ be the regular $C$-module, that is,
$Z^1=Ce_1$ is the projective simple module $S(1)$ and $Z^2=Ce_2$ is given by 
$\xymatrix{k \ar@<.5ex>[r]^-{\left[ \begin{smallmatrix} 1 \\ 0 \end{smallmatrix} \right]} 
\ar@<-.5ex>[r]_-{\left[ \begin{smallmatrix} 0 \\ 1 \end{smallmatrix} \right]} & k^2}$. Fix bases $\{z^1_1\}$ of $Z^1$ and
$\{z^2_1,\overline{z}^2_1,z^2_2\}$ of $Z^2$ such that $a_2z^2_2=z^2_1$ and $a_3z^2_2=\overline{z}^2_1$. 
The opposed endomorphism algebra $\Gamma=\End_C(Z)^{op}$ splits over its radical
$\Gamma=S'' \oplus P''$, where $S''\cong \End(Z^1) \oplus \End(Z^2)=kf''_1 \times kf''_2$ and the radical
$P'' \cong \Hom_C(Z^1,Z^2)$ is generated by a pair of morphisms $\nu,\overline{\nu}:Z^1 \to Z^2$ determined by
$z^1_1 \mapsto z^2_1$ and $z^1_1 \mapsto \overline{z}^2_1$.
Then $Z$ is an admissible rigid $C$-module. The matrix of dimension vectors of $Z$ is given by
\[
 T^Z=\left( \begin{matrix} 1&2\\0&1 \end{matrix} \right).
\]
The coactions of $Z$ in terms of the selected bases have the following form,
\begin{equation*}
\lambda(u) = \left\{
\begin{array}{l l}
\nu^* \otimes (z^1_1)^*, & \text{if $u=(z^2_1)^*$},\\
\overline{\nu}^* \otimes (z^1_1)^*, & \text{if $u=(\overline{z}^2_1)^*$},\\
0, & \text{otherwise.}
\end{array} \right.
\end{equation*}
\begin{equation*}
\rho(x) = \left\{
\begin{array}{l l}
z^2_1 \otimes \nu^* + \overline{z}^2_1 \otimes \overline{\nu}^*, & \text{if $x=z^1_{1}$},\\
0, & \text{otherwide.}
\end{array} \right.
\end{equation*}
\end{lema}
\bproof
All claims are inmediat.
\eproof

Let $B$ be the subalgebra of $A_3$ generated by the arrow $a_1$ and let
$X$, $Y$ be the $B$-modules described in lemma~\ref{(P)L:semiEje}.
We will perform reductions respect to the modules $X$, $Y$ and $Z$.
Consider the changes of basis
\[
(T^X)^tM_{K_3}T^X= \left( \begin{matrix} 1&1\\-2&-1 \end{matrix} \right) \quad \text{and} \quad 
(T^Z)^tM_{K_3}T^Z= \left( \begin{matrix} 1&2\\-1&-1 \end{matrix} \right).
\]
Their associated quivers, which will be denoted by $K^x_3$ and $K^z_3$, have respectively the form
\[
\xymatrix@R=3pc{{}_{\bullet_2} \ar@(u,l)[]_(.7){\alpha_1} \ar@(u,r)[]^(.7){\beta_1} \ar@<.5ex>@{<.}[d]^-{\gamma} 
\ar@/_4pt/[d]_-{\beta_2} \ar@/_18pt/[d]_-{\alpha_2} \\ {}_{\bullet_1} } \qquad \text{and} \qquad
\xymatrix@R=3pc{{}_{\bullet_2} \ar@(u,l)[]_(.7){a} \ar@(u,r)[]^(.7){b} \ar@{<.}[d]_-{c} \ar@<1ex>@{<.}[d]^-{d}
\ar@/_10pt/[d]_-{a'}  \\ {}_{\bullet_1} }
\]

\begin{lema} \label{(P)L:Kro3Uno}
\begin{itemize}
 \item[a)] The reduced ditalgebra $A_3^X$ is isomorphic to $(kK_3^x,\delta^x)$ 
where the nonzero values of the differential $\delta^x$ are
\[
\delta^x(\alpha_1)=\gamma \alpha_2 \quad \text{and} \quad  \delta^x(\beta_1)=\gamma \beta_2.
\]
The explicit description of the reduction functor in objects $M \mapsto F^X(M)$ has the following form
\[
\xymatrix@R=3pc{M_2 \ar@(u,l)[]_(.7){M_{\alpha_1}} \ar@(u,r)[]^(.7){M_{\beta_1}} \ar@<.5ex>@{<.}[d] 
\ar@/_4pt/[d]_-{M_{\beta_2}} \ar@/_28pt/[d]_-{M_{\alpha_2}} \\ M_1 } \qquad \mapsto \qquad 
\xymatrix@R=3pc{M_2 
\ar@<-8ex>[d]_-{\left[ \begin{smallmatrix}I\\0 \end{smallmatrix} \right]} 
\ar[d]_-{\left[ \begin{smallmatrix}M_{\alpha_1}\\M_{\alpha_2} \end{smallmatrix} \right]}
\ar@<8ex>[d]_-{\left[ \begin{smallmatrix}M_{\beta_1}\\M_{\beta_2} \end{smallmatrix} \right]} \\ M_2 \oplus M_1 }
\]
 \item[b)] The reduced ditalgebra $A_3^Z$ is isomorphic to $(kK_3^z,\delta^z)$ 
where the nonzero values of $\delta^z$ are given by
\[
 \delta^z(a)=ca' \quad \text{and} \quad \delta^z(b)=da'.
\]
The explicit description of the reduction functor in objects $N \mapsto F^Z(N)$ has the following form
\[
\xymatrix@R=3pc{N_2 
\ar@(u,l)[]_(.6){N_a} 
\ar@(u,r)[]^(.6){N_b} 
\ar@<.5ex>@{<.}[d] \ar@<1.5ex>@{<.}[d]
\ar@/_4pt/[d]_-{N_{a'}} \\ N_1 }\; \mapsto \;
\xymatrix@R=5pc{ N_2 
\ar@<-4ex>[d]_-{\left[ \begin{smallmatrix}N_a\\N_b\\N_{a'} \end{smallmatrix} \right]} 
\ar@<4ex>[d]_-{\left[ \begin{smallmatrix} 0 \\ I \\ 0 \end{smallmatrix} \right]}
\ar@<12ex>[d]_-{\left[ \begin{smallmatrix} I \\ 0 \\ 0 \end{smallmatrix} \right]} \\ 
N_2 \oplus N_2\oplus N_1. }
\]
\end{itemize}
\end{lema}
\bproof
Consider first the reduction $A_3^X$. As in the case $A_2^X$ observe that the tensor algebra 
$T_S(X^* \otimes (ka_2 \oplus ka_3) \otimes X \oplus k\sigma^*)$ has a regular associated quiver, and by lemma~\ref{(P)L:isoTens}
it is isomorphic to $kK_3^x$.
\[
 \xymatrix@R=4pc{{}_{\bullet_{\omega}} 
\ar@(u,l)[]_(.7){\txt{\scriptsize $\underline{(x^{\omega}_1)^*a_2x^{\omega}_2}$\\ \scriptsize $\alpha_1$}} 
\ar@(u,r)[]^(.7){\txt{\scriptsize $\underline{(x^{\omega}_1)^*a_3x^{\omega}_2}$\\ \scriptsize $\beta_1$}} 
\ar@<.5ex>@{<.}[d]^-{\txt{\scriptsize $\underline{\sigma^*}$\\ \scriptsize $\gamma$}} 
\ar@/_4pt/[d]_-{\txt{\scriptsize $\underline{(x^{1}_1)^*a_3x^{\omega}_2}$\\ \scriptsize $\beta_2$}} 
\ar@/_50pt/[d]_-{\txt{\scriptsize $\underline{(x^{1}_1)^*a_2x^{\omega}_2}$\\ \scriptsize $\alpha_2$}} \\ {}_{\bullet_1} }
\]
Using the coactions $\lambda$ and $\rho$ of $X$ given in lemma~\ref{(P)L:semiEje} and considering that $A_3$ has
null differential, the reduced differential $\delta^X$ is as in the statement. Hence
$A_3^X \cong (kK_3^x,\delta^x)$. 

Consider now the reduction $A_3^Z$. In a similar way we directly observe an isomorphism of tensor algebras
\[
A_3^Z=T_{S''}(Z^* \otimes ka_1 \otimes Z \oplus (P'')^*) \cong kK_3^Z,
\]
\[
 \xymatrix@R=4pc{{}_{\bullet_2} 
\ar@(u,l)[]_(.7){\txt{\scriptsize $\underline{(\overline{z}^{2}_1)^*a_1z^{2}_2}$\\ \scriptsize $a$}} 
\ar@(u,r)[]^(.7){\txt{\scriptsize $\underline{(z^{2}_1)^*a_1z^{2}_2}$\\ \scriptsize $b$}} 
\ar@<1ex>@{<.}[d]^-{\txt{\scriptsize $\underline{\overline{\nu}^*}$\\ \scriptsize $d$}} 
\ar@{<.}[d]_-{\txt{\scriptsize $\underline{\nu^*}$\\ \scriptsize $c$}} 
\ar@/_20pt/[d]_-{\txt{\scriptsize $\underline{(z^{1}_1)^*a_1z^{2}_2}$\\ \scriptsize $a'$}} \\ {}_{\bullet_1} }
\]
In order to compute the reduced differential $\delta^Z$ we use the coactions of $Z$ des\-cribed in lemma~\ref{(P)L:redZ}.
Then
\begin{eqnarray}
 \delta^Z((z^2_1)^* \otimes a_1 \otimes z^2_2) & = & \lambda((z^2_1)^*)\otimes a_1 \otimes z^2_2 
- (z^2_1)^* \otimes a_1 \otimes \rho(z^2_2) = \nonumber \\
& = & \nu^* \otimes (z^1_1)^* \otimes a_1 \otimes z^2_2, \nonumber
\end{eqnarray}
\begin{eqnarray}
 \delta^Z((\overline{z}^2_1)^* \otimes a_1 \otimes z^2_2) & = & \lambda((\overline{z}^2_1)^*)\otimes a_1 \otimes z^2_2 
- (\overline{z}^2_1)^* \otimes a_1 \otimes \rho(z^2_2) = \nonumber \\
& = & \overline{\nu}^* \otimes (z^1_1)^* \otimes a_1 \otimes z^2_2, \nonumber
\end{eqnarray}
and
\begin{eqnarray}
 \delta^Z((z^1_1)^* \otimes a_1 \otimes z^2_2) & = & \lambda((z^1_1)^*)\otimes a_1 \otimes z^2_2 
+ (z^1_1)^* \otimes a_1 \otimes \rho(z^2_2) = 0. \nonumber 
\end{eqnarray}
It is then clear that $A_3^Z \cong (kK^z_3,\delta^z)$ with $\delta^z$ as in the statement of the lemma.

\eproof

\begin{lema} \label{(P)L:Kro3Dos}
The subalgebra of $A_3^Z$ generated by the edge $a'$ is isomorphic to $B$, so it is possible to reduce $A_3^Z$ with respect to
the $B$-module $Y$. Then there is a rigid equivalence of categories $H:A_3^X\text{-mod} \to (A_3^Z)^Y\text{-mod}$.
The composition $F^YH:A^X_3\text{-mod} \to A^Z_3\text{-mod}$ has the following explicit description in objects,
\[
\xymatrix@R=3pc{M_2 \ar@(u,l)[]_(.7){M_{\alpha_1}} \ar@(u,r)[]^(.7){M_{\beta_1}} \ar@<.5ex>@{<.}[d] 
\ar@/_4pt/[d]_-{M_{\beta_2}} \ar@/_28pt/[d]_-{M_{\alpha_2}} \\ M_1 }
 \quad M \mapsto F^YH(M) \quad
\xymatrix@R=3pc{M_2 \oplus M_1 
\ar@(u,l)[]_(.6){\left[ \begin{smallmatrix} M_{\alpha_1}& 0\\M_{\alpha_2}& 0 \end{smallmatrix} \right]} 
\ar@(u,r)[]^(.6){\left[ \begin{smallmatrix}M_{\beta_1}&0\\M_{\beta_2}&0 \end{smallmatrix} \right]} 
\ar@<.5ex>@{<.}[d] \ar@<1.5ex>@{<.}[d]
\ar@/_4pt/[d]_-{\left[ \begin{smallmatrix}0 & I \end{smallmatrix} \right]} \\ M_1 }
\]
\end{lema}
\bproof
We will give equivalences of categories
\[
 \xymatrix{
A_3^X\text{-mod} \ar[r]^-{F_2} & (A_3^Z)^Y\text{-mod} \ar[r]^-{F_1} &  (A_3^Z)^Y\text{-mod},
}
\]
where $F_1$ is a change of basis and $F_2$ is the functor associated to a regularization.
Notice that the reduced tensor algebra $(A_3^Z)^Y$ can be described as follows
\[
 \xymatrix@R=4pc{*++[]{{}_{\bullet_2}} 
\ar@(u,l)[]_(.6){\txt{\scriptsize $\underline{(y^2_2)^*ay^2_2}$\\ \scriptsize $\alpha_1$}} 
\ar@(u,r)[]^(.6){\txt{\scriptsize $\underline{(y^2_2)^*by^2_2}$\\ \scriptsize $\beta_1$}} 
\ar@{<.}[d]^-{\gamma_2}
\ar@<3ex>@{<.}[d]^-{\varepsilon_2}
\ar@<-3ex>@{<.}[d]^-{\gamma}
\ar@<-4ex>@/_4pt/[d]_-{\txt{\scriptsize $\underline{(y^{\omega}_2)^*by^2_2}$\\ \scriptsize $\beta_2$}} 
\ar@<-4ex>@/_50pt/[d]_-{\txt{\scriptsize $\underline{(y^{\omega}_2)^*ay^2_2}$\\ \scriptsize $\alpha_2$}} \\ 
*+++[]{{}_{\bullet_{\omega}}}
\ar@(ld,l)[]^(.55){\txt{\scriptsize $\underline{(y^{\omega}_2)^*ay^{\omega}_2}$\\ \scriptsize $\alpha'_1$}} 
\ar@(d,ld)[]^(.55){\txt{\scriptsize $\underline{(y^{\omega}_2)^*by^{\omega}_2}$\\ \scriptsize $\beta'_1$}}  
\ar@{.>}@(d,dr)[]_(.55){\txt{\scriptsize $\underline{(y^{\omega}_2)^*cy^{\omega}_1}$\\ \scriptsize $\gamma_1$}} 
\ar@{.>}@(dr,r)[]_(.55){\txt{\scriptsize $\underline{(y^{\omega}_2)^*dy^{\omega}_1}$\\ \scriptsize $\varepsilon_1$}}  
\ar@<-5ex>@/_4pt/[u]_-{\txt{\scriptsize $\underline{(y^2_2)^*by^{\omega}_2}$\\ \scriptsize $\beta'_2$}} 
\ar@<-5ex>@/_50pt/[u]_-{\txt{\scriptsize $\underline{(y^2_2)^*ay^{\omega}_2}$\\ \scriptsize $\alpha'_2$}} \\ 
}
\]
where the dotted arrows in the center correspond to the elements $\gamma=\pi^*$, 
$\gamma_2=(y^2_2)^*\otimes c \otimes y^{\omega}_1$ and $\varepsilon_2=(y^2_2)^*\otimes d \otimes y^{\omega}_1$.
Using the coactions $\lambda'$ and $\rho'$ of $Y$ given in lemma~\ref{(P)L:semiEje} and considering the general expression
of the reduced differential $\delta^{zy}$ of $(A^Z_3)^Y$ we obtain
\[
 \delta^{zy}(u \otimes w \otimes x)=\lambda'(u)\otimes w \otimes x +\sigma_{u,x}(\delta^z(w))+(-1)^{|w|+1}u\otimes w \otimes \rho'(x),
\]
for arrows $w$ in $A_3^z$ (cf. lemma~\ref{(A)L:sigma} in the appendix for the definition of $\sigma_{u,x}$). 
We get the following values for the differential of $(A_3^Z)^Y$,
\begin{center}
 \begin{tabular}{c | l c | l}
  Arrow & Differential  & Arrow & Differential  \\
\hline
$\alpha_1 $ & $\gamma \otimes \alpha_2$, & $\gamma$ & $0$, \\
$\beta_1 $ & $\gamma \otimes \beta_2$, & $\gamma_1$ & $0$, \\
$\alpha'_1 $ & $\gamma_1 - \alpha_2 \otimes \gamma$, & $\varepsilon_1$ & $0$, \\
$\beta'_1 $ & $\varepsilon_1-\beta_2 \otimes \gamma$, & $\gamma_2$ & $\gamma \otimes \gamma_1$, \\
$\alpha_2 $ & $0$, & $\varepsilon_2$ & $\gamma \otimes \varepsilon_1$. \\
$\beta_2 $ & $0$, \\
$\alpha'_2 $ & $\gamma_2+\gamma \otimes \alpha'_1-\alpha_1\otimes \gamma$, \\
$\beta'_2 $ & $\varepsilon_2+\gamma \otimes \beta'_1-\beta_1\otimes \gamma$, \\
 \end{tabular}
\end{center}
Define una copy $\widehat{Q}$ of the quiver $Q$ associated to the reduced algebra 
$(A_3^Z)^Y$ and functions $g_1:\widehat{Q} \to kQ$ and $h_1:Q \to k\widehat{Q}$ given by
\begin{equation*}
g_1(\widehat{x}) = \left\{
\begin{array}{l l}
x, & \text{if $\widehat{x} \neq \widehat{\gamma}_1,\widehat{\gamma}_2,\widehat{\varepsilon}_1,\widehat{\varepsilon}_2$},\\
\gamma_1-\alpha_2\otimes \gamma, & \text{if $\widehat{x}=\widehat{\gamma}_1$}, \\
\varepsilon_1-\beta_2\otimes \gamma, & \text{if $\widehat{x}=\widehat{\varepsilon}_1$}, \\
\gamma_2+\gamma \otimes \alpha'_1-\alpha_1\otimes \gamma, & \text{if $\widehat{x}=\widehat{\gamma}_2$}, \\
\varepsilon_2+\gamma \otimes \beta'_1-\beta_1\otimes \gamma, & \text{if $\widehat{x}=\widehat{\varepsilon}_2$}, \\
\end{array} \right.
\end{equation*}
\begin{equation*}
h_1(x) = \left\{
\begin{array}{l l}
\widehat{x}, & \text{if $x \neq \gamma_1,\gamma_2,\varepsilon_1,\varepsilon_2$},\\
\widehat{\gamma}_1+\widehat{\alpha}_2\otimes \widehat{\gamma}, & \text{if $x=\gamma_1$}, \\
\widehat{\varepsilon}_1+\widehat{\beta}_2\otimes \widehat{\gamma}, & \text{if $x=\varepsilon_1$}, \\
\widehat{\gamma}_2-\widehat{\gamma} \otimes \widehat{\alpha}'_1+\widehat{\alpha}_1\otimes \widehat{\gamma}, & \text{if $x=\gamma_2$}, \\
\widehat{\varepsilon}_2-\widehat{\gamma} \otimes \widehat{\beta}'_1+\widehat{\beta}_1\otimes \widehat{\gamma}, & \text{if $x=\varepsilon_2$}. \\
\end{array} \right.
\end{equation*}
Then $g_1$ and $h_1$ extend to morphisms $g$ and $h$ of graded algebras which preserve degree and which are inverse from each other.
Hence $\widehat{\delta}=h\delta g$ is a differential and the ditalgebra isomorphism $h:(kQ,\delta) \to 
(k\widehat{Q},\widehat{\delta})$ induces an equivalence of categories $F_1$. The new differential has the following form
\begin{center}
 \begin{tabular}{c | l c | l}
 Arrow & Differential  & Arrow & Differential  \\
\hline
$\widehat{\alpha}_1 $ & $\widehat{\gamma} \otimes \widehat{\alpha}_2$, & $\widehat{\gamma}$ & $0$, \\
$\widehat{\beta}_1 $ & $\widehat{\gamma} \otimes \widehat{\beta}_2$, & $\widehat{\gamma_1}$ & $0$, \\
$\widehat{\alpha}'_1 $ & $\widehat{\gamma}_1$, & $\widehat{\varepsilon_1}$ & $0$, \\
$\widehat{\beta}'_1 $ & $\widehat{\varepsilon}_1$, & $\widehat{\gamma_2}$ & $0$, \\
$\widehat{\alpha}_2 $ & $0$, & $\widehat{\varepsilon_2}$ & $0$. \\
$\widehat{\beta}_2 $ & $0$, \\
$\widehat{\alpha}'_2 $ & $\widehat{\gamma}_2$, \\
$\widehat{\beta}'_2 $ & $\widehat{\varepsilon}$, \\
 \end{tabular}
\end{center}
so it is possible to regularize arrows $(\widehat{\alpha}'_1,\widehat{\gamma}_1)$, 
$(\widehat{\beta}'_1,\widehat{\varepsilon}_1)$, $(\widehat{\alpha}'_2,\widehat{\gamma}_2)$ and
$(\widehat{\beta}'_2,\widehat{\varepsilon}_2)$. We obtain in this way a ditalgebra isomorphic to $A_3^X$ and an equivalence
of categories $F_2:A_3^X\text{-mod} \to (k\widehat{Q},\widehat{\delta})\text{-mod}$. The equivalence $H$ is composition of the
functors $F_1 \circ F_2$.
\eproof

By the lemmas above there are full, faithful and rigid functors
\[
 \xymatrix@C=3pc{
A_3^X\text{-mod} \ar[d]_-{H} \ar[rr]^-{F^X} & & A_3\text{-mod} \\
(A_3^Z)^Y\text{-mod} \ar[r]_-{F^Y} & A_3^Z\text{-mod} \ar[r]_-{F^Z} & A_3\text{-mod}.
}
\]
Evaluation of these functors in an object $M$ of $A_3^X$-mod, followed by the dimension vector, has the form
(see lemma~\ref{(P)L:isoTens})
\[
 \xymatrix@C=3pc{
\vdim M \ar@{|->}[d] \ar@{|->}[rr] & & T^X\vdim M \\
\vdim M \ar@{|->}[r] & T^Y\vdim M \ar@{|->}[r] & T^ZT^Y\vdim M.
}
\]
Hence, if $\widehat{M}$ is an exceptional $A_3$-module of dimension vector $(a,b)$ and $M$ is an
$A_3^X$-module with $F^X(M)\cong \widehat{M}$, then $F^ZF^YH(M)$ is also an exceptional $A_3$-representation, 
with dimension vector
\begin{eqnarray} 
 \vdim F^ZF^YH(M) & = & T^ZT^Y(T^X)^{-1} \vdim \widehat{M}= \left( \begin{smallmatrix} 1&2\\0&1 \end{smallmatrix} \right)
\left[ \left( \begin{smallmatrix} 1&0\\1&1 \end{smallmatrix} \right)
\left( \begin{smallmatrix} 1&-1\\0&1 \end{smallmatrix} \right) \right] \vdim \widehat{M} = \nonumber \\
& = & \left( \begin{smallmatrix} 1&2\\0&1 \end{smallmatrix} \right) 
\left( \begin{smallmatrix} 1&-1\\1&0 \end{smallmatrix} \right)\vdim \widehat{M}=
\left( \begin{smallmatrix} 3&-1\\1&0 \end{smallmatrix} \right)\vdim \widehat{M}, \nonumber 
\end{eqnarray}
then the algorithm $F^X(M) \mapsto F^ZF^YH(M)$ produces an exceptional
$A_3$-repre\-sen\-ta\-tion of dimension vector $(3a-b,a)$.

In the proof of the following proposition we consider the process above to construct exceptional representations
of $A_3$. Make $a_0=0$, $a_1=1$, $a_{t+1}=3a_t-a_{t-1}$ for $t \geq 1$ and take the vectors
\[
 p^t_3=(a_{t+1},a_t), \quad \text{for $t\geq 0$}.
\]
As the quadratic form $q_{K_3}$ is given by $q_{K_3}(a,b)=a^2+b^2-3ab$, observe that
\begin{eqnarray}
q_{K_3}(a_{t+1},a_t) & = & a_{t+1}^2+a_t^2-3a_{t+1}a_t = \nonumber \\ 
& = & (3a_{t}-a_{t-1})^2+a_t^2-3(3a_{t}-a_{t-1})a_t = \nonumber \\
& = & 9a_{t}^2-6a_ta_{t-1}+a_{t-1}^2+a_{t}^2-9a_t^2+3a_{t}a_{t-1} = \nonumber \\
& = & a_{t}^2+a_{t-1}^2-3a_{t}a_{t-1} = \nonumber \\ 
& = & q_{K_3}(a_t,a_{t-1}). \nonumber
\end{eqnarray}
Since $q_{K_3}(p^0_3)=q_{K_3}(1,0)=1$ then all vectors $p^t_3$ ($t\geq 0$) are positive roots 
of the quadratic form $q_{K_3}$. We notice inductively that $a_{t+1}>2a_t$, for
\[
 a_{t+1}-2a_t=(3a_t-a_{t-1})-2a_t=a_t-a_{t-1}\geq a_t-2a_{t-1}.
\]

In the presentation of posprojective $A_3$-modules we will use the following notation.
Let $Z(a)=Z_b(a)$ be the $a\times b$ zero matrix. Let $E(a)=E_b(a)$ be the $a\times b$ matrix whose diagonal elements
are $E(a)_{ii}=1$ and all other entries are equal to zero. If $b \leq a$ then $E(a)$ has the form 
$\left[ \begin{smallmatrix}I\\0 \end{smallmatrix} \right]$. Consider the vertical concatenation of matrices 
$C_b(a)=C(a)=\left[ \begin{smallmatrix} E(a) \\ Z(a) \end{smallmatrix} \right]$ of dimension $2a\times b$.
\begin{proposicion} \label{(P)P:Kro3}
For $t \in \mathbb{N}$ denote by $P^t_3$ the $A_3$-representation given by the matrices
(with $a_{t+1}$ rows and $a_{t}$ columns)
\[
 \xymatrix@R=5pc{k^{a_{t}} 
\ar@<-12ex>[d]_-{\left[ \begin{smallmatrix} Z(a_{t-1})\\ C(a_{t-1}) \\ \vdots \\C(a_1)\\Z(2)\\E(a_{t})\end{smallmatrix} \right]} 
\ar@<4ex>[d]_-{\left[ \begin{smallmatrix} Z(a_{t}) \\ E(a_{t})\\ Z(a_{t+1}-2a_{t}) \end{smallmatrix} \right]}
\ar@<19ex>[d]_-{\left[ \begin{smallmatrix}E(a_{t})\\Z(a_{t+1}-a_{t}) \end{smallmatrix} \right]} \\ 
k^{a_{t+1}}. }
\]
If $P_3^0$ denotes the projective simple $A_3$-module, then for $t \geq 0$ the representations $P^t_3$ are exceptional
and conform a complete list of posprojective $A_3$-modules.
\end{proposicion}

\bproof
We prove that the given representations are exceptional. It is clear that $P^0_3$ is exceptional, 
se we proceed by induction over $t\geq 1$. As base case, the representation $P_3^1$ has the form
\[
 \xymatrix@R=3pc@C=4pc{k^{a_{1}} 
\ar@<-10ex>[d]_-{\left[ \begin{smallmatrix} Z(2)\\E(a_{1})\end{smallmatrix} \right]} 
\ar@<4ex>[d]_-{\left[ \begin{smallmatrix} Z(a_{1}) \\ E(a_{1})\\ Z(a_{2}-2a_{1}) \end{smallmatrix} \right]}
\ar@<17ex>[d]_-{\left[ \begin{smallmatrix}E(a_{1})\\Z(a_{2}-a_{1}) \end{smallmatrix} \right]} 
& & {} \ar@{}[d]_-{=} & k
\ar@<-5ex>[d]_-{\left[ \begin{smallmatrix} 0\\0\\1 \end{smallmatrix} \right]} 
\ar[d]_-{\left[ \begin{smallmatrix} 0\\1\\0 \end{smallmatrix} \right]} 
\ar@<5ex>[d]_-{\left[ \begin{smallmatrix} 1\\0\\0 \end{smallmatrix} \right]} \\ 
k^{a_{2}} & & {} & k^3,
}
\]
and it is easy to see that $\End_{A_3}(P^1_3) \cong k$, and thus $P^1_3$ is an exceptional module (its dimension vector
$\vdim P^1_3=p^1_3$ is a root of $q_{K_3}$). Assume that the module $P_3^t$ is exceptional and let $M(t)_1$, $M(t)_2$ and 
$M(t)_3$ be the matrices, from left to right, which conform the module $P_3^t$.
Consider the representation $\widehat{P_3^t}$, isomorphic to $P_3^t$, 
obtained by exchanging positions of the matrices $M(t)_1$ and $M(t)_3$,
\[
\widehat{P_3^t}=\xymatrix@R=7pc{k^{a_t}
\ar@<-11ex>[d]_-{\left[ \begin{smallmatrix}E(a_t)\\ \hline Z(a_{t+1}-a_t) \end{smallmatrix} \right]} 
\ar@<4ex>[d]_-{\left[ \begin{smallmatrix} Z(a_{t}) \\ \hline E(a_{t})\\ Z(a_{t+1}-2a_{t}) \end{smallmatrix} \right]}
\ar@<20ex>[d]_-{\left[ \begin{smallmatrix} Z(a_{t-1}) \\ E(a_{t-1}) \\ Z(a_{t-1}-a_{t-2})\\ \hline Z(a_{t-2})\\ C(a_{t-2}) 
\\ \vdots \\C(a_1)\\Z(2)\\E(a_{t})\end{smallmatrix} \right]} \\ 
k^{a_t} \oplus k^{2a_t-a_{t-1}}.
}
\]
In this way, using lemma~\ref{(P)L:Kro3Uno}$(a)$, we have that $\widehat{P^t_3}=F^X(M)$ where $M$ is the
$A_3^X$-representation given by
\[
\xymatrix@R=7pc{k^{a_t} \ar@(u,l)[]_(.7){Z(a_t)} 
\ar@(u,r)[]^(.7){\left[ \begin{smallmatrix} Z(a_{t-1}) \\ E(a_{t-1}) \\ Z(a_{t-1}-a_{t-2}) \end{smallmatrix} \right]} 
\ar@<-.5ex>@{<.}[d] 
\ar@/^4pt/[d]^-{\left[ \begin{smallmatrix} Z(a_{t-2})\\ C(a_{t-2}) 
\\ \vdots \\C(a_1)\\Z(2)\\E(a_{t})\end{smallmatrix} \right]} 
\ar@/_28pt/[d]_-{\left[ \begin{smallmatrix} E(a_t) \\ Z(a_{t+1}-2a_t) \end{smallmatrix} \right]} 
\\ k^{2a_t-a_{t-1}}. } 
\]
Then $M$ is exceptional for $F^X$ is a full, faithful and rigid functor. Apply now the composition $F^YH$ 
described in lemma~\ref{(P)L:Kro3Dos} to the $A_3^X$-module $M$ to obtain an exceptional $A^Z_3$-module $N=F^YH(M)$,
\[
\xymatrix@R=3pc{k^{a_{t+1}}
\ar@(u,l)[]_(.6){\left[ \begin{smallmatrix} Z(a_t)\\ E(a_t) \\ Z(a_{t+1}-2a_t) \end{smallmatrix} \right]} 
\ar@(u,r)[]^(.6){\left[ \begin{smallmatrix}Z(a_{t-1}) \\ E(a_{t-1}) \\ Z(a_{t-1}-a_{t-2})\\ \hline Z(a_{t-2})\\ C(a_{t-2}) 
\\ \vdots \\C(a_1)\\Z(2)\\E(a_{t})\end{smallmatrix} \right]} 
\ar@<.5ex>@{<.}[d] \ar@<1.5ex>@{<.}[d]
\ar@/_4pt/[d]_-{\left[ \begin{smallmatrix}0 & I \end{smallmatrix} \right]} 
\\ k^{2a_t-a_{t-1}} }
=
\xymatrix@R=3pc{k^{a_{t+1}}
\ar@(u,l)[]_(.6){\left[ \begin{smallmatrix} Z(a_t)\\ E(a_t) \\ Z(a_{t+1}-2a_t) \end{smallmatrix} \right]} 
\ar@(u,r)[]^(.6){\left[ \begin{smallmatrix}Z(a_{t-1}) \\ C(a_{t-1}) \\ C(a_{t-2}) 
\\ \vdots \\C(a_1)\\Z(2)\\E(a_{t})\end{smallmatrix} \right]} 
\ar@<.5ex>@{<.}[d] \ar@<1.5ex>@{<.}[d]
\ar@/_4pt/[d]_-{\left[ \begin{smallmatrix}0 & I \end{smallmatrix} \right]} 
\\ k^{2a_t-a_{t-1}} }
\]
Finally, using the algorithm $N \mapsto F^Z(N)$ given in lemma~\ref{(P)L:Kro3Uno}$(b)$ 
we obtain the following representation of $A_3$,
\[
 \xymatrix@R=8pc{k^{a_{t}} \oplus k^{a_{t+1}-a_t} 
\ar@<-12ex>[d]_-{\left[ \begin{smallmatrix} Z(a_t) & 0\\E(a_t)&0\\Z(a_{t+1}-2a_t)&0\\ 
\hline Z(a_{t-1})&0\\C(a_{t-1})&0\\ \vdots & \vdots \\C(a_1)&0\\Z(2)&0 \\E(a_t)&0 
\\ \hline 0&I_{a_{t+1}-a_t} \end{smallmatrix} \right]} 
\ar@<4ex>[d]_-{\left[ \begin{smallmatrix}0&0 \\ & \\ \hline 
& \\ I_{a_t}&0\\0&I_{a_{t+1}-a_t}\\ & \\ \hline & \\ 0&0 \end{smallmatrix} \right]}
\ar@<20ex>[d]_-{\left[ \begin{smallmatrix}I_{a_t}&0\\0&I_{a_{t+1}-a_t} \\ & \\ \hline & \\ 0&0 \\ & \\ \hline & \\ 0&0 \end{smallmatrix} \right]} \\ 
(k^{a_{t+1}}) \oplus (k^{a_{t+1}}) \oplus k^{a_{t+1}-a_t}. }
\]
Observe that $a_{t+1}-2a_t+a_{t-1}=3a_t-a_{t-1}-2a_t+a_{t-1}=a_t$, and hence in the left matrix the
vertical concatenation $E(a_t)|Z(a_{t+1}-2a_t)|Z(a_{t-1})$ is equal to $C(a_t)$. Moreover, the direct sum at the bottom of 
this matrix $E(a_t)\oplus I_{a_{t+1}-a_t}$ corresponds to the matrix $E(a_{t+1})$. 
Thus, this module coincides with the representation $P_3^{t+1}$ given by
\[
 \xymatrix@R=6pc{k^{a_{t+1}} 
\ar@<-13ex>[d]_-{\left[ \begin{smallmatrix} Z(a_t)\\ C(a_t)\\C(a_{t-1})\\ \vdots \\C(a_1)\\Z(2)\\E(a_{t+1})\end{smallmatrix} \right]} 
\ar@<4ex>[d]_-{\left[ \begin{smallmatrix} Z(a_{t+1}) \\ E(a_{t+1})\\ Z(a_{t+2}-2a_{t+1}) \end{smallmatrix} \right]}
\ar@<20ex>[d]_-{\left[ \begin{smallmatrix}E(a_{t+1})\\Z(a_{t+2}-a_{t+1}) \end{smallmatrix} \right]} \\ 
k^{a_{t+2}}, }
\]
and hence $P^{t+1}_3$ is also an exceptional $A_3$-module. This completes the induction step.

We compute now the Cartan matrix $C_{A_3}$ and the inverse of the Coxeter matrix $\Phi_{K_3}^{-1}$ corresponding
to the generalised Kronecker quiver $K_3$,
\[
 C_{A_3}=\left( \begin{smallmatrix} 1&3\\0&1 \end{smallmatrix} \right)
\qquad \text{and}
\qquad \Phi_{K_3}^{-1}=-C_{A_2}C_{A_2}^{-t}=\left( \begin{smallmatrix} 8&-3\\3&-1 \end{smallmatrix} \right).
\]
Notice that
\[
 \Phi_{K_3}^{-1}p^{t-1}_3=\left( \begin{smallmatrix} 8&-3\\3&-1 \end{smallmatrix} \right)
\left( \begin{smallmatrix} a_t\\a_{t-1} \end{smallmatrix} \right)=
\left( \begin{smallmatrix} 8a_t-3a_{t-1}\\3a_t-a_{t-1} \end{smallmatrix} \right)
=\left( \begin{smallmatrix} 3a_{t+1}-a_{t}\\a_{t+1} \end{smallmatrix} \right)=
\left( \begin{smallmatrix} a_{t+2}\\a_{t+1} \end{smallmatrix} \right)=p^{t+1}_3.
\]
Using lemma~\ref{(P)L:coxPP} and observing that $P_3^0$ and $P_3^1$ are indecomposable projective $A_3$-modules,
we conclude that the set $\{P^t_3\}_{t\geq 0}$ constitudes a complete list of posprojective $A_3$-modules.
\eproof

The preinjective component can be treated in a similar way, changing the module $Z$ for the direct sum of the indecomposable
injective $A_3$-modules. By lemma~\ref{(DE)L:componentes} the posprojective and preinjective
components of the Auslander-Reiten quiver $\Gamma(A_3)$ have respectively the form
\[
\xymatrix@C=1pc@R=1pc{
[P_3^0] \ar@{--}[rr] \ar@<.8ex>[rd] \ar[rd] \ar@<-.8ex>[rd] & & [P_3^2] \ar@{--}[rr] \ar@<.8ex>[rd] \ar[rd] \ar@<-.8ex>[rd] 
& & [P_3^4] \ar@<.8ex>[rd] \ar[rd] \ar@<-.8ex>[rd] & \cdots \\
& [P_3^1] \ar@{--}[rr] \ar@<.8ex>[ru] \ar[ru] \ar@<-.8ex>[ru] & & [P_3^3] \ar@{--}[rr] \ar@<.8ex>[ru] \ar[ru] \ar@<-.8ex>[ru] & & [P_3^5] & \cdots
}
\]
and
\[
\xymatrix@C=1pc@R=1pc{
\cdots & [Q_3^4] \ar@{--}[rr] \ar@<.8ex>[rd] \ar[rd] \ar@<-.8ex>[rd] & & [Q_3^2] \ar@{--}[rr] \ar@<.8ex>[rd] \ar[rd] \ar@<-.8ex>[rd] 
& & [Q_3^0] \ar@<.8ex>[rd] \ar[rd] \ar@<-.8ex>[rd] \\
& \cdots & [Q_3^5] \ar@{--}[rr] \ar@<.8ex>[ru] \ar[ru] \ar@<-.8ex>[ru] & & [Q_3^3] \ar@{--}[rr] \ar@<.8ex>[ru] \ar[ru] \ar@<-.8ex>[ru] & & [Q_3^1].
}
\]
All given presentations of posprojective $A_3$-modules correspond to those given by Ringel in proposition 3 in \cite{cmR98}. 
Coefficient quivers for the cases $P^1_3$, $P^2_3$ and $P^3_3$ can be found in \cite{cmR98} after the mentioned proposition.

\section{The Auslander-Reiten quiver of a Dynkin diagram.} \label{(DE)S:ARDynkin}
%------------------------------------------------------------------

We say that a finite solid quiver $Q$ is a (simply laced) \textbf{Dynkin quiver} if the underlying graph 
of $Q$ is isomorphic to one of the diagrams shown in table~\ref{T:dynkin}. Let $Q$ be a finite connected solid quiver.
A fundamental result in the representation theory of path algebras indicates that the following statements are equivalent.
\begin{itemize}
 \item[a)] $Q$ is a Dynkin quiver.
 \item[b)] The category $kQ$-mod is of finite representation type.
 \item[c)] The quadratic form $q_Q$ is positive definite. 
\end{itemize}
Moreover, $\vdim$ establishes a bijection between the isomorphism classes of indecomposable $A$-modules and the positive
roots of $q_Q$ (theorem 13 in \cite[section 2.4]{cmR}). 

\begin{table} [!hbt] 
\begin{center}
\renewcommand{\arraystretch}{1.2}
\begin{tabular}{c l} 
\hline 
Notation & \multicolumn{1}{c}{Valued graph} \\
\hline \\
$\mathbf{A_n} \; (n \geq 1)$ & $\xymatrix{ *+[o][F-]{{}_1} \ar@{-}[r] & {}_1 \ar@{-}[r] 
& {}_1 \ldots {}_1 \ar@{-}[r] & *+[o][F-]{{}_1}} $ \\
\raisebox{-1ex}{$\mathbf{D_n} \; (n \geq 4)$} & $\xymatrix@R=-.5pc{{}_1 \ar@{-}[rd] \\  
& {}_2 \ar@{-}[r] & {}_2 \ldots {}_2 \ar@{-}[r] & *+[o][F-]{{}_2} \ar@{-}[r] & {}_1 \\ {}_1 \ar@{-}[ru] }  $ \\
\raisebox{-2ex}{$\mathbf{E_6}$} & $\xymatrix@R=1pc{ & & *+[o][F-]{{}_2} \ar@{-}[d] \\ 
{}_1 \ar@{-}[r] & {}_2 \ar@{-}[r] & {}_3 \ar@{-}[r] & {}_2 \ar@{-}[r] & {}_1} $ \\
\raisebox{-2ex}{$\mathbf{E_7}$} & $\xymatrix@R=1pc{ & & {}_2 \ar@{-}[d] \\ 
*+[o][F-]{{}_2} \ar@{-}[r] & {}_3 \ar@{-}[r] & {}_4 \ar@{-}[r] & {}_3 \ar@{-}[r] & {}_2 \ar@{-}[r] & {}_1} $ \\
\raisebox{-2ex}{$\mathbf{E_8}$} & $\xymatrix@R=1pc{ & & {}_3 \ar@{-}[d] \\ 
{}_2 \ar@{-}[r] & {}_4 \ar@{-}[r] & {}_6 \ar@{-}[r] & {}_5 \ar@{-}[r] & {}_4 \ar@{-}[r] & {}_3 \ar@{-}[r] & *+[o][F-]{{}_2}} $ \\ [1ex]
\hline
\end{tabular}
\caption{Simply laced Dynkin diagrams $\Delta$. The integral vector indicated in the vertices is the maximal positive root 
$w_0$ of the corresponding quadratic form (respect to the partial order in $\mathbb{Z}^n$ given by $x > y$ if $x-y$ is a positive
vector). In a circle are marked the exceptional vertices of $w_0$, that is, vertices $i$ such that the partial derivative 
$(\partial/\partial_i)q_{\Delta}(w_0)$ is zero.} \label{T:dynkin}
\end{center}
\end{table}
%------------------------------------------------------------------
In this section we describe the module category of a Dynkin quiver.
Let $A_0=k\Delta$ be the path algebra of a Dynkin quiver $\Delta$ and consider the Auslander-Reiten quiver 
$\Gamma(A_0)$ of the category of $A_0$-modules.

\begin{lema} \label{(DE)L:partDyn}
There exists a partition in sections of the Auslander-Rei\-ten quiver $\Gamma(A_0)$,
\[
 \Gamma(A_0) = \mathcal{S}_0 \sqcup \mathcal{S}_1 \sqcup \ldots \sqcup \mathcal{S}_n,
\]
such that $\mathcal{S}_0$ is conformed by the isomorphism classes of indecomposable projective $A_0$-modules, 
$\tau \mathcal{S}_{i+1} \subseteq \mathcal{S}_i$ for $0 \leq i < n$
and every vertex in $\mathcal{S}_n$ is the isomorphism class of some indecomposable injective $A_0$-module. 
In particular, $\Gamma(A_0)$ can be identified with a mesh complete subquiver of $\mathbb{N}_0\Delta^{op}$.
\end{lema}
\bproof
Let $\mathcal{S}_0$ be the set of isomorphism classes of indecomposable projective $A_0$-modules.
Take as representative the modules $P(x)=Ae_x$, where $e_x$ is the trivial path of vertex $x$.\\
\underline{Step 1.} \textit{The vertex set $\mathcal{S}_0$ is a section in $\Gamma(A_0)$.}
The first condition $(i)$ is satisfied by vacuity, for $\mathcal{S}_0$ contains no nonprojective vertex.
Let $[P] \to [M]$ be an arrow with $[P]$ in $\mathcal{S}_0$ and assume that $[M]$ is not in $\mathcal{S}_0$. 
Then there is an arrow $\tau[M] \to [P]$, and since $\Gamma(A_0)$ is hereditary, $\tau[M]$
belongs to $\mathcal{S}_0$.

Assume now that we have constructed sections $\mathcal{S}_0,\ldots,\mathcal{S}_i$
such that $\tau\mathcal{S}_{\ell+1}\subseteq \mathcal{S}_{\ell}$ for each $0 \leq \ell <i$.\\
\underline{Step 2.} \textit{Defining $\mathcal{S}_{i+1}=\{ [M] \in \Ind A_0 \; | \; \tau[M] \in \mathcal{S}_i \}$
we prove that $\mathcal{S}_{i+1}$ is a section}. First, if $[M] \in \mathcal{S}_{i+1}$ then
$\tau[M] \in \mathcal{S}_i$, hence $\tau(\tau[M])\notin \mathcal{S}_i$, that is, $\tau[M]\notin \mathcal{S}_{i+1}$.
Second, let $[M] \to [X]$ be an arrow with $[M] \in \mathcal{S}_{i+1}$. Since $[M]$ is not a projective vertex
(for $\tau[M] \in \mathcal{S}_i$) then $[X]$ is not projective (because $\Gamma(A_0)$ is hereditary).
Hence we have an arrow $\tau[M] \to \tau[X]$ with $\tau[M] \in \mathcal{S}_i$, and $\mathcal{S}_i$ is a section.
Then either $\tau[X] \in \mathcal{S}_i$ (and in that case $[X] \in \mathcal{S}_{i+1}$) or
$\tau(\tau[X]) \in \mathcal{S}_i$ (in which case $\tau[X] \in \mathcal{S}_{i+1}$).

Since $\Gamma(A_0)$ has a finite number of vertices in this way can be constructed the sections
\[
 \mathcal{S}_0, \mathcal{S}_1, \ldots, \mathcal{S}_{n-1}, \mathcal{S}_n,
\]
where $\mathcal{S}_n$ consists only in isomorphism classes of indecomposable injective $A_0$-modules.\\
\underline{Step 3.} \textit{The sections $\mathcal{S}_0, \ldots, \mathcal{S}_n$ are a partition of the set
of vertices in $\Gamma(A_0)$}. We show inductively over $i$ that $\mathcal{S}_i \cap \mathcal{S}_j = \emptyset$ for
$0 \leq i < j \leq m$.  The base case $\mathcal{S}_0 \cap \mathcal{S}_j=\emptyset$
for $0<j\leq n$ is clear by definition. If $[M] \in \mathcal{S}_{i+1} \cap \mathcal{S}_{j}$
then $\tau[M] \in \mathcal{S}_i \cap \mathcal{S}_{j-1}$, which is an empty set by inductive hypothesis.

Observe finally that every $[M] \in \Ind A_0$ belongs to $\mathcal{S}_i$ for some $0 \leq i \leq n$.
For by the existence of projective covers there exist an indecomposable projective $P$ and a nonzero morphism
$P \to M$. By lemma~\ref{(P)L:KSfinita} there exists a path from $[P]$ to $[M]$ in $\Gamma(A_0)$. 
Hence $[M]$ belongs to the connected component which contains the projective modules, and by lemma~\ref{(DE)L:PosPre}, 
$[M]=\tau^{-i}[P']$ for some nonnegative integer $i$ and some indecomposable projective $P'$, that is, 
$[M] \in \mathcal{S}_i$.

\eproof

If $w$ is a vertex in the translation quiver $(\Gamma,\tau)$ denote by $\ConL(w)$ the set of vertices
$x$ in $\Gamma$ such that $x \prec w$ and by $\ConR(w)$ the set of vertices $y$ in $\Gamma$ such that 
$w \prec y$. Moreover, denote by $\dConL(w)$ the set of vertices $x$ in $\ConL(w)$ such that either $x$ is injective or 
$\tau^{-1}(x) \notin \{w\} \cup \ConL(w)$. In a similar way, denote by $\dConR(w)$ the set of vertices $y$ in $\ConR(w)$ 
such that either $y$ is projective or $\tau(y) \notin \{w\} \cup \ConR(w)$.

\begin{lema} \label{(DE)L:unicaSec}
Let $Q$ be a finite connected solid quiver with an admissible ordering of its vertices and let $\Gamma$ be a mesh complete 
connected translation subquiver of $\mathbb{Z}Q$. Let $w$ be a vertex of $\Gamma$.
\begin{itemize}
 \item[a)] The set $\mathcal{S}^L=\{w\} \cup \dConL(w)$ is the unique connected cosection contained in 
$\{w\} \cup \ConL(w)$ which contains the vertex $w$.
 \item[b)] The set $\mathcal{S}^R=\{w\} \cup \dConR(w)$ is the unique connected section contained in 
$\{w\} \cup \ConR(w)$ which contains the vertex $w$.
\end{itemize}
\end{lema}
\bproof
We show $(a)$, the proof of $(b)$ is similar. Notice that $\Gamma$ is a hereditary translation quiver 
by lemmas~\ref{(DE)L:mallas} and~\ref{(DE)L:orden}.  \\
\underline{Step 1.} \textit{We verify first that $\mathcal{S}^L$ is a cosection}.
Condition $(i)$ in the definition of cosection says that if $x \in \mathcal{S}^L$ and
$x$ is not injective, then $\tau^{-1}(x) \notin \mathcal{S}^L$. 
This condition is evident in the set $\mathcal{S}^L=\{w\} \sqcup\dConL(w)$ by its definition. 
In order to show condition $(ii)$ consider an arrow $x \to s$ with $s \in \mathcal{S}^L$. We have to show that
one of the vertices $x$ or $\tau^{-1}(x)$ belongs to $\mathcal{S}^L$.\\
\underline{Case $s=w$.} If $x$ is injective then $x \in \dConL(w)$. 
If $x$ is not injective we have that $\tau^{-1}(x) \notin \ConL(w)$ for the quiver $\Gamma$
is directed (lemmas~\ref{(DE)L:mallas} and~\ref{(DE)L:orden}). By definition of $\dConL(w)$ we have that 
$x \in \dConL(w) \subset \mathcal{S}^L$. \\
\underline{Case $s \neq w$.} By hypothesis $s \prec w$ and by transitivity $x \prec w$, that is, $x \in \ConL(w)$.
Assume that $x \notin \mathcal{S}^L$, thus $x$ is not injective and $\tau^{-1}(x) \in \ConL(w)$.
If $\tau^{-1}(x)$ is injecive then it belongs to $\mathcal{S}^L$ and we are done. If it is not injective
then $s$ is also not injective (for $\Gamma$ is hereditary) and $\tau^{-1}(s)$ does not belong to
$\ConL(w)$ (for $s \in \mathcal{S}^L$). Hence, since $\tau^{-1}(s) \prec \tau^{-2}(x)$,
we conclude that $\tau^{-2}(x) \notin \ConL(w)$ and then $\tau^{-1}(x) \in \mathcal{S}^L$.\\
\underline{Step 2.} \textit{We show connectedness}. 
If $\gamma:u \to w$ is a nontrivial path in $\Gamma$ which starts in $\dConL(w)$, then
all vertices in $\gamma$ different from $w$ are contained in the set $\dConL(w)$. Indeed, since 
$u \neq w$ the path $\gamma$ factors as $\gamma=\gamma' \alpha$, where $\alpha:u \to u'$ is an arrow and
$\gamma':u' \to w$ is a (possibly trivial) path. If $u'=w$ there is nothing to show. 
If $u' \neq w$ then $u' \in \ConL(w)$. If $u'$ is injective, all vertices in $\gamma'$ are also injective
(because $\Gamma$ is hereditary) and hence belong to $\dConL(w)$. If $u'$ is not injective, then $u$ is also not injective.
Then $\tau^{-1}(u')$ cannot belong to $\ConL(w)$, for in that case $\tau^{-1}(u') \prec w$, which is impossible because
$\tau^{-1}(u) \prec \tau^{-1}(u')$ and  $u \in \dConL(w)$. Hence $u'$ belongs to $\dConL(w)$. 
Successively applying this argument we have that every path that starts in $\dConL(w)$ and ends in $w$ 
is contained in $\dConL(w)$. In particular $\mathcal{S}^L$ is a connected cosection.\\
\underline{Step 3.} \textit{We finally show unicity}. Assume that $\mathcal{S}$ is a connected cosection in $\Gamma$ 
contained in $\{w\} \cup\ConL(w)$. By lemma~\ref{(DE)L:orbital} both $\mathcal{S}$ and $\mathcal{S}^L$ intersect each orbit 
in $\Gamma$ in exactly one vertex. By construction in $\mathcal{S}^L$, for any $s \in \mathcal{S}$ there is a nonnegative integer 
$h(s)$ such that $\tau^{-h(s)}(s) \in \mathcal{S}^L$. Let $s_w$ be the element in $\mathcal{S}$ such that
\[
 \tau^{-h(s_w)}(s_w)=w.
\]
We show by induction that $h(s) \leq h(s_w)$ for any $s \in \mathcal{S}$. 
For that purpose we order the elements $s_1^L,\ldots,s_n^L$ of $\mathcal{S}^L$ (and consequently the orbits in 
$\Gamma$ and the elements $s_1,\ldots,s_n$ in $\mathcal{S}$) in such a way that the order given to the full subquiver 
$Q_{\mathcal{S}^L}$ determined by $\mathcal{S}^L$ is admissible. Such an ordering exists because of lemmas~\ref{(DE)L:orden}
and~\ref{(DE)L:orbital}. Since $w$ is the unique sink in $\mathcal{S}^L$ we have that $w=s_1^L$ is the smallest element in
$\mathcal{S}^L$, and hence $s_1=s_w$ is the smallest element in $\mathcal{S}$ and $h(s_1)=h(s_w)$. 
Assume that the claim is valid for $s_1,\ldots,s_{\ell}$ with $\ell<n$. 
Since $\mathcal{S}^L$ is connected and $s_{\ell+1}^L$ is not a sink in $\mathcal{S}^L$, there exists $i \in \{1,\ldots,n\}$ 
such that $s^L_{\ell+1}$ is connected to $s^L_i$ by an arrow $s_{\ell+1}^L \to s_i^L$. 
Since the order in $Q_{\mathcal{S}^L}$ is admissible, $i \leq \ell$. Again by lemma~\ref{(DE)L:orbital} we have that 
$s_{\ell+1}$ is connected to $s_i$ inside of $\mathcal{S}$. \\
\underline{Case $s_{\ell+1} \to s_i$.} Since $\Gamma$ is hereditary and $\tau^{-h(s_i)}(s_i)$ is defined, 
there is an arrow $\tau^{-h(s_i)}(s_{\ell+1}) \to \tau^{-h(s_i)}(s_i)=s_i^L$. Since $s_i^L$ belongs to the cosection 
$\mathcal{S}^L$, either $\tau^{-h(s_i)}(s_{\ell+1}) \in \mathcal{S}^L$ or 
$\tau^{-h(s_i)-1}(s_{\ell+1}) \in \mathcal{S}^L$. In the second case there would be an arrow $s_i^L \to s_{\ell+1}^L$,
which is impossible for there is an arrow $s_{\ell+1}^L \to s_i^L$ and $\Gamma$ is directed. Then
$s_{\ell+1}^L=\tau^{-h(s_i)}(s_{\ell+1})$, that is, $h(s_{\ell+1})=h(s_i)$ and by inductive hypothesis
$h(s_{\ell+1})=h(s_i)\leq h(s_w)$.\\
\underline{Case $s_{i} \to s_{\ell+1}$.} In this case there is an arrow $\tau^{-h(s_{\ell+1})}(s_{i}) 
\to \tau^{-h(s_{\ell+1})}(s_{\ell+1})$, again because $\tau^{-h(s_{\ell+1})}(s_{\ell+1})=s_{\ell+1}^L$ is defined and 
$\Gamma$ is hereditary. Since $s^L_{\ell+1}$ is in the cosection $\mathcal{S}^L$, then either 
$\tau^{-h(s_{\ell+1})}(s_i) \in \mathcal{S}^L$ or $\tau^{-h(s_{\ell+1}-1)}(s_i) \in \mathcal{S}^L$.
Now, the first case is impossible for there would be an arrow $s^L_i \to s^L_{\ell+1}$.
Then $\tau^{-h(s_{\ell+1})-1}(s_i)=s_i^L$ and thus $h(s_{\ell+1})+1=h(s_i) \leq h(s_w)$, that is,
$h(s_{\ell+1}) < h(s_w)$, which completes the induction step.

Assume now that $\mathcal{S}$ is a connected cosection contained in $\{w\} \cup \ConL(w)$ which contains $w$. 
By the steps above, for any $s \in \mathcal{S}$ we have that $0 \leq h(s) \leq h(s_w)=0$, that is, 
$\mathcal{S}=\mathcal{S}^L$. This completes the proof.
\eproof

Given a set $\Lambda=\{1,\ldots,z\}$ and a function $r:\Lambda \to \mathbb{N}$ define the \textbf{star quiver}
(see Ringel~\cite{cmR}) $\mathbb{T}_r$ whose branches are indexed over $\Lambda$, the branch with index $\ell$ having length $r(\ell)$.
It is obtained from a disjoint union of copies of quivers of type $\mathbf{A}_{r(\ell)}$ ($\ell \in \Lambda$) choosing an end
point vertex for each $\mathbf{A}_{r(\ell)}$ and identifying all these points in a single vertex, the \textbf{center} of the star.
We usually denote $\mathbb{T}_r$ by $\mathbb{T}_{r(1),\ldots,r(z)}$ and call the set of integers $r(1),\ldots,r(z)$
\textbf{order} of the star. If $r(\ell)=1$ for some $\ell \in \Lambda$ then the corresponding branch
is invisible, thus we are interested in branches of length greater or equal than two.

We say that a quiver of tipe $\mathbf{A}_n$ is linearly oriented if all its arrows point in the same direction.
A \textbf{wing} of vertex $w_0$ in a translation quiver $(\Gamma,\tau)$ is a mesh complete translation subquiver $\theta(n)$
isomorphic to the Auslander-Reiten quiver $\Gamma(\mathbf{A_n})$ of a linearly oriente quiver $\mathbf{A_n}$ with $n\geq 2$, 
in which $w_0$ is the unique projective-injective vertex. Wings of order $3$ and $7$ are shown in the following figure.
%---------------------------------------------------------------------- 
\begin{displaymath}
\xy 0;/r.20pc/:
( -10,0)="X0" *{};
( -20,0)="X1" *{};
( -30,0)="X2" *{};
( -20,10)="Y0" *{};
( -17,10)="wY0" *{\scriptstyle w_0};
( -15,5)="Y1" *{};
( -25,5)="Y2" *{};
( -35,5)="theta3" *{\scriptstyle \theta(3)=};
( 0,0)="A0" *{};
(10,0)="A1" *{};
(20,0)="A2" *{};
(30,0)="A3" *{};
(40,0)="A4" *{};
(50,0)="A5" *{};
(60,0)="A6" *{};
(5,5)="B1" *{};
(10,10)="B2" *{};
(15,15)="B3" *{};
( 0,15)="theta3" *{\scriptstyle \theta(7)=};
(20,20)="B4" *{};
(25,25)="B5" *{};
(30,30)="B6" *{};
(33,30)="wB6" *{\scriptstyle w_0};
(35,25)="C1" *{};
(40,20)="C2" *{};
(45,15)="C3" *{};
(50,10)="C4" *{};
(55,5)="C5" *{};
%-----------------------
"X2";"Y0" **@{-};
"Y0";"X0" **@{-};
"Y1";"X1" **@{-};
"Y2";"X1" **@{-};
"A0";"B6" **@{-};
"A1";"C1" **@{-};
"A2";"C2" **@{-};
"A3";"C3" **@{-};
"A4";"C4" **@{-};
"A5";"C5" **@{-};
"A1";"B1" **@{-};
"A2";"B2" **@{-};
"A3";"B3" **@{-};
"A4";"B4" **@{-};
"A5";"B5" **@{-};
"A6";"B6" **@{-};
\endxy
\end{displaymath}
A vertex $x$ in $(\Gamma,\tau)$ is called \textbf{wing vertex} if every neighbor of $x$ belongs to a wing
of vertex $x$ in $(\Gamma,\tau)$.

\begin{lema} \label{(DE)L:estrella}
Let $Q$ be a finite solid connected quiver with an admissible ordering of its vertices. Then a vertex
$(x,i)$ in $\mathbb{Z}Q$ is a wing vertex if and only if $Q$ is isomorphic to a star quiver with vertex $x$.
\end{lema}

\bproof
By lemma~\ref{(DE)L:unicaSec} the subset $\mathcal{S}^L=\{(x,i)\} \cup \dConL((x,i))$ 
is a connected cosection in $\mathbb{Z}Q$. By lemma~\ref{(DE)L:orbital} the full subquiver $Q_{\mathcal{S}^L}$ 
determined by $\mathcal{S}^L$ has the same underlying graph as $Q$. If $(x,i)$ is a wing vertex then
$Q_{\mathcal{S}^L}$ is a star quiver (with all its arrows pointing towards the center of the star $(x,i)$),
hence $Q$ is also a star quiver with center $x$. For the other implication assume that $Q$ is
a star quiver $\mathbb{T}_{r(1),\ldots,r(z)}$ with center $x$. Let $i$ be an arbitrary integer. To a branch if index 
$\ell$ corresponds a full translation subquiver $\mathbb{Z}\mathbf{A_{r(\ell)}}$ in $\mathbb{Z}Q$,
corresponding to the quiver of type $\mathbf{A_{r(\ell)}}$, and corresponds hence a wing $\theta(r(\ell))$ of vertex $(x,i)$.
If $(y,j)$ is a neighbor of $(x,i)$ in $\mathbb{Z}Q$, then $y$ is neighbor of $x$ in $Q$, hence $y$ belongs
to some branch of $Q$, say of index $\ell$. Thus $(y,j)$ belongs to the wing $\theta(r(\ell))$ and
$(x,i)$ is a wing vertex.
\eproof

As can be directly observed every Dynkin diagram $\Delta$ is a star quiver. The order of the star, also called
\textbf{type} of the Dynkin diagram $\Delta$, is shown in the following table.
Only in cases $\mathbf{A}_n$ ($n\geq 2$) the type depends on the orientation of the arrows.
Fixing an end point $x$ of $\mathbf{A}_n$ one counts the number $p$ of arrows in $\mathbf{A}_n$ which
points towards $x$ and let $q$ be the number of arrows which points in the other direction. Assume that $p \geq q$.
Then, by definition, the type of $\mathbf{A}_n$ is $(p+1,q+1)$.

\begin{table} [!h] 
\begin{center}
\begin{tabular}{c | c}
Dynkin diagram & Type of the diagram \\
\hline 
 $\mathbf{A}_n$ & $\mathbb{T}_n$ \\
 $\mathbf{A}_{n_1+n_2-1}$ & $\mathbb{T}_{n_1,n_2}$  \\
 $\mathbf{D}_{n+2}$ & $\mathbb{T}_{n,2,2}$  \\
 $\mathbf{E}_6$ & $\mathbb{T}_{3,3,2}$  \\
 $\mathbf{E}_7$ & $\mathbb{T}_{4,3,2}$ \\
 $\mathbf{E}_8$ & $\mathbb{T}_{5,3,2}$  \\
\end{tabular}
\caption{Dynkin diagrams and their corresponding type. There is defined a center of star for cases
$\mathbf{D_n}$ ($n\geq 4$) and $\mathbf{E_m}$ ($m=6,7,8$). In case $\mathbf{A_n}$ any vertex can be considered as center.} 
\label{T:tipoDynkin}
\end{center}
\end{table}

Recall that an $A$-module $M$ is \textbf{generated} by the $A$-module $W$ if there exists an epimorphism
$\bigoplus_{i=1}^r W \to M$ for a finite direct sum of copies of $W$, and it is \textbf{cogenerated} by $W$ if there exists
a monomorphism $M \to \bigoplus_{i=1}^r W$.

\begin{lema} \label{(DE)L:alaDyn}
Let $A_0$ be the path algebra of a Dynkin diagram $\Delta$.
Let $W$ be an indecomposable sincere $A_0$-module, such that the isomorphism class $[W]$ belongs to the orbit in 
$\Gamma(A_0)$ corresponding to the center of star in $\Delta$ in cases $\mathbf{D_n}$ ($n\geq 4$) and $\mathbf{E_m}$
($m=6,7,8$). Then $[W]$ is a wing vertex in $\Gamma(A_0)$. 
Moreover, every $A_0$-module in $\ConL([W])$ is cogenerated by $W$ and every $A_0$-module in $\ConR([W])$ is generated by $W$. 
\end{lema}
\bproof
By lemma~\ref{(DE)L:partDyn} the quiver $\Gamma(A_0)$ is contained in $\mathbb{Z}\Delta^{op}$.
By hypothesis and lemma~\ref{(DE)L:estrella}, $[W]$ is a wing vertex in $\mathbb{Z}\Delta^{op}$.
It is then enough to prove that the wings of vertex $[W]$ in $\mathbb{Z}\Delta^{op}$ are contained in 
$\Gamma(A_0)$. This is a consequence of lemma~\ref{(P)L:KSfinita}, for $W$ is a sincere $A_0$-module and hence there exist 
nonzero morphisms from each indecomposable projective module $P$ to $W$. Thus there are paths from $[P]$ to $[W]$ 
contained in $\Gamma(A_0)$.

We want to prove now that every $A_0$-module in $\ConR([W])$ is generated by $W$ (the proof that every $A_0$-module in
$\ConL([W])$ is cogenerated by $W$ is dual). We proceed by steps. We show first that it is enough to verify the following
claim (\textasteriskcentered), then we show that the claim is valid in wing quivers and finally we use that
$[W]$ is a wing vertex to prove the claim in the general case.

\begin{itemize}
 \item[\textasteriskcentered)] Let $P$ be an indecomposable projective $A_0$-module and assume that 
there is a nonzero morphism $f:P \to N$ where $[N]$ is a vertex in $\ConR([W])$.
Then there exist morphisms $\xymatrix{P \ar[r]_-{h} & \bigoplus_{i=1}^rW \ar[r]_-{g} & N}$ such that $f=gh$.
\end{itemize} 
\underline{Step 1.} \textit{Assume that the claim (\textasteriskcentered) is true 
and that $N$ is an $A_0$-module such that $[N]$ is a vertex in $\ConR([W])$. Then $N$ is generated by $W$}. 
Consider a projective cover of $N$,
\[
\xymatrix@C=5pc{
\bigoplus_{j=1}^{z} P^j \ar[r]^-{f=\left[ \begin{smallmatrix} f_1 & \ldots & f_z \end{smallmatrix} \right] } & N,
}
\]
where each $P^j$ is indecomposable projective. By hypothesis the morphism $f_j$ factors through some direct sum 
$\bigoplus_{i=1}^{r_j}W$, say $f_j=g_jh_j$ (with $j \in \{1,\ldots,z\}$), thus if we define the following morphisms,
\[
 \xymatrix{
\bigoplus_{j=1}^{z} P^j \ar[rd]_(.3){h=\left[ \begin{smallmatrix} h_1 & & 0 \\  & {}^{\ddots} \\  0 & & h_z \end{smallmatrix} \right]} 
\ar[rr]^-{f=\left[ \begin{smallmatrix} f_1 & \ldots & f_z \end{smallmatrix} \right]} & & N \\
& \bigoplus_{j=1}^{z}\left(\bigoplus_{i=1}^{r_j} W \right), \ar[ru]_-{g=\left[ \begin{smallmatrix} g_1 & \ldots & g_z \end{smallmatrix} \right]}
}
\]
we have that $f=gh$. Since $f$ is an epimorphism, $g$ is also an epimorphism and hence $N$ is generated by $W$.\\
\underline{Step 2.} \textit{Consider a wing $\theta(n)$ of vertex $[W]$ in the Auslander-Reiten quiver $\Gamma(A_0)$. 
Then every morphism between representative of a projective and an injective vertex in $\theta(n)$, which is composition of
irreducible morphisms, factors through $W$}. In the following diagram we give a wing $\theta(n)$ of order $n$ and vertex $W=W_{1,n}$. 
For each vertex we have selected a class representative $W_{i,j}$ with $1 \leq i \leq j \leq n$. 
The difference $\ell=j-i$, which will be call level of the vertex $[W_{i,j}]$, corresponds to the $\ell$-th row in $\theta(n)$,
enumerating them from $0$ to $n-1$ from below to above. The level of a path $\gamma$ contained in $\theta(n)$
is the sum of the levels of all vertices in the path $\gamma$.
\[
 \xymatrix@!0@R=2.5pc@C=2.5pc{
& & & & W_{1,n} \ar[rd] \\
& & & W_{1,n-1} \ar[ru] \ar[rd] & & W_{2,n} \ar[rd]\\
& & W_{1,n-2} \ar[ru] & & W_{2,n-1} \ar@{}[rd]|{\cdots} \ar[ru] & & W_{3,n} \ar@{}[rd]|{\cdots} \\
& W_{1,2} \ar@{}[ru]|{\cdots} \ar[rd] & & W_{2,3} \ar@{}[rd]|{\cdots} \ar@{}[ru]|{\cdots}  & & W_{n-2,n-1} 
\ar@{}[ru]|{\cdots} \ar[rd]  & & W_{n-1,n} \ar[rd] \\
W_{1,1} \ar[ru] & & W_{2,2} \ar[ru] & & W_{i,i} \ar@{}[ru]|{\cdots} & & W_{n-1,n-1} \ar[ru] & & W_{n,n}.
}
\]
Assume that $f$ is a nonzero morphism $f:W_{1,i_0} \to W_{j_0,n}$ (for some $1 \leq i_0,j_0 \leq n$)
which is composition of irreducible maps $f_{\alpha}$ , and let $\gamma$ be the corresponding path.
Observe that if either $i_0=n$ or $j_0=1$ there is nothing to prove (for $W=W_{1,n}$), thus we can assume that 
$1 \leq i_0 < n$ and $1<j_0\leq n$. Assume that $[W_{i,j}]$ with $1<i\leq j < n$ is a vertex in $\gamma$ with minimal level.
By minimality $[W_{i,j}]$ belongs to a mesh of the form
\[
 \xymatrix@C=1pc@R=.5pc{
& [W_{i-1,j+1}] \ar[rd]^-{\beta'}  \\
[W_{i-1,j}] \ar[ru]^-{\alpha'} \ar[rd]_-{\alpha} & & [W_{i,j+1}] \\
& [W_{i,j}] \ar[ru]_-{\beta} 
}
\]
where the composition $\beta\alpha$ is a subpath of $\gamma$. Since $\theta(n)$ is a mesh complete translation subquiver
of the Auslander-Reiten quiver $\Gamma(A_0)$, to this mesh corresponds an almost split sequence in $A_0$-mod of the form,
\[
 \xymatrix{
0 \ar[r] & W_{i-1,j} \ar[r]^-{\left[ \begin{smallmatrix} f_{\alpha'} \\ f_{\alpha} \end{smallmatrix} \right]} 
& W_{i-1,j+1} \oplus W_{i,j} \ar[r]^-{\left[ \begin{smallmatrix} f_{\beta'} & -f_{\beta} \end{smallmatrix} \right]} 
& W_{i,j+1} \ar[r] & 0. 
}
\]
Since $f=f' f_{\beta}f_{\alpha}f''$ for some morphisms $f'$ and $f''$,
we have that $f=f' f_{\beta'}f_{\alpha'}f''$. Hence we can associate to $f$ a new path $\gamma'$
by exchanging the subpath $\beta \alpha$ for the subpath $\beta' \alpha'$.
Observe that the level $\mathscr{S}'$ of $\gamma'$ is strictly bigger than the level $\mathscr{S}$ of $\gamma$ 
(actually $\mathscr{S}'=\mathscr{S}+2$).

Then a path $\tilde{\gamma}:[W_{1,i_0}] \to [W_{j_0,n}]$ has maximal level among all other paths 
from $[W_{1,i_0}]$ to $[W_{j_0,n}]$ if and only if non of the vertices $[W_{i,j}]$ with 
$1 < i \leq j < n$ that belong to $\tilde{\gamma}$ have minimal level. 
This is clearly equivalent to the path $\tilde{\gamma}$ to be the only path from $[W_{1,i_0}]$
to $[W_{j_0,n}]$ which contains $[W]$. Then, by the change of factorization 
$f=f'f_{\beta}f_{\alpha}f''=f'f_{\beta'}f_{\alpha'}f''$ associated to the change of paths $\gamma \mapsto \gamma'$,
we conclude that to every morphism $f$ which is composition of irreducible maps 
corresponds a path $\tilde{\gamma}$ of maximal level. Hence $f$ factors through $W$.\\
\underline{Step 3.} \textit{We finally prove the claim (\textasteriskcentered)}.
Let $f:P \to N$ be a nonzero morphism with $P$ indecomposable projective and $[N] \in \ConR([W])$. By  
lemma~\ref{(P)L:KSfinita} there exist nonzero morphisms $f_1,\ldots,f_r:P \to N$ with $f=f_1+\cdots+f_r$ and such that
each $f_i$ is composition of irreducible maps. In particular to each $f_i$ corresponds a path $\gamma_i$ from $[P]$ to $[N]$.
%---------------------------------------------------------------------- 
\begin{displaymath}
\xy 0;/r.20pc/:
( 0,0)="A0" *{};
(20,0)="A2" *{};
(30,0)="A3" *{};
(40,0)="A4" *{};
(50,0)="A5" *{};
(60,0)="A6" *{};
(90,0)="A9" *{};
( 0,30)="B0" *{};
(30,30)="B3" *{};
(40,30)="B4" *{};
(50,30)="B5" *{};
(90,30)="B9" *{};
( 0,10)="C1" *{};
( 0,20)="C2" *{};
(90,10)="D1" *{};
(90,20)="D2" *{};
(25,5)="E1" *{};
(30,10)="E2" *{};
(35,15)="E3" *{};
(45,25)="E4" *{};
(35,25)="F1" *{};
(45,15)="F2" *{};
(50,10)="F3" *{};
(55,5)="F4" *{};
(46,20)="W" *{[W]};
(70,25)="CR" *{\ConR([W])};
(15,25)="CR" *{\ConL([W])};
(70,15)="N" *{[N]};
(-5,15)="N" *{[P]};
(40,35)="N" *{\text{Wing $\theta(n_2)$ of order $n_2=3$}};
(40,-5)="N" *{\text{Wing $\theta(n_1)$ of order $n_1=5$}};
%-----------------------
"A0";"A2" **@{-};
"A6";"A9" **@{-};
"A9";"B9" **@{-};
"A0";"B0" **@{-};
"B0";"B3" **@{-};
"B5";"B9" **@{-};
"B3";"A6" **@{-};
"B5";"A2" **@{-};
"B5";"A2" **@{-};
"E1";"A3" **@{-};
"E2";"A4" **@{-};
"E3";"A5" **@{-};
"E4";"B4" **@{-};
"F1";"B4" **@{-};
"F2";"A3" **@{-};
"F3";"A4" **@{-};
"F4";"A5" **@{-};
\endxy
\end{displaymath}
As shown in lemma~\ref{(DE)L:unicaSec} the set $\mathcal{S}^{L}=\{[W]\} \cup \dConL([W])$
is a connected cosection in $\Gamma(A_0)$ and $\mathcal{S}^{R}=\{[W]\} \cup \dConR([W])$ is a connected section in $\Gamma(A_0)$.
Since $\Gamma(A_0)$ is proper, directed and hereditary, by lemma~\ref{(DE)L:partSeccion} each path
$\gamma_i$ contains a vertex of the cosection $\mathcal{S}^L$ and a vertex of the section $\mathcal{S}^R$. 
Hence the path $\gamma_i$ goes in and out of a wing $\theta(n_i)$ of vertex $[W]$, and by step 2 
the morphism $f_i$ factors through $W$, say $f_i=g_ih_i$. Defining the morphisms $g$ and $h$ as in the following diagram 
we conclude that $gh=g_1h_1+\cdots+g_rf_r=f_1+\cdots+f_r=f$, 
\[
 \xymatrix{
P \ar[rd]_(.3){h=\left[ \begin{smallmatrix} h_1 \\{}^{\vdots}\\ h_r \end{smallmatrix} \right]} 
\ar[rr]^-{f} & & N \\
& \bigoplus_{i=1}^{r} W. \ar[ru]_-{g=\left[ \begin{smallmatrix} g_1 & \ldots & g_r \end{smallmatrix} \right]}
}
\]
This completes the proof.
\eproof

\begin{lema} \label{(DE)L:ala}
Let $W_0$ be an indecomposable $A_0$-module corresponding to the maximal positive root of the Dynkin quiver $\Delta$.
Consider the Auslander-Reiten quiver $\Gamma(A_0)$ of $A_0$-mod. Then $W_0$ is sincere
and in cases $\mathbf{D_n}$ ($n\geq 4$) and $\mathbf{E_m}$ ($m=6,7,8$) the isomorphism class
$[W_0]$ belongs to the orbit of the center of star $\Delta=\mathbb{T}_r$. 
\end{lema}

\bproof
The module $W_0$ is sincere, for its dimension vector $\vdim W_0$ is as described in table~\ref{T:dynkin}
(page~\pageref{T:dynkin}). The second claim is consequence of proposition 6.1 in Ringel~\cite{cmR}, 
where it is proved that if the orbit $[W_0]$ does not coincides with the center of star (that is, has less than three neighbors), 
then either $W_0$ is projective-injective or $W_0$ has two exceptional vertices. 
Both cases are impossible for the maximal exceptional representation $W_0$ in cases
$\mathbf{D_n}$ ($n\geq 4$) and $\mathbf{E_m}$ ($m=6,7,8$) (see table~\ref{T:dynkin}).
\eproof

Then the vertex $[W_0]$ corresponding to an indecomposable $A_0$-module $W_0$ of maximal dimension vector is a wing vertex
with at most three neighbor wings, $\theta(n_1),\ldots ,\theta(n_t)$ for $1 \leq t \leq 3$,
\[
 \xymatrix@!0@R=2.5pc@C=2.5pc{
& & & & W^s_{1,n_s} \ar[rd] \\
& & & W^s_{1,n_s-1} \ar[ru] \ar[rd] & & W^s_{2,n_s} \ar[rd]\\
& & W^s_{1,n_s-2} \ar[ru] & & W^s_{2,n_s-1} \ar@{}[rd]|{\cdots} \ar[ru] & & W^s_{3,n_s} \ar@{}[rd]|{\cdots} \\
& W^s_{1,2} \ar@{}[ru]|{\cdots} \ar[rd] & & W^s_{2,3} \ar@{}[rd]|{\cdots} \ar@{}[ru]|{\cdots}  & & W^s_{n_s-2,n_s-1} 
\ar@{}[ru]|{\cdots} \ar[rd]  & & W^s_{n_s-1,n_s} \ar[rd] \\
W^s_{1,1} \ar[ru] & & W^s_{2,2} \ar[ru] & & W^s_{i,i} \ar@{}[ru]|{\cdots} & & W^s_{n_s-1,n_s-1} \ar[ru] & & W^s_{n_s,n_s}.
}
\]

%----------------------------------------------------------------------
%----------------------------------------------------------------------
\chapter{Extended Dynkin diagrams.}
\label{Cap(DE)}
%----------------------------------------------------------------------
%----------------------------------------------------------------------

\begin{table} [!htb] 
\begin{center}
\renewcommand{\arraystretch}{1.2}
\begin{tabular}{c l} 
\hline 
Notation & \multicolumn{1}{c}{Graph} \\
\hline \\
\raisebox{2ex}{$\mathbf{\widetilde{A}_n} \; (n \geq 2)$} 
& $\xymatrix@R=-.5pc{
*+[o][F-]{{}_1} \ar@{-}[r] \ar@/^1pc/@{-}[rrrr]^-{\alpha'_0} & {}_1 \ar@{-}[r] & {}_1 
\cdots{}_1 \ar@{-}[r] & *+[o][F-]{{}_1} \ar@{-}[r]_-{\alpha_0} & {}_1 } $ \\
\raisebox{-1ex}{$\mathbf{\widetilde{D}_n} \; (n \geq 4)$} & $\xymatrix@R=-.5pc{{}_1 \ar@{-}[rd] & & & & {}_1 \\  
& {}_2 \ar@{-}[r] & {}_2 \ldots {}_2 \ar@{-}[r] & *+[o][F-]{{}_2} \ar@{-}[ru]^-{\alpha_0} \ar@{-}[rd] \\
 {}_1 \ar@{-}[ru] & & & & {}_1 }  $ \\
\raisebox{-4ex}{$\mathbf{\widetilde{E}_6}$} & $\xymatrix@R=1pc{ & & {}_1 \ar@{-}[d]^-{\alpha_0}  \\ 
& & *+[o][F-]{{}_2} \ar@{-}[d] \\
{}_1 \ar@{-}[r] & {}_2 \ar@{-}[r] & {}_3 \ar@{-}[r] & {}_2 \ar@{-}[r] & {}_1} $ \\
\raisebox{-2ex}{$\mathbf{\widetilde{E}_7}$} & $\xymatrix@R=1pc{ & & & {}_2 \ar@{-}[d] \\ 
{}_1 \ar@{-}[r]^-{\alpha_0} & *+[o][F-]{{}_2} \ar@{-}[r] & {}_3 \ar@{-}[r] & {}_4 \ar@{-}[r] & {}_3 \ar@{-}[r] & {}_2 \ar@{-}[r] & {}_1} $ \\
\raisebox{-2ex}{$\mathbf{\widetilde{E}_8}$} & $\xymatrix@R=1pc{ & & {}_3 \ar@{-}[d] \\ 
{}_2 \ar@{-}[r] & {}_4 \ar@{-}[r] & {}_6 \ar@{-}[r] & {}_5 \ar@{-}[r] & {}_4 \ar@{-}[r] & {}_3 \ar@{-}[r] 
&*+[o][F-]{{}_2} \ar@{-}[r]^-{\alpha_0} & {}_1 \\ {} } $ \\
%\vspace{.5cm}
\hline
\end{tabular}
\caption{Extended Dynkin (simply laced) quivers $\widetilde{\Delta}$. The orientation of the arrows is arbitrary. 
The Dynkin quiver $\Delta$ from which $\widetilde{\Delta}$ is extended is obtainded deleting the arrow(s) 
$\alpha_0$ (and $\alpha'_0$ in cases $\mathbf{A_n}$) and the remaining isolated vertex.
The exceptional vertices of the maximal positive root $w_0$ are marked in circle.
The integral vector given is generator of the radical of the quadratic form associated to $\widetilde{\Delta}$.} 
\label{(DE)T:dynkinE}
\end{center}
\end{table}

\section{One-point extension and reduction.} \label{(DE)S:extPunt}
%----------------------------------------------------------------------
In this section we perform reductions by admissible modules over one-point extensions,
as defined by Ringel in \cite{cmR}.

%------------------------------------------------------------------
%------------------------------------------------------------------
\textbf{Definitions and equivalences.} 
Let $A_0$ be a $k$-algebra with identity $e_0$ and $R$ an $A_0$-module, both finite dimensional over $k$.
The \textbf{one-point extension} $A_0[R]$ of $A_0$ over $R$ is given by the set of matrices of the form 
$\left( \begin{smallmatrix} a & r \\ 0 & \lambda \end{smallmatrix} \right)$
with $a \in A_0$, $r \in R$ and $\lambda \in k$, with matrix operations.
We identify the elements $a$ in $A_0$ with the matrices 
$\left( \begin{smallmatrix} a & 0 \\ 0 & 0 \end{smallmatrix} \right)$ in $A_0[R]$ and denote by $e_{\omega}$ the idempotent
$\left( \begin{smallmatrix} 0 & 0 \\ 0 & 1 \end{smallmatrix} \right)$. There are then identifications
\[
 A_0 = e_0(A_0[R])e_0, \qquad k = e_{\omega}(A_0[R])e_{\omega} \qquad  \text{and} \qquad R = e_0(A_0[R])e_{\omega}.
\]
Consider the \textbf{subspace category} $\check{\mathcal{U}}(A_0\text{-mod},|\text{-}|)$ of $A_0$-mod 
through the $k$-functor $|\text{-}|:A_0\text{-mod} \to k\text{-mod} $.
Its objects are given by triples 
\[
M=(M_0,M_{\omega},\gamma_M),
\]
where $M_0$ is an $A_0$-module, $M_{\omega}$ is a $k$-vector space and $\gamma_M$ is a linear
transformation $M_{\omega} \to |M_0|$. A morphism $f=(f_0,f_{\omega}):M \to N$ is given by a linear transformation 
$f_{\omega}:M_{\omega} \to N_{\omega}$ and a morphism of $A_0$-modules $f_0:M_0 \to N_0$ such that 
$\gamma_{N}f_{\omega}=|f_0|\gamma_M$, that is, the following diagram commutes,
\[
 \xymatrix@C=4pc{
M_{\omega} \ar[r]^-{f_{\omega}} \ar[d]_-{\gamma_M} & N_{\omega} \ar[d]^-{\gamma_N} \\
|M_0| \ar[r]_-{|f_0|} & |N_0|.
}
\]
\begin{lema} \label{(DE)L:equiv}
There is an equivalence of categories 
\[
\xymatrix{
A_0[R]\text{-mod} \ar[r]^-{G} & \check{\mathcal{U}}(A_0\text{-mod},\Hom_{A_0}(R,-)).
}
\]
\end{lema}

\bproof
If $M$ is an $A_0[R]$-module then the vector spaces $M_0=e_0M$ and $M_{\omega}=e_{\omega}M$ decompose $M$ as direct sum 
$M=M_0 \oplus M_{\omega}$. Notice that $M_0$ is an $A_0$-module by the identification $A_0 = e_0(A_0[R])e_0$.
The action of $A_0[R]$ in $M$ is given by
\[
 \left(\begin{matrix} a & r \\ 0 & \lambda \end{matrix} \right) \left(\begin{matrix} m_0 \\ m_{\omega} \end{matrix} \right)=
\left(\begin{matrix} am_0+rm_{\omega} \\ \lambda m_{\omega} \end{matrix} \right).
\]
For each $m=e_{\omega}m \in M_{\omega}$ consider the mapping $R=e_0Re_{\omega} \to M_0$ given by
\[
 \gamma_{M}(m):e_0re_{\omega} \mapsto e_0re_{\omega}m.
\]
Clearly $\gamma_{M}(m)$ is a morphism of $A_0$-modules and in fact the mapping $\gamma_M$ is $k$-linear.
We define in this way $G(M)=(M_0,M_{\omega},\gamma_M)$.
If $f:M \to N$ is a morphism of $A_0[R]$-modules then $f_0=e_0f:M_0 \to N_0$ is a morphism of $A_0$-modules and 
$f_{\omega}=e_{\omega}f:M_{\omega} \to N_{\omega}$ is a linear transformation such that the following diagram
is commutative
\[
 \xymatrix@C=4pc{
M_{\omega} \ar[r]^-{f_{\omega}} \ar[d]_-{\gamma_M} & N_{\omega} \ar[d]^-{\gamma_N} \\
\Hom_{A_0}(R,M_0) \ar[r]_-{\Hom(1,f_0)} & \Hom_{A_0}(R,N_0),
}
\]
for if $m=e_{\omega}m \in M_{\omega}$ and $r=e_0re_{\omega} \in R$ then  
\[
 \gamma_N(f_{\omega}(m))(r)=e_0re_{\omega}f(m)=e_0f(e_0re_{\omega}m)=f_0\circ \gamma_M(m)(r).
\]
Hence $G(f)=(f_0,f_{\omega})$ is a functor, since
\[
 G(f\circ g)=(e_0f\circ g,e_{\omega}f \circ g)=(e_0f,e_{\omega}f)\circ (e_0g,e_{\omega}g)=G(f)\circ G(g),
\]
and $G(I_M)=I_{G(M)}$.
The inverse of $G$ sends each object to the vector space $M=M_0 \oplus M_{\omega}$ with action of
$A_0[R]$ given by
\[
 \left( \begin{matrix}a & r \\ 0 & \lambda \end{matrix} \right)\left( \begin{matrix} m_0 \\ m_{\omega} \end{matrix} \right)=
 \left( \begin{matrix} am_0 + \gamma_M(m_{\omega})r \\ \lambda m_{\omega} \end{matrix} \right),
\]
and sends a morphism $(f_0,f_{\omega}):(M_0,M_{\omega},\gamma_{M}) \to (N_0,N_{\omega},\gamma_{N})$ to the mapping
\[
e_0f_0 \oplus e_{\omega}f_{\omega}:e_0M_0 \oplus e_{\omega}M_{\omega} \to e_0N_0 \oplus e_{\omega}N_{\omega}.
\] 
\eproof

Since the functor $|-|=\Hom_{A_0}(R,-)$ will not be changed throughout the discussion, we will omit the reference to it.
Given a module class $\mathcal{Z}$ in $A_0$-mod, we will denote by
$\check{\mathcal{U}}(\mathcal{Z})$ the full subcategory of $\check{\mathcal{U}}(A_0\text{-mod})$
whose objects $M=(M_0,M_{\omega},\gamma_M)$ satisfy that $M_0$ belongs to $\mathcal{Z}$.

Define the \textbf{tensor one-point extension} of $A_0=T_{S_0}(L'_0)$ over the $S_0$-$k$-bimodule $L''_0$ as the tensor algebra
$T_S(L)$, where $S$ is the $k$-algebra $\left( \begin{smallmatrix} S_0 & 0 \\ 0 & k \end{smallmatrix} \right)$ with 
matrix operations and $L$ is the $S$-$S$-bimodule $\left( \begin{smallmatrix} L'_0 & L''_0 \\ 0 & 0 \end{smallmatrix} \right)$ 
with actions of $S$ given by the matrix multiplication. 

\begin{lema} \label{(DE)L:prod}
Let $S$ be the $k$-algebra $\left( \begin{smallmatrix} S_0 & 0 \\ 0 & k \end{smallmatrix} \right)$ and assume 
that $A$ and $L'$ are $S_0$-$S_0$-bimo\-dules, that $B$ and $L''$ are $S_0$-$k$-bimodules
and that $A'$ is a $k$-vector space. Then there exist isomorphisms of $S$-$S$-bimodules 
\[
 \varphi: \left( \begin{matrix} A & B \\ 0 & A' \end{matrix} \right) \otimes_S 
\left( \begin{matrix} L' & L'' \\ 0 & 0 \end{matrix} \right)
 \longrightarrow \left( \begin{matrix} A \otimes_{S_0} L' & A \otimes_{S_0}L'' \\ 0 & 0 \end{matrix} \right), 
\]
\[
 \psi: \left( \begin{matrix} L' & L'' \\ 0 & 0 \end{matrix} \right) \otimes_S 
\left( \begin{matrix} A & 0 \\ 0 & A' \end{matrix} \right)
 \longrightarrow \left( \begin{matrix} L' \otimes_{S_0} A & L'' \otimes_{k}A' \\ 0 & 0 \end{matrix} \right).
\]
\end{lema}
\bproof
We observe that there are isomorphisms
\[
a) \left( \begin{smallmatrix} A & 0 \\ 0 & 0 \end{smallmatrix} \right) \otimes_S 
\left( \begin{smallmatrix} L' & 0 \\ 0 & 0 \end{smallmatrix} \right)
 \cong \left( \begin{smallmatrix} A \otimes_{S_0} L' & 0 \\ 0 & 0 \end{smallmatrix} \right), 
\]
\[
b) \left( \begin{smallmatrix} A & 0 \\ 0 & 0 \end{smallmatrix} \right) \otimes_S 
\left( \begin{smallmatrix} 0 & L'' \\ 0 & 0 \end{smallmatrix} \right)
 \cong \left( \begin{smallmatrix}  0 & A \otimes_{S_0} L'' \\ 0 & 0 \end{smallmatrix} \right), 
\]
\[
c) \left( \begin{smallmatrix} 0 & L'' \\ 0 & 0 \end{smallmatrix} \right) \otimes_S 
\left( \begin{smallmatrix} 0 & 0 \\ 0 & A' \end{smallmatrix} \right)
 \cong \left( \begin{smallmatrix}  0 & L'' \otimes_k A' \\ 0 & 0 \end{smallmatrix} \right). 
\]
We show point $(a)$ as example. We notice first that, if $\times $ denotes the catesian product,
there is a bijection between the sets
\[
 \left( \begin{smallmatrix} A & 0 \\ 0 & 0 \end{smallmatrix} \right) \times 
\left( \begin{smallmatrix} L' & 0 \\ 0 & 0 \end{smallmatrix} \right)
 \longrightarrow \left( \begin{smallmatrix} A \times L' & 0 \\ 0 & 0 \end{smallmatrix} \right)
\]
\[
 \left[ \left( \begin{smallmatrix} a & 0 \\ 0 & 0 \end{smallmatrix} \right),
\left( \begin{smallmatrix} \ell & 0 \\ 0 & 0 \end{smallmatrix} \right) \right]
 \mapsto \left( \begin{smallmatrix} (a,\ell) & 0 \\ 0 & 0 \end{smallmatrix} \right). 
\]
This bijection extends to an isomorphism between the $k$-vector spaces generated by the corresponding sets 
$F_k \left[\left( \begin{smallmatrix} A & 0 \\ 0 & 0 \end{smallmatrix} \right) \times 
\left( \begin{smallmatrix} L' & 0 \\ 0 & 0 \end{smallmatrix} \right) \right]
\longrightarrow F_k \left( \begin{smallmatrix} A \times L' & 0 \\ 0 & 0 \end{smallmatrix} \right)$.
Notice that this isomorphism establishes a correspondence between generators of the ideal for the tensor product
\[
\left[\left( \begin{smallmatrix} a & 0 \\ 0 & 0 \end{smallmatrix} \right)
\left( \begin{smallmatrix} s & 0 \\ 0 & \lambda \end{smallmatrix} \right),
\left( \begin{smallmatrix} \ell & 0 \\ 0 & 0 \end{smallmatrix} \right) \right]-
\left[\left( \begin{smallmatrix} a & 0 \\ 0 & 0 \end{smallmatrix} \right),
\left( \begin{smallmatrix} s & 0 \\ 0 & \lambda \end{smallmatrix} \right)
\left( \begin{smallmatrix} \ell & 0 \\ 0 & 0 \end{smallmatrix} \right) \right]
\mapsto \left( \begin{smallmatrix} (as,l) & 0 \\ 0 & 0 \end{smallmatrix} \right)-
\left( \begin{smallmatrix} (a,sl) & 0 \\ 0 & 0 \end{smallmatrix} \right),
\]
with $a \in A$, $\ell \in L'$, $s \in S_0$ and $\lambda \in k$. It induces an isomorphism in the quotient.
The isomorphism $\varphi$ in the statement can be established in the following way. First we use 
$(a)$, $(b)$ and distributivity of tensors to give isomorphisms
\begin{eqnarray}
 \left( \begin{smallmatrix} A & 0 \\ 0 & 0 \end{smallmatrix} \right) \otimes_S
\left( \begin{smallmatrix} L' & L'' \\ 0 & 0 \end{smallmatrix} \right) & \cong & 
\left( \begin{smallmatrix} A & 0 \\ 0 & 0 \end{smallmatrix} \right) \otimes_S
 \left[ \left( \begin{smallmatrix} L' & 0 \\ 0 & 0 \end{smallmatrix} \right) \oplus
\left( \begin{smallmatrix} 0 & L'' \\ 0 & 0 \end{smallmatrix} \right) \right] \cong \nonumber \\
& \cong & \left[ \left( \begin{smallmatrix} A & 0 \\ 0 & 0 \end{smallmatrix} \right) \otimes_S
 \left( \begin{smallmatrix} L' & 0 \\ 0 & 0 \end{smallmatrix} \right) \right] \oplus
\left[ \left( \begin{smallmatrix} A & 0 \\ 0 & 0 \end{smallmatrix} \right) \otimes_S
\left( \begin{smallmatrix} 0 & L'' \\ 0 & 0 \end{smallmatrix} \right) \right] \cong \nonumber \\
& \cong & \left( \begin{smallmatrix} A \otimes_{S_0} L' & 0 \\ 0 & 0 \end{smallmatrix} \right) \oplus
\left( \begin{smallmatrix} 0 & A \otimes_{S_0} L'' \\ 0 & 0 \end{smallmatrix} \right)= \nonumber \\
& = & \left( \begin{smallmatrix} A \otimes_{S_0} L' & A \otimes_{S_0} L'' \\ 0 & 0 \end{smallmatrix} \right). \nonumber
\end{eqnarray}
Since $\left( \begin{smallmatrix} A & B \\ 0 & A' \end{smallmatrix} \right) 
\left( \begin{smallmatrix} 1 & 0 \\ 0 & 0 \end{smallmatrix} \right)=
\left( \begin{smallmatrix} A & 0 \\ 0 & 0 \end{smallmatrix} \right)$ and
$\left( \begin{smallmatrix} 1 & 0 \\ 0 & 0 \end{smallmatrix} \right) 
\left( \begin{smallmatrix} L' & L'' \\ 0 & 0 \end{smallmatrix} \right)=
\left( \begin{smallmatrix} L' & L'' \\ 0 & 0 \end{smallmatrix} \right)$ we have that
\begin{eqnarray}
 \left( \begin{smallmatrix} A & B \\ 0 & A' \end{smallmatrix} \right) \otimes_S
\left( \begin{smallmatrix} L' & L'' \\ 0 & 0 \end{smallmatrix} \right) & = &
\left( \begin{smallmatrix} A & B \\ 0 & A' \end{smallmatrix} \right) \otimes_S
\left( \begin{smallmatrix} 1 & 0 \\ 0 & 0 \end{smallmatrix} \right)
\left( \begin{smallmatrix} L' & L'' \\ 0 & 0 \end{smallmatrix} \right) = \nonumber \\
& = & \left( \begin{smallmatrix} A & B \\ 0 & A' \end{smallmatrix} \right)
\left( \begin{smallmatrix} 1 & 0 \\ 0 & 0 \end{smallmatrix} \right) \otimes_S
\left( \begin{smallmatrix} L' & L'' \\ 0 & 0 \end{smallmatrix} \right)= \nonumber \\
& = & \left( \begin{smallmatrix} A & 0 \\ 0 & 0 \end{smallmatrix} \right) \otimes_S
\left( \begin{smallmatrix} L' & L'' \\ 0 & 0 \end{smallmatrix} \right) \cong  \nonumber \\
& \cong & \left( \begin{smallmatrix} A \otimes_{S_0} L' & A \otimes_{S_0} L'' \\ 0 & 0 \end{smallmatrix} \right). \nonumber
\end{eqnarray}
In a similar way we get the isomorphism $\psi$.
\eproof 

\begin{lema} \label{(DE)L:extTens}
If $A_0$ is a tensor algebra of the form $T_{S_0}(L'_0)$ then the one-point extension 
$A_0[R]$ by the projective $A_0$-m\'odulo $R=A_0 \otimes_{S_0} L''_0$ is isomorphic to the tensor algebra 
$T_{\left( \begin{smallmatrix} S_0 & 0 \\ 0 & k \end{smallmatrix} \right)}\left( \begin{matrix} L'_0 & L''_0 
\\ 0 & 0 \end{matrix} \right)$.
\end{lema}

\bproof
Take $S=\left( \begin{smallmatrix} S_0 & 0 \\ 0 & k \end{smallmatrix} \right)$ and consider the $S$-$S$-bimodule
\[
L=\left( \begin{matrix} L'_0 & L''_0 \\ 0 & 0 \end{matrix} \right)
\cong \left( \begin{matrix} L'_0 & S_0 \otimes_{S_0}L''_0 \\ 0 & 0 \end{matrix} \right).
\]
By the isomorphism $\varphi$ in the last lemma we have that
\[
 L^{\otimes 2}=L\otimes_S L \cong \left( \begin{matrix} L'_0\otimes_{S_0} L'_0 & L'_0 
\otimes_{S_0}  L''_0 \\ 0 & 0 \end{matrix} \right),
\]
and inductively for $m \geq 1$
\begin{eqnarray}
 L^{\otimes (m+1)} & = & L \otimes_S L^{\otimes m} \cong
\left( \begin{matrix} L'_0 &  L''_0 \\ 0 & 0 \end{matrix} \right) \otimes_S
\left( \begin{matrix} {L'_0}^{\otimes m} & {L'_0}^{\otimes (m-1)} \otimes_{S_0}  L''_0 \\ 0 & 0 \end{matrix} \right) \cong \nonumber \\
& \cong & 
\left( \begin{matrix} {L'_0}^{\otimes (m+1)} & {L'_0}^{\otimes m} \otimes_{S_0}  L''_0 \\ 0 & 0 \end{matrix} \right). \nonumber
\end{eqnarray}
In this way, considering that ${L'_0}^{\otimes 0}=S_0$, we have
\begin{eqnarray}
 T_S(L) & \cong & S \oplus \bigoplus_{m \geq 1}L^{\otimes m} \cong  \nonumber \\
& \cong & \left( \begin{matrix} S_0 & 0 \\ 0 & k \end{matrix} \right) \oplus
\left( \begin{matrix} \bigoplus_{m \geq 1}{L'_0}^{\otimes m} 
& \left[ \bigoplus_{m \geq 1}{L'}_0^{\otimes (m-1)} \right] \otimes_{S_0} L''_0 \\ 0 & 0 \end{matrix} \right) = \nonumber \\
& = & \left( \begin{matrix} T_{S_0}(L'_0) & \left[ T_{S_0}(L'_0) \right] \otimes_{S_0} L''_0 
\\ 0 & k \end{matrix} \right) = \nonumber \\
& = & A_0[R]. \nonumber
\end{eqnarray}
\eproof

The proof of the following lemma is direct.
\begin{lema} \label{(DE)L:epsilonBimod}
Let $A_0$ and $R$ be as in the last lemma and assume that $X_0$ is an admissible $A_0$-modulo,
with splitting $(Z_0,P_0)$ of its opposed endomorphism algebra.
Then there is a morphism $\varepsilon$ of right $Z_0$-modules,
\[
 \xymatrix@R=.5pc{
\Hom_{A_0}(R,X_0) \otimes_{Z_0} P_0 \ar[r]^-{\varepsilon} & \Hom_{A_0}(R,X_0) \\
h\otimes p \ar@{|->}[r] & h p=p\circ h.
}
\]
Its right $Z_0$-dual $\varepsilon^*$ is given by
\[
 \xymatrix@R=.5pc{
\Hom_{A_0}(R,X_0)^*  \ar[r]^-{\varepsilon^*} & (\Hom_{A_0}(R,X_0) \otimes_{Z_0} P_0)^* \\
F \ar@{|->}[r] & [h \otimes p \mapsto F(h p)].
}
\]
\end{lema}

\begin{lema} \label{(DE)L:isoDual}
With the notation of the last lemma, if $R$ is a finitly generated projective $A_0$-module, 
then the mapping
\[
 \xymatrix@R=.5pc{
X_0^* \otimes_{A_0} R \ar[r]^-{\eta} & \Hom_{A_0}(R,X_0)^* \\
g \otimes r  \ar@{|->}[r] & [f \mapsto g(f(r))],
}
\]
is an isomorphism of left $Z_0$-modules.
\end{lema}

\bproof
Clearly the last mapping is a well defined morphism.
We will show the claim for the case $R=A_0e$. The general case $R=\bigoplus_id_i(A_0e_i)$ follows from additivity.
Observe first that the dual of the isomorphism 
\[
 \xymatrix@R=.5pc{
\Hom_{A_0}(A_0e,X_0) \ar[r] & eX_0 \\
f \ar@{|->}[r] & f(e),
}
\]
is given by
\[
 \xymatrix@R=.5pc{
(eX_0)^* \ar[r] & \Hom_{A_0}(A_0e,X_0)^* \\
g \ar@{|->}[r] & [f \mapsto g(f(e))].
}
\]
We follow the isomorphisms
\[
X_0^* \otimes_{A_0} R = X_0^* \otimes_{A_0} A_0e \cong X_0^*e \cong (eX_0)^* \cong \Hom_{A_0}(R,X_0)^*. 
\]
Taking $r=e$ we have
\[
 g \otimes e = g \otimes e \mapsto g \cdot e \mapsto [f \mapsto g(f(e))],
\]
that is, the morphism given in the statement is an isomorphism.
\eproof

Let $\psi:P_0^* \otimes_{Z_0} \Hom_{A_0}(R,X_0)^* \longrightarrow (\Hom_{A_0}(R,X_0) \otimes_{Z_0} P_0)^*$ 
be the natural isomorphism given by $\psi(\gamma \otimes F)[h \otimes p]=\gamma(F(h)p)$. 
See lemma~\ref{(A)L:psiBimod} in the appendix for a description of the inverse $\psi^{-1}$
in terms of the dual bases of $\Hom_{A_0}(R,X_0)$ and $P_0$.

\begin{proposicion} \label{(DE)P:extPuntRed}
Let $A_0$ be a tensor algebra of the form $T_{S_0}(L'_0)$ and $R$ a projective $A_0$-module $R=A_0 \otimes_{S_0} L''_0$.
Let $A$ be the one-point extension $A_0[R]$ and take the subalgebra $B$ of $A$ given by $B=A_0[0]$.
Consider the ditalgebra $\mathcal{A}=(A,0)$ with zero differential and assume that $X_0$ is
an admissible $A_0$-module, with splitting $(Z_0,P_0)$ of its opposed endomorphism algebra.
Let $S(\omega)$ be the simple $B$-module corresponding to the extension vertex.
Then the $B$-module $X=X_0 \oplus S(\omega)$ is admissible and the reduced ditalgebra $\mathcal{A}^X$
is isomorphic to the ditalgebra $(A_0^{X_0}[R^{X_0}],\delta^x)$, where $R^{X_0}$ is the $A_0^{X_0}$-module given by
$A_0^{X_0}\otimes_{Z_0} \Hom_{A_0}(R,X_0)^*$ and the differential $\delta^x$ is determined by the differential
$\delta^{X_0}$ in $A_0^{X_0}$ and the transformation
\[
 \widehat{\delta}=\psi^{-1}\varepsilon^*: \Hom_{A_0}(R,X_0)^* \longrightarrow P_0^* \otimes_{Z_0} \Hom_{A_0}(R,X_0)^*.
\]
Assume that there are finite dual bases $\{(p_j,\gamma_j)\}_{j \in J}$ of $P_0$ and $\{(a_i,\lambda_i)\}_{i \in I}$
of $\Hom_{A_0}(R,X_0)$. For an element $H \in \Hom_{A_0}(R,X_0)^*$ the transformation $\widehat{\delta}$ 
has the explicit form
\[
 \widehat{\delta}(H)= \sum_{i \in I,j \in J}H(a_ip_j)\gamma_j \otimes \lambda_i.
\]
\end{proposicion}
\bproof
Observe that there is a decomposition of the opposed endomorphism algebra
\[
 \End_B(X)^{op}= \left( \begin{matrix} \End_{A_0}(X_0)^{op} & 0\\0 &\End_k(S(\omega))^{op} \end{matrix} \right) \cong
\left( \begin{matrix} Z_0 & 0\\0 & k \end{matrix} \right) \oplus \left( \begin{matrix} P_0 & 0\\0 &0 \end{matrix} \right),
\]
which determines an splitting $(Z=\left( \begin{smallmatrix} Z_0 & 0\\0 & k \end{smallmatrix} \right),
P=\left( \begin{smallmatrix} P_0 & 0\\0 & 0 \end{smallmatrix} \right))$ of the opposed endomorphism algebra of $X$.
Observe that the $Z$-module $X_Z$ has a finite dual basis, for $(X_0)_{Z_0}$ has one as 
$Z_0$-module. Hence $X$ is an admissible module. Moreover, the module $Z$ is semi-simple if $Z_0$ is semi-simple and $Z$
is trivial if $Z_0$ is trivial.

As shown in the proof of lemma~\ref{(DE)L:extTens}, there exists an isomorphism
\[
 A_0[R] \cong T_S(L), \qquad \text{with $L=L' \oplus L''$},
\]
where $S=\left( \begin{smallmatrix} S_0 & 0\\0 & k \end{smallmatrix} \right)$,
$L'=\left( \begin{smallmatrix} L'_0 & 0\\0 & 0 \end{smallmatrix} \right)$ and 
$L''=\left( \begin{smallmatrix} 0 & L''_0\\0 & 0 \end{smallmatrix} \right)$.
It is also clear that $B=T_S(L')$. Following definitions~\ref{(A)D:subditalg} and~\ref{(A):tensRed} 
in the appendix we have that
\[
 \underline{L}=BL''B=\left( \begin{matrix} A_0 & 0\\0 & k \end{matrix} \right)
\left( \begin{matrix} 0 & L''_0\\0 & k \end{matrix} \right)
\left( \begin{matrix} A_0 & 0\\0 & k \end{matrix} \right) \cong
\left( \begin{matrix} 0 & R\\0 & 0 \end{matrix} \right),
\]
and using lemma~\ref{(DE)L:isoDual},
\begin{eqnarray}
 L_0^X & = & X^*\otimes_B \underline{L} \otimes_B X\cong \left( \begin{matrix} X_0^* & 0\\0 & S(\omega)^* \end{matrix} \right)
\otimes_B \left( \begin{matrix} 0 & R\\0 & 0 \end{matrix} \right) \otimes_B
\left( \begin{matrix} X_0 & 0\\0 & S(\omega) \end{matrix} \right) \cong \nonumber \\
& \cong & \left( \begin{matrix} 0 & X_0^* \otimes_{A_0} R \otimes_k S(\omega)\\0 & 0 \end{matrix} \right) 
\cong \left( \begin{matrix} 0 & \Hom_{A_0}(R,X_0)^*\\0 & 0 \end{matrix} \right), \nonumber
\end{eqnarray}
and $L_1^X=P^* \cong \left( \begin{matrix} P_0^* & 0\\0 & 0 \end{matrix} \right)$. The reduced tensor algebra
$A^X$ is defined as
\[
 A^X=T_Z(L_0^X\oplus L_1^X)\cong T_{\left( \begin{smallmatrix} Z_0 & 0\\0 & k \end{smallmatrix} \right)}
\left( \begin{matrix} P_0^* & \Hom_{A_0}(R,X_0)^* \\ 0& 0 \end{matrix} \right),
\]
and again by lemma~\ref{(DE)L:extTens}, noticing that $A_0^{X_0}=T_{Z_0}(P_0^*)$ and recalling the notation
$R^{X_0}=A_0^{X_0} \otimes_{Z_0} \Hom_{A_0}(R,X_0)^*$, we have that $A^X \cong A_0^{X_0}[R^{X_0}]$.

Now, the reduced differential $\delta^X$ is determined by the morphisms of $S$-$S$-bimodules
\[
 [X^* \otimes_B \underline{L} \otimes_B X] \oplus P^* \cong 
\left( \begin{matrix} P_0^* & X_0^* \otimes_{A_0} R \otimes_k S(\omega) \\ 0& 0 \end{matrix} \right) \longrightarrow A^X,
\]
given by the coproduct in $P$ and the mapping established in the appendix~\ref{(A)L:diferen},
\[
\delta^X(g \otimes r \otimes x^{\omega})= \lambda(g) \otimes r \otimes w 
+\sigma_{g,x^{\omega}}(\delta(r))+(-1)^{|r|+1}g \otimes r \otimes \rho(x^{\omega}),
\]
where $x^{\omega}$ is a generator of the simple $S(\omega)$, $r \in R$, $g \in X_0^*$ and $\lambda$, $\rho$ are the left and
right coactions of $X$ (definition~\ref{(A)D:coacciones}). Observe that the morphism $m_r$ in the definition of
$\rho$ is null in $S(\omega)$, since for any $p$ in $P_0$,
\[
 m_r\left( \begin{matrix} 0 & 0\\0 & x^{\omega} \end{matrix} \right)
\left[ \left( \begin{matrix} p & 0\\0 & 0 \end{matrix} \right) \right] =
\left( \begin{matrix} p & 0\\0 & 0 \end{matrix} \right)
\left( \begin{matrix} 0 & 0\\0 & x^{\omega} \end{matrix} \right)=
\left( \begin{matrix} 0 & 0\\0 & 0 \end{matrix} \right).
\]
Then, since the differential $\delta$ of the ditalgebra $\mathcal{A}$ is zero, the reduced differential $\delta^X$
has the form
\begin{eqnarray}
\delta^X(g \otimes r \otimes x^{\omega}) & = & \lambda(g) \otimes r \otimes x^{\omega}. \nonumber
\end{eqnarray}
In particular the image of $\delta^X$ is contained in $P_0^* \otimes_{Z_0}[X^* \otimes_{A_0}R \otimes_k S(\omega)]$.
Considering the natural isomorphism $F:X_0^* \otimes_{A_0} R \otimes_k S(\omega) \to X_0^* \otimes_{A_0} R$
and the isomorphism $\eta$ in lemma~\ref{(DE)L:isoDual}, we want to show that the following diagram is commutative
\[
 \xymatrix@C=3pc{
X_0^* \otimes_{A_0} R \otimes_k S(\omega) \ar[r]^-{\eta \circ F} \ar[dd]_-{\delta^X} 
& \Hom_{A_0}(R,X_0)^* \ar[d]^-{\varepsilon^*} \ar@<5ex>@/^45pt/[dd]^-{\widehat{\delta}=\psi^{-1}\varepsilon^*} \\
& (\Hom_{A_0}(R,X_0) \otimes_{Z_0}P_0 )^* \\
P_0^* \otimes_{Z_0} [X_0^* \otimes _{A_0} R \otimes_k S(\omega)] \ar[r]_-{Id_{P_0^*} \otimes (\eta \circ F)} & 
P_0^* \otimes_{Z_0} \Hom_{A_0}(R,X_0)^* \ar[u]_-{\psi}.
}
\]
In order to compute $\psi \circ [Id_{P_0^*} \otimes (\eta \circ F)] \circ \delta^X$ we fix right dual bases
$(x_i,\nu_i)_i$ and $(p_j,\gamma_j)_j$ of $(X_0)_{Z_0}$ and $(P_0)_{Z_0}$ respectively.
Then, for an element $g \otimes r \otimes x^{\omega}$ of $X^* \otimes_{A_0} R \otimes_k S(\omega)$
we have that $[Id_{P_0^*}\otimes (\eta \circ F)](\delta^X(g \otimes r \otimes x^{\omega}))$ is given by
\[ 
[Id_{P_0^*}\otimes (\eta \circ F)]\left( \sum_{i,j}g(x_ip_j)\gamma_j \otimes (\nu_i \otimes r \otimes x^{\omega}) \right) 
= \sum_{i,j}g(x_ip_j)\gamma_j \otimes F^r_i, 
\]
where $F^r_i(f)=\nu_i(f(r))$. We evaluate $\psi$ in the last expression,
\[
 \psi \left( \sum_{i,j}g(x_ip_j)\gamma_j \otimes F^r_i \right) =[h \otimes p \mapsto \sum_{i,j}g(x_ip_j)\gamma_j(F^r_i(h) p)].
\]
Hence, for an element $h \otimes p$ of $\Hom_{A_0}(R,X_0) \otimes_{Z_0} P_0$ we have
\begin{eqnarray}
\psi[Id_{P_0^*}\otimes (\eta \circ F)]\delta^X(g \otimes r \otimes x^{\omega})(h \otimes p)
& = & \sum_{i,j}g(x_ip_j)\gamma_j(F^r_i(h) p) = \nonumber \\
& = & \sum_{i,j}g(x_ip_j)\gamma_j(\nu_i[h(r)] p) = \nonumber \\
& = & \sum_ig \left( x_i \sum_jp_j\gamma_j(\nu_i[h(r)] p) \right)  =\nonumber \\
& = & \sum_ig(x_i\nu_i[h(r)] p) = \nonumber \\
& = & g \left(\sum_ix_i\nu_i[h(r)]p \right)= \nonumber \\
& = & g(h(r)p)=g(p(h(r))). \nonumber
\end{eqnarray}
On the other hand, by the expressions of $\varepsilon^*$ (lemma~\ref{(DE)L:epsilonBimod}) and 
$\eta$ (lemma~\ref{(DE)L:isoDual}) we notice that
\[
 \varepsilon^*[(\eta \circ F)(g \otimes r \otimes x^{\omega})](h \otimes p)= g(p(h(r))).
\]
Then the diagram above is commutative, which shows that the differential $\delta^x$ is isomorphic
to the reduced differential $\delta^X$. Finally we observe that the explicit form of $\widehat{\delta}$ 
is obtained from the expression given for $\varepsilon^*$ in lemma~\ref{(DE)L:epsilonBimod} and that of the inverse 
$\psi^{-1}$ in lemma~\ref{(A)L:psiBimod} in the appendix.
\eproof

The ditalgebra $\mathcal{A}^X$ of the result above has as layer the pair 
\[
\left( \left( \begin{matrix} Z_0&0\\0&k \end{matrix} \right),
\left( \begin{matrix}P_0^*&\Hom_{A_0}(R,X_0)^*\\0&0 \end{matrix} \right)  \right),
\]
where the elements of $\Hom_{A_0}(R,X_0)^*$ have degree zero and those in $P_0^*$ have degree one.

%------------------------------------------------------------------
%------------------------------------------------------------------
\textbf{Choice of extension modules.} 
Let $\Delta$ be a Dynkin quiver with arbitrary orientation of its arrows and let $\widetilde{\Delta}$ 
be its extension (at the level of graphs) such that the added vertex is a source.
The purpose of the rest of this section is to determine a projective $\Delta$-module $R$ such that 
$k\widetilde{\Delta}$ is isomorphic to the one-point extension $k\Delta[R]$, as described by Ringel in \cite[3.6(4)]{cmR}.
Assume that $A_0$ is the path $k$-algebra of a finite solid quiver with an admissible ardering of its vertices
$Q$ and let $A=A_0[R]$ be the one-point extension of $A_0$ by the projective $A_0$-module $R$.
Let $C_0$ be the Cartan matrix of $A_0$, that is, the columns of $C_0$ are given by the dimension vectors of the 
indecomposable projective $P(i)=Ae_i$.  Since $A \cong Ae_0 \oplus Ae_{\omega}$ and $Ae_{\omega} \cong R \oplus k$, 
The Cartan matrix of $A$ has the form
\[
C=\left( \begin{matrix} C_0 & r \\ 0 & 1 \end{matrix} \right),
\]
where $r=\vdim R$. Since $C_0$ is invertible, $C$ has as inverse the matrix
\[
C^{-1}=\left( \begin{matrix} C_0^{-1} & -C_0^{-1}r \\ 0 & 1 \end{matrix} \right).
\]
Then the Coxeter matrix of $A$ is 
\[
 \Phi=-C^tC^{-1}= \left( \begin{matrix} -C_0^tC_0^{-1} & C_0^tC_0^{-1}r \\ -r^tC_0^{-1} & r^tC_0^{-1}r-1 \end{matrix} \right)
=\left( \begin{matrix} \Phi_0 & -\Phi_0 r \\ -r^tC_0^{-1} & q_0(r)-1 \end{matrix} \right),
\]
where $\Phi_0$ is the Coxeter matrix of $A_0$. By lemma \ref{(DE)L:extTens}, $A=A_0[R]$ is the path algebra of a quiver $\widetilde{Q}$.
Let $\langle x,y \rangle$ and $\langle x,y \rangle_0$ be the bilinear forms associated to the quivers $\widetilde{Q}$ and $Q$
respectively. Recall that $\langle x,y \rangle=x^t M_{\widetilde{Q}} y$ where $M_{\widetilde{Q}}$ is the matrix associated to the
quiver $\widetilde{Q}$. By lemma~\ref{(P)L:invers}, $M_{\widetilde{Q}}=C^{-t}$, for the ordering of the vertices in
$\widetilde{Q}$ is admissible. Let $(x,y)=\frac{1}{2}(\langle x,y \rangle+\langle y,x \rangle)$ be the symmetrization of
$\langle x,y \rangle$ and $q(x)=\langle x,x \rangle$ its associated quadratic form. In a similar way define
$\langle \cdot,\cdot \rangle_0$, $(\cdot,\cdot)_0$ and $q_0$ respect to the quiver $Q$ and its associated matrix $M_Q=C_0^{-t}$.
Recall that an integral vector $x$ is in the radical of the quadratic form $q$ if $(x,y)=0$ for any integral vector $y$.

\begin{lema} \label{(DE)L:vctRad}
 For an integral vector $w_0$ the following conditions are equivalent.
\begin{itemize}
 \item[i)] $w_0+\mathbf{e}_{\omega}$ is a vector in the radical of $q$.
 \item[ii)] $r=(I-\Phi_0^{-1})w_0$ and $q_0(w_0)=1$.
 \item[iii)] The linear forms $\langle r,- \rangle_0$ and $2(w_0,-)_0$ coincide and $q_0(w_0)=1$. 
\end{itemize}
\end{lema}

\bproof
We transcript the proof of Ringel \cite[2.5 point (11)]{cmR}. To show the equivalence of $(ii)$ and $(iii)$ 
notice that the forms $\langle r,- \rangle_0$ and $2(w_0,-)_0$ corres\-pond to multiplication by the vectors 
$\langle r,- \rangle_0=r^tC_0^{-t}$ and $2(w_0,-)_0=w_0^t(C_0^{-t}+C_0^{-1})$. Hence the forms $\langle r,- \rangle_0$ 
and $2(w_0,-)_0$ coincide if and only if 
\[
r^tC_0^{-t}=w_0^t(C_0^{-t}+C_0^{-1}),
\]
and transposing
\[
C_0^{-1}r=(C_0^{-t}+C_0^{-1})w_0.
\]
This happens if and only if $r=(I+C_0C_0^{-t})w_0=(I-\Phi_0^{-1})w_0$.
To prove the equivalence of $(i)$ and $(iii)$ we compute the linear form $2(-,w_0+\mathbf{e}_{\omega})_0$,
\[
 \left( \begin{matrix} C_0^{-1}+C_0^{-t} & -C_0^{-1}r \\ -r^tC_0^{-t} & 2 \end{matrix} \right)
\left( \begin{matrix} w_0 \\ 1 \end{matrix} \right)= 
\left( \begin{matrix} (C_0^{-1}+C_0^{-t})w_0 - C_0^{-1}r \\ -r^tC_0^{-t}w_0+2 \end{matrix} \right).
\]
Thus $w_0+\mathbf{e}_{\omega}$ is in the radical of $q$ if and only if the forms 
$\langle r,- \rangle_0$ and $2(w_0,-)_0$ coincide and $\langle r , w_0 \rangle_0=2(w_0,w_0)_0=2$.
\eproof

Assume now that $Q=\Delta$ is a Dynkin quiver and take $A_0=k\Delta$. Let $W_0$ be an indecomposable $A_0$-module
such that $\vdim W_0=w_0$ is the maximal positive root of $\Delta$. If $d_i$ is the partial derivative of the quadratic 
form $q_{\Delta}$ evaluated in the maximal root, 
\[
d_i=\frac{\partial}{\partial x_i} q_{\Delta}(w_0),
\]
then $d_i\geq 0$ for all $i=1,\ldots,n$. Indeed, since $d_i=2(\mathbf{e}_i,w_0)$, the simple reflexion $\sigma_i(w_0)$ has
the form
\[
 \sigma_i(w_0)=w_0-d_i,
\]
and by maximality $d_i \geq 0$.

\begin{lema} \label{(DE):extEucl}
Let $\Delta$ be a Dynkin quiver and $W_0$ an indecomposable $\Delta$-modulo whose dimension vector 
$w_0=\vdim W_0$ is the maximal root of $\Delta$. Let $\widetilde{\Delta}$ be the extended Dynkin quiver
where the extended vertex is a source and let $R$ be the projective $\Delta$-module $R=\bigoplus_{i} d_iP(i)$
where $P(i)$ is projective cover of the simple module of vertex $i$ and $d_i=\frac{\partial}{\partial x_i} q_{\Delta}(w_0)$. 
Then $k\widetilde{\Delta} \cong k\Delta[R]$. 
\end{lema}

\bproof
Take $A_0=k\Delta$. Since $\widetilde{\Delta}$ is obtained from $\Delta$ by adding the vertex $\omega$ and arrows
from $\omega$ to certain vertices in $\Delta$ called exceptional vertices (marked with a circle in table~\ref{(DE)T:dynkinE}), 
it is clear that $k\widetilde{\Delta}=A_0[R']$ for some projective $\Delta$-modulo $R'$. 
As is well known $w_0+\mathbf{e}_{\omega}$ is generator of the radical of $q_{\widetilde{\Delta}}$, hence the equivalence
of $(i)$ and $(ii)$ in the lemma above implies that $r'=\vdim R'=(I-\Phi_0^{-1})w_0$.

On the other hand, if $r=\vdim R=\sum_id_ip_i$ (with $p_i=\vdim P(i)$) we show that the bilinear forms $\langle r,- \rangle_0$ and 
$2(w_0,-)_0$ coincide. It is enough to see that these forms coincide in the canonical basis
\[
 \langle r,\mathbf{e}_j \rangle_0 = \sum_i d_i \langle p_i,\mathbf{e}_j \rangle_0 d_j=d_j = 2(w_0,\mathbf{e}_j)_0,
\]
for $\langle p_i,\mathbf{e}_j \rangle_0 = p_i^tC_0^{-t}\mathbf{e}_j=(C_0\mathbf{e}_i)^tC_0^{-t}\mathbf{e}_j=
\mathbf{e}_i^t(C_0^tC_0^{-t})\mathbf{e}_j=\mathbf{e}_i^t\mathbf{e}_j=\delta_{i,j}$. Then the forms $\langle r,- \rangle_0$ 
and $2(w_0,-)_0$ are the same, and by the equivalence of $(ii)$ and $(iii)$ in the lemma above
$r=(I-\Phi_0^{-1})w_0$. Then $r=r'$. Since $R'=\bigoplus_{i} d'_iP(i)$ for some nonnegative integers
$d'_i$ and the set of vectors $\{p_i\}$ is linearly independent we have that $d_i=d'_i$ for all vertices $i$. 
Hence $R'=R$ and $k\widetilde{\Delta}\cong A_0[R]$.
\eproof

We end this section with the following lemma, that can be found in Ringel \cite[section 3.4(4)]{cmR}
and that is fundamental for the description of the reductions we will perform at the end of this chapter.
\begin{lema} \label{(DE)L:valores}
Let $W_0$ be an indecomposable $A_0$-module corresponding to the ma\-xi\-mal root of $\Delta$.
Consider the Auslander-Reiten quiver $\Gamma(A_0)$ of $A_0$-mod and the following subsets of vertices determined by $[W_0]$,
\[
 \mathcal{X}^0=\dConL([W_0])=\{ [M] \in \Ind A_0 \; | \; [M] \prec [W_0] \text{ but } \tau^{-1}[M] \npreceq [W_0]\},
\]
\[
 \mathcal{Y}^0=\dConR([W_0])=\{ [M] \in \Ind A_0 \; | \; [W_0] \prec [M] \text{ but } [W_0] \npreceq \tau[M] \},
\]
\[
 \mathcal{W}^0=\{ [M] \in \Ind A_0 \; | \; [M] \npreceq [W_0] \npreceq [M] \}.
\]
Then for any indecomposable $A_0$-module $M$ we have
 \begin{equation*}
\dimk_k \Hom_{A_0}(R,M) = \left\{
\begin{array}{rl}
2, & \text{if } [M] \cong W_0,\\
1, & \text{if } [M] \in \mathcal{X}^0 \cup \mathcal{Y}^0,\\
0, & \text{if } [M] \in \mathcal{W}^0.
\end{array} \right.
\end{equation*}
\end{lema}

\bproof
Since $R$ is projective we have that $\Ext^1_{A_0}(R,M)=0$, thus using lemma~\ref{(P)L:Euler}
and the point $(iii)$ in lemma~\ref{(DE)L:vctRad} we have
\[
 \dimk_k \Hom_{A_0}(R,M)=\langle \vdim R,\vdim M \rangle_0=2(\vdim W_0,\vdim M)_0.
\]
Hence $\dimk_k \Hom_{A_0}(R,W_0)=2$.
Assume now that $[M] \in \mathcal{W}_0$, $\mathcal{X}_0$ or $\mathcal{Y}_0$.
Using the Auslander-Reiten formulas in the hereditary case and since there are neither
paths from $M$ to $\tau W_0$ nor from $\tau^{-1}W_0$ to $M$, we have that
\[
\Ext^1_{A_0}(W_0,M) \cong D\Hom_{A_0}(M,\tau W_0)=0, \quad \text{and}
\]
\[
\Ext^1_{A_0}(M,W_0) \cong D\Hom_{A_0}(\tau^{-1}W_0,M)=0.
\]
Then
\[
 2(\vdim W_0,\vdim M)_0= \dimk_k \Hom_{A_0}(W_0,M)+\dimk_k \Hom_{A_0}(M,W_0),
\]
and $[M] \in \mathcal{W}_0$ implies that $\dimk_k \Hom_{A_0}(R,M)=0$.

Recall that $[W_0]$ is a wing vertex in $\Gamma(A_0)$. The notation for the wings $\theta(n_s)$ 
of $[W_0]$ given at the end of section~\ref{(DE)S:ARDynkin} is, for $1 \leq s \leq t$,
\[
 \xymatrix@!0@R=2.5pc@C=2.5pc{
& & & & W^s_{1,n_s} \ar[rd] \\
& & & W^s_{1,n_s-1} \ar[ru] \ar[rd] & & W^s_{2,n_s} \ar[rd]\\
& & W^s_{1,n_s-2} \ar[ru] & & W^s_{2,n_s-1} \ar@{}[rd]|{\cdots} \ar[ru] & & W^s_{3,n_s} \ar@{}[rd]|{\cdots} \\
& W^s_{1,2} \ar@{}[ru]|{\cdots} \ar[rd] & & W^s_{2,3} \ar@{}[rd]|{\cdots} \ar@{}[ru]|{\cdots}  & & W^s_{n_s-2,n_s-1} 
\ar@{}[ru]|{\cdots} \ar[rd]  & & W^s_{n_s-1,n_s} \ar[rd] \\
W^s_{1,1} \ar[ru] & & W^s_{2,2} \ar[ru] & & W^s_{i,i} \ar@{}[ru]|{\cdots} & & W^s_{n_s-1,n_s-1} \ar[ru] & & W^s_{n_s,n_s}.
}
\]
Then the set $\mathcal{X}_0=\dConL([W_0])$ consists in the isomorphism classes of the modules $W_{1,i}^{s}$ 
with $i=1,\ldots,n_s-1$ for $s=1,\ldots,t$ and $\mathcal{Y}_0=\dConR([W_0])$ is formed by the set $\{ [W_{i,1}^{s}] \}$ 
with $i=1,\ldots,n_s-1$ and $s=1,\ldots,t$. 
Moreover, $\mathcal{W}_0$ is given by the elements of the form $[W_{i,j}^{s}]$ for $1<i \leq j < n_s,$ and $s=1,\ldots,t$.
Assume that $[M] \in \mathcal{X}_0$. If $M \cong W_{1,n_s-1}^{(s)}$, then $\rad^2(M,W_0)=0$ and
\[
 \dimk_k \Hom_{A_0}(M,W_0) =  \dimk_k (\rad (M,W_0)/\rad^2(M,W_0))=1.
\]
By induction, if $M \cong W_{1,i-1}^{(s)}$ and we assume that $\dimk_k \Hom_{A_0}(W_{1,i}^{(s)},W_0)=1$, then 
\[
\dimk_k\Hom_{A_0}(W_{1,i-1}^{(s)},W_{1,i}^{(s)})=1
\] 
(since $\rad^2(W_{1,i-1}^{(s)},W_{1,i}^{(s)})=0$), 
hence $\Hom_{A_0}(W_{1,i-1}^{(s)},W_0)=1$ or $0$. The last case is impossible for $W_{1,i-1}^{(s)}$ is cogenerated by $W_0$
(lemma~\ref{(DE)L:alaDyn}). A similar argument for $[M] \in \mathcal{Y}_0$ completes the proof.
\eproof

\section{Lifting of exact sequences.} \label{(DE)S:levantaSuc}
%----------------------------------------------------------------------
%------------------------------------------------------------------
Let $A_0$ be a finite dimensional  $k$-algebra and $R$ a projective $A_0$-module.
We consider the category of $A_0$-modules included in the category of subspaces $\check{\mathcal{U}}(A_0\text{-mod})$ 
respect to the functor $\Hom_{A_0}(R,-)$ in two essentially different ways,
\[
A_0\text{-mod} \longrightarrow \check{\mathcal{U}}(A_0\text{-mod}),
\]
through the inclusion functor $M_0 \mapsto M_0=(M_0,0,0)$ and induction functor
$M_0 \mapsto \overline{M_0}=(M_0,|M_0|,I_{|M_0|})$. The \textbf{liftings} $M_0$ and $\overline{M_0}$ coincide as long
as $|M_0|=0$. The following result (Ringel, \cite{cmR}2.5(6)) relates the almost split sequences in $A_0$-mod and 
$\check{\mathcal{U}}(A_0\text{-mod})$.

\begin{lema} \label{(DE)L:levAR}
\begin{itemize}
 \item[a)] If $\varepsilon_0=\xymatrix{0 \ar[r] & X_0 \ar[r]^-{f} & Y_0 \ar[r]^-{g} & Z_0 \ar[r] & 0}$ is an almost split sequence 
of $A_0$-modules, then 
\[
\overline{\varepsilon_0}=\xymatrix{0 \ar[r] & \overline{X_0} \ar[r]^-{(f,|I_{X_0}|)} & (Y_0,|X_0|,|f|) \ar[r]^-{(g,0)} & Z_0 \ar[r] & 0}
\]
is an almost split sequence in $\check{\mathcal{U}}(A_0\text{-mod})$, called lifting of $\varepsilon_0$.

 \item[b)] Let $\varepsilon=\xymatrix{0 \ar[r] & X \ar[r]^-{f} & Y \ar[r]^-{g} & Z \ar[r] & 0}$ be an almost split sequence
in the category $\check{\mathcal{U}}(A_0\text{-mod})$. If $Z_{\omega}\neq 0$ then the restriction $\varepsilon|_{A_0}$ is trivial. 
If $Z_{\omega} = 0 $ then $Z_0$ is not a projective $A_0$-module and $\varepsilon|_{A_0}$ is an almost split sequence. 
\end{itemize}
\end{lema}

\bproof
For $(a)$ we follow Ringel \cite{cmR}2.5(lemmas 5 and 6). Observe first that $\overline{\varepsilon_0}$ 
is an exact sequence, for $|-|=\Hom_{A_0}(R,-)$ is an exact functor ($R$ is projective), 
and the following diagram commutes
\[
 \xymatrix{
0 \ar[r] & |X_0| \ar@{=}[r] \ar@{=}[d] & |X_0| \ar[d]^-{|f|} \ar[r] & 0 \ar[r] \ar[d] & 0 \\
0 \ar[r] & |X_0| \ar[r]_{|f|} & |Y_0| \ar[r]_-{|g|} & |Z_0| \ar[r] & 0.
}
\]
On the other hand, $(f,|I_{X_0}|)$ is not a section for $f$ is not a section. Assume that $v=(v_0,v_{\omega}):\overline{X_0}\to V $ 
is a morphism which is not a section. If there exists $v_0':V_0 \to X_0$ such that $v_0'v_0=I_{X_0}$, taking $v'=(v_0',|v_0'|\gamma_V)$
we have that $v'v=I_{\overline{X_0}}$, hence $v_0$ is not a section. Since $f$ is left almost split, there exists 
$\eta_0:Y_0 \to V_0$ such that $v_0=\eta_0f$.
\[
 \xymatrix{
\overline{X_0} \ar[d]_-{(v_0,v_{\omega})} \ar[r]^-{(f,|I_{X_0}|)} & (Y_0,|X_0|,|f|) \ar[ld]^-{(\eta_0,v_{\omega})} 
& X_0 \ar[d]_-{v_0} \ar[r]^-{f} & Y_0 \ar[ld]^-{\eta_0} \\
(V_0,V_{\omega},\gamma_V) & & V_0
}
\]
Notice that $(\eta_0,v_{\omega})(f,|I_{X_0}|)=(v_0,v_{\omega})$.
Finally, if $\xi=(\xi_0,\xi_{\omega})$ is an element in $\End(Y_0,|X_0|,|f|)$ that satisfies $\xi(f,|I_{X_0}|)=(f,|I_{X_0}|)$ then
$\xi_0f=f$ and $\xi_{\omega}=|I_{X_0}|$. Thus $\xi$ is an automorphism and $(f,|I_{X_0}|)$ is a left minimal almost split
morphism. By lemma~\ref{(DE)L:ARequiv}, the sequence $\overline{\varepsilon_0}$ is almost split.

We show now $(b)$. Since the extended vertex $\omega$ is a source, $Z_0$ is a submodule of $Z$ in $\check{\mathcal{U}}(A_0\text{-mod})$.
If $Z_{\omega} \neq 0$ then the inclusion $i:Z_0 \hookrightarrow Z$ is not a retraction, and since $g$ is right almost split,
there exists $t:Z_0 \to Y$ such that $i=gt$. Then $I_{Z_0}=g_0t_0$ and $g_0$ is a retraction.
On the other hand, if $Z_{\omega}=0$ then $Z_0$ cannot be a projective $A_0$-module, for in the contrary
the epimorphism $g_0:Y_0 \to Z_0$ would be a retraction and hence $g=(g_0,0):Y \to Z=Z_0$ would be also a retraction
(impossible since $g$ is right almost split).
In this way, if $\varepsilon_0$ is the almost split sequence in $A_0$-mod that ends in $Z_0$, then
by point $(a)$, $\varepsilon$ and $\overline{\varepsilon_0}$ are almost split sequences that start in $Z$.
By uniqueness these are isomorphic sequences, and hence $\varepsilon|_{A_0} \cong \overline{\varepsilon_0}|_{A_0}=\varepsilon_0$
is an almost split sequence in $A_0$-mod.

\eproof

We return to the case in which $A_0$ is the path algebra of a Dynkin quiver $\Delta$, $W_0$
is an indecomposable $A_0$-module with dimension vector the maximal root and $R$ is the projective $A_0$-module 
described in lemma~\ref{(DE):extEucl}. For any morphism $\rho \in \Hom_{A_0}(R,W_0)$ we denote by $W_0(\rho)$ 
the object in $\check{\mathcal{U}}(A_0\text{-mod})$ given by
\[
 W_0(\rho)=(W_0,k,\rho),
\]
that is, $\gamma_{W_0(\rho)}$ is the transformation which sends $1 \in k$ to the morphism $\rho$. 
Since the endomorphism algebra of $W_0$ is isomorphic to the field $k$, the morphisms
$f:W_0(\rho) \to W_0(\rho')$ are given by pairs of scalars $f=(aId_k,bId_{W_0})$ which make the following diagram
commutative
\[
 \xymatrix{
k \ar[r]^-{a} \ar[d]_-{\rho} & k \ar[d]^-{\rho'}\\
|W_0| \ar[r]_{|bI_{W_0}|} & |W_0|.
}
\]
Hence there is a nonzero morphism $f:W_0(\rho) \to W_0(\rho')$ if and only if $a,b\neq 0$ and $\rho'=(b/a)\rho$. 
In this case $f$ is an isomorphism. Clearly $W_0(\rho)$ is indecomposable if and only if $\rho \neq 0$.
Then there is a $\mathbb{P}\Hom_{A_0}(R,W_0)$-parametric family of nonisomorphic indecomposable $A$-modules $W_0(\rho)$. 
All of them share as dimension vector the positive generator $w_0 + \mathbf{e}_{\omega}$ of the radical of the quadratic form
associated to $A$.

\begin{lema} \label{(DE)L:ciclo}
For any nonzero morphism $\rho \in \Hom_{A_0}(R,W_0)$ the indecomposable $A$-module $W_0(\rho)$ is neither
posprojective nor preinjective.
\end{lema}
\bproof
Assume that $W_0(\rho)$ is a posprojective $A$-module. In particular $W_0(\rho)$ is exceptional.
By lemmas~\ref{(DE)L:radCero} and~\ref{(P)L:coxPP}$(b)$, since $\vdim W_0(\rho)$ is in the radical of the quadratic form
of $A$, we have that
\[
 \vdim \tau^{-1}W_0(\rho)=\vdim W_0(\rho).
\]
Since $\tau^{-1}W_0(\rho)$ is also exceptional, by lemma~\ref{(P)L:determ} there is an isomorphism 
$W_0(\rho) \cong \tau^{-1}W_0(\rho)$. This is a contradiction, for the posprojective component of
$A$-mod has no periodic orbits (lemmas~\ref{(DE)L:componentes} and~\ref{(DE)L:dirigido}).
In a similar way it can be shown that $W_0(\rho)$ is not a preinjective $A$-module.
\eproof

Consider a wing $\theta(n_s)$ of vertex $[W_0]$ in the Auslander-Reiten quiver $\Gamma(A_0)$,
as in section~\ref{(DE)S:ARDynkin} and recall the notation $[W_{i,j}^s]$ for the vertices in
$\theta(n_s)$ given in that section. Fix a nonzero morphism $\rho_s':R \to W_{1,n_s-1}^{s}$ 
(the group of morphisms $|W_{1,n_s-1}^{s}|=\Hom_{A_0}(R,W_{1,n_s-1}) \cong k$ is one dimensional).
Fixing an irreducible mapping $f:W_{1,n_s-1}^s \to W_0$ and defining $\rho_s=f\rho'_s$ we have
a morphism $\xymatrix{\overline{W_{1,n_s-1}^{s}} \ar[r] & W_0(\rho_s)}$ in the category $A$-mod as shown in the following
figure,
\[
\xymatrix@R=1pc{
 & W_0 & |W_{1,n_s-1}^{s}| \ar@{=}[d] \ar[r]^-{\rho_s' \mapsto 1} & k \ar[d]\\
R \ar[ru]^-{\rho_s} \ar[r]_-{\rho_s'} & W_{1,n_s-1}^{s} \ar[u]_-{f} & |W_{1,n_s-1}^{s}| \ar[r]_-{\rho_s' \mapsto \rho_s} & |W_0|.
}
\]
In this way Ringel shows~\cite[3.4 point (5)]{cmR} that the wing $\theta(n_s)$ lifts in a completely to
a subquiver of $\Gamma(A)$, as shown in the following figure. The projective-inyective vertex corresponds to the module $W_0(\rho_s)$, 
\[
 \xymatrix@!0@R=2.5pc@C=2.5pc{
& & & & W_0(\rho_s) \ar[rd] \\
& & & \overline{W_{1,n_s-1}^{s}} \ar[ru] \ar[rd] & & W_{2,n_s}^{s} \ar[rd]\\
& & \overline{W_{1,n_s-2}^{s}} \ar[ru] & & W_{2,n_s-1}^{s} \ar@{}[rd]|{\cdots} \ar[ru] & & W_{3,n_s}^{s} \ar@{}[rd]|{\cdots} \\
& \overline{W_{1,2}^{s}} \ar@{}[ru]|{\cdots} \ar[rd] & & W_{2,3}^{s} \ar@{}[rd]|{\cdots} \ar@{}[ru]|{\cdots}  
& & W_{n_s-2,n_s-1}^{s} \ar@{}[ru]|{\cdots} \ar[rd]  & & W_{n_s-1,n_s}^{s} \ar[rd] \\
\overline{W_{1,1}^{s}} \ar[ru] & & W_{2,2}^{s} \ar[ru] & & W_{i,i}^{s} \ar@{}[ru]|{\cdots} & & W_{n_s-1,n_s-1}^{s} \ar[ru] 
& & W_{n_s,n_s}^{s}.
}
\]

As consequence of the lemmas above we have the following result, which is fundamental to our objectives.
For any subsets $\mathcal{L}_1$ and $\mathcal{L}_2$ of $\Ind A_0$ denote by $\mathcal{L}_1 \vee \mathcal{L}_2$
the smallest class of objects in $A_0$-mod which contains the class representatives in $\mathcal{L}_1$ and
$\mathcal{L}_2$. 

\begin{proposicion} \label{(DE)C:ciclo}
Consider the subsets in $\Ind A_0$ given by 
\[
\mathcal{X}=\ConL([W_0]) \quad \text{and} \quad \mathcal{Y}=\ConR([W_0]).
\]
\begin{itemize}
 \item[a)] Any posproyective $A$-module lies in the full subcategory
$\check{\mathcal{U}}(\mathcal{X} \vee [W_0])$ of $A$-mod.
 \item[b)] Any preinjective $A$-module lies in the full subcategory
$\check{\mathcal{U}}([W_0] \vee \mathcal{Y})$ of $A$-mod.
\end{itemize} 
\end{proposicion}

\bproof
We start with the proof of $(a)$.
Recall that the posprojective component of $\Gamma(A)$ admites an admissible ordering of its vertices
(lemma~\ref{(DE)L:orden}). To avoid the use of negative indexes, we order representatives of the indecomposable 
posprojective $A$-modules $\{N_{i}\}$ in such a way that if $f:N_j \to N_i$ is an irreducible morphism then $j < i$
(opposed admissible order). We will prove the following auxiliar statement.
\begin{itemize}
 \item[$i)$] If $Z_0$ is an indecomposable nonprojective $A_0$-module which is direct summand of the restriction 
$N|_{A_0}$ for some indecomposable posprojective $A$-module $N$, then $Z_0$ and $\overline{\tau_0(Z_0)}$ are posprojective 
$A$-modules, smaller or equal to $N$ in the ordering given above $\{N_i\}_{i \in \mathbb{N}}$.
\end{itemize}
We prove the claim inductively. Since the first modules in the list $N_1,N_2,\ldots$ are projective,
the base case holds by vacuity. Assume that the condition $(i)$ is valid for the modules that belong to the restriction
$N_i|_{A_0}$ of some posprojective module $N_i$ with $1\leq i \leq \ell$ and assume that $Z_0$ is direct summand of 
$N_{\ell+1}|_{A_0}$. If $Z_0$ is in the restriction of some of the modules $N_1,\ldots,N_{\ell}$, then we have
the result by induction hypothesis. So we can assume that $Z_0$ is not direct summando of any $N_{i}|_{A_0}$ for 
$i=1,\ldots,\ell$. Consider the almost split sequence which ends in $N_{\ell+1}$,
\[
\varepsilon: \xymatrix{
0 \ar[r] & \tau(N_{\ell+1}) \ar[r] & E \ar[r] & N_{\ell+1} \ar[r] & 0.
}
\]
Notice that the restriction $\varepsilon|_{A_0}$ cannot be trivial, for $Z_0$ is direct summand of $N_{\ell+1}|_{A_0}$ but not 
of $E|_{A_0}$. Then by lemma~\ref{(DE)L:levAR}$(b)$ we have that $(N_{\ell+1})_{\omega}=0$
thus $Z_0 \cong N_{\ell+1}$. Consider now the almost split sequence in $A_0$-mod which ends in $Z_0$,
\[
\varepsilon_0: \xymatrix{
0 \ar[r] & \tau_0(Z_0) \ar[r] & F \ar[r] & Z_0 \ar[r] & 0.
}
\]
By lemma~\ref{(DE)L:levAR}$(a)$ the lifting $\overline{\varepsilon_0}$ and $\varepsilon$ are both Auslander-Reiten sequences
in $A$-mod which end in $N_{\ell+1}$, hence $\overline{\tau_0(Z_0)} \cong \tau(N_{\ell+1})$.
That is, $Z_0$ and $\overline{\tau_0(Z_0)}$ are posprojective $A$-modules, smaller or equal to $N_{\ell+1}$.

We prove now in steps the statement of the proposition.\\
\underline{Step 1.} \textit{Assume that $Z_0$ is an indecomposable $A_0$-module such that $\tau_0^i(Z_0) \in \mathcal{X}_0$ for
some $i \geq 1$. Then $Z_0$ is not a module in the restriction to $A_0$ of the posprojective component of $A$-mod.}
Those liftings to $A$-mod of elements in $\mathcal{X}_0$ are connected by irreducible morphisms to the module $W_0(\rho_s)$, 
for the morphism $\rho_s$ corresponding to the lifting of some wing $\theta(n_s)$ of vertex $[W_0]$. 
Then by lemma~\ref{(DE)L:ciclo} the liftings of elements in 
$\mathcal{X}_0$ are connected in $\Gamma(A)$ to an object that is not posprojective, and hence cannot be posprojective.
Applying successively the point $(i)$ above it is clear that if $Z_0$ is direct summand of $N|_{A_0}$ 
for some posprojective $A$-modulo $N$, then the lifting $\overline{\tau_0^{i}(Z_0)}$ belongs to the posprojective
component, which is impossible since $\tau_0^{i}(Z_0)\in \mathcal{X}_0$. This shows that $Z_0$ cannot be direct summand 
of the restriction to $A_0$ of any posprojective module.\\
\underline{Step 2.} \textit{Assume that $Z_0$ is an indecomposable $A_0$-module such that $\tau_0^{i}(Z_0) \cong W_0$ 
for some $i\geq 1$. Then $Z_0$ is not a module in the restriction to $A_0$ of the posprojective component of $A$-mod.}
Again by successive applications of $(i)$ it is enough to show that $\overline{W_0}$ is not a posprojective $A$-module.
If $\overline{W_0}$ is posprojective, then it is not an injective $A$-module. By lemma~\ref{(DE)L:levAR} the almost split
sequence that starts in $\overline{W_0}$ is lifting to $A$-mod of an almost split sequence in 
$A_0$-mod that starts in $W_0$ (since we assume that $W_0 \cong \tau_0^{i}(Z_0)$, the module $W_0$ is not injective in $A_0$-mod),
\[
 \xymatrix{0 \ar[r] & W_0 \ar[r] & E \ar[r] & \tau_0^{-1}W_0 \ar[r] & 0}.
\]
Since there is an irreducible mapping $W_0 \to W^1_{2,n_1}$ in
$A_0$-mod, we have that $W^1_{2,n_1}$ is direct summand of the restriction to $A_0$-mod of some posprojective $A$-module.
This contradicts step 1, for $\tau_0^{-1}[W^1_{2,n_1}]\in \mathcal{X}_0$.\\
\underline{Step 3.} We finally conclude that if $Z_0$ is isomorphic to a direct summand of the restriction to $A_0$ of a
posprojective $A$-module $N$, by the steps above $\tau_0^{i}[Z_0]$ does not belong to $\mathcal{X}_0$ nor it is equal to the 
vertex $[W_0]$ for any $i \geq 1$. Hence $[Z_0] \in \mathcal{X} \cup \{[W_0]\}$, that is, 
$N$ belongs to $\check{\mathcal{U}}(\mathcal{X} \vee [W_0])$.

The proof of $(b)$ is similar, we give the details next.
We give an admissible ordering to the vertices in the preinjective component of $A$-mod  $\{M_i\}_{i \in \mathbb{N}}$
(lemma~\ref{(DE)L:orden}), and prove the following auxiliar claim.
\begin{itemize}
 \item[$i')$] If $Z_0$ is an indecomposable noninjective $A_0$-module that is direct summand of the restriction 
$M|_{A_0}$ for some indecomposable preinjective $A$-module $M$, then $\overline{Z_0}$ and $\tau_0^{-1}(Z_0)$ are 
preinjective $A$-modules, smaller or equal to $M$ in the order given above.
\end{itemize}
We give an inductive proof. Since the first modules in the list $M_1,M_2,\ldots$ are injective,
the base case holds by vacuity. Assume that the condition $(i')$ is valid for modules in the restriction $M_i|_{A_0}$
of some preinjective $M_i$ with $1\leq i \leq \ell$ and assume that $Z_0$ is direct summand of $M_{\ell+1}|_{A_0}$. 
If $Z_0$ is in the restriction of some of the modules $M_1,\ldots,M_{\ell}$, by induction hypothesis we have the result. 
Then we can assume that $Z_0$ is not direct summand of any $M_{i}|_{A_0}$ for $i=1,\ldots,\ell$. Consi\-der the almost split
sequence that starts in $M_{\ell+1}$,
\[
\varepsilon: \xymatrix{
0 \ar[r] & M_{\ell+1} \ar[r] & E \ar[r] & \tau^{-1}(M_{\ell+1}) \ar[r] & 0.
}
\]
Notice that the restriction $\varepsilon|_{A_0}$ cannot be trivial, for $Z_0$ is direct summand of $M_{\ell+1}|_{A_0}$ but not of 
$E|_{A_0}$. Then by lemma~\ref{(DE)L:levAR}$(b)$ we have that $(\tau^{-1}(M_{\ell+1}))_{\omega}=0$ and $\varepsilon|_{A_0}$ 
is an almost split sequence in $A_0$-mod. In particular $M_{\ell+1}|_{A_0}$ is indecomposable and thus 
$Z_0 \cong M_{\ell+1}|_{A_0}$. Then $\tau_0^{-1}Z_0 \cong \tau^{-1}(M_{\ell+1})|_{A_0}=\tau^{-1}(M_{\ell+1})$ and 
$M_{\ell+1} \cong \overline{Z_0}$, which are indecomposable preinjectives smaller or equal to $M_{\ell+1}$. 
This completes the induction step.\\
\underline{Step 1'.} \textit{Assume that $Z_0$ is an indecomposable $A_0$-module such that $\tau_0^{-i}(Z_0) \in \mathcal{Y}_0$
for some $i \geq 1$. Then $Z_0$ is not a module in the restriction to $A_0$ of the preinjective component of $A$-mod.}
The elements in $\mathcal{Y}_0$ (considered as $A$-modules) are connected by irreducible morphisms to the 
$A$-module $W_0(\rho_s)$ for the morphism $\rho_s$ corresponding to the lifting of some wing 
$\theta(n_s)$ of vertex $W_0$. Then by lemma~\ref{(DE)L:ciclo} the elements in 
$\mathcal{Y}_0$ are connected in $\Gamma(A)$ to a nonpreinjective vertex, and hence cannot be preinjective.
Applying successively the claim $(i')$ above it is clear that if $Z_0$ is direct summand of $M|_{A_0}$ for some
preinjective $M$, then $\tau_0^{-i}(Z_0)\in \mathcal{Y}_0$ belongs to the preinjective component, 
which is impossible. This shows that $Z_0$ cannot be direct summand of the restriction to $A_0$ of a preinjective module.\\
\underline{Step 2'.} \textit{Assume that $Z_0$ is an indecomposable $A_0$-module such that $\tau_0^{-i}(Z_0) \cong W_0$ 
for some $i\geq 1$. Then $Z_0$ is not a module in the restriction to $A_0$ of the preinjective component of $A$-mod.}   
Again by successive applications of $(i')$ it is enough to show that $W_0$ is not a preinjective $A$-module.
For if $W_0$ is preinjective then $W_0$ is not projective and the almost split sequence $\varepsilon$ which ends in $W_0$ 
(which is contained in $\mathcal{I}$) restricts $\varepsilon|_{A_0}$ to an almost split sequence 
in $A_0$-mod (lemma~\ref{(DE)L:levAR}(b)). Since there is an irreducible morphism $W^1_{1,n_1-1} \to W_0$ in
$A_0$-mod, we have that $W^1_{1,n_1-1}$ is direct summand of the restriction to $A_0$-mod of a preinjective $A$-module.
This contradicts step 1', since $\tau^{-1}(W^1_{1,n_1-1})\in \mathcal{Y}_0$.\\
\underline{Step 3'.} From steps 1' y 2' it follows that no element at the left of the section 
$\{[W_0]\} \cup \mathcal{Y}_0$ appears as direct summand of the restriction to $A_0$ of a preinjective $A$-module, which implies
that any preinjective $A$-module is an element in $\check{\mathcal{U}}([W_0] \vee \mathcal{Y})$.
\eproof

\section{Lifting of functors.} \label{(DE)S:levantaFun}
%-------------------------------------
%-------------------------------------
%-------------------------------------
We now study situations in which functors in $A_0$-mod can be lifted to the extended category $A_0[R]$-mod. 
The following preliminar result will be necessary.
\begin{lema} \label{(DE)L:isoNat}
Let $\mathcal{Z}$ be an additive $k$-category and $F,G,H:\mathcal{Z} \to \mathcal{Z}$ additive functors.
For a natrual transformation of functors $\eta: F \to G$ denote by $\eta\cdot H: F \circ H \to G \circ H$ 
the natural transformation given by $(\eta\cdot H)^M=\eta^{H(M)}$ for any object $M$ in $\mathcal{Z}$.
The following conditions are equivalent,
\begin{itemize}
 \item[a)] the natural transformation $\eta \cdot H$ is an isomorphism of functors,
 \item[b)] if $M \in \mathcal{Z}$ is isomorphic to an object of the form $H(N)$ then 
           $\eta^{M}:F(M) \to G(M)$ is an isomorphism.
\end{itemize}
\end{lema}
\bproof
If $f:M \to N$ is a morphism in $\mathcal{Z}$, by naturality of $\eta$ the following diagram is commutative,
\[
 \xymatrix@R=1.5pc{
F(H(M)) \ar[r]^-{\eta^{F(H(M))}} \ar[d]_-{H(f)} & G(H(M)) \ar[d]^-{G(H(f))} \\
F(H(N)) \ar[r]_-{\eta^{H(N)}} & G(H(N)),
}
\]
thus $\eta\cdot H$ is a natural transformation.
On the other hand, if we assume that $\eta\cdot H$ is a natural isomorphism and that $s:M \to H(N)$ is an isomorphism
for some $N$ in $\mathcal{Z}$, then by naturality of $\eta$ the following diagram commutes,
\[
 \xymatrix@R=1.5pc{
F(M) \ar[r]^{\eta^{M}} \ar[d]_-{F(s)} & G(M) \ar[d]^-{G(s)} \\
F(H(N)) \ar[r]_-{\eta^{H(N)}} & G(H(N)),   
}
\]
that is, $\eta^{M}=G(s)^{-1}\eta^{H(N)}F(s)$ is composition of isomorphisms and hence an isomorphism. 
This shows that $(a)$ implies $(b)$, the other implication is evident.

\eproof

\begin{lema} \label{(DE)L:levFunUno}
Let $\mathcal{Z}\subset A_0$-mod be a class of $A_0$-modules. Assume there is an additive functor
$F_0:\mathcal{Z} \longrightarrow \mathcal{Z}$ and a natural transformation $\eta_0:F_0 \to Id_{\mathcal{Z}}$.
Observe that the pair $(F_0,\eta_0)$ determines two natural transformations
\[
 \eta_0\cdot F_0,F_0\cdot \eta_0 : F^2_0 \longrightarrow F_0,
\]
where $\eta_0\cdot F_0$ is obtained as in the lemma above and the transformation $F_0 \cdot \eta_0$ is given by
$(F_0\cdot \eta_0)^{M_0}=F_0(\eta_0^{M_0})$ for an object $M_0$ en $\mathcal{Z}$.
Then there are a functor $F=\overline{F_0}:\check{\mathcal{U}}(\mathcal{Z}) \to \check{\mathcal{U}}(\mathcal{Z})$ 
and a natural transformation
$\eta=\overline{\eta_0}:F \to Id_{\check{\mathcal{U}}(\mathcal{Z})}$ such that $F|_{\mathcal{Z}}=F_0$ and $\eta|_{\mathcal{Z}}=\eta_0$.
Moreover, 
\begin{itemize}
 \item[a)]  if $\eta_0\cdot F_0: F_0^2 \longrightarrow F_0$ is an isomorphism of functors then $\eta \cdot F$ is also a
natural isomorphism;
 \item[b)]  if $\eta_0\cdot F_0=F_0 \cdot \eta_0$ then $\eta\cdot F=F \cdot \eta$;
 \item[c)]  if $\eta_0^{M_0}:F_0(M_0) \to M_0$ is an epimorphism for any $M_0$ in $\mathcal{Z}$ then $F$ is a faithful functor; 
 \item[d)]  if in addition to $(c)$ the functor $F_0$ preserves exact sequences in $\mathcal{Z}$ then $F$ preserves 
exact sequences in $\check{\mathcal{U}}(\mathcal{Z})$.
\end{itemize}
\end{lema}

\bproof
We define first the functor $F$. If $M=(M_0,M_{\omega},\gamma_M)$ is an object in $\check{\mathcal{U}}(\mathcal{Z})$
take $F(M)_{\omega}$ as the pull-back of $|\eta^{M_0}_0|$ and $\gamma_M$
\[
 \xymatrix@C=3pc{
F(M)_{\omega} \ar[r]^-{\eta^M_{\omega}} \ar[d]_-{\gamma_{F(M)}} & M_{\omega} \ar[d]^-{\gamma_M} \\
|F_0(M_0)| \ar[r]_-{|\eta^{M_0}_0|} & |M_0|.
}
\]
Take then $F(M)=(F_0(M_0),F(M)_{\omega},\gamma_{F(M)})$
and
\[
\eta^M=(\eta^{M_0}_0,\eta^M_{\omega}):F(M) \to M.
\]
Assume that $f:M \to N$ is a morphism in $\check{\mathcal{U}}(\mathcal{Z})$.
Since $\eta_0$ is a natural transformation, the external paths from $F(M)_{\omega}$ to $|N_0|$ in the following diagram commute 
($\gamma_N f_{\omega} \eta^M_{\omega}=|f_0|\gamma_M \eta^M_{\omega}=|\eta^{N_0}_0||F_0(f_0)|\gamma_{F(M)}$),
hence there exists $F(f)_{\omega}:F(M)_{\omega} \to F(N)_{\omega}$ such that the complete diagram is commutative
\[
 \xymatrix@C=1pc@R=1pc{
& M_{\omega} \ar[rr]^-{f_{\omega}} \ar[dd]^(.3){\gamma_{M}} &  & N_{\omega} \ar[dd]^-{\gamma_N} \\
F(M)_{\omega} \ar@{.>}[rr]_(.7){F(f)_{\omega}} \ar[dd]_-{\gamma_{F(M)}} \ar[ru]^-{\eta^M_{\omega}} & 
& F(N)_{\omega} \ar[dd] \ar[ru]_-{\eta^N_{\omega}} \\
& |M_0| \ar[rr]_(.3){|f_0|} & & |N_0| \\
|F_0(M_0)| \ar[ru]^-{|\eta^{M_0}_0|} \ar[rr]_-{|F_0(f_0)|} & & |F_0(N_0)| \ar[ru]_-{|\eta^{N_0}_0|}
}
\]
Denote by $F(f)$ the morphism between the objects $F(M)$ and $F(N)$ given by
$(F_0(f_0),F(f)_{\omega})$. Assume now that $\xymatrix{L \ar[r]^-{g} & M \ar[r]^-{f} & N}$ are composable
morphisms. By uniqueness in the universal property of the pull-back we have that
\[
F(fg)_{\omega}=F(f)_{\omega}F(g)_{\omega}.
\]
\[
 \xymatrix@C=1pc@R=1pc{
& L_{\omega} \ar[dd] \ar[rr]^-{g_{\omega}} & & M_{\omega} \ar[rr]^-{f_{\omega}} \ar[dd]_(.3){\gamma_{M}} &  
& N_{\omega} \ar[dd]^-{\gamma_N} \\
F(L)_{\omega} \ar@{.>}[rr]^(.7){F(g)_{\omega}} \ar@{.>}@/_.7pc/[rrrr] \ar[ru] \ar[dd] & 
& F(M)_{\omega} \ar@{.>}[rr]^(.7){F(f)_{\omega}} \ar[dd] 
\ar[ru]^-{\eta^M_{\omega}} & & F(N)_{\omega} \ar[dd] \ar[ru]_-{\eta^N_{\omega}} \\
& |L_0| \ar[rr]_(.3){|g_0|} & & |M_0| \ar[rr]_(.3){|f_0|} & & |N_0| \\
|F_0(L_0)| \ar[ru]^-{|\eta^{L_0}_0|} \ar[rr]_-{|F_0(g_0)|} & 
& |F_0(M_0)| \ar[ru]_-{|\eta^{M_0}_0|} \ar[rr]_-{|F_0(f_0)|} & & |F_0(N_0)| \ar[ru]_-{|\eta^{N_0}_0|}
}
\]
It is also clear that $F(I_M)=I_{F(M)}$, thus $F$ is an additive functor. 
By construction we have that $f \eta^M=\eta^N F(f)$, that is, $\eta:F \to Id_{\check{\mathcal{U}}(\mathcal{Z})}$ 
is a natural transformation. It is evident that for any $M_0$ in $\mathcal{Z}$ we have $F(M_0)|_{\mathcal{Z}}=F_0(M_0)$ and 
$\eta^{M_0}|_{\mathcal{Z}}=\eta_0^{M_0}$. 

To show $(a)$ notice first that in a pull-back diagram of $k$-vector spaces,
\[
 \xymatrix{
A \ar[r]^-{x} \ar[d]_-{y} & B \ar[d]^-{g} \\
C \ar[r]_-{f}  & D
}
\]
if $f$ is an isomorphism then $x$ is an isomorphism. Indeed, since $f$ is surjective, the following sequence is exact
\[
 \xymatrix{0 \ar[r]
& A \ar[r]^-{\left[ \begin{smallmatrix} x \\ y \end{smallmatrix} \right]}
& B \oplus C \ar[r]^-{\left[ \begin{smallmatrix} g & -f \end{smallmatrix} \right]} & D \ar[r] & 0.
}
\]
Hence, if $a\in A$ is in the kernel of $x$ then $0=g(x(a))=f(y(a))=0$ and $y(a)=0$ since $f$ is
injective. Then $a$ is in the kernel of $\left[ \begin{smallmatrix} x \\ y \end{smallmatrix} \right]$,
and $a=0$. Since $\dimk_k(B \oplus C)=\dimk_k A +\dimk_kD$ and $\dimk_k C=\dimk_k D$ for $f$ is an isomorphism,
we conclude that $\dimk_k A=\dimk_k B$ and then the injective transformation $x$ is an isomorphism.

Using the equivalence proved in lemma~\ref{(DE)L:isoNat}, we must show that if $M=(M_0,M_{\omega},\gamma_{M})$
is an object in $\check{\mathcal{U}}(\mathcal{Z})$ isomorphic to $F(N)$ for some $N$ in
$\check{\mathcal{U}}(\mathcal{Z})$, then $\eta^M$ is an isomorphism.
Since $M_0 \cong (F(N))_0=F_0(N_0)$, by hypothesis $\eta_0^{M_0}$ is an isomorphism, hence
$|\eta_0^{M_0}|$ is also an isomorphism. Since the following is a pull-back diagram,
\[
 \xymatrix@C=3pc{
F(M)_{\omega} \ar[r]^-{\eta^M_{\omega}} \ar[d]_-{\gamma_{F(M)}} & M_{\omega} \ar[d]^-{\gamma_M} \\
|F_0(M_0)| \ar[r]_-{|\eta^{M_0}_0|} & |M_0|,
}
\]
and $|\eta_0^{M_0}|$ is an isomorphism, by the above $\eta^M_{\omega}$, and thus $\eta^M$, are isomorphisms.
The point $(b)$ is immediate from the definition of $F(\eta^M)=(F_0(\eta^{M_0}_0),F(\eta^M)_{\omega})$, as shown 
in the following diagram
\[
 \xymatrix@C=1pc@R=1pc{
& F(M)_{\omega} \ar[rr]^-{\eta^M_{\omega}} \ar[dd]^(.3){\gamma_{F(M)}} &  & M_{\omega} \ar[dd]^-{\gamma_M} \\
F(F(M))_{\omega} \ar@{.>}[rr]_(.7){F(\eta^M)_{\omega}} \ar[dd]_-{\gamma_{F(F(M))}} \ar[ru]^-{\eta^{F(M)}_{\omega}} & 
& F(M)_{\omega} \ar[dd] \ar[ru]_-{\eta^M_{\omega}} \\
& |F_0(M_0)| \ar[rr]_(.3){|\eta_0^{M_0}|} & & |M_0| \\
|F_0(F_0(M_0))| \ar[ru]^-{|\eta^{F_0(M_0)}_0|} \ar[rr]_-{|F_0(\eta^{M_0}_0)|} & & |F_0(M_0)|. \ar[ru]_-{|\eta^{M_0}_0|}
}
\]
For $F(\eta^M)_{\omega}$ is the only transformation (dotted arrow) that makes the diagram above commutative
and since by hypothesis $F_0(\eta_0^{M_0})=\eta_0^{F_0(M_0)}$ (and hence
$|F_0(\eta_0^{M_0})|=|\eta_0^{F_0(M_0)}|$), we have that $\eta^{F(M)}_{\omega}$ also makes the diagram commutative.
From the uniqueness in the pull-back it follows that $\eta^{F(M)}_{\omega}=F(\eta^M)_{\omega}$ and
\[
 F(\eta^M)=(F_0(\eta_0^{M_0}),F(\eta^M)_{\omega})=(\eta_0^{F_0(M_0)},\eta_{\omega}^{F(M)})=\eta^{F(M)}.
\]

We prove $(c)$. Notice that in a pull-back diagram of $k$-vector spaces if $f$ is an epimorphism then $x$ is an epimorphism. 
\[
 \xymatrix{
k \ar@/^1pc/[rrd]^-{b} \ar@/_1pc/[ddr]_-{c} \ar@{.>}[rd]^(.6){a} \\
& A \ar[r]^-{x} \ar[d]_-{y} & B \ar[d]^-{g} \\
& C \ar[r]_-{f}  & D
}
\]
Indeed, for an arbitrary $b \in B$ take $c \in C$ such that $f(c)=g(b)$ (this is possible for $f$ is an epimorphism).
Define the transformations $k \to B$ and $k \to C$ given by $1 \mapsto b$ and $1 \mapsto c$. Then
the external square in the diagram above is commutative, thus there exists $a:k \to A$ such that $x(a)=b$. Then $x$ 
is an epimorphism. By definition of $\eta^M=(\eta_0^{M_0},\eta_{\omega}^M)$ as pull-back and since $|\eta_0^{M_0}|$ is
an epimorphism ($|-|$ is an exact functor) then $\eta_{\omega}^M$ is also an epimorphism.
Since $f \eta^M=\eta^N F(f)$ and $\eta^M$ is an epimorphism we have that $F$ is a faithful functor.
 
To verify $(d)$ assume that we have en exact sequence 
\[
\varepsilon:\xymatrix{0 \ar[r] & L \ar[r]^-{u} & M \ar[r]^-{v} & N \ar[r] & 0}
\] 
in $\check{\mathcal{U}}(\mathcal{Z})$. By naturality of $\eta$ the following diagram is commutative
\[
\xymatrix{
0 \ar[r] & L \ar[r]^-{u} & M \ar[r]^-{v} & N \ar[r] & 0\\
0 \ar[r] & F(L) \ar[r]_-{F(u)} \ar[u]^-{\eta^L} & F(M) \ar[r]_-{F(v)} \ar[u]^-{\eta^M} & F(N) \ar[r] \ar[u]^-{\eta^N} & 0.
}
\]
Going to kernels we obtain the following diagram where all columns are exact, as well as the rows in the bottom grid.
In the upper grid the first row is exact and the third one is isomorphic to the third row in the bottom, 
by definition of $F$ through pull-backs.  By the nine lemma, the second upper row is also exact.
This shows that $F(\varepsilon)$ is an exact sequence.
\[
\xymatrix@!0@C=2.5pc@R=2.5pc{
& & & & & 0 & & 0 & & 0 \\
& & & & & 0 & & 0 & & 0 \\
& & 0 \ar[rr] & & L_{\omega} \ar[d] \ar[rr] \ar[ruu] & & M_{\omega} \ar[d] \ar[rr] \ar[ruu] & & N_{\omega} \ar[d] \ar[rr] \ar[ruu] & & 0 \\
& & 0 \ar[rr] & & |L_0| \ar[rr] \ar[ruu] & & |M_0| \ar[rr] \ar[ruu] & & |N_0| \ar[rr] \ar[ruu] & & 0 \\
& 0 \ar[rr] & & F(L)_{\omega} \ar[d] \ar[rr] \ar[ruu] & & F(M)_{\omega} \ar[d] 
\ar[rr] \ar[ruu] & & F(N)_{\omega} \ar[d] \ar[rr] \ar[ruu] & & 0 \\
& 0 \ar[rr] & & |F_0(L_0)| \ar[rr] \ar[ruu] & & |F_0(M_0)| \ar[rr] \ar[ruu] & & |F_0(N_0)| \ar[rr] \ar[ruu] & & 0 \\
0 \ar[rr] & & K^L \ar@{=}[d] \ar[rr] \ar[ruu] & & K^M \ar@{=}[d] \ar[rr] \ar[ruu] & & K^N \ar@{=}[d] \ar[rr] \ar[ruu] & & 0 \\
0 \ar[rr] & & K^L \ar[rr] \ar[ruu] & & K^M \ar[rr] \ar[ruu] & & K^N \ar[rr] \ar[ruu] & & 0 \\
& 0 \ar[ruu] & & 0 \ar[ruu] & & 0 \ar[ruu]  \\
& 0 \ar[ruu] & & 0 \ar[ruu] & & 0 \ar[ruu] 
}
\]
\eproof

\begin{lema} \label{(DE)L:levFunDos}
Let $\mathcal{Z} \subseteq A_0$-mod be a class of $A_0$-modules.
Assume we have an additive functor
$G_0:\mathcal{Z} \longrightarrow \mathcal{Z}$ and a natural transformation of functors $\eta_0:Id_{\mathcal{Z}} \to G_0$.
There are two natural transformations
\[
 \eta_0\cdot G_0,G_0 \cdot \eta_0: G_0 \longrightarrow G^2_0.
\]
Then there exist an additive functor $G=\overline{G_0}:\check{\mathcal{U}}(\mathcal{Z}) \to \check{\mathcal{U}}(\mathcal{Z})$ 
and a natural transformation $\eta=\overline{\eta_0}:Id_{\check{\mathcal{U}}(\mathcal{Z})} \to G$ such that $G|_{\mathcal{Z}}=G_0$ 
and $\eta|_{\mathcal{Z}}=\eta_0$. Moreover, 
\begin{itemize}
 \item[a)] if $\eta_0 \cdot G_0: G_0 \longrightarrow G_0^2$ is a natural isomorphism then $\eta \cdot G$ is
also an isomorphism of functors;
\item[b)] if $\eta_0 \cdot G_0=G_0 \cdot \eta_0$ then $\eta \cdot G=G\cdot \eta$;
 \item[c)] if $\eta_0^{M_0}:M_0 \to G_0(M_0)$ is a monomorphism for each $M_0$ in $\mathcal{Z}$ then $G$ is a faithful functor;
 \item[d)] if $G_0$ preserves exact sequences in $\mathcal{Z}$ then $G$ preserves exact sequences in $\check{\mathcal{U}}(\mathcal{Z})$.
\end{itemize}
\end{lema}

\bproof
For an object $M=(M_0,M_{\omega},\gamma_{M})$ in $\check{\mathcal{U}}(\mathcal{Z})$ define 
\[
G(M)=(G_0(M_0),M_{\omega},\gamma_{G(M)})
\]
where $\gamma_{G(M)}=|\eta_o^{M_0}|\gamma_M$. For a morphism
$g=(g_0,g_{\omega}):M \to N$ take $G(g)=(G_0(g_0),g_{\omega})$, see the following diagram,
\[
 \xymatrix{
& M_{\omega} \ar[rr]^-{g_{\omega}} \ar[d]^-{\gamma_M} \ar@/_1pc/[ddl]_(.4){\gamma_{G(M)}} &  
& N_{\omega} \ar[d]^-{\gamma_N} \ar@/_1pc/[ddl]_(.4){\gamma_{G(N)}} \\
& |M_0| \ar[rr]_(.4){|g_0|} \ar[dl]^-{|\eta_0^{M_0}|} & & |N_0| \ar[dl]^-{|\eta_0^{N_0}|} \\
|G_0(M_0)| \ar[rr]_-{|G_0(g_0)|} & & |G_0(N_0)|.
}
\]
We observe directly that $G$ is an additive functor. On the other hand, the mapping $\eta^M=(\eta_0^{M_0},Id_{M_{\omega}})$ 
satisfies $G(g) \eta^M=\eta^N g$, therefore $\eta$ is a natural transformation. 

To prove $(a)$, by lemma~\ref{(DE)L:isoNat} it is enough to show that if $M$ is isomorphic to an object of the form
$G(N)$ then $\eta^M$ is an isomorphism. Since $\eta^M=(\eta_0^{M_0},Id_{M_{\omega}})$, it is then enough to verify that
$\eta_0^{M_0}$ is an isomorphism, which is direct consequence of the hypothesis for $M_0 \cong G_0(N_0)$.

The proof of $(b)$ is also direct consequence of the definition of $G(\eta^M)$, for by the hypothesis
we have (see the diagram below),
\[
 G(\eta^M)=G(\eta^{M_0}_0,Id_{M_{\omega}})=(G_0(\eta^{M_0}_0),Id_{M_{\omega}})
=(\eta^{G_0(M_0)}_0,Id_{M_{\omega}})=\eta^{G(M)},
\]
\[
 \xymatrix{
& M_{\omega} \ar[rr]^-{Id_{M_{\omega}}} \ar[d]^-{\gamma_M} \ar@/_1pc/[ddl]_(.4){\gamma_{G(M)}} &  
& G(M)_{\omega}=M_{\omega} \ar[d]^-{\gamma_{G(M)}=|\eta_0^{M_0}|\gamma_M} \ar@/_1pc/[ddl]_(.4){\gamma_{G(G(N))}} \\
& |M_0| \ar[rr]_(.4){|\eta^{M_0}_0|} \ar[dl]^-{|\eta_0^{M_0}|} & & |G_0(M_0)| \ar[dl]^-{|\eta_0^{G_0(M_0)}|} \\
|G_0(M_0)| \ar[rr]_-{|G_0(\eta^{M_0}_0)|} & & |G_0(G_0(M_0))|.
}
\]
To show $(c)$, it is clear that if $\eta_0^{M_0}$ is a monomorphism then $\eta^M=(\eta_0^{M_0},Id_{M_{\omega}})$ is a monomorphism.
Since $\eta^M g=G(g) \eta^N $ and $\eta^M$ is a monomorphism we have that $G$ is a faithful functor. 

To verify $(d)$, if $\varepsilon:\xymatrix{0 \ar[r] & L \ar[r] & M \ar[r] & N \ar[r] & 0}$ is an exact sequence in 
$\check{\mathcal{U}}(\mathcal{Z})$, its image under $G$ has the form
\[
 \xymatrix{
0 \ar[r] & L_{\omega} \ar[r] \ar[d]^-{\gamma_{G(L)}} & M_{\omega} \ar[r] \ar[d]^-{\gamma_{G(M)}} 
& N_{\omega} \ar[r] \ar[d]^-{\gamma_{G(N)}} & 0 \\
0 \ar[r] & |G_0(L_0)| \ar[r] & |G_0(M_0)| \ar[r] & |G_0(N_0)| \ar[r] & 0.
}
\]
By definition of exactness in $\varepsilon$ the upper row in the diagram is exact.
Since $G_0$ preserves exact sequences and $|-|$ is an exact functor, the lower row is also exact,
hence $G(\varepsilon)$ is an exact sequence in $\check{\mathcal{U}}(\mathcal{Z})$.
\eproof

\section{Universal mappings and rank.} \label{(DE)S:rank}
%----------------------------------------------------------------------
%-------------------------------------
Let $A_0=T_{S_0}(L_0)$ be the path algebra of a Dynkin diagram $\Delta$ 
(where $L_0$ is the $S_0$-$S_0$-bimodule of arrows in $\Delta$)
and $A=A_0[R]\cong k\widetilde{\Delta}$ is the one-point extension of $A_0$ by the module $R$ (as in lemma~\ref{(DE):extEucl}).
Let $W_0$ be an indecomposable $A_0$-module with maximal dimension vector.
For an $A_0$-module $M_0$ denote by $\sharp M_0$ the number of indecomposable summands of $M_0$.
Let $m_{M_0}$ be the multiplicity of $W_0$ in $M_0$ and $e_{M_0}$ the number of indecomposable summands 
of $M_0$ not isomorphic to $W_0$ (so that $\sharp M_0=m_{M_0}+e_{M_0}$). Define the \textbf{rank} of an indecomposable 
$A$-module $M=(M_0,M_{\omega},\gamma_{M})$ as
\[
 \Rnk(M)=\dimk_k M_{\omega} - m_{M_0}.
\]
One of our main interests is to determine the behavior of the rank with respect to the Auslander-Reiten translations
in $A$-mod.  We will achieve that by means of the universal mappings defined by Ringel in
\cite[section 3.4(17)]{cmR} and their liftings to the category $A$-mod.

%----------------------------------------
The following lemma and corollary~\ref{(DE)C:ContComp} correspond to theorem~11.5 in Ga\-briel and Roiter \cite{GR97}.

\begin{lema} \label{(DE)L:subKro}
\textbf{Kronecker subcategory.}
Let $B$ be the subalgebra of $A=A_0[R]$ given by $B=A_0[0]$ and
$W$ the $B$-module $W=W_0\oplus S(\omega)$, where $W_0$ is an indecomposable $A_0$-module with maximal dimension vector 
and $S(\omega)$ is the simple module of vertex $\omega$. Consider the ditalgebra $\mathcal{A}=(A,0)$ with zero differential. 
Then the reduced ditalgebra $\mathcal{A}^W$ is isomorphic to the classical Kronecker algebra $A_2$ (with zero differential). 
The image of the functor associated to the reduction $F^W$ corresponds to the subalgebra $\check{\mathcal{U}}(W_0)$ of $A$-mod.
The images of the posprojective and preinjective components of $\mathcal{A}^W$-mod under $F^W$ will be denoted by 
$\mathcal{P}_0$ and $\mathcal{I}_0$ respectively.
\end{lema}
\bproof
Since $W_0$ is exceptional and $Z_0=\End_{A_0}(W_0)^{op}\cong k$ (lemma~\ref{(P)L:endExc}), the $A_0$-module $W_0$ is admissible.
The reduced algebra $A_0^{W_0}=T_{Z_0}(0)$ is isomorphic to the field $k$ (with zero differential).  
By proposition~\ref{(DE)P:extPuntRed} and lemma~\ref{(DE)L:extTens} the reduced tensor algebra $A^W$ is isomorphic to 
\[
A_0^{W_0}[R^{W_0}] \cong k[R^{W_0}] =T_{\left( \begin{smallmatrix} k&0\\0&k \end{smallmatrix} \right)} 
\left( \begin{matrix} 0& \Hom_{A_0}(R,W_0)^* \\0&0 \end{matrix} \right). 
\]
Now, by lemma~\ref{(DE)L:valores}, the $k$-vector space $\Hom_{A_0}(R,W_0)^*$ has dimension two,
therefore there exists an isomorphism between the extension $A_0^{W_0}[R^{W_0}]$ and the classical Kronecker algebra 
$\widetilde{\mathbf{A_1}}=A_2$. Finally observe that the reduced differential $\delta^W$ is zero, due to 
proposition~\ref{(DE)P:extPuntRed} and since the morphism $\varepsilon$ in lemma~\ref{(DE)L:epsilonBimod} is zero
(for $P_0=0$).
\eproof

The \textbf{universal mappings} to be considered are the functors
\[
F_0= W_0 \otimes_k \Hom_{A_0}(W_0,-): A_0\text{-mod} \longrightarrow A_0 \text{-mod},
\]
\[
G_0= W_0 \otimes_k D\Hom_{A_0}(-,W_0): A_0\text{-mod} \longrightarrow A_0 \text{-mod}.
\]
Their liftings, which are projections to the Kronecker subcategory $\check{\mathcal{U}}(W_0)$, determine
the behavior of the rank in the preinjective and posprojective components in $A$-mod.

%-------------------------------------
\begin{lema} \label{(DE)L:FunUno}
Consider the functor
\[
F_0= W_0 \otimes_k \Hom_{A_0}(W_0,-): A_0\text{-mod} \longrightarrow A_0 \text{-mod}.
\]
Given a module $M_0$ in the category $A_0$-mod define the \textbf{right universal mapping}
$\eta^{M_0}_0:F_0(M_0) \to M_0$ as the evaluation $w \otimes f \mapsto f(w)$. 
Then $\eta_0:F_0 \to Id_{A_0\text{-mod}}$ is a natural transformation and
\begin{itemize}
 \item[a)] $\eta_0\cdot F_0$ is an isomorphism of functors,
 \item[b)] $\eta_0\cdot F_0=F_0\cdot \eta_0$.
\end{itemize} 
Restrict the functor $F_0$ and the transformation $\eta_0$ to the subcategory 
$\mathcal{Z}=W_0 \vee \mathcal{Y}$ of $A_0$-mod, where $\mathcal{Y}$ is the class of $A_0$-modules
generated by representatives of the set $\ConR([W_0])=\{ [M] \in \Ind A_0 \; | \; [W_0] \prec [M] \}$.
Then
\begin{itemize}
 \item[c)] $\eta_0^{M_0}$ is an epimorphism for each $M_0$ in $\mathcal{Z}$, and
 \item[d)] $F_0$ preserves exact sequences in $\mathcal{Z}$.
\end{itemize}
\end{lema}
\bproof
Observe that $\eta_0$ is a natural transformation, for
\[
\eta_0^{N_0}(F_0(f)(w \otimes h))=\eta_0^{N_0}(w \otimes fh)=f(h(w))=f(\eta_0^{M_0}(w \otimes h)),
\]
that is, the following is a commutative diagram,
\[
 \xymatrix{
F_0(M_0) \ar[d]_-{\eta_0^{M_0}} \ar[r]^-{F_0(f)} & F_0(N_0) \ar[d]^-{\eta_0^{N_0}} \\
M_0 \ar[r]_-{f} & N_0. 
}
\]
Consider now the following transformations
\[
 \xymatrix{
\End_{A_0}(W_0) \otimes_k \Hom_{A_0}(W_0,M_0) \ar[dd]_-{\Phi} \ar[rd]^-{\Psi} & \\
& \Hom_{A_0}(W_0,M_0) \\
\Hom_{A_0}(W_0,W_0\otimes_k \Hom_{A_0}(W_0,M_0))
\ar[ru]_(.6){\Hom(W_0,\eta_0^{M_0})} 
}
\]
given for morphisms $f\in \End_{A_0}(W_0)$, $g \in \Hom_{A_0}(W_0,M_0)$ and an element $w$ in $W_0$ in the following way:
$\Phi(f\otimes g)(w)=f(w)\otimes g$ and $\Psi(f\otimes g)=gf$. Notice that the diagram above is commutative, 
for by definition of $\eta_0^{M_0}$ we have $\eta_0^{M_0}(\Phi(f\otimes g)(w))=\eta_0^{M_0}(f(w)\otimes g)=g(f(w))
=\Psi(f \otimes g)(w)$. In this way, applying the functor $W_0\otimes -$ to the diagram above we have the 
following commutative diagram
\[
 \xymatrix{
W_0 \otimes_k \End_{A_0}(W_0) \otimes_k \Hom_{A_0}(W_0,M_0) \ar[dd]_-{Id_{W_0} \otimes \Phi} \ar[rd]^-{Id_{W_0} \otimes \Psi} & \\
& W_0 \otimes_k \Hom_{A_0}(W_0,M_0) \\
W_0 \otimes_k \Hom_{A_0}(W_0,W_0\otimes_k \Hom_{A_0}(W_0,M_0))
\ar[ru]_(.6){F_0(\eta_0^{M_0})} 
}
\]
since $F_0(\eta_0^{M_0})=Id_{W_0}\otimes \Hom_{A_0}(W_0,\eta_0^{M_0})$. We observe that the diagram above is also 
commutative when substituting $F_0(\eta_0^{M_0})$ by the morphism $\eta_0^{F_0(M_0)}$.
Indeed, on the one hand, using the notation $f$, $g$ and $w$ given above, we have
\[
 \eta_0^{F_0(M_0)}([Id_{W_0} \otimes \Phi](w\otimes f\otimes g))=\eta_0^{F_0(M_0)}(w\otimes \Phi(f\otimes g))
=\Phi(f\otimes g)(w)=f(w)\otimes g.
\]
On the other hand, $[Id_{W_0}\otimes \Psi](w\otimes f \otimes g)=w \otimes \Psi(f\otimes g)=w\otimes gf$,
and since the endomorphism algebra $\End_{A_0}(W_0)$ is isomorphic to the field $k$, the equality
$f(w)\otimes g=w\otimes gf$ is evident. Then $\eta_0^{F_0(M_0)}=[Id_{W_0}\otimes \Psi][Id_{W_0}\otimes \Phi]^{-1}$
is an isomorphism and $\eta_0^{F_0(M_0)}=F_0(\eta_0^{M_0})$. This proves $(a)$ and $(b)$.

Clearly $\eta_0^{M_0}$ is an epimorphism for modules $M_0$ in $\mathcal{Z}= W_0 \vee \mathcal{Y}$ (point $(c)$)
for the modules in $\mathcal{Y}$ are generated by $W_0$ (lemma~\ref{(DE)L:alaDyn}).
 
To show that $F_0$ preserves exact sequences assume that $\varepsilon$ is an exact sequence in $\mathcal{Z}$,
\[
\varepsilon= \xymatrix{
0 \ar[r] & M_0 \ar[r] & E_0 \ar[r] & N_0 \ar[r] & 0.
}
\]
We use the Auslander-Reiten formulas in the hereditary case 
\[\Ext_{A_0}^1(W_0,M_0)\cong D\Hom_{A_0}(M_0,\tau W_0)\] and the fact that $\Hom_{A_0}(M_0,\tau W_0)=0$ since 
$M_0 \in \mathcal{Y}$ and $\tau W_0 \in \mathcal{X}$. Then $\Hom(W_0,\varepsilon)$ is an exact sequence
\[
\xymatrix@C=1pc{
0 \ar[r] & \Hom_{A_0}(W_0,M_0) \ar[r] & \Hom_{A_0}(W_0,E_0) \ar[r] & \Hom_{A_0}(W_0,N_0) \ar[r] & 0.
}
\]
Since $Id_{W_0}\otimes -$ es an exact functor we conclude that $F_0(\varepsilon)=Id_{W_0} \otimes \Hom(W_0,\varepsilon)$ 
is an exact sequence, so $(d)$ holds.
\eproof

Denote by $(\overleftarrow{},\eta)$ the liftings of $(F_0,\eta_0)$ to $\check{\mathcal{U}}(W_0 \vee \mathcal{Y})$. 
By lemma~\ref{(DE)L:levFunUno} the functor 
\[
F=\overline{W_0 \otimes_k \Hom_{A_0}(W_0,-)}:\check{\mathcal{U}}(W_0 \vee \mathcal{Y}) 
\longrightarrow \check{\mathcal{U}}(W_0)
\] 
is additive, faithful, preserves exact sequences and its restriction to $\check{\mathcal{U}}(W_0)$ is full and dense.
\begin{lema} \label{(DE)L:ARUno}
Assume that $\xymatrix{ 0 \ar[r] & \tau Q \ar[r]^-{g} & E \ar[r]^-{f} & Q \ar[r] & 0 }$
is an almost split sequence in $\check{\mathcal{U}}(W_0 \vee \mathcal{Y})$. 
If $Q \in \check{\mathcal{U}}(W_0)$ then $\overleftarrow{f}$ is a right almost split morphism in 
$\check{\mathcal{U}}(W_0)$. On the contrary $\overleftarrow{f}$ is a retraction.
\end{lema}

\bproof
Consider the universal projection
\[
 \xymatrix{ 
0 \ar[r] &  \overleftarrow{\tau Q} \ar[r]^-{\overleftarrow{g}} \ar[d]^-{\eta^{\tau Q}} 
& \overleftarrow{E} \ar[d]^-{\eta^E} \ar[r]^-{\overleftarrow{f}} 
& \overleftarrow{Q} \ar@{.>}[dl]^-{t} \ar[d]^-{\eta^Q} \ar[r] & 0 \\
0 \ar[r] & \tau Q \ar[r]_-{g} & E \ar[r]_-{f} & Q \ar[r] & 0. 
}
\]
If $Q$ is not in $\check{\mathcal{U}}(W_0)$ then $\eta^Q$ is not a retraction, hence there exists $t:\overleftarrow{Q} \to E$ 
such that $ft=\eta^Q$. Then $\overleftarrow{\eta^Q}=\overleftarrow{f}\overleftarrow{t}$. 
By the points $(a)$ and $(b)$ in lemmas~\ref{(DE)L:FunUno} and~\ref{(DE)L:levFunUno} we have that 
$\overleftarrow{\eta^Q}=\eta^{\overleftarrow{Q}}$ is an isomorphism. Hence $f$ is a retraction. 

Assume now that $Q$ is in $\check{\mathcal{U}}(W_0)$.
Then $\eta^Q$ is an isomorphism and $(\eta^Q)^{-1}f$ is not a retraction, for $f$ is not a retraction
and $(\eta^Q)^{-1}fs=Id_{\overleftarrow{Q}}$ if and only if $fs(\eta^Q)^{-1}=Id_{Q}$. 
Therefore $\overleftarrow{f}=(\eta^Q)^{-1} f\eta^E$ is not a retraction.
Assume on the other hand that $V$ is an object in $\check{\mathcal{U}}(W_0)$ and that 
$v: V \to \overleftarrow{Q}$ is not a retraction,
\[
 \xymatrix{
& \overleftarrow{\overleftarrow{Q}} \ar[d]_-{\eta^{\overleftarrow{Q}}} & \overleftarrow{V} \ar[d]^-{\eta^V} \ar[l]_-{\overleftarrow{v}} \\
\overleftarrow{E} \ar[r]_-{\overleftarrow{f}} \ar[d]_-{\eta^W} & \overleftarrow{Q} \ar[d]^-{\eta^Q} & V \ar[l]^-{v} \\
E \ar[r]_-{f} & Q.
}
\]
Since $\eta^Q v$ is not a retraction and $f$ is right almost split, there exists $s:V \to E$ such that 
$\eta^Q v=fs$. Then, using again that $\eta^{\overleftarrow{Q}}=\overleftarrow{\eta^Q}$ (lemma~\ref{(DE)L:levFunUno}$(b)$) 
and the commutativity of the diagram above, we have that 
$v\eta^V=\eta^{\overleftarrow{Q}}\overleftarrow{v}=\overleftarrow{\eta^Q}\overleftarrow{v}=\overleftarrow{f}\overleftarrow{s}$, 
thus $v=\overleftarrow{f}\overleftarrow{s}(\eta^V)^{-1}$. That is, $\overleftarrow{f}$ is right almost split.
\eproof

For the following proposition denote by $J^i$ the indecomposable injective $A$-modules of the form $D(e_iA)$. 
Observe that the indecomposable injectives in the Kronecker subcategory
$\check{\mathcal{U}}(W_0)$ are given by $I^{\omega}=(0,k,0)$ and $I=\overline{W_0}$.

\begin{lema} \label{(DE)P:Uno}
Let $w_0=\vdim W_0$ be the dimension vector of $W_0$. Then
\begin{itemize}
 \item[a)] $\overleftarrow{J^{\omega}} \cong I^{\omega}$ and $\overleftarrow{J^i} \cong (w_0)_iI$ for $i=1,\ldots,n$.
 \item[b)] For any $M \in \mathcal{I}$ we have $\overleftarrow{M} \in \mathcal{I}_0$.
 \item[c)] The function $\sharp(\overleftarrow{\cdot}):\mathcal{I} \to \mathbb{N}$ is additive.
\end{itemize} 
\end{lema}

\bproof
$(a)$ Since $J^{\omega}=E_{\omega}=I^{\omega}$ is the simple module of vertex $\omega$, we have
$\overleftarrow{J^{\omega}}=I^{\omega}$.
For $i=1,\ldots , n$ there exist isomorphisms $D(W_0)e_i \cong D(e_iW_0)$. Thus
\begin{eqnarray}
 \dimk_k \Hom_{A_0}(W_0,J^i_0) & = & \dimk_k \Hom(W_0,D(e_iA_0)) = \nonumber \\
& = & \dimk_k \Hom(e_iA_0,D(W_0)) = \nonumber \\
& = & \dimk_k D(W_0)e_i = \dimk_k e_iW_0 = \nonumber \\
& = & (w_0)_i. \nonumber 
\end{eqnarray}
Then the multiplicity of $W_0$ in $\overleftarrow{(J^i)}_0$ is $m_{\overleftarrow{(J^i)}_0}=(w_0)_i$ 
and using lemma~\ref{(DE)L:valores} we have 
$\dimk_k |\overleftarrow{(J^i)}_0|=\dimk_k \Hom_{A_0}(R,\overleftarrow{(J^i)}_0)=2(w_0)_i$.
Since $J^i$ is injective, $\overleftarrow{J^i}$ is obtained through a pull-back of the form 
\[
 \xymatrix{
(\overleftarrow{J^i})_{\omega} \ar[r] \ar[d]_-{\gamma_{\overleftarrow{J^i}}} & J^i_{\omega} \ar@{=}[d] \\
|(\overleftarrow{J^i})_0| \ar[r]_-{|\eta_0^{J^i_0}|} & |J^i_0|.
}
\]
Since $|\eta_0^{J^i_0}|$ is an epimorphism, by the argument given in the proof of point $(c)$ in lemma~\ref{(DE)L:levFunUno},
we have that $\gamma_{\overleftarrow{J^i}}$ is an isomorphism. Hence $\overleftarrow{J^i} \cong (w_0)_iI$.
The point $(b)$ is inductively obtained from the lemma above as follows. Consider the ordering of indecomposable 
preinjective modules given in lemma~\ref{(DE)L:orden}, $M_1,M_2,M_3,\ldots$
As base case, by the point ($a$) the claim is valid for the indecomposable injective $A$-modules
(which are the first ones in the list). Assume that the statement is valid for the preinjectives $M_1,\ldots,M_i$ and let us show it
for $M_{i+1}$. By the properties of this ordering, the almost split sequence that starts in $M_{i+1}$,
\[
\varepsilon : \xymatrix{0 \ar[r] & M_{i+1} \ar[r]^-{g} & E \ar[r]^-{f} & M_{\ell}  \ar[r] & 0},
\]
satisfies $\ell<i+1$ and for any $M_j$ isomorphic to a direct summand of $E$ we have $j<i+1$.
By lemma~\ref{(DE)L:ARUno}, if $M_{\ell}$ is not in the Kronecker subcategory $\check{\mathcal{U}}(W_0)$ 
then $\overleftarrow{\varepsilon}$ is a trivial sequence. Therefore $\overleftarrow{M_{i+1}}\in \mathcal{I}_0$ since 
by induction hypothesis $\overleftarrow{E} \in \mathcal{I}_0$.
Assume now that $M_{\ell}$ is an object in the subcategory $\check{\mathcal{U}}(W_0)$. By lemmas~\ref{(DE)L:ARUno} 
and~\ref{(DE)L:minimal}, $\overleftarrow{\varepsilon}$ is isomorphic to a sequence of the form
\[
 \overleftarrow{\varepsilon} \cong \xymatrix@C=3.5pc{
0 \ar[r] & P \oplus E'' \ar[r]^-{\left[ \begin{smallmatrix} \Ker f' & 0 \\0 &  I_{E''} \end{smallmatrix} \right]} & E' \oplus E'' 
\ar[r]^-{[f'\; 0]} & M_{\ell} \ar[r] & 0,
}
\]
where $\xymatrix@C=1.7pc{0 \ar[r] & P \ar[r] & E' \ar[r]^-{f'} & M_{\ell} \ar[r] & 0}$ is an almost split sequence in
$\check{\mathcal{U}}(W_0)$, thus $P$ is preinjective. By induction $E''$ is preinjective in $\check{\mathcal{U}}(W_0)$,
hence $\overleftarrow{M_{i+1}} \cong P \oplus E'' \in \mathcal{I}_0$.

The point $(c)$ is evident for almost split sequences $\varepsilon$ such that $\overleftarrow{\varepsilon}$ is trivial.
Whenever $\overleftarrow{\varepsilon}$ is not a trivial sequence, the final term $M_{\ell}$ in $\varepsilon$ is isomorphic to 
$\overleftarrow{M_{\ell}}$ and again by lemma~\ref{(DE)L:minimal},
$\overleftarrow{\varepsilon}$ is isomorphic to a sequence as above.
But the Auslander-Reiten sequences $\xymatrix{0 \ar[r] & P \ar[r] & E' \ar[r]^-{f'} & M_{\ell} \ar[r] & 0}$ 
in the preinjective component of the Kronecker subcategory satisfy $\sharp M_{\ell} + \sharp P= \sharp E'$
(see description of the preinjective component $\mathcal{I}_0$ at the end of section~\ref{(P)S:KroDos}). 
Since $\overleftarrow{}$ is an additive functor,
\[
\sharp M_{i+1} + \sharp M_{\ell}=\sharp E'' + \sharp P + \sharp M_{\ell} 
= \sharp E'' + \sharp E'=\sharp E,
\]
that is, $\sharp{(\overleftarrow{\cdot}})$ is an additive function in $\mathcal{I}$.
\eproof

%-------------------------------------
\begin{lema} \label{(DE)L:FunDos}
Consider now the functor 
\[
G_0=W_0 \otimes_k D\Hom_{A_0}(-,W_0):A_0\text{-mod} \longrightarrow A_0\text{-mod}.
\]
For a module $M_0$ in the category $A_0$-mod fix a basis $\{\alpha_1 ,\ldots,\alpha_u\}$ of the $k$-vector space
$\Hom_{A_0}(M_0,W_0)$ and define $\eta_0^{M_0}:M_0 \to G_0(M_0)$ as $m \mapsto \sum_{i=1}^u\alpha_i(m)\otimes \alpha_i^*$. 
Then $\eta_0:Id_{A_0\text{-mod}}\to G_0$ is a well defined natural transformation and
\begin{itemize}
 \item[a)] $\eta_0\cdot G_0$ is an isomorphism of functors,
 \item[b)] $\eta_0\cdot G_0=G_0\cdot \eta_0$.
\end{itemize}
Restrict the functor $G_0$ and the transformation $\eta_0$ to the subcategory 
$\mathcal{Z}=\mathcal{X} \vee W_0$ of $A_0$-mod, where $\mathcal{X}$ is the class of $A_0$-modules
generated by representatives of the set $\ConL([W_0])=\{ [M] \in \Ind A_0 \; | \; [M] \prec [W_0] \}$.
Then
\begin{itemize}
 \item[c)] $\eta_0^{M_0}$ is a monomorphism for any $M_0$ in $\mathcal{Z}$, and
 \item[d)] $G_0$ preserves exact sequences in $\mathcal{Z}$.
\end{itemize}
\end{lema}
\bproof
Let $f_0:M_0 \to N_0$ be a morphism and fix a basis $\{ \beta_1,\ldots,\beta_v \}$ of the space $\Hom_{A_0}(N_0,W_0)$. Since
$\beta_jf_0 \in \Hom_{A_0}(M_0,W_0)$ for all $j$, there exist $\lambda_{ij} \in k$ such that 
$\beta_jf_0=\sum_{i=1}^u\lambda_{ij} \alpha_i$. On the other hand, the morphism
\[
 D\Hom_{A_0}(f_0,W_0):D\Hom_{A_0}(M_0,W_0) \to D\Hom_{A_0}(N_0,W_0),
\]
is given in the basis $\{\alpha_i^*\}_{i=1}^u$ of $D\Hom_{A_0}(M_0,W_0)$ by
$\alpha_i^* \mapsto [g \mapsto \alpha_i^*(gf_0)]$. Notice that
$D\Hom_{A_0}(f_0,W_0)[\alpha_i^*]=\sum_{j=1}^v\lambda_{ij}\beta_j^*$, for when evaluating in the basis $\{\beta_j\}_{j=1}^v$ we have
\[
 D\Hom_{A_0}(f_0,W_0)[\alpha_i^*](\beta_j)=\alpha_i^*(\beta_jf_0)=\alpha_i^*(\sum_{\ell=1}^u\lambda_{\ell j}\alpha_{\ell})
=\lambda_{ij}.
\]
We want to verify that the following diagram commutes
\[
\xymatrix{
 M_0 \ar[d]_-{\eta_0^{M_0}} \ar[r]^-{f_0} & N_0 \ar[d]^-{\eta_0^{N_0}} \\
G_0(M_0) \ar[r]_-{G_0(f_0)} & G_0(N_0).
}
\]
For a $m \in M_0$ we have
\begin{eqnarray}
 G_0(f_0)\eta_0^{M_0}(m) & = & [Id_{W_0} \otimes D\Hom(f_0,W_0)] \left( \sum_{i=1}^u\alpha_i(m)\otimes \alpha_i^* \right)= \nonumber \\
& = & \sum_{i=1}^u \alpha_i(m) \otimes \left( \sum_{j=1}^v \lambda_{ij}\beta_j^* \right) 
= \sum_{j=1}^v\left( \sum_{i=1}^u \lambda_{ij} \alpha_i(m) \right) \otimes \beta_j^* = \nonumber \\
& = & \sum_{j=1}^v \beta_j(f_0(m)) \otimes \beta_j^* =\eta_0^{N_0}(f_0(m)). \nonumber
\end{eqnarray}
In particular the morphism $\eta_0^{M_0}$ does not depend on the choice of basis.
Then the mapping $\eta_0:Id_{A_0\text{-mod}} \to G_0$ is a well defined natural transformation.
To show the points $(a)$ and $(b)$ consider the morphisms
\[
 \xymatrix{
\End_{A_0}(W_0) \otimes_k \Hom_{A_0}(M_0,W_0) \ar[dd]_-{\widetilde{\Phi}} \ar[rd]^-{\widetilde{\Psi}} & \\
& \Hom_{A_0}(M_0,W_0) \\
\Hom_{A_0}(W_0\otimes_k D\Hom_{A_0}(M_0,W_0),W_0)
\ar[ru]_(.6){\Hom(\eta_0^{M_0},W_0)} 
}
\]
defined in the following way. For an element $f$ in $\End_{A_0}(W_0)$ and a morphism $g$ in $\Hom_{A_0}(M_0,W_0)$ define
$\widetilde{\Psi}(f \otimes g)=fg$ and $\widetilde{\Phi}(f \otimes g)(w \otimes h)=h(fg)w$ for $w \in W_0$ and 
$h \in D\Hom_{A_0}(M_0,W_0)$. Since the algebra $\End_{A_0}(W_0)$ is isomorphic to the field $k$, both $\widetilde{\Psi}$ and
$\widetilde{\Phi}$ are isomorphisms. We verify that the diagram above is commutative, since for $m \in M_0$ we have
\begin{eqnarray}
\Hom_{A_0}(\eta_0^{M_0},W_0)(\widetilde{\Phi}(f \otimes g))(m) 
& = & \widetilde{\Phi}(f \otimes g)\left( \sum_{i=1}^u \alpha_i(m) \otimes \alpha_i^* \right) = \nonumber \\
& = & \sum_{i=1}^u \alpha_i^*(fg) \alpha_i(m) = \nonumber \\
& = & f(g(m))=\widetilde{\Psi}(f \otimes g)(m). \nonumber
\end{eqnarray}
Applying the functor $W_0 \otimes_k D(-)$ to the diagram above we obtain the following commutative diagram,
\[
 \xymatrix@!0@C=12pc@R=4pc{
W_0 \otimes_k D[\End_{A_0}(W_0) \otimes_k \Hom_{A_0}(M_0,W_0)] 
\ar@{<-}[dd]_(.5){Id_{W_0} \otimes D\widetilde{\Phi}} \ar@{<-}[rd]^(.6){Id_{W_0} \otimes D\widetilde{\Psi}} & \\
& W_0 \otimes_k D\Hom_{A_0}(M_0,W_0) \\
W_0 \otimes_k D\Hom_{A_0}(W_0\otimes_k D\Hom_{A_0}(M_0,W_0),W_0)
\ar@{<-}[ru]_(.6){G_0(\eta_0^{M_0})}
}
\]
We will show that when changing $G_0(\eta_0^{M_0})$ by $\eta_0^{G_0(M_0)}$ in the diagram above we still have 
commutativity. Take a dual basis $(f_i,f_i^*)_{i=1}^u$ of the $k$-vector space 
$\Hom_{A_0}(W_0 \otimes_k D\Hom_{A_0}(M_0,W_0),W_0)$. On the one hand, for an element $w \in W_0$ and a function  
$h \in D\Hom_{A_0}(M_0,W_0)$ we have that $[Id_{W_0} \otimes D\widetilde{\Psi}](w \otimes h)=w \otimes (h \circ \widetilde{\Psi})$. 
On the other hand,
\begin{eqnarray}
 [Id_{W_0} \otimes D\widetilde{\Phi}]\eta_0^{G_0(M_0)}(w \otimes h) 
& = & [Id_{W_0} \otimes D\widetilde{\Phi}] \left( \sum_{i=1}^u f_i(w\otimes h) \otimes f_i^* \right)= \nonumber \\
& = & \sum_{i=1}^u f_i(w\otimes h) \otimes (f_i^* \circ \widetilde{\Phi}).  \nonumber
\end{eqnarray}
We must show then that the expression $w \otimes (h \circ \widetilde{\Psi})$ equals the expression
$\sum_{i=1}^u f_i(w\otimes h) \otimes (f_i^* \circ \widetilde{\Phi})$. Consider the isomorphism
$K$ given by
\[
 \xymatrix{
W_0 \otimes_k D[\End_{A_0}(W_0) \otimes_k \Hom_{A_0}(M_0,W_0)], \ar[d]_-{K} & w\otimes H \ar@{|->}[d]^-{K} \\
\Hom_k([\End_{A_0}(W_0) \otimes_k \Hom_{A_0}(M_0,W_0)],W_0), & [f \otimes g \mapsto H(f \otimes g)w].
}
\]
For $f \in \End_{A_0}(W_0)$ and $g \in \Hom_{A_0}(M_0,W_0)$ we have that
\begin{eqnarray}
 K(\sum_{i=1}^uf_i(w \otimes h)\otimes (f_i^*\circ \widetilde{\Phi})) (f \otimes g) & = & 
\sum_{i=1}^u f_i^* (\widetilde{\Phi}(f\otimes g))  f_i(w \otimes h) = \nonumber \\
& = & \widetilde{\Phi}(f \otimes g)(w \otimes h)=h(fg)w, \nonumber
\end{eqnarray}
while
\[
 K(w \otimes (h \circ \widetilde{\Psi}))(f \otimes g)= h(\widetilde{\Psi}(f \otimes g))w=h(fg)w.
\]
Therefore the diagram above (substituting $G_0(\eta_0^{M_0})$ by $\eta_0^{G_0(M_0)}$) is commutative and 
we have $(b)$. Since $Id_{W_0}\otimes D\widetilde{\Psi}$ and $Id_{W_0}\otimes D\widetilde{\Phi}$ are both isomorphisms, 
we have that $G_0(\eta_0^{M_0})$ and $\eta_0^{G_0(M_0)}$ are also isomorphisms, which proves $(a)$.

We restrict now to the class of $A_0$-modules $\mathcal{X} \vee W_0$. In there any module $M_0$ is cogenerated by $W_0$
(lemma~\ref{(DE)L:alaDyn}) thus $\eta_0^{M_0}$ is a monomorphism (point $(c)$).
Finally we prove that $G_0$ preserves exact sequences. Assume that $\varepsilon$ is an exact sequence in $\mathcal{Z}$,
\[
\varepsilon: \xymatrix{
0 \ar[r] & M_0 \ar[r] & E_0 \ar[r] & N_0 \ar[r] & 0.
}
\]
We use the Auslander-Reiten formulas in the hereditary case \[\Ext_{A_0}^1(N_0,W_0)\cong D\Hom_{A_0}(\tau^{-1}W_0,M_0)\] 
and the fact that $\Hom_{A_0}(\tau^{-1}W_0,M_0)=0$ for $M_0 \in \mathcal{X}$ and $\tau^{-1} W_0 \in \mathcal{Y}$. 
Then $\Hom_{A_0}(\varepsilon,W_0)$ is an exact sequence
\[
\xymatrix@C=1pc{
0 \ar[r] & \Hom_{A_0}(N_0,W_0) \ar[r] & \Hom_{A_0}(E_0,W_0) \ar[r] & \Hom_{A_0}(M_0,W_0) \ar[r] & 0.
}
\]
Since $Id_{W_0}\otimes D(-)$ is an exact functor we conclude that $G_0(\varepsilon)=Id_{W_0} \otimes D\Hom_{A_0}(\varepsilon,W_0)$ 
is an exact sequence, and we have $(d)$.
\eproof

Let $(\overrightarrow{},\eta)$ be the liftings of $(G_0,\eta_0)$. By lemma~\ref{(DE)L:levFunDos} we have an additive faithful functor
which preserves exact sequences
\[
 G=\overline{W_0 \otimes D\Hom_{A_0}(-,W_0)}:\check{\mathcal{U}}(\mathcal{X} \vee W_0) \to \check{\mathcal{U}}(W_0),
\]
and whose restriction to $\check{\mathcal{U}}(W_0)$ is full and dense.

\begin{lema} \label{(DE)L:ARDos}
Assume that $\xymatrix{ 0 \ar[r] & P \ar[r]^-{g} & E \ar[r]^-{f} & \overline{\tau} P \ar[r] & 0}$
is an almost split sequence in $\check{\mathcal{U}}(\mathcal{X} \vee W_0)$. 
If $P \in \check{\mathcal{U}}(W_0)$ then $\overrightarrow{g}$ is a left almost split morphism 
in $\check{\mathcal{U}}(W_0)$. On the contrary $\overrightarrow{g}$ is a section.
\end{lema}

\bproof
Consider the diagram given by the universal mappings
\[
 \xymatrix{ 
0 \ar[r] & P \ar[d]_-{\eta^P} \ar[r]^-{g} & E \ar@{.>}[dl]_-{t} \ar[d]^-{\eta^E} \ar[r]^-{f} 
& \overline{\tau} P \ar[d]^-{\eta^{\overline{\tau}P}} \ar[r] & 0 \\
0 \ar[r] & \overrightarrow{P} \ar[r]_-{\overrightarrow{g}} & \overrightarrow{E} \ar[r]_-{\overrightarrow{f}} 
& \overrightarrow{\overline{\tau} P} \ar[r] & 0.
}
\]
If $P$ is not an object in $\check{\mathcal{U}}(W_0)$ then $\eta^P$ is not a section, 
hence there exists $t:E \to \overrightarrow{P}$ such that
$\eta^P=tg$. Then $\overrightarrow{\eta^P}=\overrightarrow{f}\overrightarrow{g}$, and since $\overrightarrow{\eta^P}
=\eta^{\overrightarrow{P}}$ is an isomorphism (points $(a)$ and $(b)$ in lemma~\ref{(DE)L:levFunDos}), 
$\overrightarrow{g}$ is a section.

On the other hand, assume that $P \in \check{\mathcal{U}}(W_0)$. Then $\eta^P$ is an isomorphism and 
$\overrightarrow{g} = \eta^E g(\eta^P)^{-1}$ is not a section, for $g(\eta^P)^{-1}$ is not a section.  
If $U \in \check{\mathcal{U}}(W_0)$ and $u:\overrightarrow{P} \to U$ is not a section,
then $u\eta^P$ is not a section, thus there exists $s:E \to U$ such that $u\eta^P=sg$. 
\[
 \xymatrix{
& P \ar[r]^-{g} \ar[d]_-{\eta^P} & E \ar[d]^-{\eta^E} \\
U  \ar[d]_-{\eta^U} & \overrightarrow{P} \ar[d]^-{\eta^{\overrightarrow{P}}} \ar[r]^-{\overrightarrow{g}} \ar[l]_-{u} & \overrightarrow{E} \\
\overrightarrow{U}  & \overrightarrow{\overrightarrow{P}} \ar[l]^-{\overrightarrow{u}}
}
\]
Then by lemma~\ref{(DE)L:levFunDos}$(b)$ we have that 
$\eta^U u=\overrightarrow{u}\eta^{\overrightarrow{P}}=\overrightarrow{u}\overrightarrow{\eta^P}=\overrightarrow{s}\overrightarrow{g}$.
Since $\eta^U$ is invertible, we conclude that $u=(\eta^U)^{-1}\overrightarrow{s}\overrightarrow{g}$,
that is, $\overrightarrow{g}$ is a left almost split morphism.
\eproof

For the following proposition denote by $P^i$ the indecomposable projective $A$-modules $Ae_i$ for $i=1,\ldots,n,\omega$. 
Observe that the simple projective module in the Kronecker subcategory $\check{\mathcal{U}}(W_0)$ is given by
$Q=W_0=(W_0,0,0)$. The other indecomposable projective module $Q^{\omega}=(W_0 \otimes D\Hom_{A_0}(R,W_0),k,\gamma_{Q^{\omega}})$ 
is given by
\[
 \gamma_{Q^{\omega}}: 1 \mapsto \eta^R.
\]
This can be justify noticing that the endomorphism algebra $\End_A(Q^{\omega})$ is isomorphic to the field $k$
and that $\vdim Q^{\omega}=(2,1)$ is the corresponding projective root.

\begin{lema} \label{(DE)P:Dos}
Let $w_0=\vdim W_0$ be the dimension vector of $W_0$. Then
\begin{itemize}
 \item[a)] $\overrightarrow{P^{\omega}} \cong Q^{\omega}$ and $\overrightarrow{P^i} \cong (w_0)_iQ$ for $i=1,\ldots,n$.
 \item[b)] For any $N \in \mathcal{P}$ we have $\overrightarrow{N} \in \mathcal{P}_0$.
 \item[c)] The function $\sharp(\overrightarrow{\cdot}):\mathcal{P} \to \mathbb{N}$ is additive.
\end{itemize} 
\end{lema}

\bproof
We show $(a)$. For $i=1,\ldots,n$ it is clear that $Ae_i=A_0e_i=P^i_0$ and since $\omega$ is a source vertex, $P^i_{\omega}=0$.
Since $\dimk_k \Hom_{A_0}(A_0e_i,W_0) = \dimk_k (e_iW_0)=(w_0)_i$ we have that
\[
\overrightarrow{P^i}=(W_0\otimes_k D\Hom_{A_0}(P^i_0,W_0),0,0) \cong  (w_0)_iQ.
\]
On the other hand, $P^{\omega}=Ae_{\omega}=\left( \begin{smallmatrix} 0 & R \\ 0 & k \end{smallmatrix} \right)$, hence
$P^{\omega}=(R,k,\gamma_{P^{\omega}})$ where $\gamma_{P^{\omega}}$ is the transformation $1 \mapsto Id_R$.
By the presentation given for the projective Kronecker module $Q^{\omega}$ it is clear that $\overrightarrow{P^{\omega}} = Q^{\omega}$.
The claims $(b)$ and $(c)$ can be shown in a similar way as in lemma~\ref{(DE)P:Uno}.
\eproof

\begin{proposicion} \label{(DE)C:invar}
Let $\mathcal{I}$ and $\mathcal{P}$ be the preinjective and posprojective components of $A$-mod respectively. 
Then the functions 
$\sharp (\overleftarrow{\cdot}):\mathcal{I} \to \mathbb{Z}$ and $\sharp (\overrightarrow{\cdot}):\mathcal{P} \to \mathbb{Z}$ 
are invariant under translation.
\end{proposicion}
\bproof
We will show the result in the preinjective component $\mathcal{I}$.
As we have shown in lemma~\ref{(DE)P:Uno}$(c)$ the function $f=\sharp(\overleftarrow{\cdot})$ is additive.
By lemma~\ref{(DE)L:componentes} the component $\mathcal{I}$ is a hereditary directed proper translation quiver
with no projective vertices and with a finite number of orbits.
Let $\mathcal{S}_0$ be the cosection in $\mathcal{I}$ which consists in the isomorphism classes of indecomposable 
injective $A$-modules. Then the full subquiver $Q_{\mathcal{S}_0}$ of $\mathcal{I}$
determined by $\mathcal{S}_0$ is isomorphic to $\widetilde{\Delta}^{op}$.
By the point $(b)$ in lemma~\ref{(DE)L:aditCox} we have that
\[
 f(\tau \mathcal{S}_0)=\Phi_{\widetilde{\Delta}^{op}} f(\mathcal{S}_0),
\]
where $\Phi_{\widetilde{\Delta}^{op}}$ is the Coxeter matrix associated to the quiver $\widetilde{\Delta}^{op}$.

Now, by lemma~\ref{(DE)P:Uno}$(a)$ we have the equality $f(\mathcal{S}_0)=w_0+\mathbf{e}_{\omega}$.
By construction of the algebra $A=k\widetilde{\Delta}$ and lemma~\ref{(DE)L:vctRad}, 
the vector $w_0+\mathbf{e}_{\omega}$ belongs to the radical of the quadratic form $q_{\widetilde{\Delta}^{op}}$.
From the equivalence of $(a)$ and $(c)$ in lemma~\ref{(DE)L:radCero} we have that 
$\Phi_{\widetilde{\Delta}}(w_0+\mathbf{e}_{\omega})=w_0+\mathbf{e}_{\omega}$, hence
\[
 f(\tau \mathcal{S}_0)=\Phi_{\widetilde{\Delta}^{op}} f(\mathcal{S}_0)=f(\mathcal{S}_0).
\]
The posprojective case can be shown in a similar way.
\eproof

Recall that $\mathcal{A}^W$ is the reduced Kronecker algebra given in lemma~\ref{(DE)L:subKro}.

\begin{corolario} \label{(DE)C:ContComp}
 For any preinjective (posprojective) $\mathcal{A}^W$-module $\widetilde{M}$, the $A$-module $F^W(\widetilde{M})$ is preinjective
(posprojective respectively).
\end{corolario}
\bproof
We prove the case of preinjective $\mathcal{A}^W$-modules, the posprojective case is similar. 
Let $\widetilde{M}_{i}$ and $M_j$ be admissible orderings ($i,j \in \mathbb{N}$) of representatives of
indecomposable preinjective $\mathcal{A}^W$ and $A$-modules respectively, as in lemma~\ref{(DE)L:orden}.
We will show that for any $i \in \mathbb{N}$ there exists $j_i \in \mathbb{N}$ such that 
$F^W(\widetilde{M}_i) \cong M_{j_i}$. Consider the set
\[
 \mathcal{N}=\{i\in \mathbb{N} \; | \; \text{there exists $j_i \in \mathbb{N}$ such that $F^W(\widetilde{M}_i)\cong M_{j_i}$}\}.
\]
Recall the notation of $A$-modules $I^{\omega}=(0,k,0)$,
 $I=\overline{W_0}$ which are injective in the Kronecker subcategory $\check{\mathcal{U}}(W_0)$.
Let $\widetilde{I}$ and $\widetilde{I^{\omega}}$ be indecomposable injective modules in $\mathcal{A}^W$-mod such that
$F^W(\widetilde{I})=I$ and $F^W(\widetilde{I^{\omega}})=I^{\omega}$.\\
\underline{Step 1.} \textit{For any $j \in \mathbb{N}$ there exist $g(j)$, $g^{\omega}(j)$, $d(j)$, $d^{\omega}(j)\geq 0$ such that
\[
 \overleftarrow{M_{j}}\cong d(j)F^W(\tau_{K}^{g(j)}\widetilde{I})\oplus d^{\omega}(j)F^W(\tau_K^{g^{\omega}(j)}\widetilde{I^{\omega}}),
\]
where $\tau_K$ is the Auslander-Reiten translation in the classical Kronecker category.} 
We will show the claim inductively. As base case the claim is valid on the indecomposable injective $A$-modules
because of point $(a)$ in lemma~\ref{(DE)P:Uno}. Assume it is valid for $M_1,\ldots,M_{\ell}$ and
consider the exact sequence in $A$-mod that starts in $M_{\ell+1}$,
\[
\varepsilon:\xymatrix{0 \ar[r] & M_{\ell+1} \ar[r] & E \ar[r] & \tau^{-1}(M_{\ell+1}) \ar[r] & 0}. 
\]
By lemma~\ref{(DE)L:ARUno}, if $\tau^{-1}(M_{\ell+1})$ is not an object in $\check{\mathcal{U}}(W_0)$ then 
$\overleftarrow{\varepsilon}$ is a trivial sequence,  and the result follows by induction hypothesis.
Assume then that $\tau^{-1}(M_{\ell+1}) \cong M_{\ell_0}$ is in $\check{\mathcal{U}}(W_0)$.
Again by lemma~\ref{(DE)L:ARUno} and proposition~\ref{(DE)C:invar}, the sequence $\overleftarrow{\varepsilon}$ is almost
split in $\check{\mathcal{U}}(W_0)$,
\[
\overleftarrow{\varepsilon}=\xymatrix{0 \ar[r] & \overleftarrow{M_{\ell+1}} \ar[r] & \overleftarrow{E} \ar[r] 
& \overleftarrow{M_{\ell_0}} \ar[r] & 0}. 
\]
Since $\ell_0\leq \ell$, by induction hypothesis we have that $\overleftarrow{M_{\ell_0}}\cong N$ with
\[
  N= d(\ell_0)\tau_{K}^{g(\ell_0)}\widetilde{I}\oplus 
d^{\omega}(\ell_0)\tau_K^{g^{\omega}(\ell_0)}\widetilde{I^{\omega}},
\]
Hence $\overleftarrow{M_{\ell+1}}\cong F^W(\tau_K(N))$, which completes the induction step.\\
\underline{Step 2.} \textit{The set $\mathcal{N}$ is infinite}. By proposition~\ref{(DE)C:invar}
the set $\{d(j)+d^{\omega}(j)\}_{j \in \mathbb{N}}$ is bounded. Since the dimension of $\overleftarrow{M_j}$ over the field $k$
is unbounded when $j$ tends to infinity, then (at least) one of the sets
$\{g(j)\}_{j\in \mathbb{N}}$ or $\{g^{\omega}(j)\}_{j\in \mathbb{N}}$ is unbounded.
By lemma~\ref{(DE)L:ARUno} this implies that $\mathcal{N}$ is an infinite set.\\
\underline{Step 3.} \textit{The set $\mathcal{N}$ is equal to $\mathbb{N}$.} 
Let $i \in \mathbb{N}$ be arbitrary. Then there exists $i_0$ such that for any integer $i'\geq i_0$
there is a nonzero morphism $f:\widetilde{M_{i'}} \to \widetilde{M_{i}}$ in the Kronecker category $\mathcal{A}^W$-mod. 
By step 2, there exists $i_1 \geq i_0$ with $i_1 \in \mathcal{N}$. Then there exists $j_{i_1}$ such that 
$F^W(\widetilde{M_{i_1}})\cong M_{j_{i_1}}$, and therefore
there is a nonzero morphism $F^W(f):M_{j_{i_1}} \to F^W(\widetilde{M_i})$. Since $M_{j_{i_1}}$ is preinjective,
we conclude that $F^W(\widetilde{M_i})$ is an indecomposable preinjective $A$-module. Then $i \in \mathcal{N}$, 
which completes the proof.
\eproof

The reductions in the next section will describe the subcategory 
$\check{\mathcal{U}}(\mathcal{X}_0 \vee W_0 \vee \mathcal{Y}_0)$ of $A$-mod, where
\[
 \mathcal{X}_0=\dConL([W_0])=\{ [N_0] \in \Ind A_0 \; | \; [N_0] \prec [W_0] \text{ but } \tau^{-1}_0[N_0] \npreceq [W_0]\}, \text{ and}
\]
\[
 \mathcal{Y}_0=\dConR([W_0])=\{ [M_0] \in \Ind A_0 \; | \; [W_0] \prec [M_0] \text{ but } [W_0] \npreceq \tau_0[M_0]\}.
\]
This subcategory containes almost complete the preinjective and posprojective components in $A$-mod. 

\begin{teorema} \label{(DE)T:finito}
Consider the subsets of $\Ind A_0$ given by
\[
\mathcal{X}_0=\dConL([W_0]) \quad \text{and} \quad \mathcal{Y}_0=\dConR([W_0]).
\]
\begin{itemize}
 \item[a)] The number of isomorphism classes of indecomposable posprojective $A$-modules which are not in the 
subcategory $\check{\mathcal{U}}(\mathcal{X}_0 \vee [W_0])$ is finite.
 \item[b)] The number of isomorphism classes of indecomposable preinjective $A$-modules which are not in the 
subcategory $\check{\mathcal{U}}([W_0] \vee \mathcal{Y}_0)$ is finite. 
\end{itemize}
\end{teorema}
\bproof
We show $(a)$ by steps. Fix an ordering contrary to the admissible $\{[N_{i}]\}_{i \in \mathbb{N}}$
in representatives of elements of the posprojective component $\mathcal{P}$, that is, if there is an irreducible
morphism $f:N_i \to N_j$ then $i < j$. Consider the subsets of $\Ind A_0$ given by
\[
 G^{\mathcal{P}}= \{ [Z_0] \in \Ind A_0 \; | \; \text{ there exists $h([Z_0])\in \mathbb{N}$ such that $Z_0 \cong N_{h([Z_0])}$} \},
\]
\[
 \overline{G^{\mathcal{P}}}= \{ [Z_0] \in \Ind A_0 \; | \; \text{ there exists $\overline{h}([Z_0])\in \mathbb{N}$ 
such that $\overline{Z_0} \cong N_{\overline{h}([Z_0])}$} \}.
\]
\underline{Step 1.} \textit{We have that $\overline{G^{\mathcal{P}}} \subset G^{\mathcal{P}}$ and 
$\overline{G^{\mathcal{P}}} = \tau_0 G^{\mathcal{P}}$}. For the first inclusion observe that
if $[Z_0] \in \overline{G^{\mathcal{P}}}$, then $Z_0$ is isomorphic to a direct summand of the restriction to
$A_0$-mod of $N_{\overline{h}([Z_0])}$. By the point $(i)$ in the proof of proposition~\ref{(DE)C:ciclo}, 
there exists $h([Z_0])\leq \overline{h}([Z_0])$ such that $N_{h([Z_0])} \cong Z_0$, that is, $[Z_0] \in G^{\mathcal{P}}$. 
To verify the equality 
$\overline{G^{\mathcal{P}}}=\tau_0 G^{\mathcal{P}}$ observe first that if $Z_0$ is an injective $A_0$-module,
then $\overline{Z_0}$ is an injective $A$-module. Indeed, if $u:L \to M$ is a monomorphism in $A$-mod and 
$g:L \to \overline{Z_0}$ is any morphism as in the following diagram,
\[
\xymatrix@R=1pc@C=1.5pc{
& L_{\omega} \ar[dl]_-{g_{\omega}} \ar[r]^-{u_{\omega}} \ar[dd]_-{\gamma_L} & M_{\omega} \ar[dd]^-{\gamma_M} \ar@{-->}[dll] \\
|Z_0| \ar@{=}[dd] \\
& |L_0| \ar[dl]_-{|g_0|} \ar[r]^-{|u_0|} & |M_0| \ar[dll]^-{|g'_0|} \\
|Z_0|
}
\]
then there exists $g'_0:M_0 \to Z_0$ such that the lower triangle commutes (for $u_0$ is a monomorphism and $Z_0$ 
is injective in $A_0$-mod), thus if we define $g'=(g'_0,|g'_0|\gamma_M):M \to \overline{Z_0}$ then $g=g'u$, and therefore
$\overline{Z_0}$ is injective.

In this way, if $Z_0 \in \overline{G^{\mathcal{P}}}$ then the $A_0$-module $Z_0$ is not injective,
for $\overline{Z_0}$ is posprojective and $A$-mod is not of finite representation type (lemma~\ref{(DE)L:componentes}).
Then there exists an almost split sequence in $A$-mod of the form
\[
\overline{\varepsilon_0}=\xymatrix@C=3pc{0 \ar[r] & \overline{Z_0} \ar[r]^-{(f,|I_{Z_0}|)} 
& (Y_0,|Z_0|,|f|) \ar[r]^-{(g,0)} & \tau^{-1}_0Z_0 \ar[r] & 0},
\]
and hence $\tau^{-1}_0Z_0$ is posprojective in $A$-mod, that is, 
$\overline{G^{\mathcal{P}}} \subseteq \tau_0G^{\mathcal{P}}$. If on the other hand $\tau^{-1}_0Z_0$ belongs
to $G^{\mathcal{P}}$, then by the sequence above $\overline{\varepsilon_0}$ we have that
$\overline{Z_0}$ is posprojective, thus $\overline{G^{\mathcal{P}}} \supseteq \tau_0G^{\mathcal{P}}$.

In particular, if $[Z_0] \in G^{\mathcal{P}}$, then $Z_0$ is not injective in $A_0$-mod.\\
\underline{Step 2.} \textit{The set $\mathcal{S}^{\mathcal{P}}=G^{\mathcal{P}}-\overline{G^{\mathcal{P}}}$
is a connected cosection in $\Gamma(A_0)$}. The first condition for a cosection is clearly satisfied, since 
$\overline{G^{\mathcal{P}}} = \tau_0 G^{\mathcal{P}}$. We verify the second condition, that is,
if $x \to s$ is an arrow in $\Gamma(A_0)$ with $s \in \mathcal{S}^{\mathcal{P}}$, then either 
$x \in \mathcal{S}^{\mathcal{P}}$ or $\tau^{-1}_0(x) \in \mathcal{S}^{\mathcal{P}}$.
Observe that $x \in G^{\mathcal{P}}$, for if $s=[Z_0]$ then $x$ corresponds to the isomorphism class of a direct summand
of the $A_0$-module $Y_0$ in the exact sequence above $\overline{\varepsilon_0}$. Hence, by the
point $(i)$ in the proof of proposition~\ref{(DE)C:ciclo}, $x \in G^{\mathcal{P}}$. Now, if 
$x \notin \overline{G^{\mathcal{P}}}$ then $x \in \mathcal{S}^{\mathcal{P}}$. 

It is left to show that if $x \in \overline{G^{\mathcal{P}}}$, then $\tau^1_0(x) \in \mathcal{S}^{\mathcal{P}}$.
By the equality $\overline{G^{\mathcal{P}}}=\tau_0G^{\mathcal{P}}$, since $x \in \overline{G^{\mathcal{P}}}$
then $\tau_0^{-1}(x) \in G^{\mathcal{P}}$. Thus it is enough to see that $\tau_0^{-1}(x) \notin G^{\mathcal{P}}$.
Assume that $\tau_0^{-1}(x) \in \overline{G^{\mathcal{P}}}$. Then $\tau_0^{-1}(x)$ is not injective
and $\tau_0^{-2}(x) \in G^{\mathcal{P}}$. If $x=[X_0]$, then there exists an Auslander-Reiten sequence in
the posprojective component of $A$-mod of the form
\[
\xymatrix@C=3pc{0 \ar[r] & \overline{\tau_0^{-1}X_0} \ar[r]& (W_0,|\tau_0^{-1}X_0|,|f|) \ar[r] & \tau_0^{-2}X_0 \ar[r] & 0},
\]
where one of the direct summands of $W_0$ belong to the class $\tau_0^{-1}(s)=[\tau_0^{-1}Z_0]$ again by the point 
$(i)$ in proposition~\ref{(DE)C:ciclo}, we have that $\tau_0^{-1}(s)\in G^{\mathcal{P}}$. 
This contradicts that $s \in G^{\mathcal{P}}-\overline{G^{\mathcal{P}}}$. 
Then $\tau_0^{-1}(x) \in \mathcal{S}^{\mathcal{P}}$ and hence $\mathcal{S}^{\mathcal{P}}$ is a cosection.
The connectedness is due to lemma~\ref{(DE)L:orbital}, for by the equality $\overline{G^{\mathcal{P}}} = \tau_0 G^{\mathcal{P}}$ the
set $\mathcal{S}^{\mathcal{P}}$ contains at most one element in each orbit. This completes the step 2.

For any $[Z_0] \in \Ind A_0$ define the set
\[
 \mathcal{N}([Z_0])=\{ i \in \mathbb{N} \; | \; \text{$Z_0$ is isomorphic to a direct summand of $N_{i}|_{A_0}$} \}.
\]
\underline{Step 3.} \textit{$\mathcal{S}^{\mathcal{P}}=\{[Z_0] \in \Ind A_0 \; | \; 
\text{$\mathcal{N}([Z_0])$ is infinite} \}$. } Observe first that 
$[Z_0]$ belongs to $G^{\mathcal{P}}$ if and only if $\mathcal{N}([Z_0]) \neq \emptyset$, and in that case
\[
 h([Z_0])=\textbf{min} \mathcal{N}([Z_0]).
\]
This is consequence of the point $(i)$ in the proof of proposition~\ref{(DE)C:ciclo}.
On the other hand, we show that $[Z_0] \in \overline{G^{\mathcal{P}}}$ if and only if $\mathcal{N}([Z_0])$ is
a finite set, and in that case
\[
 \overline{h}([Z_0])=\textbf{max}\mathcal{N}([Z_0]).
\]
Indeed, if $[Z_0] \in \overline{G^{\mathcal{P}}}$ then for any $N_{i}$ such that $Z_0$ is isomorphic to a direct
summand of $(N_{i})_0$, there exists a morphism of $A_0$-modules $f_0:(N_{i})_0 \to Z_0$ and therefore
a morphism of $A$-modules $f:N_{i} \to \overline{Z_0}$ given by
\[
 \xymatrix{
(N_{i})_{\omega} \ar[r]^-{|f_0|\gamma_{N_{i}}} \ar[d]_-{\gamma_{N_{i}}} & |Z_0| \ar@{=}[d] \\
|(N_{i})_0| \ar[r]_-{|f_0|} & |Z_0|.
}
\]
Hence $i \leq \overline{h}([Z_0])$ and $\mathcal{N}([Z_0])$ is a finite set. 
Assume now that $\mathcal{N}([Z_0])$ is finite and let
$j_0$ be its maximal element. Then the restriction to $A_0$ of the almost split sequence that
starts in $N_{j_0}$ is not trivial, and by lemma~\ref{(DE)L:levAR} we have that $N_{j_0} \cong \overline{Z_0}$.
In this way $[Z_0] \in \overline{G^{\mathcal{P}}}$ and $\overline{h}([Z_0])=\textbf{max} \mathcal{N}([Z_0])$.
As consequence $[Z_0] \in G^{\mathcal{P}}$ if and only if $\mathcal{N}([Z_0]) \neq \emptyset$, and in that case $[Z_0]$ 
is an element in $\overline{G^{\mathcal{P}}}$ if and only if $\mathcal{N}([Z_0])$ is a finite set. This completes
step~3.

We finally show the point $(a)$ in the theorem. By step 2 the set  
$\mathcal{S}^{\mathcal{P}}=G^{\mathcal{P}}-\overline{G^{\mathcal{P}}}$ is a connected cosection, which is contained
in $\{[W_0]\} \cup \ConL([W_0])$ by the proposition~\ref{(DE)C:ciclo}. Moreover, $\mathcal{S}^{\mathcal{P}}$ contains 
$[W_0]$ by corollary~\ref{(DE)C:ContComp}. From uniqueness in lemma~\ref{(DE)L:unicaSec} it follows that the cosection
$\mathcal{S}^{\mathcal{P}}$ is equal to the cosection $\{[W_0]\} \cup \dConL([W_0])$, that is,
\[
 \mathcal{S}^{\mathcal{P}}=\{[W_0]\} \cup \mathcal{X}_0.
\]
By step 3 the following subset of $\mathbb{N}$ is finite,
\[
 \mathcal{N}=\bigcup_{[Z_0] \in \overline{G^{\mathcal{P}}}}\mathcal{N}([Z_0]).
\]
Assume that $N_j$ is a posprojective $A$-module $(j \in \mathbb{N})$ which is not in 
$\check{\mathcal{U}}(\mathcal{X}_0 \vee W_0)$. Then the restriction $N_j|_{A_0}$ has as direct summand an
$A_0$-module $Z_0$ which is in $\overline{G^{\mathcal{P}}}$, and hence $j \in \mathcal{N}$ which is a finite set.
This proves $(a)$, the proof of $(b)$ is similar.
\eproof

The functions $\sharp (\overleftarrow{\cdot})$ and 
$\sharp (\overrightarrow{\cdot})$ described above are related to the rank inside the classes of 
$A$-modules in the theorem. Recall that, as defined at the beginning of this section, 
the rank of an indecomposable $A$-module $M=(M_0,M_{\omega},\gamma_{M_0})$ is given by
\[
 \Rnk(M)=\dimk_k M_{\omega} - m_{M_0},
\]
where $m_{M_0}$ is the number of indecomposable summands of $M_0$ isomorphic to $W_0$.

\begin{proposicion} \label{(DE)P:rank}
Let $\mathcal{I}'$ and $\mathcal{P}'$ be the subsets of the preinjective component $\mathcal{I}$ and the posprojective
component $\mathcal{P}$ of $A$-mod determined by the indecomposable modules that lie in $\check{\mathcal{U}}(W_0 \vee \mathcal{Y}_0)$
and $\check{\mathcal{U}}(\mathcal{X}_0 \vee W_0)$ respectively. 
\begin{itemize}
 \item[a)] If $[M] \in \mathcal{I}'$ then $\Rnk(M)= \sharp \overleftarrow{M}$.
 \item[b)] If $[N] \in \mathcal{P}'$ then $\Rnk(N)+ \sharp \overrightarrow{N}=e_{N_0}$.
\end{itemize}
\end{proposicion}

\bproof
First we show $(a)$. Notices that in a pull-back diagram with surjective horizontal arrows
\[
\xymatrix{
 \overleftarrow{M}_{\omega} \ar[r] \ar[d] & M_{\omega} \ar[d] \\
 |\overleftarrow{M}_0| \ar[r] & |M_0|,
}
\]
the equality $\dimk_k \overleftarrow{M}_{\omega}-\dimk_k |\overleftarrow{M}_0|=\dimk_k M_{\omega}-\dimk_k |M_0|$ holds.
Since $[M_0] \in W_0 \vee \mathcal{Y}_0$ and $\dimk_k\Hom_{A_0}(W_0,Y)=1$ for any $[Y]$ in $\mathcal{Y}_0$, we have
\begin{eqnarray}
m_{\overleftarrow{M}_0} & = & m_{M_0}+e_{M_0}. \nonumber
\end{eqnarray}
Since $\dimk_k|W_0|=2$ and $\dimk_k |Y|=1$ (lemma~\ref{(DE)L:valores}) we have that 
$\dimk_k |\overleftarrow{M}_0|=2m_{\overleftarrow{M}_0}$ and $\dimk_k |M_0|=2m_{M_0}+e_{M_0}$. Hence
 \begin{eqnarray}
\dimk_k \overleftarrow{M}_{\omega}-\dimk_k |\overleftarrow{M}_0| & = & \dimk_k M_{\omega}-\dimk_k |M_0|, \nonumber \\
\dimk_k \overleftarrow{M}_{\omega}-2m_{\overleftarrow{M}_0} & = & \dimk_k M_{\omega}-(2m_{M_0}+e_{M_0}), \nonumber \\
\dimk_k \overleftarrow{M}_{\omega}-m_{\overleftarrow{M}_0} -(m_{\overleftarrow{M}_0}) 
& = & \dimk_k M_{\omega}-m_{M_0}-(m_{M_0}+e_{M_0}), \nonumber \\
\Rnk(\overleftarrow{M})-(m_{\overleftarrow{M}_0}) & = & \Rnk(M)-(m_{M_0}+e_{M_0}). \nonumber
 \end{eqnarray}
Then $\Rnk(\overleftarrow{M})=\Rnk(M)$ for $m_{\overleftarrow{M}_0}  =  m_{M_0}+e_{M_0}$.
To complete the proof of $(a)$ notice that $\sharp \overleftarrow{M}=\Rnk(\overleftarrow{M})$, 
for $\overleftarrow{M}$ is preinjective in the Kronecker subcategory $\check{\mathcal{U}}(W_0)$
(lemma~\ref{(DE)P:Uno}$(b)$). Indeed, if $N=(N_1,N_2,N_{\alpha_1},N_{\alpha_2})$ is a preinjective module in the
Kronecker category $\mathcal{A}^W$-mod then $\sharp N=\dimk_k N_2-\dimk_k N_1$. Under the functor
$F^W$ given in lemma~\ref{(DE)L:subKro} we have
\begin{eqnarray}
 \sharp F^W(N)=\sharp (W \otimes N) & = & \dimk_k (S(\omega)\otimes N_2)-\sharp (W_0 \otimes N_1)= \nonumber \\
& = & \dimk_k N_2-\sharp (W_0 \otimes N_1)=\Rnk(F^W(N)). \nonumber
\end{eqnarray}
We turn now to the proof of $(b)$. Observe that if $[N_0] \in \mathcal{X}_0 \vee W_0$ then 
\[
 m_{\overrightarrow{N_0}}=m_{N_0}+e_{N_0},
\]
for if $[X]$ belongs to $\mathcal{X}_0$ then $\dimk_k \Hom_{A_0}(X,W_0)=1$. In this way
\begin{eqnarray}
\Rnk(N)-\Rnk(\overrightarrow{N}) & = & (\dimk_k N_{\omega}-m_{N_0})-(\dimk_k N_{\omega}-m_{\overrightarrow{N_0}}) \nonumber \\
& = &  e_{N_0}. \nonumber
\end{eqnarray}
We show that $\sharp \overrightarrow{N}=-\Rnk(\overrightarrow{N})$ since $\overrightarrow{N}$ is posprojective in the
Kronecker sub\-ca\-te\-gory $\check{\mathcal{U}}(W_0)$. Indeed, if $M=(M_1,M_2,M_{\alpha_1},M_{\alpha_2})$ is a posprojective module
in the Kronecker category $\mathcal{A}^W$-mod then $\sharp M=\dimk_k M_1-\dimk_k M_2$. Again under the functor
$F^W$ we have
\begin{eqnarray}
 \sharp F^W(M)=\sharp (W \otimes M) & = &\sharp (W_0 \otimes M_1)- \dimk_k (S(\omega)\otimes M_2)= \nonumber \\
& = & \sharp (W_0 \otimes M_1) - \dimk_k M_2=-\Rnk(F^W(M)). \nonumber
\end{eqnarray}
We conclude that
\[
 \Rnk(N)+ \sharp \overrightarrow{N}=e_{N_0}.
\]
\eproof

Let $\mathcal{S}$ be a connected section in the preinjective component of $A$-mod contained in $\mathcal{I}'$. 
By lemma~\ref{(DE)L:orbital}, $\mathcal{S}$ intersects each orbit of $\mathcal{I}'$ in exactly one vertex,
hence by proposition~\ref{(DE)C:invar}, $\Rnk(\mathcal{S})=w_0 + \mathbf{e}_{\omega}$ 
is the radical vector of the quadratic form $q$ associated to $A$. On the other hand, there exists a minimal integer
$p_A\geq 1$ such that for any connected section $\mathcal{S}$ in $\mathcal{P}'$ we have 
$\Rnk(\mathcal{S})=\Rnk(\tau^{-p_A}\mathcal{S})$. In other words, there exists an integer $m_A$ such that
\[
 \sum_{i=0}^{p_A-1}\Rnk(\tau^{-i} \mathcal{S})=m_A(w_0 + \mathbf{e}_{\omega}).
\]
The integers $p_A$ and $m_A$, which will be called \textbf{period} and \textbf{multiplicity} of the posprojective component, 
depend only on the type of the Dynkin diagram $\Delta$. Their values are shown in the following table.

\begin{table} [!hbt] 
\begin{center}
 \begin{tabular}{c c c c c} 
\hline 
Diagram & & Tubular type & Period & Multiplicity \\
\hline \\
$\widetilde{\mathbf{A}}_{p,q}$ &  & $(p,q)$ & $\mathbf{m.c.m.}(p,q)$ & $\frac{pq-(p+q)}{\mathbf{m.c.d.}(p,q)}$ \\
$\widetilde{\mathbf{D}}_{n}$ & ${}_{n\text{ odd}}$ & $(n-2,2,2)$ & $2(n-2)$ & $2(n-3)$  \\
$\widetilde{\mathbf{D}}_{n}$ & ${}_{n\text{ even}}$ & $(n-2,2,2)$ &$n-2$ & $n-3$ \\
$\widetilde{\mathbf{E}}_{6}$ & & $(3,3,2)$ & 6 & 5   \\
$\widetilde{\mathbf{E}}_{7}$ & & $(4,3,2)$ & 12 & 11  \\
$\widetilde{\mathbf{E}}_{8}$ & & $(5,3,2)$ & 30 & 29  \\
\hline
\end{tabular}
\end{center}
\end{table}

\section{Reduction of extended Dynkin algebras.} \label{(DE)S:algRed}
%----------------------------------------------------------------------

Let $\Delta$ be a Dynkin quiver and $\widetilde{\Delta}$ the corresponding extended Dynkin quiver such that 
the extended vertex is a source. Let $A_0$ and $A$ be the path algebras of $\Delta$ y 
$\widetilde{\Delta}$ respectively and consider the ditalgebra $\mathcal{A}=(A,0)$ with zero differential.
In this section we study reduced ditalgebras $\mathcal{A}^X$ and $\mathcal{A}^Y$ by complete rigid admissible modules 
$X$, $Y$ that are defined from the $A_0$-module $W_0$ and representatives of the vertices in $\mathcal{X}_0=\dConL([W_0])$ and
$\mathcal{Y}_0=\dConR([W_0])$ respectively. 

%---------------------------------------------------------------------- 
\begin{displaymath}
\xy 0;/r.20pc/:
( 0,0)="A0" *{};
(20,0)="A2" *{{}_{\bullet}};
(30,0)="A3" *{};
(40,0)="A4" *{};
(50,0)="A5" *{};
(60,0)="A6" *{{}_{\bullet}};
(90,0)="A9" *{};
( 0,30)="B0" *{};
(30,30)="B3" *{{}_{\bullet}};
(40,30)="B4" *{};
(50,30)="B5" *{{}_{\bullet}};
(90,30)="B9" *{};
( 0,10)="C1" *{};
( 0,20)="C2" *{};
(90,10)="D1" *{};
(90,20)="D2" *{};
(25,5)="E1" *{{}_{\bullet}};
(30,10)="E2" *{{}_{\bullet}};
(35,15)="E3" *{{}_{\bullet}};
(45,25)="E4" *{{}_{\bullet}};
(35,25)="F1" *{{}_{\bullet}};
(45,15)="F2" *{{}_{\bullet}};
(50,10)="F3" *{{}_{\bullet}};
(55,5)="F4" *{{}_{\bullet}};
(40,20)="WB" *{{}_{\bullet}};
(46,20)="W" *{\scriptstyle [W_0]};
(52,33)="CR0" *{\scriptstyle \dConR([W_0])};
(25,33)="CL0" *{\scriptstyle \dConL([W_0])};
(70,15)="CR" *{\scriptstyle \ConR([W_0])};
(15,15)="CL" *{\scriptstyle \ConL([W_0])};
%-----------------------
"A0";"A2" **@{-};
"A6";"A9" **@{-};
"A9";"B9" **@{-};
"A0";"B0" **@{-};
"B0";"B3" **@{-};
"B5";"B9" **@{-};
"B3";"A6" **@{-};
"B5";"A2" **@{-};
"B5";"A2" **@{-};
"E1";"A3" **@{-};
"E2";"A4" **@{-};
"E3";"A5" **@{-};
"E4";"B4" **@{-};
"F1";"B4" **@{-};
"F2";"A3" **@{-};
"F3";"A4" **@{-};
"F4";"A5" **@{-};
\endxy
\end{displaymath}
One of our objectives is to describe exceptional representations for the reductions $\mathcal{A}^X$ and $\mathcal{A}^Y$, 
in such a way that through the use of the functors $F^X$ and $F^Y$ we obtain explicit exceptional representations 
of the extended Dynkin quiver $\widetilde{\Delta}$.

For a set of vertices $\mathcal{S}$ in $\Gamma(A_0)$ denote by $\Delta_{\mathcal{S}}$ the full subquiver 
of $\Gamma(A_0)$ determined by $\mathcal{S}$. Recall that $\mathcal{X}_0=\dConL([W_0])=\{ W^s_{1,i} \}$ with $i=1,\ldots, n_s$ and 
$s=1,\ldots, t$ (the integer $t$ takes the values
$1$, $2$ or $3$ depending on the tubular type of $\Delta$, cf. section~\ref{(DE)S:ARDynkin}). 
By simplicity we switch to the notation $W^s_{1,i}=X^s_i$. With this notation the modules $W_0$ and $X^s_{n_s}$ ($s=1,\ldots,t$) 
coincide. See the following figure for an example with two wings of order 3 and 5.

%---------------------------------------------------------------------- 
\begin{displaymath}
\xy 0;/r.20pc/:
( 0,0)="A0" *{};
(40,20)="A1" *{{}_{\bullet}};
(53,20)="A1e" *{\scriptstyle [W_0]=[X^1_5]=[X^2_3]};
(20,0)="A2" *{{}_{\bullet}};
(25,0)="A2e" *{\scriptstyle [X^1_1]};
( 0,30)="B0" *{};
(30,30)="B3" *{{}_{\bullet}};
(35,30)="B3e" *{\scriptstyle [X^2_1]};
(25,5)="E1" *{{}_{\bullet}};
(30,5)="E1e" *{\scriptstyle [X^1_2]};
(30,10)="E2" *{{}_{\bullet}};
(35,10)="E2e" *{\scriptstyle [X^1_3]};
(35,15)="E3" *{{}_{\bullet}};
(40,15)="E3e" *{\scriptstyle [X^1_4]};
(35,25)="F1" *{{}_{\bullet}};
(40,25)="F12" *{\scriptstyle [X^2_2]};
%-----------------------
"A0";"A2" **@{-};
"A0";"B0" **@{-};
"B0";"B3" **@{-};
"B3";"A1" **@{-};
"A1";"A2" **@{-};
%"B5";"A2" **@{-};
\endxy
\end{displaymath}
Consider the $A_0$-module
\[
X_0=W_0 \oplus \bigoplus_{\substack{1\leq i < n_s \\ s=1,\ldots,t}}X^s_i.
\]
Denote by $\Delta_{X_0}$ the quiver $\Delta_{\{[W_0]\}\cup \mathcal{X}_0}$.
Let $M_{\Delta_{X_0}}$ be the matrix associated to $\Delta_{X_0}$ and let $\Delta_{X_0}^{-1}$ be the quiver whose
associated matrix is the inverse $M_{\Delta_{X_0}}^{-1}$.

\begin{lema} \label{(DE)L:AXcero}
The $A_0$-module $X_0$ given above is admissible and complete. Denote by $(Z_0,P_0)$ a
splitting of the opposed endomorphism algebra of $X_0$. The $k$-algebra $Z_0$ is trivial.
The reduced tensor algebra $A_0^{X_0}$, given by $T_{Z_0}(P_0^*)$, is isomorphic to the path algebra
of the quiver $\Delta_{X_0}^{-1}$. Moreover, there exist finite dual bases of legible elements 
\[
 \{(p^s_{i,j},\gamma^s_{i,j})\}_{\substack{1 \leq i < j \leq n_s \\ s=1,\ldots,t}} \qquad \text{and} 
\qquad \{(a^s_i,\xi^s_i),(a,\xi),(a',\xi')\}_{\substack{1 \leq i < n_s \\ s=1,\ldots,t}},
\]
of $P_0$ and $\Hom_{A_0}(R,X_0)$ respectively, where the restrictions of $(a,\xi)$ and $(a',\xi')$ 
to $\Hom_{A_0}(R,W_0)$ are a dual basis. We can choose the basis of $P_0$ in such a way that for $1 \leq i < j \leq n_s$, 
$1 \leq k < \ell \leq n_r$ and $1 \leq r,s \leq t$ the following equations hold,
\begin{equation*}
p^s_{i,j}p^r_{k,\ell} = \left\{
\begin{array}{ll}
p^s_{i,\ell}, & \text{ if } j=k, s=r,\\
0, & \text{ otherwise.}
\end{array} \right.
\end{equation*} 
Defining $a^s_{n_s}=a^s_{n_s-1}p^s_{n_s-1,n_s}$, it is possible to choose a basis of $\Hom_{A_0}(R,X_0)$ such that for  
$1 \leq i < n_s$, $1 \leq k < \ell \leq n_r$ and $1 \leq r,s \leq t$ the following equations are satisfied,
\begin{equation*} 
a^s_jp^r_{k,\ell} = \left\{
\begin{array}{ll}
a^s_{\ell}, & \text{ if } j=k, s=r,\\
0, & \text{ otherwise.}
\end{array} \right.
\end{equation*} 
Denote by $\rho_s$ the restriction to $\Hom_{A_0}(R,W_0)$ of the morphisms $a^s_{n_s}$ (that
belong to $\Hom_{A_0}(R,X_0)$). These morphisms determine nonisomorphic indecomposable representations
 $W_0(\rho_s)$, for $s=1,\ldots,t$.
\end{lema}
\bproof
Observe that there exists a nonzero morphism in $\Hom_{A_0}(X^r_i,X^s_j)$ if and only if $r=s$, $i \leq j$.
Consider 
\[
 Z_0 = \End_{A_0}(W_0)^{op} \oplus  \bigoplus_{\substack{1\leq i < n_s\\ s=1,\ldots,t}}\End_{A_0}(X^s_i)^{op},
\]
\[
 P_0 =\bigoplus_{\substack{1\leq i < j \leq  n_s\\ s=1,\ldots,t}}\Hom_{A_0}(X^s_iX^s_j)^{op}.
\]
Then the opposed endomorphism algebra of $X_0$, considered as matrix algebra, is triangular
with diagonal isomorphic to $Z_0$. By lemma~\ref{(P)L:endExc} the $k$-algebra $Z_0$ is trivial, for all indecomposable
direct summands of $X_0$ are exceptional. It is clear that there is a decomposition of $Z_0$-$Z_0$-bimodules
\[
 \Gamma_0= \End_{A_0}(X_0)^{op} \cong Z_0 \oplus P_0,
\]
thus $X_0$ is admissible and complete. Now, by definition
\[
 A_0^{X_0}=(T_{A_0}(0))^{X_0}= T_{Z_0}(P_0^*).
\]
Fix a dual basis of the right $Z_0$-module $P_0$. For any $i=1,\ldots,n_s-1$ 
denote by $p_{i,i+1}^s$ a nonzero morphism between $X^s_i$ and $X^s_{i+1}$.
All morphisms $p^s_{i,i+1}$ are irreducible and monomorphisms, for each $X^s_i$ is cogenerated 
by $W_0$ (lemma~\ref{(DE)L:alaDyn}),
\[
 \xymatrix@C=3pc{
X^s_1 \ar[r]^-{p^s_{1,2}} & X^s_2 \ar[r]^-{p^s_{2,3}} & X^s_3  \cdots X^s_{n_s-2} \ar[r]^-{p^s_{n_s-2,n_s-1}} & X^s_{n_s-1} 
\ar[r]^-{p^s_{n_s-1,n_s}} & W_0.
}
\]
Define for any $1 \leq i < j \leq n_s$ the composition
\[
 p^s_{i,j}=p^s_{i,i+1}p^s_{i+1,i+2}\ldots p^s_{j-1,j}.
\]
Then $\{(p^s_{i,j},\gamma^s_{i,j})\}$ is a dual basis of legible elements of $P_0$,
which by construction satisfy the equation in the statement.

Recall that the path algebra $k\Delta_{X_0}^{-1}$ is isomorphic to the tensor algebra
$T_{Z_1}(W_{X_0})$ where $Z_1$ is the trivial $k$-algebra corresponding to the idempotent paths $e_i^s$ of 
$\Delta_{X_0}^{-1}$ and $W_{X_0}$ is the arrow bimodule of the quiver $\Delta_{X_0}^{-1}$.
Clearly we have an isomorphism of algebras $\Phi_0:Z_1 \to Z_0$ which assigns each idempotent 
$e_i^s$ to the identity morphism $f^s_i$ of the module $X_i^s$.

By lemma~\ref{(P)L:invers} to each arrow $\alpha:i \to j$ of $\Delta_{X_0}^{-1}$ corresponds a path
$p^s_{i,j}$ from $i$ to $j$ in $\Delta_{X_0}$. This path is contained in a wing $\theta(n_s)$ of $\Gamma(k\Delta)$
(for some $s \in \{1,\ldots,t \}$) for $[W_0]$ is a sink in $\Delta_{X_0}$. 
Define $\Phi_1(\alpha)=\gamma^s_{i,j}$. Since $\alpha=\alpha e_i$ and
$\alpha=e_j \alpha$, observe that
\[
 \Phi_1(\alpha)\Phi_0(e_i)=\gamma^s_{i,j}f^s_i=\gamma^s_{i,j}=\Phi_1(\alpha e_i),
\]
\[
 \Phi_0(e_j)\Phi_1(\alpha)=f^s_j\gamma^s_{i,j}=\gamma^s_{i,j}=\Phi_1(e_j\alpha).
\]
In this way the mapping $\Phi_1$ extends to an isomorphism of $Z_1$-$Z_1$-bimo\-du\-les through $\Phi_0$,
\[
 \Phi_1:W_{X_0} \to P_0^*.
\]
Hence the isomorphisms $\Phi_0$ and $\Phi_1$ extend in turn to an isomorphism of tensor algebras
\[
 \Phi:T_{Z_1}(W_{X_0}) \to T_{Z_0}(P_0^*),
\]
whose restrictions to $T$ and $W_{X_0}$ are $\Phi_0$ and $\Phi_1$ respectively.

Fix now a basis of $\Hom_{A_0}(R,X_0)$ starting from the chosen basis for $P_0$. Take first a nonzero element 
$a^s_1 \in \Hom_{A_0}(R,X^s_1)$ for each $s=1,\ldots, t$. Recursively define $a^s_{i}=a^s_{i-1}p^s_{i-1,i}$ for
$1<i<n_s$ and take a basis $\{a,a' \}$ of $\Hom_{A_0}(R,W_0)$. Consider all these morphisms as elements in
\[
 \Hom_{A_0}(R,X_0) \cong \Hom_{A_0}(R,W_0) \oplus \bigoplus_{\substack{1 \leq i < n_s \\ s=1,\ldots,t}} \Hom_{A_0}(R,X^s_{i}).
\]
Then $\{(a^s_i,\xi^s_i),(a,\xi),(a',\xi')\}$ (with $i=1,\ldots,n_s-1$, $s=1,\ldots,t$) is a dual basis 
of legible elements in $\Hom_{A_0}(R,X_0)$ which satisfies, by construction, the formula of the statement.
\eproof

\begin{lema} \label{(DE)L:AX}
Take $A=A_0[R]$ and let $B$ be the subalgebra of $A$ given by $B=A_0[0]$.
Denote by $\mathcal{A}$ the ditalgebra $\mathcal{A}=(A,0)$ with zero differential.
Then the $B$-module $X=X_0 \oplus S(\omega)$ is complete admissible and the reduced ditalgebra
$\mathcal{A}^X$ is isomorphic to $(k\Delta_{X_0}^{-1}[R^{X_0}],\delta^x)$, where
$R^{X_0}=\Hom_{A_0}(R,X_0)$. The differential $\delta^x$ has the following explicit form
in solid arrows, in terms of the bases in lemma~\ref{(DE)L:AXcero},
\begin{center}
 \begin{tabular}{c | l}
Vector in $\Hom_{A_0}(R,X_0)^*$ & Its differential $\delta^x$ \\
\hline
$\xi_j^s$ & $\sum_{1 \leq i < j}\gamma^s_{i,j} \otimes \xi^s_i$,  \\
$\xi$       & $\sum_{s=1,\ldots,t} \beta_s  \sum_{1 \leq i < n_s}\gamma^s_{i,n_s} \otimes \xi^s_i$, \\
$\xi'$   & $\sum_{s=1,\ldots,t} \beta_s' \sum_{1 \leq i < n_s}\gamma^s_{i,n_s} \otimes \xi^s_i$, \\
 \end{tabular}
\end{center}
where $\rho_s=\beta_s a+\beta'_sa'$ for each $s=1,\ldots,t$.
\end{lema}
\bproof
The proof is a direct application of proposition~\ref{(DE)P:extPuntRed} and the results of lemma~\ref{(DE)L:AXcero}.
For the computation of the reduced differential use the properties of the bases given in lemma~\ref{(DE)L:AXcero}.
\eproof

Recall that $\dimk_k\Hom_{A_0}(R,N_0)=1$ for any $N_0 \in \mathcal{X}_0$ and that the space
$\Hom_{A_0}(R,W_0)$ has dimension two (lemma~\ref{(DE)L:valores}). Then the one-point extension
$A_0^{X_0}[R^{X}]$ can be formed by adding two solid arrows from the extended vertex
$\omega$ to the vertex $m$ co\-rresponding to the maximal module in $\Delta_{X_0}^{-1}$, and a solid arrow from $\omega$
to the rest of the vertices.  The following figure shows the quiver of the reduced ditalgebra $\mathcal{A}^X$ 
when $A$ is the path algebra of the extension of a linearly ordered Dynkin diagram $\mathbf{A}_{\ell}$
(we omit the superscript $s=1$, for there is only one wing).
%------------------------------------------------ 
\begin{equation} \label{(DE)F:brazoX}
 \xymatrix@C=3pc{
& & & & *+++[]{\bullet_{\omega}} \\
\\
\\
\bullet_1 \ar@{.>}[r]^(.7){\gamma_{1,2}} \ar@{.>}@/^1.5pc/[rrr]^(.35){\gamma_{1,\ell-1}} \ar@{.>}@/_2pc/[rrrr]_(.5){\gamma_{1,\ell}} 
\ar@{<-}@/^2pc/[rrrruuu]^(.27){\xi_1} 
 & \bullet_2 \ar@{.>}@/_.5pc/[rr]_(.5){\gamma_{2,\ell-1}} \ar@{.>}@/^1.5pc/[rrr]^(.5){\gamma_{2,\ell}} 
\ar@{<-}@/^1pc/[rrruuu]^(.35){\xi_2} 
& \ldots & \bullet_{\ell-1} \ar@{.>}[r]^(.4){\gamma_{\ell-1,\ell}} \ar@{<-}@/^.5pc/[ruuu]^(.42){\xi_{\ell-1}}
& *++[]{\bullet_{\ell}} \ar@{<-}[uuu]^-{\xi} \ar@{<-}@<-1ex>[uuu]_-{\xi'} \\
}
\end{equation}

%-------------------------------------
The case of the reduction $\mathcal{A}^Y$ is similar to the case above. We give the corresponding details.
Write $Y^s_i$ instead of the notation $W^s_{i,n_s}$ ($1<i\leq n_s$ and $s=1,\ldots,t$) for the representatives of 
elements in $\mathcal{Y}_0$. In this case all $A_0$-modules $W_0$ and $Y^s_1$ coincide ($s=1,\ldots,t$), 
see for instance the following figure,
%---------------------------------------------------------------------- 
\begin{displaymath}
\xy 0;/r.20pc/:
(40,20)="A1" *{{}_{\bullet}};
(26,20)="A1e" *{\scriptstyle [Y^1_1]=[Y^2_1]=[W_0]};
(60,0)="A6" *{{}_{\bullet}};
(55,0)="A6e" *{\scriptstyle [Y^1_5]};
(90,0)="A9" *{};
(50,30)="B5" *{{}_{\bullet}};
(45,30)="B5e" *{\scriptstyle [Y^2_3]};
(90,30)="B9" *{};
(45,25)="E4" *{{}_{\bullet}};
(40,25)="E4e" *{\scriptstyle [Y^2_2]};
(45,15)="F2" *{{}_{\bullet}};
(40,15)="F2e" *{\scriptstyle [Y^1_2]};
(50,10)="F3" *{{}_{\bullet}};
(45,10)="F3e" *{\scriptstyle [Y^1_3]};
(55,5)="F4" *{{}_{\bullet}};
(50,5)="F4e" *{\scriptstyle [Y^1_4]};
%-----------------------
"A6";"A9" **@{-};
"A9";"B9" **@{-};
"B5";"B9" **@{-};
"A1";"A6" **@{-};
"B5";"A1" **@{-};
\endxy
\end{displaymath}
Define the $A_0$-module
\[
Y_0=W_0 \oplus \bigoplus_{\substack{1< i \leq n_s \\ s=1,\ldots,t}}Y^s_i,
\]
and denote by $\Delta_{Y_0}$ the quiver $\Delta_{\{[W_0]\}\cup \mathcal{Y}_0}$.
Let $M_{\Delta_{Y_0}}$ be the matrix associated to $\Delta_{Y_0}$ and let $\Delta_{Y_0}^{-1}$ be the quiver whose 
associated matrix is the inverse $M_{\Delta_{Y_0}}^{-1}$.

\begin{lema} \label{(DE)L:AYcero}
The $A_0$-module $Y_0$ given above is admissible and complete. Denote by $(Z'_0,Q_0)$ the splitting
of the opposed endomorphism algebra $Y_0$. The $k$-algebra $Z'_0$ is trivial.
The reduced tensor algebra $A_0^{Y_0}$, given by $T_{Z'_0}(Q_0^*)$, is isomorphic to the path algebra of
the quiver $\Delta_{Y_0}^{-1}$. Moreover, there exist finite dual bases of legible elements 
\[
 \{(q^s_{i,j},\varepsilon^s_{i,j})\}_{\substack{1 \leq i < j \leq n_s \\ s=1,\ldots,t}} \qquad \text{and} 
\qquad \{(b^s_i,\theta^s_i),(b,\theta),(b',\theta')\}_{\substack{1 < i \leq n_s \\ s=1,\ldots,t}},
\]
of $Q_0$ and $\Hom_{A_0}(R,Y_0)$ respectively, where the restrictions of $(b,\theta)$ and $(b',\theta')$ 
to $\Hom_{A_0}(R,W_0)$ are a dual basis. We can choose the basis of $Q_0$ in such a way that for $1 \leq i < j \leq n_s$, 
$1 \leq k < \ell \leq n_r$ and $1 \leq r,s \leq t$ the following equations are satisfied,
\begin{equation*}
q^s_{i,j}q^r_{k,\ell} = \left\{
\begin{array}{ll}
q^s_{i,\ell}, & \text{ if } j=k, s=r,\\
0, & \text{ otherwise.}
\end{array} \right.
\end{equation*} 
Let $b^s_1$ be a nonzero morphism in $\Hom_{A_0}(R,W_0)$, different to the morphism $\rho_s$
given in lemma~\ref{(DE)L:AXcero}, and consider $b^s_1$ as an element of $\Hom_{A_0}(R,Y_0)$.
Then it is possible to choose a basis in $\Hom_{A_0}(R,Y_0)$ such that for $1 \leq i \leq n_s$, $1 \leq k < \ell \leq n_r$ 
and $1 \leq r,s \leq t$ the following equations hold,
\begin{equation*} 
b^s_jq^r_{k,\ell} = \left\{
\begin{array}{ll}
b^s_{\ell}, & \text{ if } j=k, s=r,\\
0, & \text{ otherwise.}
\end{array} \right.
\end{equation*} 
\end{lema}
\bproof
Observe that there is a nonzero morphism $\Hom_{A_0}(Y^r_i,Y^s_j)$ if and only if $r=s$, $i \leq j$. 
Consider
\[
Z'_0 = \End_{A_0}(W_0)^{op} \oplus  \bigoplus_{\substack{1 < i \leq n_s\\ s=1,\ldots,t}}\End_{A_0}(Y^s_i)^{op},
\]
\[
Q_0 = \bigoplus_{\substack{1 \leq i < j \leq n_s\\ s=1,\ldots,t}}\End_{A_0}(Y^s_i,Y^s_j)^{op}.
\]
Then the opposed endomorphism algebra of $Y_0$, considered as matrix algebra, is triangular
with diagonal isomorphic to $Z'_0$. Again by lemma~\ref{(P)L:endExc} the $k$-algebra $Z'_0$ is trivial.
There is clearly a decomposition of $Z'_0$-$Z'_0$-bimodules
\[
 \Gamma'_0= \End_{A_0}(Y_0)^{op} \cong Z'_0 \oplus Q_0,
\]
thus $Y_0$ is complete admissible. By definition
\[
 A_0^{Y_0}=(T_{A_0}(0))^{Y_0}= T_{Z'_0}(Q_0^*).
\]
As in lemma~\ref{(DE)L:AXcero} there is an isomorphism of graded algebras $T_{Z'_0}(Q_0^*) \cong k\Delta_{Y_0}^{-1}$.
Fix a dual basis of the right $Z'_0$-module $Q_0$. For any $i=1,\ldots,n_s-1$ 
denote by $q_{i,i+1}^s$ a nonzero morphism between $Y^s_i$ and $Y^s_{i+1}$ (recall that every $Y^s_1$ is a
copy of $W_0$). All morphisms $q^s_{i,i+1}$ are irreducible and epimorphisms, since any $Y^s_i$ is generated by 
$W_0$ (lemma~\ref{(DE)L:alaDyn}),
\[
 \xymatrix@C=3pc{
W_0 \ar[r]^-{q^s_{1,2}} & Y^s_2 \ar[r]^-{q^s_{2,3}} & Y^s_3  \cdots Y^s_{n_s-2} \ar[r]^-{q^s_{n_s-2,n_s-1}} & Y^s_{n_s-1} 
\ar[r]^-{q^s_{n_s-1,n_s}} & Y^s_{n_s}.
}
\]
For any $1 \leq i < j \leq n_s$ take the composition
\[
 q^s_{i,j}=q^s_{i,i+1}q^s_{i+1,i+2}\ldots q^s_{j-1,j}.
\]
Then $\{(q^s_{i,j},\varepsilon^s_{i,j})\}$ is a dual basis of legible elements of $Q_0$,
which by construction satisfy the equation in the statement.

We want to give a basis of $\Hom_{A_0}(R,Y_0)$ starting from the chosen basis of $Q_0$.
Notice first that for any $s=1,\ldots,t$ there is a commutative diagram in $A_0$-mod (cf. lemma~\ref{(DE)L:ala}) 
given by
\[
\xymatrix@R=1pc{
& W_0 \ar[rd]^-{q^s_{1,2}} \\
X^s_{n_s-1} \ar[ru]^-{p^s_{n_s-1,n_s}} \ar[rd] & & Y^x_2 \\
& W^s_{2,n_s-1}. \ar[ru]
}
\]
Apply the functor $\Hom_{A_0}(R,-)$ to the diagram above and recall that for any $s$ we have that 
$\Hom_{A_0}(R,W^s_{2,n_s-1})=0$ (lemma~\ref{(DE)L:valores}). We obtain thus an exact sequence
\[
 \xymatrix@C=1.5pc{
0 \ar[r] & \Hom_{A_0}(R,X^s_{n_s-1}) \ar[r]^-{\Phi} & \Hom_{A_0}(R,W_0) \ar[r]^-{\Psi} & \Hom_{A_0}(R,Y^s_{2}) \ar[r] & 0,
}
\]
where $\Phi=\Hom(R,p^s_{n_s-1,n_s})$ and $\Psi=\Hom(R,q^s_{1,2})$.
Indeed, since $\Phi$ is not zero and $\Hom_{A_0}(R,X^s_{n_s-1})$ is one dimensional (lemma~\ref{(DE)L:valores})
we have that $\Phi$ is injective. Moreover, $\Psi$ is surjective, for it is nonzero and the vector space $\Hom_{A_0}(R,Y^s_{2})$
is also one dimensional (lemma~\ref{(DE)L:valores}).
Since the space $\Hom_{A_0}(R,W_0)$ has dimension two, we conclude that the sequence is exact, for $\Psi \circ \Phi=0$.
Then we have that
\[
0= \Psi(\Phi(a^s_{n_s-1}))=\Psi(a^s_{n_s-1}p^s_{n_s-1,n_s})=\Psi(\rho_s)=\rho_s q^s_{1,2}.
\]
In this way, starting for each $s=1,\ldots,t$ with a $b^s_1 \in \Hom_{A_0}(R,W_0)$ different from zero and from $\rho_s$,
and considering these morphisms as elements of $\Hom_{A_0}(R,W_0)$, we recursively obtain morphisms
$b^s_{i+1}=b^s_{i}q^s_{i,i+1}$ for $1 \leq i < n_s$.
Then $\{(b^s_i,\theta^s_i),(b,\theta),(b',\theta')\}$ (with $i=2,\ldots,n_s$, $s=1,\ldots,t$) is a dual basis
of legible elements of $\Hom_{A_0}(R,Y_0)$, where the restriction of $b$ and $b'$ to $\Hom_{A_0}(R,W_0)$ determine a basis,
and which by construction satisfies the formula in the statement.
\eproof

\begin{lema} \label{(DE)L:AY}
Take $A=A_0[R]$ and let $B$ be tha subalgebra of $A$ given by $B=A_0[0]$.
Denote by $\mathcal{A}$ the ditalgebra $\mathcal{A}=(A,0)$ with zero differential.
Then the $B$-module $Y=Y_0 \oplus S(\omega)$ is complete admissible and the reduced ditalgebra
$\mathcal{A}^Y$ is isomorphic to $(k\Delta_{Y_0}^{-1}[R^{Y_0}],\delta^y)$, where
$R^{Y_0}=\Hom_{A_0}(R,Y_0)$. The differential $\delta^y$ has the following explicit form
is solid arrows, in terms of the bases of lemma~\ref{(DE)L:AYcero},
\begin{center}
\begin{tabular}{c | l}
Vector in $\Hom_{A_0}(R,Y_0)^*$ & Its differential $\delta^y$ \\
\hline
$\theta_j^s$ & $\alpha_s\varepsilon^s_{1,j} \otimes \theta + \alpha_s'\varepsilon^s_{1,j} \otimes \theta' 
+ \sum_{1 < i < j}\varepsilon^s_{i,j} \otimes \theta^s_i$, \\ 
$\theta$ & $0$, \\
$\theta'$ & $0$, \\
\end{tabular}
\end{center}
where $b^s_1=\alpha_s b + \alpha'_s b'$, for $s=1,\ldots, t$.
\end{lema}
\bproof
Also consequence of proposition~\ref{(DE)P:extPuntRed} and lemma~\ref{(DE)L:AYcero}.
\eproof

Again by lemma~\ref{(DE)L:valores} we have that $\dimk_k\Hom_{A_0}(R,M_0)=1$ for any $[M_0] \in \mathcal{Y}_0$ and
$\dimk_k\Hom_{A_0}(R,W_0)=2$. Then, the quiver of the extension $A_0^{Y_0}[R^{Y_0}]$ can be formed by adding two solid 
arrows from the vertex $\omega$ to the vertex $m$ corresponding to the maximal module in $\Delta_{Y_0}^{-1}$, 
and a solid arrow from $\omega$ to the rest of vertices. 
The following figure shows the quiver associated to the reduced ditalgebra $\mathcal{A}^Y$ 
when $A$ is the path algebra of the extension of a linearly ordered Dynkin diagram $\mathbf{A}_{\ell}$.

%------------------------------------------------ 
\begin{equation} \label{(DE)F:brazoY}
 \xymatrix@C=3pc{
*+++[]{\bullet_{\omega}} 
\ar@/^2pc/[rrrrddd]^(.77){\theta_{\ell}}
\ar@/^1pc/[rrrddd]^(.7){\theta_{\ell-1}}
\ar@/^.5pc/[rddd]^(.59){\theta_{2}}
\ar[ddd]^-{\theta} \ar@<-1ex>[ddd]_-{\theta'}\\
\\
\\
\bullet_1 \ar@{.>}[r]^(.65){\varepsilon_{1,2}} \ar@{.>}@/^1.5pc/[rrr]^(.5){\varepsilon_{1,\ell-1}} 
\ar@{.>}@/_2pc/[rrrr]_(.5){\varepsilon_{1,\ell}}  
& \bullet_2 \ar@{.>}@/_.5pc/[rr]_(.5){\varepsilon_{2,\ell-1}} \ar@{.>}@/^1.5pc/[rrr]^(.65){\varepsilon_{2,\ell}}  
& \ldots & \bullet_{\ell-1} \ar@{.>}[r]^(.3){\varepsilon_{\ell-1,\ell}} 
& *++[]{\bullet_{\ell}} \\
}
\end{equation} 

\begin{lema} \label{(DE)L:XYRig}
The $A_0$-modules $X_0$ and $Y_0$ are rigid.
\end{lema}
\bproof
We prove the case $X_0$. By the Auslander-Reiten formulas in the he\-re\-di\-tary case,
for any pair of $A_0$-modules $M$, $N$ such that $[M],[N] \in \{[W_0]\} \cup \mathcal{X}_0$ we have that
\[
 \Ext^1_{A_0}(M,N) \cong D\Hom_{A_0}(\tau^{-1}N,M)=0,
\]
for $M \preceq W_0$ in $\Gamma(A_0)$ and $\tau^{-1}N \npreceq W_0$.

The case $Y_0$ is similar. If $M$, $N$ are $A_0$-modules such that $[M],[N] \in \{[W_0]\} \cup \mathcal{Y}_0$ we have that
\[
 \Ext^1_{A_0}(M,N) \cong D\Hom_{A_0}(N,\tau M)=0,
\]
for $W_0 \preceq N$ and $W_0 \npreceq \tau M$.
\eproof

%-------------------------------------
%-------------------------------------
Denote by $\Delta^X$ and $\Delta^Y$ the quivers determined by the tensor algebras $A^X$ and $A^Y$.
Let $\mathfrak{C}^X$ be the subquiver of $\Delta^X$ obtained by deleting the extension vertex 
$\omega$ and the vertex $m$ corresponding to the maximal module, and all arrows that start or end on them. 
Since $\Delta^Y$ is obtained from $\Delta^X$ by changing the orientation of the dotted arrows,
$\mathfrak{C}^Y=(\mathfrak{C}^X)^{op}$. 
We recover $\Delta^X$ from $\mathfrak{C}^X$ adding a copy of the Kronecker quiver 
$\xymatrix{\omega \ar@<-.5ex>[r] \ar@<.5ex>[r] & m}$, a solid arrow from $\omega$
to every vertex in $\mathfrak{C}^X$ and a dotted arrow from $m$ to every vertex in $\mathfrak{C}^X$.
In a similar way, we recover $\Delta^Y$ from $\mathfrak{C}^Y$.
The quadratic forms of $\Delta^X$ and $\Delta^Y$ coincide, and will be denoted by $q^{xy}$.
Let $\underline{d}=(d_1,\ldots,d_n,d_{\omega})$ be a positive integral vector in $\mathbb{Z}^{n+1}$.
Observe that the value $q^{xy}(\underline{d})$ does not depend on the integers $d_m$, $d_{\omega}$, 
only on their difference $r=d_{\omega}-d_m$, which will be called \textbf{rank of the vector} $\underline{d}$. 
Indeed, if we denote by $\mathfrak{C}$ the underlying graph of $\mathfrak{C}^X$ (and $\mathfrak{C}^Y$), we have
\begin{eqnarray} 
q^{xy}(\underline{d}) & = & d_{\omega}^2+\sum_{i=1}^n d_i^2-2d_md_{\omega} - \sum_{\substack{i=1 \\ i \neq m}}^n d_id_{\omega} 
 +\sum_{\substack{i=1 \\ i \neq m}}^nd_id_m
+ \sum_{\substack{\gamma \in \mathfrak{C} \\ \gamma:i\to j}}d_id_j = \nonumber \\
 & = & (d_{\omega} -d_m)^2 +\sum_{\substack{i=1 \\ i \neq m}}^nd_i^2 -\sum_{\substack{i=1 \\ i \neq m}}^nd_i(d_{\omega}-d_m) 
+ \sum_{\substack{\gamma \in \mathfrak{C} \\ \gamma:i\to j}}d_id_j= \nonumber \\
& = & r^2 + \sum_{\substack{i=1 \\ i \neq m}}^n d_i(d_i-r) + 
\sum_{\substack{\gamma \in \mathfrak{C} \\ \gamma:i\to j}}d_id_j.\label{(DE)E:raiz} 
\end{eqnarray} 
Therefore, if $\underline{d}$ is a root of $q^{xy}$, then for any nonnegative integer $u$ the vector
\[
 (d_1,\ldots,d_{m-1},d_m+u,d_{m+1},\ldots,d_n,d_{\omega}+u)
\]
is also a root of $q^{xy}$. Then  the roots of the reductions $\mathcal{A}^X$ and $\mathcal{A}^Y$ 
are given in countable families for each integer $r=(d_{\omega}-d_m)$. 
\begin{proposicion} \label{(DE)P:rango6}
Let $\mathcal{A}^X=(A^X,\delta^x)$ be the ditalgebra given in lemma~\ref{(DE)L:AX} and let
$\mathcal{A}^Y=(A^Y,\delta^y)$ be the ditalgebra of lemma~\ref{(DE)L:AY}.
Then the endomorphism algebras of any exceptional $\mathcal{A}^X$-module and any exceptional $\mathcal{A}^Y$-module 
are isomorphic to the field $k$. Assume that $w_0$ is the maximal root of the quadratic form $q_0$ and let
$\ell$ be the biggest integer that appears as an entry of $w_0$.
\begin{itemize}
 \item[a)] If $\underline{d}=(d_1,\ldots,d_n,d_{\omega})$ is the dimension vector of an exceptional $\mathcal{A}^Y$-module 
$\widetilde{M}$ such that $F^Y(\widetilde{M})$ is a preinjective $A$-module, then 
$r=(d_{\omega}-d_m)$ satisfies $1\leq r \leq \ell$.
 \item[b)] If $\underline{d}=(d_1,\ldots,d_n,d_{\omega})$ is the dimension vector of an exceptional $\mathcal{A}^X$-module 
$\widetilde{N}$ such that $F^X(\widetilde{N})$ is a posprojective $A$-module, then 
$r=(d_{\omega}-d_m)$ satisfies $-1\leq r < \sum_{\substack{i=1 \\ i \neq m}}^n d_i$.
\end{itemize}
\end{proposicion}
\bproof
By lemma~\ref{(DE)L:XYRig} it is clear that the $B$-modules $X$ and $Y$ are rigid.
The first claim is consequence of lemma~\ref{(P)L:endExc}, for the functors $F^X$ and $F^Y$ are rigid.

We show first the point $(a)$. Recall that the restriction of $F^Y(\widetilde{M})$ to $B$ is given by 
$Y \otimes_{Z_0'} \widetilde{M}$, hence $e_{\omega}F^Y(\widetilde{M}) \cong d_{\omega} S(\omega)$. Moreover we have
\[
 F^Y(\widetilde{M})|_{A_0}\cong Y_0 \otimes_{Z_0'} \widetilde{M} \cong \bigoplus_{i=1}^n d_iY_i, 
\]
thus $m_{(F^Y(\widetilde{M}))_0}=d_m$. In this way
\[
r=d_{\omega}-d_m=\dimk_k (F^Y(\widetilde{M}))_{\omega}-m_{(F^Y(\widetilde{M}))_0}=\Rnk(F^Y(\widetilde{M})). 
\]
By proposition~\ref{(DE)P:rank} and since $F^Y(\widetilde{M})$ is a preinjective in $A$-mod que which lies in 
$\check{\mathcal{U}}(W_0 \vee \mathcal{Y}_0)$, we have $r=\Rnk(F^Y(\widetilde{M}))=\sharp \overleftarrow{F^Y(\widetilde{M})}$. 
Due to proposition~\ref{(DE)C:invar} we conclude that $1\leq r \leq \ell$.

We show now $(b)$.
Since $\underline{d}=\vdim \widetilde{N}$ is the dimension vector of an exceptional representation $\widetilde{N}$ we have that 
$\underline{d}$ is root of $q^{xy}$. By formula~(\ref{(DE)E:raiz}) it is clear that $q^{xy}(\underline{d})=1$ implies $-1 \leq r$. 
Similarly as in the case $(a)$ it can be shown that $r=\Rnk(F^X(\widetilde{N}))$.
By proposition~\ref{(DE)P:rank} and since $F^X(\widetilde{N})$ is a posprojective in $A$-mod which lies in 
$\check{\mathcal{U}}(\mathcal{X}_0 \vee W_0)$, we have the following equation,
\[
r+\sharp \overrightarrow{F^X(\widetilde{N})}=\Rnk(F^X(\widetilde{N}))+\sharp \overrightarrow{F^X(\widetilde{N})}=
e_{(F^X(\widetilde{N}))_0}.
\]
The result follows observing that $e_{(F^X(\widetilde{N}))_0}=\sum_{\substack{i=1 \\ i \neq m}}^n d_i$.
\eproof

Tables~\ref{(DE)T:raicesMenor} and~\ref{(DE)T:raicesMayor} below show families of positive roots of the quadratic form $q^{xy}$ 
associated to the reduced ditalgebras $\mathcal{A}^X$ and $\mathcal{A}^Y$. 
In the notation $\mathcal{R}^{(r)}_{\ell_1,\ell_2,\ell_3}$ we specify the rank $r$ and the tubular type $\ell_1,\ell_2,\ell_3$ 
of the subgraph where the root is sincere. The presence of an exponent in the tubular typer $\ell_s^2$, $\ell_s^3$ or 
$\ell_s^{22}$ gives an account of the integers $d_i$ different from one outside the vertices $\omega$ and $m$. 

%-----------------------------------------------------------------
\begin{table} [!ht] 
\begin{center}
\begin{tabular}{c l l} 
\hline 
Rank & \multicolumn{1}{c}{Roots} & Sincere subgraph $\mathfrak{C}$ \\
\hline \\
-1 & $\mathcal{R}_{1,1,1}^{(-1)}$ & $\emptyset$\\
0  & ${}^{\dag}\mathcal{R}_{2,1,1}^{(0)}$ & $\xymatrix@!0{{}_1}$ \\
1  & $\mathcal{R}_{1,1,1}^{(1)}$; ${}^{\ddag}\mathcal{R}_{2,1,1}^{(1)}$; $\mathcal{R}_{2,2,1}^{(1)}$; $\mathcal{R}_{2,2,2}^{(1)}$ & 
$\emptyset $; $\xymatrix@!0{{}_1}$; $\xymatrix@!0@R=.6pc{{}_1 \\ {}_1}$; 
$\xymatrix@!0@R=.6pc{{}_1 \\ {}_1 \\ {}_1}$\\
2  & $\mathcal{R}_{2,2,2}^{(2)}$; $\mathcal{R}_{3,2,2}^{(2)}$; $\mathcal{R}_{3,3,2}^{(2)}$ 
& $\xymatrix@!0@R=.6pc{{}_1 \\ {}_1 \\ {}_1}$; 
$\xymatrix@!0@R=.6pc@C=1.6pc{{}_1 \ar@{.}[r] & {}_1 \\ {}_1 \\ {}_1}$; 
$\xymatrix@!0@R=.6pc@C=1.6pc{{}_1 \ar@{.}[r] & {}_1 \\ {}_1 \ar@{.}[r] & {}_1 \\ {}_1}$ \\
\hline
\end{tabular}
\caption{
Reduced positive roots of minor rank. The roots marked with $\dag$ (rank zero) are not dimension vector of any
exceptional $\mathcal{A}^Y$-module, while those marked with $\ddag$ do not correspond to exeptional $\mathcal{A}^X$-modules. 
} 
\label{(DE)T:raicesMenor}
\end{center}
\end{table}

\begin{lema} \label{(DE)L:tablaComp}
Let $q^{xy}$ be the quadratic form associated to the reduced ditalgebras $\mathcal{A}^X$ and $\mathcal{A}^Y$
and let $\underline{d}=(d_1,\ldots,d_n,d_{\omega})$ be a positive root of $q^{xy}$. If the rank $r=d_{\omega}-d_m$
of the vector $\underline{d}$ is smaller or equal than two, then the restriction of $\underline{d}$ to its support
is one of the vectors in table~\ref{(DE)T:raicesMenor}.
\end{lema}
\bproof
Recall that
\[
 q^{xy}(\underline{d})=r^2 + \sum_{\substack{i=1 \\ i \neq m}}^n d_i(d_i-r) + 
\sum_{\substack{\gamma \in \mathfrak{C} \\ \gamma:i\to j}}d_id_j.
\]
If $r \leq 1$ then all three summand on the right of the equation are nonnegative.\\
\underline{Case $r<0$}. In this case $q^{xy}(\underline{d})=1$ implies that $r=-1$ and hence
\[
 0=\sum_{\substack{i=1 \\ i \neq m}}^n d_i(d_i+1) + 
\sum_{\substack{\gamma \in \mathfrak{C} \\ \gamma:i\to j}}d_id_j.
\]
Then $d_i=0$ for $i \in \{1,\ldots,\widehat{m},\ldots,n\}$ thus $\underline{d}$ has the form $\mathcal{R}^{(-1)}_{1,1,1}$.\\
\underline{Case $r=0$}. Now we have
\[
 1=q^{xy}(\underline{d})=\sum_{\substack{i=1 \\ i \neq m}}^n d_i^2 + 
\sum_{\substack{\gamma \in \mathfrak{C} \\ \gamma:i\to j}}d_id_j,
\]
thus there exists a unique $i \in \{1,\ldots,\widehat{m},\ldots,n\}$ with $d_i$ different than zero and moreover $d_i=1$. 
In this case the vector $\underline{d}$ has the form $\mathcal{R}^{(0)}_{2,1,1}$. \\
\underline{Case $r=1$}. For this case we have the equality
\[
 0=\sum_{\substack{i=1 \\ i \neq m}}^n d_i(d_i-1) + 
\sum_{\substack{\gamma \in \mathfrak{C} \\ \gamma:i\to j}}d_id_j,
\]
hence if there is an entry $d_i$ with $i \in \{1,\ldots,\widehat{m},\ldots,n\}$ different than zero, then
$d_i=1$. Moreover, two vertices $i$, $j$ such that $d_i=d_j=1$ cannot be in the same branch of the reduced ditalgebras quivers
$\Delta^X$ or $\Delta^Y$, for in that case there exists a dotted arrow $\gamma \in \mathfrak{C}$
between the vertices $i$ and $j$ (see figures~(\ref{(DE)F:brazoX}) and~(\ref{(DE)F:brazoY})), and thus the summand
$\sum_{\substack{\gamma \in \mathfrak{C} \\ \gamma:i\to j}}d_id_j$ is greater than zero, which is impossible.
Hence, if there exist vertices $i$ such that $d_i>0$, then $d_i=1$ and each vertex $i$ lies in different branches.
This corresponds to the four rank one roots $\mathcal{R}^{(1)}_{1,1,1}$,
$\mathcal{R}^{(1)}_{2,1,1}$, $\mathcal{R}^{(1)}_{2,2,1}$ and $\mathcal{R}^{(1)}_{2,2,2}$.\\
\underline{Case $r=2$}. We have now an equality of the form
\[
-3= \sum_{\substack{i=1 \\ i \neq m}}^n d_i(d_i-2) + 
\sum_{\substack{\gamma \in \mathfrak{C} \\ \gamma:i\to j}}d_id_j.
\]
Assume there exist indices $i_1,\ldots,i_{\ell}$ in the same branch of $\Delta^X$ (or $\Delta^Y$) such that
$d_{i_j}>0$. The contribution of the entries $d_{i_j}$ in the equation above is given by
\[
 \sum_{j_1}^{\ell} d_{i_j}(d_{i_j}-2) + 
\sum_{\gamma \in \mathfrak{C}'}d_{s(\gamma)}d_{t(\gamma)},
\]
where $\mathfrak{C}'$ is the set of (dotted) arrows between the vertices $i_1,\ldots,i_{\ell}$. Since all these vertices
belong to a same branch, $\mathfrak{C}'$ has exactly $\ell(\ell-1)/2$ elements. Hence
\[
 \sum_{\gamma \in \mathfrak{C}'}d_{s(\gamma)}d_{t(\gamma)} \geq \frac{\ell(\ell-1)}{2},
\]
and it is clear that $\sum_{j_1}^{\ell} d_{i_j}(d_{i_j}-2) \geq -\ell$. Then
\[
 \sum_{j_1}^{\ell} d_{i_j}(d_{i_j}-2) + 
\sum_{\gamma \in \mathfrak{C}'}d_{s(\gamma)}d_{t(\gamma)} \geq \frac{\ell(\ell-1)}{2}-\ell=\frac{\ell^2-3\ell}{2}.
\]
In this way the contribution of each of the branches of $\Delta^X$ (or $\Delta^Y$) is
\begin{itemize}
 \item[a)] positive, if in the branch there are more than two vertices $i$ with $d_i>0$,
 \item[b)] $-1$, if in the branch there are exactly two vertices $i$, $j$ with $d_i,d_j>0$ and moreover $d_i=d_j=1$,
 \item[c)] $-1$, if in the branch there is exactly one vertex $i$ with $d_i>0$ and moreover $d_i=1$. 
\end{itemize}

To reach the equality $q^{xy}(\underline{d})=1$ we must substract three unities to the integer $r^2=4$. 
This is possible only when the quiver
$\Delta^X$ (or $\Delta^Y$) has three branches and all of them are in the situations $(b)$ or $(c)$ above.
This corresponds to the cases $\mathcal{R}^{(2)}_{2,2,2}$, $\mathcal{R}^{(2)}_{3,2,2}$ and $\mathcal{R}^{(2)}_{3,3,2}$
(the case $\mathcal{R}^{(2)}_{3,3,3}$ does not appear since $3,3,3$ is not the type of any Dynkin diagram, 
see table~\ref{T:tipoDynkin} in the previous chapter).

\eproof

%-----------------------------------------------------------------
\begin{table} [!hbt] 
\begin{center}
\renewcommand{\arraystretch}{1.2}
\begin{tabular}{c l l} 
\hline 
Rank & \multicolumn{1}{c}{Roots} &  Sincere subgraph $\mathfrak{C}$ \\
\hline \\
3  & $\mathcal{R}_{3,3,2}^{(3)}$; $\mathcal{R}_{4,3,2}^{(3)}$ 
& $\xymatrix@!0@R=.6pc@C=1.6pc{{}_1 \ar@{.}[r] & {}_1 \\ {}_1 \ar@{.}[r] & {}_1 \\ {}_1}$;
$\xymatrix@!0@R=.6pc@C=1.6pc{{}_1 \ar@{.}@/^3pt/[rr] \ar@{.}[r] & {}_1 \ar@{.}[r] & {}_1 \\ {}_1 \ar@{.}[r] & {}_1 \\ {}_1}$ \\  
   & $\mathcal{R}_{3,3,2^2}^{(3)}$; $\mathcal{R}_{4,3,2^2}^{(3)}$ 
& $\xymatrix@!0@R=.6pc@C=1.6pc{{}_1 \ar@{.}[r] & {}_1 \\ {}_1 \ar@{.}[r] & {}_1 \\ {}_2}$;
$\xymatrix@!0@R=.6pc@C=1.6pc{{}_1 \ar@{.}@/^3pt/[rr] \ar@{.}[r] & {}_1 \ar@{.}[r] & {}_1 \\ {}_1 \ar@{.}[r] & {}_1 \\ {}_2}$ \\ 
4  & $\mathcal{R}_{4,3,2^2}^{(4)}$; $\mathcal{R}_{5,3,2^2}^{(4)}$ 
& $\xymatrix@!0@R=.6pc@C=1.6pc{{}_1 \ar@{.}@/^3pt/[rr] \ar@{.}[r] & {}_1 \ar@{.}[r] & {}_1 \\ {}_1 \ar@{.}[r] & {}_1 \\ {}_2}$;
$\xymatrix@!0@R=.6pc@C=1.6pc{{}_1 \ar@{.}@/_4pt/[rrr] \ar@{.}@/^3pt/[rr] \ar@{.}[r] 
& {}_1 \ar@{.}@/^3pt/[rr] \ar@{.}[r] & {}_1 \ar@{.}[r] & {}_1  \\ {}_1 \ar@{.}[r] & {}_1 \\ {}_2}$ \\ 
   &  $\mathcal{R}_{4,3^2,2^2}^{(4)}$; $\mathcal{R}_{5,3^2,2^2}^{(4)}$ 
& $\xymatrix@!0@R=.6pc@C=1.6pc{{}_1 \ar@{.}@/^3pt/[rr] \ar@{.}[r] & {}_1 \ar@{.}[r] & {}_1 \\ {}_i \ar@{.}[r] & {}_j \\ {}_2}$;
$\xymatrix@!0@R=.6pc@C=1.6pc{{}_1 \ar@{.}@/_4pt/[rrr] \ar@{.}@/^3pt/[rr] \ar@{.}[r] 
& {}_1 \ar@{.}@/^3pt/[rr] \ar@{.}[r] & {}_1 \ar@{.}[r] & {}_1  \\ {}_i \ar@{.}[r] & {}_j \\ {}_2}$ \\ 
5  & $\mathcal{R}_{5,3^2,2^2}^{(5)}$; $\mathcal{R}_{5,3^2,2^3}^{(5)}$ 
& $\xymatrix@!0@R=.6pc@C=1.6pc{{}_1 \ar@{.}@/_4pt/[rrr] \ar@{.}@/^3pt/[rr] \ar@{.}[r] 
& {}_1 \ar@{.}@/^3pt/[rr] \ar@{.}[r] & {}_1 \ar@{.}[r] & {}_1  \\ {}_i \ar@{.}[r] & {}_j \\ {}_2}$;
$\xymatrix@!0@R=.6pc@C=1.6pc{{}_1 \ar@{.}@/_4pt/[rrr] \ar@{.}@/^3pt/[rr] \ar@{.}[r] 
& {}_1 \ar@{.}@/^3pt/[rr] \ar@{.}[r] & {}_1 \ar@{.}[r] & {}_1  \\ {}_i \ar@{.}[r] & {}_j \\ {}_3}$ \\
   & $\mathcal{R}_{5,3^{22},2^2}^{(5)}$; $\mathcal{R}_{5,3^{22},2^3}^{(5)}$ 
& $\xymatrix@!0@R=.6pc@C=1.6pc{{}_1 \ar@{.}@/_4pt/[rrr] \ar@{.}@/^3pt/[rr] \ar@{.}[r] 
& {}_1 \ar@{.}@/^3pt/[rr] \ar@{.}[r] & {}_1 \ar@{.}[r] & {}_1  \\ {}_2 \ar@{.}[r] & {}_2 \\ {}_2}$;
$\xymatrix@!0@R=.6pc@C=1.6pc{{}_1 \ar@{.}@/_4pt/[rrr] \ar@{.}@/^3pt/[rr] \ar@{.}[r] 
& {}_1 \ar@{.}@/^3pt/[rr] \ar@{.}[r] & {}_1 \ar@{.}[r] & {}_1  \\ {}_2 \ar@{.}[r] & {}_2 \\ {}_3}$ \\
6  & $\mathcal{R}_{5,3^{22},2^3}^{(6)}$; $\mathcal{R}_{5^2,3^{22},2^3}^{(6)}$ 
& $\xymatrix@!0@R=.6pc@C=1.6pc{{}_1 \ar@{.}@/_4pt/[rrr] \ar@{.}@/^3pt/[rr] \ar@{.}[r] 
& {}_1 \ar@{.}@/^3pt/[rr] \ar@{.}[r] & {}_1 \ar@{.}[r] & {}_1  \\ {}_2 \ar@{.}[r] & {}_2 \\ {}_3}$;
$\xymatrix@!0@R=.6pc@C=1.6pc{{}_a \ar@{.}@/_4pt/[rrr] \ar@{.}@/^3pt/[rr] \ar@{.}[r] 
& {}_b \ar@{.}@/^3pt/[rr] \ar@{.}[r] & {}_c \ar@{.}[r] & {}_d  \\ {}_2 \ar@{.}[r] & {}_2 \\ {}_3}$ \\
\hline
\end{tabular}
\caption{
Reduced positive roots of major rank. In the roots of rank 4 and 5 the pairs $(i,j)$ take the values
$(1,2)$ and $(2,1)$. In the last case the integers $a,b,c,d$ are greater or equal to one and their sum is five.
} 
\label{(DE)T:raicesMayor}
\end{center}
\end{table}

\section{Families of $\mathcal{A}^Y$-modules.} \label{S:AYEjem}
%----------------------------------------------------------------------

We present now some families of exceptional $\mathcal{A}^Y$-modules corresponding to the roots of minor rank
(table~\ref{(DE)T:raicesMenor}). Observe that the quiver of the reduced ditalgebra
$\mathcal{A}^Y$ containes a Kronecker subquiver $K_2$ (which consists in the arrows with zero differential
$\theta$ and $\theta'$). Hence, if $\widetilde{M}$ is an indecomposable $\mathcal{A}^Y$-module
with dimension vector $\underline{d}$ and rank zero, then the restriction $\widetilde{M}|_{K_2}$ is not rigid,
for $d_{\omega}-d_m=0$ (see proposition~\ref{(P)P:Kro2}). By lemma~\ref{(P)L:sub}, the representation
$\widetilde{M}$ cannot be exceptional, that is, there are no exceptional $\mathcal{A}^Y$-modules of rank zero.

We describe next the series of sincere roots of ranks $-1$, $0$ and $1$. 
\begin{equation*}
\xy
( 0,0)="R" *{\RaizM};
( 25,0)="R" *{\RaizC};
( 55,0)="R" *{\RaizUc};
( 85,0)="R" *{\RaizUu};
( 20,-30)="R" *{\RaizUd};
( 60,-30)="R" *{\RaizUt};
\endxy 
\vspace{1.5cm}
\end{equation*} 
For the families of roots $\mathcal{R}$ of rank different than zero given above we give series of exceptional 
$\mathcal{A}^Y$-modules $\widetilde{M}$ with corresponding dimension vectors. Recall that the matrices $I^{\downarrow}$, 
$I^{\uparrow}$ are formed from $I$ adding a row of zeros at the end and beginning of $I$ respectively. 
%---------------------------------------------------------------------- 
\begin{equation} \label{(DE)F:Rm}
\xy
( -4,-15)="A0" *{};
( 0,-15)="A1" *{};
( 4,-15)="A2" *{};
( 8,-15)="A3" *{};
( 12,-15)="A4" *{};
( 16,-15)="A5" *{};
( -2,-2)="B0" *{};
( 2,-2)="B1" *{};
( 6,-2)="B2" *{};
( 10,-2)="B3" *{};
( 14,-2)="B4" *{};
( 22,-2)="B5" *{};
( -30,0)="R" *{\widetilde{M}^{(-1)}_{1,1,1}[\ell]=};
(-8,0)="Gr" *{\xymatrix@R=3pc@C=1pc{k^{\ell-1} \ar@<-.5ex>[d]_-{I^{\downarrow}} \ar@<.5ex>[d]^-{I^{\uparrow}} \\ k^{\ell}}};
%-----------------------
"A0";"B0" **@{-};
"A1";"B1" **@{-};
"A2";"B2" **@{-};
"A3";"B3" **@{-};
"A4";"B4" **@{-};
"B0";"A1" **@{-};
"B1";"A2" **@{-};
"B2";"A3" **@{-};
"B3";"A4" **@{-};
"B4";"A5" **@{-};
\endxy
\vspace{.3cm}
\end{equation} 
This series of rank minus one roots corresponds to the posprojective $K_2$-re\-pre\-sen\-ta\-tions (proposition~\ref{(P)P:Kro2}),
and their images under $F^Y$ lie in the posprojective component of $\Gamma(A)$ (by corollary~\ref{(DE)C:ContComp}).
To the right we show the corresponding coefficient quiver for the case $\ell=6$. On the other hand, the
family of roots $\mathcal{R}_{1,1,1}^{(1)}$ of rank one correspond to the preinjective $K_2$-modules,
thus their images under $F^Y$ are preinjective $A$-modules. Recall also that the matrices $I^{\leftarrow}$, 
$I^{\rightarrow}$ are formed from $I$ by adding a column of zeros at the left and right of $I$ respectively.
\begin{equation} \label{(DE)F:R1}
 \xy
( -4,-15)="A0" *{};
( 0,-15)="A1" *{};
( 4,-15)="A2" *{};
( 8,-15)="A3" *{};
( 12,-15)="A4" *{};
( 16,-15)="A5" *{};
( -2,-2)="B0" *{};
%( -2,0)="B0F" *{\leftarrow};
( 2,-2)="B1" *{};
( 6,-2)="B2" *{};
%( 6,0)="B2F" *{\leftrightarrow};
( 10,-2)="B3" *{};
( 14,-2)="B4" *{};
%( 14,0)="B4F" *{\rightarrow};
( -30,0)="R" *{\widetilde{M}_{1,1,1}^{(1)}[\ell]=};
(-8,0)="Gr" *{\xymatrix@R=3pc@C=1pc{k^{\ell+1} \ar@<-.5ex>[d]_-{I^{\leftarrow}} \ar@<.5ex>[d]^-{I^{\rightarrow}} \\ k^{\ell}}};
%-----------------------
"A1";"B1" **@{-};
"A2";"B2" **@{-};
"A3";"B3" **@{-};
"A4";"B4" **@{-};
"B0";"A1" **@{-};
"B1";"A2" **@{-};
"B2";"A3" **@{-};
"B3";"A4" **@{-};
\endxy
\vspace{.3cm}
\end{equation} 
In the figure above we show the coefficient quiver for the case $\ell=4$. For the rank one roots of the form 
$\mathcal{R}_{2,1,1}^{(1)}$ we distinguish three cases.
Let $\widetilde{M}_{\theta}=I^{\leftarrow}$ and $\widetilde{M}_{\theta'}=I^{\rightarrow}$.
%---------------------------------------------------------------------- 

%---------------------------------------------------------------------- 
\begin{equation} \label{(DE)F:R1a}
\xy
(9,6)="R1" *{\scriptstyle \text{Case $\delta^y(\theta^s)=x(\varepsilon^s \otimes \theta)$}};
(9,3)="R1Et" *{\scriptstyle \widetilde{M}_{\theta^s}=[1,0,\ldots,0]};
(60,6)="R2" *{\scriptstyle \text{Case $\delta^y(\theta^s)=x'(\varepsilon^s \otimes \theta')$}};
(60,3)="R2Et" *{\scriptstyle \widetilde{M}_{\theta^s}=[0,\ldots,0,1]};
(35.3,-34)="R3" *{\scriptstyle \text{Case $\delta^y(\theta^s)=x(\varepsilon^s \otimes \theta)+x'(\varepsilon^s \otimes \theta')$}};
(35.3,-37)="R3Et" *{\scriptstyle \widetilde{M}_{\theta^s}=[0,\ldots 0,1,0\ldots,0]};
%----
( -30,0)="R" *{\widetilde{M}^{(1)}_{2,1,1}[\ell]=};
(-8,0)="Gr" *{\xymatrix@!0@R=2pc@C=1pc{k^{\ell+1} \ar@<-.5ex>[dd]_(.7){I^{\leftarrow}} 
\ar@<.5ex>[dd]^(.7){I^{\rightarrow}} 
\ar@/^15pt/[rrdddd]^(.8){\widetilde{M}_{\theta^s}}
\\ \\ k^{\ell} \ar@{.>}[rrdd] \\  \\ & & k }};
%-----
( 4,-15)="A1" *{};
( 8,-15)="A2" *{};
( 12,-15)="A3" *{};
( 16,-15)="A4" *{};
( 20,-15)="A5" *{};
( 2,-2)="B0" *{{}_{\bullet}};
( 2,-32)="C" *{};
( 2,0)="B0F" *{\leftarrow};
( 6,-2)="B1" *{};
( 10,-2)="B2" *{};
( 14,-2)="B3" *{};
( 18,-2)="B4" *{};
%----
( 29,-15)="A'1" *{};
( 33,-15)="A'2" *{};
( 37,-15)="A'3" *{};
( 41,-15)="A'4" *{};
( 45,-15)="A'5" *{};
( 27,-2)="B'0" *{};
( 35,-32)="C'" *{};
( 31,-2)="B'1" *{};
( 35,-2)="B'2" *{{}_{\bullet}};
( 35,0)="B'2F" *{\leftrightarrow};
( 39,-2)="B'3" *{};
( 43,-2)="B'4" *{};
%----
( 54,-15)="A''1" *{};
( 58,-15)="A''2" *{};
( 62,-15)="A''3" *{};
( 66,-15)="A''4" *{};
( 70,-15)="A''5" *{};
( 52,-2)="B''0" *{};
( 68,-32)="C''" *{};
( 56,-2)="B''1" *{};
( 60,-2)="B''2" *{};
( 64,-2)="B''3" *{};
( 68,-2)="B''4" *{{}_{\bullet}};
( 68,0)="B''4F" *{\rightarrow};
%-----------------------
"A1";"B1" **@{-};
"A2";"B2" **@{-};
"A3";"B3" **@{-};
"A4";"B4" **@{-};
"B0";"A1" **@{-};
"B1";"A2" **@{-};
"B2";"A3" **@{-};
"B3";"A4" **@{-};
"B0";"C" **@{-};
%-----
"A'1";"B'1" **@{-};
"A'2";"B'2" **@{-};
"A'3";"B'3" **@{-};
"A'4";"B'4" **@{-};
"B'0";"A'1" **@{-};
"B'1";"A'2" **@{-};
"B'2";"A'3" **@{-};
"B'3";"A'4" **@{-};
"B'2";"C'" **@{-};
%-----
"A''1";"B''1" **@{-};
"A''2";"B''2" **@{-};
"A''3";"B''3" **@{-};
"A''4";"B''4" **@{-};
"B''0";"A''1" **@{-};
"B''1";"A''2" **@{-};
"B''2";"A''3" **@{-};
"B''3";"A''4" **@{-};
"B''4";"C''" **@{-};
\endxy
\end{equation}
%\vspace{.3cm}
For the case $\delta^y(\theta^s)=x(\varepsilon^s \otimes \theta)$ take $\widetilde{M}_{\theta^s}=[1,0,\ldots,0]$ and whenever
$\delta^y(\theta^s)=x'(\varepsilon^s \otimes \theta')$ take $\widetilde{M}_{\theta^s}=[0,\ldots,0,1]$ ($x,x' \neq 0$). In the case
$\delta^y(\theta^s)=x(\varepsilon^s \otimes \theta) + x'(\varepsilon^s \otimes \theta')$ we can place the 1 in the 
matrix $\widetilde{M}_{\theta^s}$ in any position, $\widetilde{M}_{\theta^s}=[0,\ldots 0,1,0\ldots,0]$ ($s \in \{1,\ldots,t \}$). 
We give corresponding coefficient quivers in any case for $\ell=4$. To complete the rank one case observe that we can give
coefficient quivers for $\widetilde{M}^{(1)}_{2,1,1}[\ell]$ and $\widetilde{M}^{(1)}_{2,2,1}[\ell]$ which determine exceptional
representations for any form of the transformation $\delta^y$.
%---------------------------------------------------------------------- 
\begin{equation} \label{(DE)F:R1b}
\xy
( -30,0)="R" *{};
(0,5)="M" *{\widetilde{M}^{(1)}_{2,2,1}[\ell]};
(-8,0)="Gr" *{\xymatrix@!0@R=2pc@C=1pc{
k^{\ell+1} \ar@<-.5ex>[dd]_(.7){I^{\leftarrow}} \ar@<.5ex>[dd]^(.7){I^{\rightarrow}} 
\ar@/^10pt/[rrddd]|(.7){\widetilde{M}_{\theta^r}} \ar@/^15pt/[rrrdddd]|(.85){\widetilde{M}_{\theta^s}}
\\ \\ k^{\ell} \ar@{.>}[rrd] \ar@/_6pt/@{.>}[rrrdd]  \\ & &  k \\ & & & k }};
(5,-43)="R1Et" *{\scriptstyle \widetilde{M}_{\theta^r}=[1,0,\ldots,0]};
(5,-47)="R2Et" *{\scriptstyle \widetilde{M}_{\theta^s}=[0,\ldots ,0,1]};
%-----
( 12,-15)="A1" *{};
( 16,-15)="A2" *{};
( 20,-15)="A3" *{};
( 24,-15)="A4" *{};
( 28,-15)="A5" *{};
( 10,-2)="B0" *{};
( 14,-2)="B1" *{};
( 18,-2)="B2" *{};
( 22,-2)="B3" *{};
( 26,-2)="B4" *{};
( 10,-23)="C0" *{};
( 26,-31.5)="C4" *{};
%-----------------------
"A1";"B1" **@{-};
"A2";"B2" **@{-};
"A3";"B3" **@{-};
"A4";"B4" **@{-};
"B0";"A1" **@{-};
"B1";"A2" **@{-};
"B2";"A3" **@{-};
"B3";"A4" **@{-};
"B0";"C0" **@{-};
%"B2";"C2" **@{-};
"B4";"C4" **@{-};
\endxy
%%%%%%%%%%%%%%%%%%%%
\xy
( -30,0)="R" *{};
(0,5)="M" *{\widetilde{M}^{(1)}_{2,2,2}[\ell]};
(-8,0)="Gr" *{\xymatrix@!0@R=2pc@C=1pc{
k^{\ell+1} \ar@<-.5ex>[dd]_(.7){I^{\leftarrow}} \ar@<.5ex>[dd]^(.7){I^{\rightarrow}} 
\ar@/^10pt/[rrddd]_(.7){\widetilde{M}_{\theta^r}} \ar@/^15pt/[rrrdddd]|(.8){\widetilde{M}_{\theta^s}} 
\ar@/^20pt/[rrrrddddd]^(.9){\widetilde{M}_{\theta^u}}
\\ \\ k^{\ell} \ar@{.>}[rrd] \ar@/_3pt/@{.>}[rrrdd] \ar@/_6pt/@{.>}[rrrrddd]  \\ & &  k \\ & & & k \\ & & & & k }};
(5,-47)="R3Et" *{\scriptstyle \widetilde{M}_{\theta^u}=[0,\ldots 0,1,0\ldots,0]};
%-----
( 12,-15)="A1" *{};
( 16,-15)="A2" *{};
( 20,-15)="A3" *{};
( 24,-15)="A4" *{};
( 28,-15)="A5" *{};
( 10,-2)="B0" *{};
( 14,-2)="B1" *{};
( 18,-2)="B2" *{};
( 22,-2)="B3" *{};
( 26,-2)="B4" *{};
( 10,-23)="C0" *{};
( 18,-40)="C2" *{};
( 26,-31.5)="C4" *{};
%-----------------------
"A1";"B1" **@{-};
"A2";"B2" **@{-};
"A3";"B3" **@{-};
"A4";"B4" **@{-};
"B0";"A1" **@{-};
"B1";"A2" **@{-};
"B2";"A3" **@{-};
"B3";"A4" **@{-};
"B0";"C0" **@{-};
"B2";"C2" **@{-};
"B4";"C4" **@{-};
\endxy
\vspace{.5cm}
\end{equation}
Indeed, if there exists $i \in \{1,\ldots,t\}$ such that $\delta^y(\theta^i) \in k(\varepsilon^i \otimes \theta)$ 
(recall that $\widetilde{M}_{\theta}=I^{\leftarrow}$), then take 
$r=i$. Similarly, if there exists $j \in \{1,\ldots,t\}$ such that $\delta^y(\theta^j) \in k(\varepsilon^j \otimes \theta')$ take 
$s=j$ ($\widetilde{M}_{\theta}=I^{\leftarrow}$). Take also $\widetilde{M}_{\theta^r}=[1,0,\ldots,0]$, 
$\widetilde{M}_{\theta^s}=[0,\ldots,0,1]$ and $\widetilde{M}_{\theta^u}=[0,\ldots 0,1,0\ldots,0]$

\begin{lema} \label{(DE)L:excRedRnkUno}
The series of rank minus one $\mathcal{A}^Y$-modules $\widetilde{M}^{(-1)}_{1,1,1}[\ell]$ ($\ell \geq 1$) and
rank one $\mathcal{A}^Y$-modules $\widetilde{M}^{(1)}_{1,1,1}[\ell]$, 
$\widetilde{M}^{(1)}_{2,1,1}[\ell]$, $\widetilde{M}^{(1)}_{2,2,1}[\ell]$ and $\widetilde{M}^{(1)}_{2,2,2}[\ell]$ ($\ell \geq 0$) 
given above consist in exceptional representations.
\end{lema}
\bproof
As shown in proposition~\ref{(P)P:Kro2} the representations $\widetilde{M}^{(-1)}_{1,1,1}[\ell]$ 
and $\widetilde{M}^{(1)}_{1,1,1}[\ell]$ are exceptional. Since the dimension vectors of the given representations
are roots of $q^{xy}$, by lemma~\ref{(P)L:Euler} it is enough to show that their endomorphism algebras
have dimension one.

Consider first the family $\widetilde{M}^{(1)}_{2,1,1}[\ell]$ with
$\delta^y(\theta^s)=x(\varepsilon^s \otimes \theta)$ ($x \neq 0$).
In the following diagram we give an endomorphism $f=(f^0,f^1)$ through the matrices $A^0$, $B^0$ and $C^0=[c^0]$ (part $f^0$)  
and $A^1=[a^1_1\ldots,a^1_{\ell}]$ (part $f^1$). 
%---------------------------------------------------------------------- 
\begin{equation*}
 \xymatrix@!0@R=2pc@C=1.5pc{
k^{\ell+1} \ar@<-.5ex>[dd]_(.6){I^{\leftarrow}} \ar@<.5ex>[dd]^(.6){I^{\rightarrow}} 
\ar@/^15pt/[rrrdddd]|(.85){[1,0 \ldots 0,0]} \ar[rrrrrrrrrr]^-{A^0} & & & & & & & & & & 
k^{\ell+1} \ar@<-.5ex>[dd]_(.6){I^{\leftarrow}} \ar@<.5ex>[dd]^(.6){I^{\rightarrow}} 
\ar@/^15pt/[rrrdddd]|(.85){[1,0 \ldots 0,0]}
\\ \\ k^{\ell} \ar@/^3pt/@{.>}[rrrrrrrrrrrrrdd]|-{[a^1_1\ldots a^1_{\ell}]} \ar[rrrrrrrrrr]^-{B^0} & & & & & & & & & & k^{\ell}
\\ \\ & & & k \ar[rrrrrrrrrr]_-{[c^0]} & & & & & & & & & & k }
\end{equation*}
Since the arrows $\theta$ and $\theta'$ determine a Kronecker quiver 
we have $A^0=[aI]$ and $B^0=[aI]$ for some scalar $a$. Moreover, by the definition of morphism~(\ref{(A)EQ:hom}) we have
\[
 [1,0 \ldots 0,0]A^0=C^0[1,0 \ldots 0,0]+xA^1I^{\leftarrow},
\]
that is,
\[
 [a-c^0,0 \ldots 0,0]=x[0,a^1_1 \ldots a^1_{\ell}],
\]
and hence $a=c^0$ and $A^1=0$. Then $\End_{\mathcal{A}^Y}(\widetilde{M}^{(1)}_{2,1,1}[\ell])=aId$.
The case $\delta^y(\theta^s)=x'(\varepsilon^s \otimes \theta')$ ($x' \neq 0$) can be treated in a similar way and 
for the case $\delta^y(\theta^s)=x(\varepsilon^s \otimes \theta)+x'(\varepsilon^s \otimes \theta')$ 
(figure~\ref{(DE)F:R1a}) observe that the equation above has the form
\[
 [\ldots 0,a-c^0,0 \ldots]=x[0,a^1_1 \ldots a^1_{\ell}]+x'[a^1_1 \ldots a^1_{\ell},0],
\]
therefore we conclude again that $a=c^0$ and $a^1_1=\ldots=a^1_{\ell}=0$.

We compute now the endomorphism algebra of $\widetilde{M}^{(1)}_{2,2,2}[\ell]$ 
(the case $\widetilde{M}^{(1)}_{2,2,1}[\ell]$ is similar) with help of the following diagram,
%---------------------------------------------------------------------- 
\begin{equation*}
 \xymatrix@!0@R=2pc@C=2pc{
k^{\ell+1} \ar@<-.5ex>[dd]_(.6){I^{\leftarrow}} \ar@<.5ex>[dd]^(.6){I^{\rightarrow}} 
\ar@/^10pt/[rrddd]|(.83){[1,0 \ldots 0,0]} \ar@/^15pt/[rrrdddd]|(.86){[0,0 \ldots 0,1]} 
\ar@/^20pt/[rrrrddddd]|(.9){[ \ldots 0,1,0 \ldots]} \ar[rrrrrrrr]^-{A^0} & & & & & & & & 
k^{\ell+1} \ar@<-.5ex>[dd]_(.6){I^{\leftarrow}} \ar@<.5ex>[dd]^(.6){I^{\rightarrow}} 
\ar@/^10pt/[rrddd]|(.83){[1,0 \ldots 0,0]} \ar@/^15pt/[rrrdddd]|(.86){[0,0 \ldots 0,1]} 
\ar@/^20pt/[rrrrddddd]|(.9){[ \ldots 0,1,0 \ldots]}
\\ \\ k^{\ell} \ar@/^3pt/@{.>}[rrrrrrrrrrd]|(.7){[a^1_1\ldots a^1_{\ell}]} \ar@/^3pt/@{.>}[rrrrrrrrrrrdd]|(.75){[b^1_1\ldots b^1_{\ell}]}
\ar@/^3pt/@{.>}[rrrrrrrrrrrrddd]|(.8){[c^1_1\ldots c^1_{\ell}]} \ar[rrrrrrrr]^-{B^0} & & & & & & & & k^{\ell}
\\ & &  k \ar[rrrrrrrr]|(.5){[c^0]} & & & & & & & & k
\\ & & & k \ar[rrrrrrrr]|(.5){[d^0]} & & & & & & & & k 
\\ & & & & k \ar[rrrrrrrr]|(.5){[e^0]} & & & & & & & & k }
\end{equation*}
Since the central square is commutative we have $A^0=[aI]$ and $B^0=[aI]$ (for some $a \in k$). 
The following equations are obtained from the definition of morphism and the conditions 
$\delta^y(\theta^r)=\varepsilon^r \otimes \theta$, $\delta^y(\theta^s)=\varepsilon^s \otimes \theta'$ and 
$\delta^y(\theta^u)=x(\varepsilon^u \otimes \theta)+x'(\varepsilon^u \otimes \theta')$ ($x,x' \neq 0$),
\[
 [a-c^0,0 \ldots 0,0]=[0,a^1_1 \ldots a^1_{\ell}],
\]
\[
 [0,0 \ldots 0,a-d^0]=[0,b^1_1 \ldots b^1_{\ell}],
\]
\[
 [\ldots 0,a-e^0,0 \ldots]=x[0,c^1_1 \ldots c^1_{\ell}]+x'[c^1_1 \ldots c^1_{\ell},0],
\]
then $a=c^0=d^0=e^0$ and $a^1_1=\ldots=a^1_{\ell}=b^1_1=\ldots=b^1_{\ell}=c^1_1=\ldots=c^1_{\ell}=0$. 
Hence the endomorphism algebra of $\widetilde{M}^{(1)}_{2,2,2}[\ell]$ is one dimensional,
which completes the proof.
\eproof

The case of rank two $\mathcal{A}^Y$-representations is considered by cases depending on the parity of the parameter
$\ell$. In figure~(\ref{(DE)F:R2a}) we explicitly show a series of modules for $\ell$ even
(representations $\widetilde{N}^{(2)}_{3,3,2}[\ell]$) corresponding to the roots of type $\mathcal{R}^{(2)}_{3,3,2}$ 
together with coefficient quivers for tha case $\ell=8$. Representations corresponding to roots of type
$\mathcal{R}^{(2)}_{3,2,2}$ and type $\mathcal{R}^{(2)}_{2,2,2}$ can be obtained deleting from the module 
$\widetilde{N}^{(2)}_{3,3,2}[\ell]$ the matrices $\widetilde{N}_{\theta^v}$ and/or $\widetilde{N}_{\theta^w}$ 
(dotted arrows in the coefficient quiver), depending on the value of the differential $\delta^y$.
The case $\ell$ odd is described in figure~(\ref{(DE)F:R2b}).
%---------------------------------------------------------------------- 
\begin{equation} \label{(DE)F:R2a}
\xy
(21,5)="Et" *{\widetilde{N}^{(2)}_{3,3,2}[\ell], \scriptstyle \text{ case $\ell$ even}};
(-18,0)="Gr" *{\xymatrix@!0@R=2.4pc@C=1.5pc{
k^{\ell+2} \ar@<-.5ex>[dd]_(.6){\widetilde{N}} \ar@<.5ex>[dd]^(.6){\widetilde{N}'} 
\ar@/^10pt/[rrddd]|(.6){\widetilde{N}_{\theta^r}} \ar@/^15pt/[rrrdddd]|(.85){\widetilde{N}_{\theta^s}} 
\ar@/^20pt/[rrrrddddd]|(.88){\widetilde{N}_{\theta^u}} \ar@/^10pt/[rrrrrddd]|(.65){\widetilde{N}_{\theta^v}}
\ar@/^13pt/[rrrrrdddd]|(.65){\widetilde{N}_{\theta^w}}
\\ \\ k^{\ell} \ar@{.>}[rrd] \ar@/_3pt/@{.>}[rrrdd] \ar@/_6pt/@{.>}[rrrrddd] 
\ar@/^5pt/@{.>}[rrrrrd] \ar@/^7pt/@{.>}[rrrrrdd]
\\ & &  k \ar@{.>}[rrr] & & & k \\ & & & k \ar@{.>}[rr] & & k \\ & & & & k }};
(-30,-43)="R1Et" *{\scriptstyle \widetilde{N}=I_{\ell/2}^{\leftarrow} \oplus I_{\ell/2}^{\leftarrow}};
(-30,-47)="R2Et" *{\scriptstyle \widetilde{N}'=I_{\ell/2}^{\rightarrow} \oplus I_{\ell/2}^{\rightarrow}};
(61,-16)="S1Et" *{\scriptstyle \widetilde{N}_{\theta^r}=[1,0 \ldots | \ldots 0,0]};
(61,-20)="S2Et" *{\scriptstyle \widetilde{N}_{\theta^s}=[0,0 \ldots | \ldots 0,1]};
(61,-24)="S3Et" *{\scriptstyle \widetilde{N}_{\theta^u}=[\ldots 0,1,0 \ldots | \ldots 0,1,0 \ldots]};
(61,-30)="S4Et" *{\scriptstyle \widetilde{N}_{\theta^v}=[0,0 \ldots |1,0 \ldots ]};
(61,-34)="S5Et" *{\scriptstyle \widetilde{N}_{\theta^w}=[ \ldots 0,1 | \ldots 0,0 ]};
%-----
( 4,-20)="A1" *{};
( 8,-20)="A2" *{};
( 12,-20)="A3" *{};
( 16,-20)="A4" *{};
( 20,-20)="A5" *{};
( 2,-2)="B0" *{};
( 2,-30)="C" *{};
( 6,-2)="B1" *{};
( 10,-2)="B2" *{};
( 14,-2)="B3" *{};
( 18,-2)="B4" *{{}_{\bullet}};
( 18,0)="B4F" *{\rightarrow};
( 18,-40)="E" *{};
%----
( 29,-20)="A'1" *{};
( 33,-20)="A'2" *{};
( 37,-20)="A'3" *{};
( 41,-20)="A'4" *{};
( 45,-20)="A'5" *{};
( 27,-2)="B'0" *{{}_{\bullet}};
( 27,0)="B'0F" *{\leftarrow};
( 43,-40)="C'" *{};
( 31,-2)="B'1" *{};
( 35,-2)="B'2" *{};
( 39,-2)="B'3" *{};
( 43,-2)="B'4" *{};
( 23,-52)="D" *{};
( 27,-30)="E'" *{};
%-----------------------
"A1";"B1" **@{-};
"A2";"B2" **@{-};
"A3";"B3" **@{-};
"A4";"B4" **@{-};
"B0";"A1" **@{-};
"B1";"A2" **@{-};
"B2";"A3" **@{-};
"B3";"A4" **@{-};
"B0";"C" **@{-};
"B4";"E" **@{--};
%-----
"A'1";"B'1" **@{-};
"A'2";"B'2" **@{-};
"A'3";"B'3" **@{-};
"A'4";"B'4" **@{-};
"B'0";"A'1" **@{-};
"B'1";"A'2" **@{-};
"B'2";"A'3" **@{-};
"B'3";"A'4" **@{-};
"B'4";"C'" **@{-};
"B'0";"E'" **@{--};
%-----
"B2"; "D"  **\crv{(9,-52) & (20,-38)};
"B'2"; "D" **\crv{(37,-52) & (26,-38)};
%-----
\endxy
\vspace{1cm}
\end{equation}
To give the transformations $\widetilde{N}_{\theta^r}$, $\widetilde{N}_{\theta^s}$ and $\widetilde{N}_{\theta^u}$ we
can assume, without lost of generality, that the arrows $\theta^{r}$, $\theta^{s}$ and $\theta^{u}$ satisfy
\[
 \delta^y(\theta^r)=\varepsilon^r \otimes \theta, \qquad \delta^y(\theta^s)=\varepsilon^s \otimes \theta' 
\quad \text{and} \quad \delta^y(\theta^u)=x(\varepsilon^u \otimes \theta) + x'(\varepsilon^u \otimes \theta'),
\]
where $x$ and $x'$ are nonzero scalars. Take $\widetilde{N}_{\theta^r}=[1,0 \ldots | \ldots 0,0]$,
$\widetilde{N}_{\theta^s}=[0,0 \ldots | \ldots 0,1]$ and $\widetilde{N}_{\theta^u}=[\ldots 0,1,0 \ldots | \ldots 0,1,0 \ldots]$,
where the positions of the 1's in the last matrix are arbitrary inside the blocks separated by the central line.

For the odd case we make the same assumption in the arrows $\theta^r$, $\theta^s$ and $\theta^u$
and show as example the coefficient quiver for the case $\ell=7$.
%---------------------------------------------------------------------- 
\begin{equation} \label{(DE)F:R2b}
\xy
(21,5)="Et" *{\widetilde{M}^{(2)}_{3,3,2}[\ell], \scriptstyle \text{ case $\ell$ odd}};
(-18,0)="Gr" *{\xymatrix@!0@R=2.4pc@C=1.5pc{
k^{\ell+2} \ar@<-.5ex>[dd]_(.6){\widetilde{M}} \ar@<.5ex>[dd]^(.6){\widetilde{M}'} 
\ar@/^10pt/[rrddd]|(.6){\widetilde{M}_{\theta^r}} \ar@/^15pt/[rrrdddd]|(.85){\widetilde{M}_{\theta^s}} 
\ar@/^20pt/[rrrrddddd]|(.88){\widetilde{M}_{\theta^u}} \ar@/^10pt/[rrrrrddd]|(.65){\widetilde{M}_{\theta^v}}
\ar@/^13pt/[rrrrrdddd]|(.65){\widetilde{M}_{\theta^w}}
\\ \\ k^{\ell} \ar@{.>}[rrd] \ar@/_3pt/@{.>}[rrrdd] \ar@/_6pt/@{.>}[rrrrddd] 
\ar@/^5pt/@{.>}[rrrrrd] \ar@/^7pt/@{.>}[rrrrrdd]
\\ & &  k \ar@{.>}[rrr] & & & k \\ & & & k \ar@{.>}[rr] & & k \\ & & & & k }};
(-30,-43)="R1Et" *{\scriptstyle \widetilde{M}=I_{\ell+1/2}^{\leftarrow} \oplus I_{\ell-1/2}^{\leftarrow}};
(-30,-47)="R2Et" *{\scriptstyle \widetilde{M}'=I_{\ell+1/2}^{\rightarrow} \oplus I_{\ell-1/2}^{\rightarrow}};
(61,-16)="S1Et"  *{\scriptstyle \widetilde{M}_{\theta^r}=[1,0 \ldots 0,0|1,0 \ldots 0]};
(61,-20)="S2Et"  *{\scriptstyle \widetilde{M}_{\theta^s}=[0,0 \ldots 0,0|0 \ldots 0,1]};
(61,-24)="S3Et"  *{\scriptstyle \widetilde{M}_{\theta^u}=[0,0 \ldots 0,0| \ldots 0,1,0 \ldots]};
(61,-30)="S4Et"  *{\scriptstyle \widetilde{M}_{\theta^v}=[0,0 \ldots 0,0|1,0 \ldots 0]};
(61,-34)="S5Et"  *{\scriptstyle \widetilde{M}_{\theta^w}=[0,0 \ldots 0,1 | \ldots 0,0 ]};
%-----
( 4,-20)="A1" *{};
( 8,-20)="A2" *{};
( 12,-20)="A3" *{};
( 16,-20)="A4" *{};
( 20,-20)="A5" *{};
( 2,-2)="B0" *{};
(31,-52)="C" *{};
( 6,-2)="B1" *{};
( 10,-2)="B2" *{};
( 14,-2)="B3" *{};
( 18,-2)="B4" *{{}_{\bullet}};
( 18,0)="B4F" *{\rightarrow};
( 18,-40)="E" *{};
%----
( 29,-20)="A'1" *{};
( 33,-20)="A'2" *{};
( 37,-20)="A'3" *{};
%( 41,-20)="A'4" *{};
( 45,-20)="A'5" *{};
( 27,-2)="B'0" *{{}_{\bullet}};
( 27,0)="B'0F" *{\leftarrow};
( 39,-40)="C'" *{};
( 31,-2)="B'1" *{};
( 35,-2)="B'2" *{};
( 39,-2)="B'3" *{};
%( 43,-2)="B'4" *{};
( 2,-30)="D" *{};
( 27,-30)="E'" *{};
%-----------------------
"A1";"B1" **@{-};
"A2";"B2" **@{-};
"A3";"B3" **@{-};
"A4";"B4" **@{-};
"B0";"A1" **@{-};
"B1";"A2" **@{-};
"B2";"A3" **@{-};
"B3";"A4" **@{-};
"B'1";"C" **@{-};
"B4";"E" **@{--};
%-----
"A'1";"B'1" **@{-};
"A'2";"B'2" **@{-};
"A'3";"B'3" **@{-};
%"A'4";"B'4" **@{-};
"B'0";"A'1" **@{-};
"B'1";"A'2" **@{-};
"B'2";"A'3" **@{-};
%"B'3";"A'4" **@{-};
"B'3";"C'" **@{-};
"B'0";"E'" **@{--};
%-----
"B0"; "D"  **@{-};
"B'0"; "D" **\crv{(25,-32) & (14,-18)};
%-----
\endxy
\vspace{1cm}
\end{equation}

\begin{lema} \label{(DE)L:excRedRnkDos}
The series of rank two representations $\widetilde{N}^{(2)}_{3,3,2}[\ell]$ ($\ell \geq 0$ even) and
$\widetilde{M}^{(2)}_{3,3,2}[\ell]$ ($\ell \geq 1$ odd) given above 
consist in exceptional $\mathcal{A}^Y$-modules. Moreover, the representations
$\widetilde{N}^{(2)}_{3,2,2}[\ell]$, $\widetilde{N}^{(2)}_{2,2,2}[\ell]$ and 
$\widetilde{M}^{(2)}_{3,2,2}[\ell]$, $\widetilde{M}^{(2)}_{2,2,2}[\ell]$ obtained by re\-mo\-ving the matrices
corresponding to the arrows $\theta^r$ and/or $\theta^s$ are also exceptional.
\end{lema}
\bproof
We compute again endomorphism algebras and use lemma~\ref{(P)L:Euler}.
Consider first the case $\ell$ even (representations $\widetilde{N}^{(2)}_{3,3,2}[\ell]$). We give an endomorphism $(f^0,f^1)$
of $\widetilde{N}^{(2)}_{3,3,2}[\ell]$ through the matrices $A^0-G^0$, which conform the part $f^0$ and $A^1-G^1$ for the part
$f^1$ (which are shown below in separated diagrams for clarity).
%---------------------------------------------------------------------- 
\begin{equation*}
 \xymatrix@!0@R=1.7pc@C=1.5pc{
k^{\ell+2} \ar@<-.5ex>[dd]_(.6){Q} \ar@<.5ex>[dd]^(.6){Q'} 
\ar@/^10pt/[rrddd]|(.8){R} \ar@/^15pt/[rrrrrdddd]|(.78){V} \ar@/^20pt/[rrrddddd]|(.8){S} 
\ar@/^30pt/[rrrrddddddd]|(.85){U} \ar@/^30pt/[rrrrrrdddddd]|(.83){W}
\ar[rrrrrrrrrr]^-{A^0} & & & & & & & & & & 
k^{\ell+2} \ar@<-.5ex>[dd]_(.6){Q} \ar@<.5ex>[dd]^(.6){Q'} 
\ar@/^10pt/[rrddd]|(.8){R} \ar@/^15pt/[rrrrrdddd]|(.78){V} \ar@/^20pt/[rrrddddd]|(.8){S} 
\ar@/^30pt/[rrrrddddddd]|(.85){U} \ar@/^30pt/[rrrrrrdddddd]|(.83){W}
\\ \\ k^{\ell} \ar[rrrrrrrrrr]^-{B^0} & & & & & & & & & & k^{\ell}
\\ & &  k \ar[rrrrrrrrrr]|(.5){C^0} & & & & & & & & & & k
\\ & & & & & k \ar[rrrrrrrrrr]|(.5){F^0} & & & & & & & & & & k
\\ & & & k \ar[rrrrrrrrrr]|(.5){D^0} & & & & & & & & & & k 
\\ & & & & & & k \ar[rrrrrrrrrr]|(.5){G^0} & & & & & & & & & & k
\\ & & & & k \ar[rrrrrrrrrr]|(.5){E^0} & & & & & & & & & & k }
\end{equation*}
In this case we have $Q=\widetilde{N}=I^{\leftarrow}_{\ell/2}\oplus I^{\leftarrow}_{\ell/2}$ and 
$Q'=\widetilde{N}'=I^{\rightarrow}_{\ell/2}\oplus I^{\rightarrow}_{\ell/2}$. In addition we have
$R=\widetilde{N}_{\theta^r}=[1,0 \ldots | \ldots 0,0]$,
$S=\widetilde{N}_{\theta^s}=[0,0 \ldots | \ldots 0,1]$, $U=\widetilde{N}_{\theta^u}=[\ldots 0,1,0 \ldots | \ldots 0,1,0 \ldots]$,
$V=\widetilde{N}_{\theta^v}=[0,0 \ldots |1,0 \ldots ]$ and $W=\widetilde{N}_{\theta^w}=[ \ldots 0,1 | \ldots 0,0 ]$
(figure~\ref{(DE)F:R2a}).
%---------------------------------------------------------------------- 
\begin{equation*}
 \xymatrix@!0@R=1.7pc@C=1.5pc{
k^{\ell+2} \ar@<-.5ex>[dd]_(.6){Q} \ar@<.5ex>[dd]^(.6){Q'} 
\ar@/^10pt/[rrddd]|(.8){R} \ar@/^15pt/[rrrrrdddd]|(.78){V} \ar@/^20pt/[rrrddddd]|(.8){S} 
\ar@/^30pt/[rrrrddddddd]|(.85){U} \ar@/^30pt/[rrrrrrdddddd]|(.83){W} & & & & & & & & & & 
k^{\ell+2} \ar@<-.5ex>[dd]_(.6){Q} \ar@<.5ex>[dd]^(.6){Q'} 
\ar@/^10pt/[rrddd]|(.8){R} \ar@/^15pt/[rrrrrdddd]|(.78){V} \ar@/^20pt/[rrrddddd]|(.8){S} 
\ar@/^30pt/[rrrrddddddd]|(.85){U} \ar@/^30pt/[rrrrrrdddddd]|(.83){W}
\\ \\ k^{\ell} \ar@/^10pt/@{.>}[rrrrrrrrrrrrd]|(.6){A^1} \ar@/^15pt/@{.>}[rrrrrrrrrrrrrddd]|(.6){B^1}
\ar@/^20pt/@{.>}[rrrrrrrrrrrrrrddddd]|(.6){C^1} \ar@/_40pt/@{.>}[rrrrrrrrrrrrrrrdd]|(.6){D^1} 
\ar@/_40pt/@{.>}[rrrrrrrrrrrrrrrrdddd]|(.6){E^1}
& & & & & & & & & & k^{\ell} 
\\ & & k \ar@{.>}[rrrrrrrrrrrrrd]|(.5){F^1} & & & & & & & & & & k
\\ & & & & & k & & & & & & & & & & k
\\ & & & k \ar@{.>}[rrrrrrrrrrrrrd]|(.5){G^1} & & & & & & & & & & k 
\\ & & & & & & k & & & & & & & & & & k
\\ & & & & k & & & & & & & & & & k }
\end{equation*}
The matrices $A^0-G^0$ and $A^1-G^1$ are subject to the following equations which are obtained
from the definition of morphism~(\ref{(A)EQ:hom}). Recall we have assumed that $\delta^y(\theta)=\delta^y(\theta')=0$, 
$\delta^y(\theta^r)=\varepsilon^r \otimes \theta$, $\delta^y(\theta^s)=\varepsilon^s \otimes \theta'$ and 
$\delta^y(\theta^u)=x(\varepsilon^u \otimes \theta)+x'(\varepsilon^u \otimes \theta')$,
\begin{equation} \label{(DE)EQ:prueba1}
\left\{
\begin{array}{l}
\scriptstyle QA^0=B^0Q,\\
\scriptstyle Q'A^0=B^0=Q',
\end{array} \right.
\left\{
\begin{array}{l}
\scriptstyle RA^0=C^0R+A^1Q,\\
\scriptstyle SA^0=D^0S+B^1Q',\\
\scriptstyle UA^0=E^0U+xC^1Q+x'C^1Q',
\end{array} \right.
\left\{
\begin{array}{l}
\scriptstyle VA^0=F^0V+D^1Q+F^1R,\\
\scriptstyle WA^0=G^0W+E^1Q'+G^1S.
\end{array} \right.
\end{equation}
From the pair of equations on the left it follows that the matrices $A^0$ and $B^0$ have the form
\[
A^0=\left[ \begin{matrix} aI & b'I \\ a'I & bI \end{matrix} \right], \qquad
B^0=\left[ \begin{matrix} aI & b'I \\ a'I & bI \end{matrix} \right] 
\]
for scalars $a,b,a',b'\in k$. The first two equations in the middle~(\ref{(DE)EQ:prueba1}) can be rewritten as follows,
\begin{eqnarray}
 \left[ a-c^0,0 \ldots 0,0|b',0 \ldots 0,0 \right] & = & \left[ 0,a^1_1 \ldots a^1_{\ell/2}|0,a^1_{\ell/2+1} \ldots a^1_{\ell} \right], \nonumber \\
 \left[ 0,0 \ldots 0,a'|0,0 \ldots 0,b-d^0 \right] & = & \left[ b^1_1 \ldots b^1_{\ell/2},0|b^1_{\ell/2+1} \ldots b^1_{\ell},0 \right], \nonumber
\end{eqnarray}
in particular we have that $a'=b'=0$, and hence the third equation in the middle has the form
\begin{eqnarray}
\left[ \ldots 0,a-e^0,0 \ldots| \ldots 0,b-e^0,0 \ldots \right] & = & 
x \left[ 0,c^1_1 \ldots c^1_{\ell/2}|0,c^1_{\ell/2+1} \ldots c^1_{\ell} \right] + \nonumber \\
& & +x'\left[ c^1_1 \ldots c^1_{\ell/2},0|c^1_{\ell/2+1} \ldots c^1_{\ell},0 \right]. \nonumber
\end{eqnarray}
Then we have $a=b=c^0=d^0=e^0$ and $A^1=B^1=C^1=0$ (in particular the representations $\widetilde{N}^{(2)}_{2,2,2}[\ell]$ are
exceptional). On the other hand, the equations on the right are
\begin{eqnarray}
\left[ 0,0 \ldots 0,0|a-f^0,0 \ldots 0,0 \right] & = & \left[ 0,d^1_1 \ldots d^1_{\ell/2}|0,d^1_{\ell/2+1} \ldots d^1_{\ell} \right]+ \nonumber\\
& & + \left[ f^1,0 \ldots 0,0|0,0 \ldots 0,0\right], \nonumber \\
\left[ 0,0 \ldots 0,a-g^0|0,0 \ldots 0,0\right] & = & \left[ e^1_1 \ldots e^1_{\ell/2},0|e^1_{\ell/2+1} \ldots e^1_{\ell},0 \right]+ \nonumber \\
& & + \left[ 0,0 \ldots 0,0|0,0 \ldots 0,g^1\right] , \nonumber
\end{eqnarray}
and it follows that $a=f^0=g^0$, $D^1=E^1=0$ and $F^1=G^1=0$. In this way $f^0=aId$ and $f^1=0$, that is, 
the endomorphism algebra of $\widetilde{N}^{(2)}_{3,3,2}[\ell]$ is one dimensional and $\widetilde{N}^{(2)}_{3,3,2}[\ell]$
is exceptional (as well as the modules $\widetilde{N}^{(2)}_{3,2,2}[\ell]$).

Consider now the case $\widetilde{M}^{(2)}_{3,3,2}[\ell]$ given in figure~(\ref{(DE)F:R2b}).
In this case we have $Q=\widetilde{M}=I^{\leftarrow}_{\ell-1/2}\oplus I^{\leftarrow}_{\ell+1/2}$ and 
$Q'=\widetilde{M}'=I^{\rightarrow}_{\ell-1/2}\oplus I^{\rightarrow}_{\ell+1/2}$. In addition we have
$R=\widetilde{M}_{\theta^r}=[1,0 \ldots |1,0 \ldots 0]$,
$S=\widetilde{M}_{\theta^s}=[0 \ldots 0|0 \ldots 0,1]$, $U=\widetilde{M}_{\theta^u}=[0\ldots 0| \ldots 0,1,0 \ldots]$,
$V=\widetilde{M}_{\theta^v}=[0 \ldots 0|1,0 \ldots 0]$ and $W=\widetilde{M}_{\theta^w}=[0 \ldots 0,1 |0 \ldots 0]$.
From commutativity of the inner square, which corresponds to the Kronecker subquiver, we have
\[
 A^0=\left[ \begin{matrix} aI & 0 \\ a'I^{\leftarrow}+b'I^{\rightarrow} & bI \end{matrix} \right], \qquad
B^0=\left[ \begin{matrix} aI & 0 \\ a'I^{\leftarrow}+b'I^{\rightarrow} & bI \end{matrix} \right],
\] 
for some $a,b,a',b'\in k$. From the first two equations in the middle of~(\ref{(DE)EQ:prueba1}) we obtain the equations
\begin{eqnarray} 
\left[ a+a'-c^0,b',0 \ldots 0|b-c^0,0 \ldots 0 \right] = 
\left[ 0,a^1_1 \ldots a^1_{\ell+1/2}|0,a^1_{\ell+1/2+1} \ldots a^1_{\ell} \right] \label{(DE)EQ:prueba2} \\
\left[ 0 \ldots 0,a',b'|0 \ldots 0,b-d^0\right] = 
\left[ b^1_1 \ldots b^1_{\ell+1/2},0|b^1_{\ell+1/2+1} \ldots b^1_{\ell},0 \right], \nonumber
\end{eqnarray}
from where it follows that $b'=0$, $b=d^0$ and $B^1=0$. Then the third equation is given by
\begin{eqnarray}
\left[ \ldots 0,a',0 \ldots| \ldots 0,b-e^0,0 \ldots \right] & = & 
x \left[ 0,c^1_1 \ldots c^1_{\ell+1/2}|0,c^1_{\ell+1/2+1} \ldots c^1_{\ell}\right] + \nonumber \\
& & +x' \left[ c^1_1 \ldots c^1_{\ell+1/2},0|c^1_{\ell+1/2+1} \ldots c^1_{\ell},0 \right], \nonumber
\end{eqnarray}
and hence $a'=0$, $b=e^0$ and $C^1=0$. Since $a'=0$, from equation~(\ref{(DE)EQ:prueba2})
it follows that $a=c^0=b$, and since $b'=0$ we have $A^1=0$ (in particular the representations $\widetilde{M}^{(2)}_{2,2,2}[\ell]$
are exceptional). Finally, from the pair of equations on the right of~(\ref{(DE)EQ:prueba1})
we obtain the equalities
\begin{eqnarray}
 \left[0,0 \ldots 0,0|a-f^0 \ldots 0 \right] & = & \left[ 0,d^1_1 \ldots d^1_{\ell+1/2}|0,d^1_{\ell+1/2+1} \ldots d^1_{\ell}\right]+ \nonumber \\
& & +\left[ f^1,0 \ldots 0,0|f^1 \ldots 0\right], \nonumber \\
\left[ 0,0 \ldots 0,a-g^0|0 \ldots 0\right] & = & \left[ e^1_1 \ldots e^1_{\ell+1/2},0|e^1_{\ell+1/2+1} \ldots e^1_{\ell},0\right]+ \nonumber \\
& & +\left[ 0,0 \ldots 0,0|0 \ldots g^1\right]. \nonumber
\end{eqnarray}
Therefore $f^1=0$ and $a=f^0=g^0$, $D^1=E^1=0$ and $F^1=G^1=0$, that is, the endomorphism algebra of 
$\widetilde{M}^{(2)}_{3,3,2}[\ell]$ is one dimensional and $\widetilde{M}^{(2)}_{3,3,2}[\ell]$ is exceptional
(as well as the modules $\widetilde{M}^{(2)}_{3,2,2}$). This completes the proof.
\eproof

%----------------------------------------------------------------------
%----------------------------------------------------------------------
%----------------------------------------------------------------------
\chapter{Construction of exceptional modules for $\widetilde{\Delta}$.}
\label{Cap(B)}

In this chapter we use the reduced exceptional representations given at the end of chapter~\ref{Cap(DE)} and the
reduction functors $F^X$ and $F^Y$ to construct series of exceptional $\widetilde{\Delta}$-modules, where
$\widetilde{\Delta}$ is an extended Dynkin quiver such that its extension vertex is a source.
This construction can be performed directly in coefficient quivers, for which we need the explicit knowledge of
some representation of the Dynkin diagram $\Delta$. The coefficient quiver corresponding to the maximal exceptional 
$\Delta$-module $W_0$ is indicated with the figure $\CoefMax$. The marked vertices constitude a basis of the vector space
$\Hom_{A_0}(R,W_0)$. All other needed $\Delta$-modules, whose isomorphism classes conform the sets $\mathcal{X}_0$ and $\mathcal{Y}_0$, 
are indicated with the figures $\CoefAlaU{i}$. The index $i$ determines the corresponding representation $X_i$ or $Y_i$,
while the orientation of the triangle symbolized the different wings of vertex $[W_0]$ in the Auslander-Reiten quiver of $A_0$.

%------------------------------
\begin{displaymath}
\xy
( 0,-15)="E1" *{\CoefAlaU{i}};
( 0,15)="A1" *{};
( 6,15)="A2" *{};
(12,15)="A3" *{};
(18,15)="A4" *{};
(24,15)="A5" *{};
( 2, 0)="B1" *{};
( 8, 0)="B2" *{};
(14, 0)="B3" *{};
(20, 0)="B4" *{};
(26, 0)="B5" *{};
( 4, 0)="C1" *{};
(10, 0)="C2" *{};
(16, 0)="C3" *{};
(22, 0)="C4" *{};
(28, 0)="C5" *{};
( 3, 0)="D1"  *{\CoefMax};
( 9, 0)="D2"  *{\CoefMax};
( 15, 0)="D3" *{\CoefMax};
( 21, 0)="D3" *{\CoefMax};
%-----------------------
"A1";"B1" **@{-};
"A2";"B2" **@{-};
"A3";"B3" **@{-};
"A4";"B4" **@{-};
"A2";"C1" **@{-};
"A3";"C2" **@{-};
"A4";"C3" **@{-};
"A5";"C4" **@{-};
"A1";"E1" **@{-};
\endxy
%%%%%%%%%%%%%%%%%%%%%%%%%%%
\xy
( 0,-12)="E1"  *{\CoefAlaU{i}};
( 24,-18)="E5" *{\CoefAlaV{j}};
( 0,15)="A1" *{};
( 6,15)="A2" *{};
(12,15)="A3" *{};
(18,15)="A4" *{};
(24,15)="A5" *{};
( 2, 0)="B1" *{};
( 8, 0)="B2" *{};
(14, 0)="B3" *{};
(20, 0)="B4" *{};
(26, 0)="B5" *{};
( 4, 0)="C1" *{};
(10, 0)="C2" *{};
(16, 0)="C3" *{};
(22, 0)="C4" *{};
(28, 0)="C5" *{};
( 3, 0)="D1"  *{\CoefMax};
( 9, 0)="D2"  *{\CoefMax};
( 15, 0)="D3" *{\CoefMax};
( 21, 0)="D3" *{\CoefMax};
%-----------------------
"A1";"B1" **@{-};
"A2";"B2" **@{-};
"A3";"B3" **@{-};
"A4";"B4" **@{-};
"A2";"C1" **@{-};
"A3";"C2" **@{-};
"A4";"C3" **@{-};
"A5";"C4" **@{-};
"A1";"E1" **@{-};
"A5";"E5" **@{-};
\endxy
%%%%%%%%%%%%%%%%%
\xy
( 0,-9)="E1"   *{\CoefAlaU{i}};
( 12,-15)="E3" *{\CoefAlaW{k}};
( 24,-21)="E5" *{\CoefAlaV{j}};
( 0,15)="A1" *{};
( 6,15)="A2" *{};
(12,15)="A3" *{};
(18,15)="A4" *{};
(24,15)="A5" *{};
( 2, 0)="B1" *{};
( 8, 0)="B2" *{};
(14, 0)="B3" *{};
(20, 0)="B4" *{};
(26, 0)="B5" *{};
( 4, 0)="C1" *{};
(10, 0)="C2" *{};
(16, 0)="C3" *{};
(22, 0)="C4" *{};
(28, 0)="C5" *{};
( 3, 0)="D1"  *{\CoefMax};
( 9, 0)="D2"  *{\CoefMax};
( 15, 0)="D3" *{\CoefMax};
( 21, 0)="D3" *{\CoefMax};
%-----------------------
"A1";"B1" **@{-};
"A2";"B2" **@{-};
"A3";"B3" **@{-};
"A4";"B4" **@{-};
"A2";"C1" **@{-};
"A3";"C2" **@{-};
"A4";"C3" **@{-};
"A5";"C4" **@{-};
"A1";"E1" **@{-};
"A3";"E3" **@{-};
"A5";"E5" **@{-};
\endxy
\end{displaymath}

%----------------------------------------------------------------------
%----------------------------------------------------------------------
\section{Representations of $\widetilde{\mathbf{A_n}}$.} \label{(B)S:1}
%--------------------------------------------------------------------------------------
%--------------------------------------------------------------------------------------

Consider first the extended Dynkin quivers $\widetilde{\mathbf{A_n}}$ ($n \geq 2$) with linear orientation of its arrows.
In this case the Auslander-Reiten quiver of the Dynkin diagram $\Gamma(k\mathbf{A_n})$ consists in a single wing of order $n$.
%------------------
\begin{displaymath}
 \xymatrix@!0@C=24pt@R=24pt{
& & & & {\vctACi{1}{1}{\cdots}{1}{1}} \ar[rd] \\
& & & {\vctACi{0}{1}{\cdots}{1}{1}} \ar[ru] \ar[rd] & & {\vctACi{1}{1}{\cdots}{1}{0}} \ar[rd] \\
& & {} \ar[ru] \ar@{}[rd] & & {} \ar[ru] \ar@{}[rd] & & {} \ar@{}[rd]|-{\cdots}  \\
& {\vctACi{0}{0}{\cdots}{1}{1}} \ar@{}[ru]|-{\cdots} \ar[rd] & & {} \ar@{}[ru] 
\ar@{}[rd] & & {} \ar@{}[ru] \ar[rd] & & {\vctACi{1}{1}{\cdots}{0}{0}} \ar[rd] \\
{\vctACi{0}{0}{\cdots}{0}{1}} \ar[ru] & & {} \ar@{}[rrrr]|-{\cdots} \ar[ru] & & {} 
\ar@{}[ru] & & {} \ar[ru] & & {\vctACi{1}{0}{\cdots}{0}{0}} \\
}
\end{displaymath}
According to lemmas~\ref{(DE)L:AX} and~\ref{(DE)L:AY} we compute the reduced ditalgebras $\mathcal{A}^X$ and $\mathcal{A}^Y$,
where $X_0$ is direct sum of representatives of the projective $A_0$-modules and $Y_0$ is direct sum of representatives of
the injective $A_0$-modules (ordered from left to right in both cases). Fix bases of the vector spaces $\Hom_{A_0}(R,X_0)$ and 
$\Hom_{A_0}(R,Y_0)$ as in lemmas~\ref{(DE)L:AXcero} and~\ref{(DE)L:AYcero}. Choose as basis
of $\Hom_{A_0}(R,W_0)$ the vectors $a=\rho_1$ and $a'$ arbitrary, and $b_1^1=b=a'$, $b'=\rho_1$.
The quivers corresponding to the reductions $\mathcal{A}^X$ and $\mathcal{A}^Y$ are shown next.
%------------------------------------------------ 
\begin{displaymath}
 \xymatrix@C=1.3pc{
& & & & *+++[]{\bullet_{\omega}} \ar[ddd]^-{\xi'} \ar@<-1ex>[ddd]_-{\xi} \\
\\
\\
\bullet_1 \ar@{.>}[r]^(.7){\gamma^1_{1,2}} \ar@{.>}@/^1.5pc/[rrr]^(.35){\gamma^1_{1,n-1}} \ar@{.>}@/_2pc/[rrrr]_(.5){\gamma^1_{1,n}} 
\ar@{<-}@/^2pc/[rrrruuu]^(.27){\xi^1_1} 
 & \bullet_2 \ar@{.>}@/_.5pc/[rr]_(.5){\gamma^1_{2,n-1}} \ar@{.>}@/^1.5pc/[rrr]^(.5){\gamma^1_{2,n}} 
\ar@{<-}@/^1pc/[rrruuu]^(.35){\xi^1_2} 
& \ldots & \bullet_{n-1} \ar@{.>}[r]^(.4){\gamma^1_{n-1,n}} \ar@{<-}@/^.5pc/[ruuu]^(.42){\xi^1_{n-1}}
& *++[]{\bullet_{n}} \\
}
\xymatrix@C=1.3pc{
*+++[]{\bullet_{\omega}} 
\ar@/^2pc/[rrrrddd]^(.74){\theta^1_{n}}
\ar@/^1pc/[rrrddd]^(.67){\theta^1_{n-1}}
\ar@/^.5pc/[rddd]^(.59){\theta^1_{2}}
\ar[ddd]^-{\theta} \ar@<-1ex>[ddd]_-{\theta'}  \\
\\
\\
\bullet_1 \ar@{.>}[r]^(.65){\varepsilon^1_{1,2}} \ar@{.>}@/^1.5pc/[rrr]^(.5){\varepsilon^1_{1,n-1}} 
\ar@{.>}@/_2pc/[rrrr]_(.5){\varepsilon^1_{1,n}}  
& \bullet_2 \ar@{.>}@/_.5pc/[rr]_(.5){\varepsilon^1_{2,n-1}} \ar@{.>}@/^1.5pc/[rrr]^(.65){\varepsilon^1_{2,n}}  
& \ldots & \bullet_{n-1} \ar@{.>}[r]^(.3){\varepsilon^1_{n-1,n}} 
& *++[]{\bullet_{n}} \\
}
\end{displaymath}

With respect to the selected bases the differentials $\delta^x$ and $\delta^y$ have the following form.

%---------------------------------------------
\begin{center}
 \begin{tabular}{c | l  c |  l}
Arrow & Its differential $\delta^x$ & Arrow & Its differential $\delta^y$\\
\hline
$\xi^1_j $ & $\sum_{1 \leq i<j}\gamma^1_{i,j} \otimes \xi^1_i$, & $\theta^1_j $ & $\varepsilon^1_{1,j} \otimes 
	      \theta+\sum_{1<i<j}\varepsilon^1_{i,j} \otimes \theta^1_i$, \\
$\xi$ & $\sum_{1\leq i<n_1}\gamma^1_{i,n_1} \otimes \xi^1_i$, & $\theta $ & $0$, \\
$\xi' $ & $0$, & $\theta' $ & $0$. \\
 \end{tabular}
\end{center}

The coefficient quiver of an exceptional representation corresponding to the maximal root of $\mathbf{A_n}$, 
indicated with the figure $\CoefMax$, is substituted
in the coefficient quivers of the $\mathcal{A}^Y$-modules with rank minus one $\widetilde{M}^{(-1)}_{1,1,1}[\ell]$ 
and rank one $\widetilde{M}^{(1)}_{1,1,1}[\ell]$, given in figures~(\ref{(DE)F:Rm}) and~(\ref{(DE)F:R1}) in section~\ref{S:AYEjem}. 
In this way we directly determine coefficient quivers for $F^Y(\widetilde{M}^{(-1)}_{1,1,1}[\ell])$
and $F^Y(\widetilde{M}^{(1)}_{1,1,1}[\ell])$. We show for instance the cases $\ell=4$.
\vspace{1cm}

Rank minus one representations $M^{(-1)}_{1,1,1}[\ell]=F^Y(\widetilde{M}^{(1)}_{1,1,1}[\ell])$ ($\ell>0$)
and rank one representations $M^{(1)}_{1,1,1}[\ell]=F^Y(\widetilde{M}^{(1)}_{1,1,1}[\ell])$ ($\ell \geq 0$).
\begin{equation} \label{(C)F:AnO1}
\xy
(-15,10)="Et1" *{\scriptstyle \mathcal{C}(F^Y(\widetilde{M}^{(-1)}_{1,1,1}[4]))};
(40,10)="Et2" *{\scriptstyle \mathcal{C}(F^Y(\widetilde{M}^{(1)}_{1,1,1}[4]))};
(10,10)="Et3" *{\scriptstyle \mathcal{C}(W_0)};
(-19,3)="Col" *{\xymatrix@!0@R=18pt{{}_{k^{\ell-1}} \ar[d]^-{I^{\uparrow}} \ar@/_13pt/[ddddd]_-{I^{\downarrow}} 
\\{}_{k^{\ell}} \ar[d]^-{=} \\ {}_{k^{\ell}} \ar@{}[d]|{\cdots}
\\ {}_{k^{\ell}} \ar[d]^-{=} \\ {}_{k^{\ell}} \ar[d]^-{=} \\ {}_{k^{\ell}}} };
( -13,3.5)="Col1" *{
\xymatrix@!0@R=18pt@C=18pt{ {} & {} \ar@{-}[lddddd] \ar@{-}[d] 
& {} \ar@{-}[lddddd] \ar@{-}[d] & {} \ar@{-}[lddddd] \ar@{-}[d]
\\ {}_{\bullet} \ar@{-}[d] & {}_{\bullet} \ar@{-}[d] & {}_{\bullet} \ar@{-}[d] & {}_{\bullet} \ar@{-}[d] 
\\ {} \ar@{}[d]|{\cdots} & {} \ar@{}[d]|{\cdots} & {} \ar@{}[d]|{\cdots} & {} \ar@{}[d]|{\cdots}
\\ {} \ar@{-}[d] & {} \ar@{-}[d] & {} \ar@{-}[d] & {} \ar@{-}[d]
\\ {} \ar@{-}[d] & {} \ar@{-}[d] & {} \ar@{-}[d] & {} \ar@{-}[d]
\\ {}_{\bullet} & {}_{\bullet} & {}_{\bullet} & {}_{\bullet}  }};
( 10,15)="A1" *{\CoefMax};
( 31,3)="Col2" *{\xymatrix@!0@R=18pt{{}_{k^{\ell+1}} \ar[d]_-{I^{\rightarrow}} \ar@/^13pt/[ddddd]^-{I^{\leftarrow}} 
\\{}_{k^{\ell}} \ar[d]_-{=} \\ {}_{k^{\ell}} \ar@{}[d]|{\cdots}
\\ {}_{k^{\ell}} \ar[d]_-{=} \\ {}_{k^{\ell}} \ar[d]_-{=} \\ {}_{k^{\ell}}} };
( 5,0)="Col" *{\xymatrix@!0@R=18pt{{}_{\bullet} \ar@{-}[d] \\ {} \ar@{}[d]|{\cdots}\\ {} \ar@{-}[d]\\ {} \ar@{-}[d]\\ {}_{\bullet}} };
( 13,3.5)="Col" *{
\xymatrix@!0@R=18pt@C=18pt{ {} \ar@{-}[d] & {} \ar@{-}[lddddd] \ar@{-}[d] 
& {} \ar@{-}[lddddd] \ar@{-}[d] & {} \ar@{-}[lddddd] \ar@{-}[d] & {} \ar@{-}[lddddd]
\\ {}_{\bullet} \ar@{-}[d] & {}_{\bullet} \ar@{-}[d] & {}_{\bullet} \ar@{-}[d] & {}_{\bullet} \ar@{-}[d] 
\\ {} \ar@{}[d]|{\cdots} & {} \ar@{}[d]|{\cdots} & {} \ar@{}[d]|{\cdots} & {} \ar@{}[d]|{\cdots}
\\ {} \ar@{-}[d] & {} \ar@{-}[d] & {} \ar@{-}[d] & {} \ar@{-}[d]
\\ {} \ar@{-}[d] & {} \ar@{-}[d] & {} \ar@{-}[d] & {} \ar@{-}[d]
\\ {}_{\bullet} & {}_{\bullet} & {}_{\bullet} & {}_{\bullet}  }};
%--------------
\endxy 
\vspace{2.5cm}
\end{equation}
Since in the linearly ordered case in $\mathbf{A_n}$ there is only one wing of vertex $[W_0]$ in the
Auslander-Reiten quiver, we observe that, besides the given representations $\widetilde{M}^{(1)}_{1,1,1}[\ell]$, there are 
only reduced representations of rank one of the form $\widetilde{M}^{(1)}_{2,1,1}[\ell]$ (point~(\ref{(DE)F:R1a})). 
In this case the reduced representation ${}_i\widetilde{M}^{(1)}_{2,1,1}[\ell]$ depends on an index $i$
corresponding to the representation $Y^1_i$ ($1<i\leq n$), for which we also give a coefficient quiver.

Rank one representations ${}_iM^{(1)}_{2,1,1}[\ell]=F^Y({}_i\widetilde{M}^{(1)}_{2,1,1}[\ell])$ ($\ell\geq 0$).
\begin{equation} \label{(C)F:AnO2}
\xy
(45,10)="Et" *{\scriptstyle \mathcal{C}({}_iM^{(1)}_{2,1,1}[4])};
( -15,3)="Col" *{\xymatrix@!0@R=18pt{{}_{k^{\ell+1}} \ar[d]^-{I^{\rightarrow}} \ar@/_17pt/[ddddd]_-{=} 
\\{}_{k^{\ell}} \ar@{}[d]|{\cdots} \\ {}_{k^{\ell}} \ar[d]^-{I^{\uparrow}} \\ {}_{k^{\ell+1}} \ar@{}[d]|{\cdots} 
\\ {}_{k^{\ell+1}} \ar[d]^-{=}\\ {}_{k^{\ell+1}}} };
(-10,15)="B1" *{\CoefAlaU{i}};
(-10,10)="B1Et" *{\scriptstyle \mathcal{C}(Y^1_i)};
(-5,-3.3)="Col1" *{\xymatrix@!0@R=18pt{{} \\ {}_{i} \ar@{}[d]|{\cdots}\\ {}_{{}_2} \ar@{-}[d]\\ {}_{\bullet_1}} };
( 10,15)="A1" *{\CoefMax};
( 10,10)="A1" *{\scriptstyle \mathcal{C}(W_0)};
( 5,0)="Col" *{\xymatrix@!0@R=18pt{{}_{\bullet} \ar@{-}[d] \\ {} \ar@{}[d]|{\cdots}\\ {} \ar@{-}[d]\\ {} \ar@{-}[d]\\ {}_{\bullet}} };
( 15,3.5)="Col" *{
\xymatrix@!0@R=18pt@C=18pt{{} & {} \ar@{-}[lddddd] \ar@{-}[d] & {} \ar@{-}[lddddd] \ar@{-}[d] 
& {} \ar@{-}[lddddd] \ar@{-}[d] & {} \ar@{-}[lddddd] \ar@{-}[d] & {} \ar@{-}[lddddd]
\\{} & {}_{\bullet} \ar@{-}[d] & {}_{\bullet} \ar@{-}[d] & {}_{\bullet} \ar@{-}[d] & {}_{\bullet} \ar@{-}[d] 
\\ {} & {} \ar@{}[d]|{\cdots} & {} \ar@{}[d]|{\cdots} & {} \ar@{}[d]|{\cdots} & {} \ar@{}[d]|{\cdots}
\\ {}_i \ar@{}[d]|{\cdots} & {} \ar@{-}[d] & {} \ar@{-}[d] & {} \ar@{-}[d] & {} \ar@{-}[d]
\\ {} \ar@{-}[d] & {} \ar@{-}[d] & {} \ar@{-}[d] & {} \ar@{-}[d] & {} \ar@{-}[d]
\\ {}_{\bullet} & {}_{\bullet} & {}_{\bullet} & {}_{\bullet} & {}_{\bullet}  }};
%--------------
%"A1";"A2" **@{-};
\endxy 
\vspace{2.5cm}
\end{equation}

\begin{proposicion} \label{(C)L:An}
Assume that $\widetilde{\mathbf{A_n}}$ has a linear orientation of its arrows.
Then every indecomposable preinjective $\widetilde{\mathbf{A_n}}$-representation is isomorphic to one of the
modules $M^{(1)}_{1,1,1}[\ell]$ or ${}_iM^{(1)}_{2,1,1}[\ell]$ (with $1<i\leq n$)
given in figures~(\ref{(C)F:AnO1}) and~(\ref{(C)F:AnO2}).
\end{proposicion}
\bproof
By theorem~\ref{(DE)T:finito}($b$) and the construction of the reduction module $Y$, almost every indecomposable
preinjective module is isomorphic to an object in the image of the functor $F^Y$. Since $Y$ is direct sum of
representatives of the indecomposable injective $A_0$-modules, in fact every indecomposable preinjective module
is isomorphic to an object in the image of $F^Y$. 

By proposition~\ref{(DE)P:rango6}($a$), if $\widetilde{M}$ is an indecomposable $\mathcal{A}^Y$-module such that
$F^Y(\widetilde{M})$ is preinjective, then the dimension vector $\vdim \widetilde{M}$ has rank one.
By lemma~\ref{(DE)L:tablaComp} and since the quiver $\Delta^Y$ associated to the ditalgebra $\mathcal{A}^Y$ has
a single branch, $\vdim \widetilde{M}$ is one of the roots $\mathcal{R}^{(1)}_{1,1,1}[\ell]$ or 
$\mathcal{R}^{(1)}_{2,1,1}[\ell]$.
Hence every indecomposable preinjective module is isomorphic to a representation of the form $F^Y(\widetilde{M})$
where $\widetilde{M}$ is one of the representations $\widetilde{M}^{(1)}_{1,1,1}[\ell]$
and $\widetilde{M}^{(1)}_{2,1,1}[\ell]_i$ for $1<i\leq n$ given in lemma~\ref{(DE)L:excRedRnkUno}.
\eproof

We turn now to the case $\mathbf{A_n}$ with a nonlinear ordering of its arrows, thus there exist two wings of vertex $[W_0]$
in the Auslander-Reiten quiver $\Gamma(k\mathbf{A_n})$.
\begin{displaymath}
 \xymatrix@!0@C=24pt@R=24pt{
{\scriptstyle \cdots} \ar[rd] & & {{}_{x^2_1}} \ar[rd] & & {} \ar@{}[rd] \ar@{}[rrd]|-{\cdots} 
& & {} \ar[rd] & & {{}_{y^2_{n_2}}} \ar[rd] & & {\scriptstyle \cdots} \\
& {} \ar@{}[rd] \ar[ru] & & {{}_{x^2_2}} \ar@{}[rd] \ar[ru] & & {} \ar@{}[rd] \ar@{}[ru] 
& {} & {{}_{y^2_{n_2-1}}} \ar@{}[rd] \ar[ru] & & {} \ar@{}[rd] \ar[ru] \\
{} \ar[rd] \ar@{}[ru] & & {} \ar[rd] \ar@{}[ru] & & {{}_{x^2_{n_2-1}}} \ar[rd] \ar@{}[ru] 
& & {{}_{y^2_2}} \ar[rd] \ar@{}[ru] & & {} \ar[rd] \ar@{}[ru] & & {} \\
{\scriptstyle \cdots} & {} \ar[rd] \ar[ru] & & {} \ar[rd] \ar[ru] 
& & {{}_{w_0}} \ar[rd] \ar[ru] & & {} \ar[rd] \ar[ru] & & {} \ar[rd] \ar[ru] & {\scriptstyle \cdots} \\
{} \ar@{}[rd] \ar[ru] & & {} \ar@{}[rd] \ar[ru] & & {{}_{x^1_{n_1-1}}} \ar@{}[rd] \ar[ru] 
& & {{}_{y^1_2}} \ar@{}[rd] \ar[ru] & & {} \ar@{}[rd] \ar[ru] & & {} \\
& {} \ar[rd] \ar@{}[ru] & & {{}_{x^1_{3}}} \ar[rd] \ar@{}[ru] & & {} \ar@{}[rd] \ar@{}[ru] 
& & {{}_{y^1_{n_1-2}}} \ar[rd] \ar@{}[ru] & & {} \ar[rd] \ar@{}[ru] \\
{\scriptstyle \cdots} \ar[rd] \ar[ru] & & {{}_{x^1_2}} \ar[rd] \ar[ru] & & {} \ar@{}[rd] \ar@{}[ru] \ar@{}[rr]|-{\cdots} 
& & {} \ar[rd] \ar[ru] & & {{}_{y^1_{n_1-1}}} \ar[rd] \ar[ru] & & {\scriptstyle \cdots} \\
& {{}_{x^1_1}} \ar[ru] & & {} \ar[ru] & & {} \ar@{}[ru] & & {} \ar[ru] & & {{}_{y^1_{n_1}}} \ar[ru] \\
}
\end{displaymath}
Observe that there are exactly $2n-2$ indecomposable representations $Z_0$ of $\mathbf{A_n}$ such that 
$\dimk_k \Hom_{A_0}(R,Z_0)=1$. By lemma~\ref{(DE)L:valores} this are the elements which conform the sets $\mathcal{X}_0$ 
and $\mathcal{Y}_0$. There are rearrangements of these roots as shown below, where as usual $w_0$ is the maximal root
of $\mathbf{A_n}$,
\[
 \{ x^1_1,\ldots,x^1_{n_1-1},w_0,y^2_2,\ldots,y^2_{n_2} \}=\left\{\vctACi{0}{0}{\cdots}{0}{1},\vctACi{0}{0}{\cdots}{1}{1},\ldots,
\vctACi{0}{1}{\cdots}{1}{1},\vctACi{1}{1}{\cdots}{1}{1} \right\},
\]
\[
 \{ x^2_1,\ldots,x^2_{n_2-1},w_0,y^1_2,\ldots,y^1_{n_1} \}=\left\{ \vctACi{1}{0}{\cdots}{0}{0},\vctACi{1}{1}{\cdots}{0}{0},\ldots,
\vctACi{1}{1}{\cdots}{1}{0},\vctACi{1}{1}{\cdots}{1}{1} \right\}.
\]
Moreover, since there are nonzero morphisms $W_0 \to Y^1_i$ and $W_0 \to Y^2_j$ for $1<i\leq n_1$ and $1<j\leq n_2$,
we notice that the elements of $\mathcal{Y}_0$ correspond to those roots such that the arrow joining the last entry one
with the first entry zero has the orientation $1 \to 0$. In a similar way the elements of $\mathcal{X}_0$ correspond
to those arrows where the arrow joining the series of ones with the series of zeros has the orientation $0 \to 1$.
For the choice of the bases of $\Hom_{A_0}(R,X_0)$ and $\Hom_{A_0}(R,Y_0)$ take $a=\rho_1$, $a'=\rho_2$, 
$b_1^1=b=a'$ and $b_1^2=b'=a$ (cf. lemmas~\ref{(DE)L:AXcero} and~\ref{(DE)L:AYcero}). In this way the differentials
$\delta^x$ and $\delta^y$ have the forms given in the following table.
%---------------------------------------------
\begin{center}
 \begin{tabular}{c | l  c |  l}
Arrow & Its differential $\delta^x$ & Arrow & Its differential $\delta^y$\\
\hline
$\xi^1_j $ & $\sum_{1 \leq i<j}\gamma^1_{i,j} \otimes \xi^1_i$, & $\theta^1_j $ & $\varepsilon^1_{1,j} \otimes 
	      \theta+\sum_{1<i<j}\varepsilon^1_{i,j} \otimes \theta^1_i$, \\
$\xi^2_j $ & $\sum_{1 \leq i<j}\gamma^2_{i,j} \otimes \xi^2_i$, & $\theta^2_j $ & $\varepsilon^2_{1,j} \otimes 
	      \theta'+\sum_{1<i<j}\varepsilon^2_{i,j} \otimes \theta^2_i$, \\
$\xi$ & $\sum_{1\leq i<n_1}\gamma^1_{i,n_1} \otimes \xi^1_i$, & $\theta $ & $0$, \\
$\xi' $ & $\sum_{1\leq i<n_2}\gamma^2_{i,n_2} \otimes \xi^2_i$, & $\theta' $ & $0$. \\
 \end{tabular}
\end{center}
To supplement the list of rank one modules that appear in cases $\mathbf{A_n}$ consider reduced representations 
of the form ${}_i\widetilde{M}^{(1)}_{2,2,1}[\ell]_j$ (point~(\ref{(DE)F:R1b})). We describe the prepresentations 
$F^Y({}_i\widetilde{M}^{(1)}_{2,2,1}[\ell]_j))$. The quivers indicated with $\CoefAlaU{i}$ and $\CoefAlaV{j}$ correspond
to coefficient quivers of $Y^1_i$ and $Y^2_j$ respectively. We show for instance the cases $\ell=4$.
\vspace{1cm}

Rank one representations ${}_iM^{(1)}_{2,2,1}[\ell]_j=F^Y({}_i\widetilde{M}^{(1)}_{2,2,1}[\ell]_j)$ ($\ell\geq 0$).
\begin{displaymath}
\xy
( -17,3)="Col" *{\xymatrix@!0@R=18pt{{}_{\omega} \ar[d] \ar@/_13pt/[ddddd] \\{}_{n} \ar@{-}[d] \\ {}_{n-1} \ar@{}[d]|{\cdots}
\\ {}_{3} \ar@{-}[d]\\ {}_{2} \ar@{-}[d]\\ {}_{1}} };
(45,10)="Et" *{\scriptstyle \mathcal{C}({}_iM^{(1)}_{2,2,1}[4]_j)};
(-20,15)="B1" *{\CoefAlaU{i}};
(-20,10)="B1Et" *{\scriptstyle \mathcal{C}(Y^1_i)};
(-10,-3.3)="Col1" *{\xymatrix@!0@R=18pt{{}_i \ar@{-}[d]\\ {} \ar@{}[d]|{\cdots}\\ {} \ar@{-}[d]\\ {}_{\bullet}} };
( 10,15)="A1" *{\CoefMax};
( 10,10)="A1" *{\scriptstyle \mathcal{C}(W_0)};
( 5,0)="Col" *{\xymatrix@!0@R=18pt{{}_{\bullet} \ar@{-}[d] \\ {} \ar@{}[d]|{\cdots}\\ {} \ar@{-}[d]\\ {} \ar@{-}[d]\\ {}_{\bullet}} };
(-5,15)="B2" *{\CoefAlaV{j}};
(-5,10)="B2" *{\scriptstyle \mathcal{C}(Y^2_j)};
(-2.5,0)="Col2" *{\xymatrix@!0@R=18pt{{}_{\bullet} \ar@{-}[d]\\ {} \ar@{}[d]|{\cdots}\\ {} \ar@{-}[d]\\ {}_{j}} };
( 15,3.5)="Col" *{
\xymatrix@!0@R=18pt@C=18pt{{} & {} \ar@{-}[lddddd] \ar@{-}[d] & {} \ar@{-}[lddddd] \ar@{-}[d] 
& {} \ar@{-}[lddddd] \ar@{-}[d] & {} \ar@{-}[lddddd] \ar@{-}[d] & {} \ar@{-}[lddddd] \ar@{-}[d]
\\{} & {}_{\bullet} \ar@{-}[d] & {}_{\bullet} \ar@{-}[d] & {}_{\bullet} \ar@{-}[d] & {}_{\bullet} \ar@{-}[d] & {}_{\bullet} \ar@{-}[d]
\\ {}_i \ar@{-}[d] & {} \ar@{}[d]|{\cdots} & {} \ar@{}[d]|{\cdots} & {} \ar@{}[d]|{\cdots} & {} \ar@{}[d]|{\cdots} & {} \ar@{}[d]|{\cdots}
\\ {} \ar@{}[d]|{\cdots} & {} \ar@{-}[d] & {} \ar@{-}[d] & {} \ar@{-}[d] & {} \ar@{-}[d] & {} \ar@{-}[d]
\\ {} \ar@{-}[d] & {} \ar@{-}[d] & {} \ar@{-}[d] & {} \ar@{-}[d] & {} \ar@{-}[d] & {}_{j}
\\ {}_{\bullet} & {}_{\bullet} & {}_{\bullet} & {}_{\bullet} & {}_{\bullet}  }};
%--------------
%"A1";"A2" **@{-};
\endxy 
\vspace{2.5cm}
\end{displaymath}
These $\widetilde{\mathbf{A_n}}$- modules correspond to those modules $N_{\pi}$ determined
by (counterclockwise) walks $\pi$ given by Gabriel and Roiter \cite[section~11.3]{GR97}.

\section{Representations of $\widetilde{\mathbf{D_n}}$.} \label{(B)S:2}
%--------------------------------------------------------------------------------------
%--------------------------------------------------------------------------------------
%--------------------------------------------------------------------------------------

For the extended Dynkin quivers of type $\widetilde{\mathbf{D_n}}$ we fix the following orientation of arrows,
\begin{displaymath}
 \xymatrix@R=1pc{
{}_{\bullet_{2}} \ar@{<-}[rd] &  &  & & {}_{\bullet_{\omega}} \\
        & {}_{\bullet_{3}} \ar@{<-}[r] & {}_{\bullet_{4}} \cdots {}_{{}_{n-2}\bullet} \ar@{<-}[r] 
& {}_{{}_{n-1}\bullet} \ar@{<-}[ru] \ar@{<-}[rd] \\ {}_{\bullet_1} \ar@{<-}[ru] &  &  &  & {}_{\bullet_n} \\
}
\end{displaymath}
%------------------------------------------------
%------------------------------------------------
\begin{sidewaysfigure}
\begin{center}
$ \xymatrix@!0@C=28pt@R=28pt{
& & & & & & {\vctDn{1}{1}{1}{1}{1}{1}{1}} \ar[rd] & & {} \ar[rd] & & {} \ar[rd] & & {} \ar@{}[rd] \ar@{}[rr]|-{\cdots}
& & {} \ar[rd] & & {\vctDn{0}{0}{0}{0}{0}{1}{0}} \ar[rd] & & {} \\
& & & & & {} \ar[ru] \ar[rd] & & {\vctDn{1}{1}{2}{1}{1}{1}{1}} \ar[ru] \ar[rd] & & {} \ar[ru] \ar[rd] & & {} \ar[ru] \ar@{}[rd] 
& & {} \ar[ru] \ar[rd] & & {\vctDn{0}{0}{0}{0}{1}{1}{0}} \ar[ru] \ar[rd] & & {} \ar[ru] \\
& & & & {} \ar[ru] \ar@{}[rd] & & {} \ar[ru] \ar@{}[rd] & & {\vctDn{1}{1}{2}{2}{1}{1}{1}} \ar[ru] \ar@{}[rd] & & {} \ar@{}[rr]|-{\cdots} \ar[ru] \ar@{}[rd] 
& & {} \ar[ru] \ar@{}[rd] & & {\vctDn{0}{0}{0}{0}{1}{1}{0}} \ar[ru] \ar@{}[rd] & & {} \ar[ru] \\
& & & {} \ar@{}[ru]|-{\cdots} \ar[rd] & & {} \ar@{}[ru] \ar[rd] & & {} \ar@{}[rr]|-{\cdots} \ar@{}[ru]|-{\cdots} \ar@{}[rd] 
& & {\vctDn{1}{1}{2}{2}{1}{1}{1}} \ar@{}[ru]|-{\cdots} \ar[rd] 
& & {} \ar@{}[ru] \ar[rd] & & {\vctDn{0}{0}{0}{1}{1}{1}{0}} \ar@{}[ru] \ar[rd] & & {} \ar@{}[ru]|-{\cdots} \\
& & {} \ar[ru] \ar[rd] & & {} \ar[ru] \ar[rd] & & {} \ar[ru] \ar@{}[rd] & & {} \ar[ru] \ar[rd] 
& & {\vctDn{1}{1}{2}{2}{2}{1}{1}} \ar[ru] \ar[rd] & & {\vctDn{0}{0}{1}{1}{1}{1}{0}} \ar[ru] \ar[rd] & & {} \ar[ru] \\
{} \ar[r] & {} \ar[ru] \ar[rd] \ar[r] & {} \ar[r] & {} \ar[ru] \ar[rd] \ar[r] & {} \ar[r] 
& {} \ar[ru] \ar@{}[rd] & & {} \ar[ru] \ar[rd] \ar[r] & {} \ar[r] & {} \ar[ru] \ar[rd] \ar[r] & {\vctDn{1}{0}{1}{1}{1}{1}{0}} \ar[r] 
& {\vctDn{1}{1}{2}{2}{2}{2}{1}} \ar[ru] \ar[rd] \ar[r] & {\vctDn{0}{1}{1}{1}{1}{1}{1}} \ar[r] & {} \ar[ru] \\
{} \ar[ru] & & {} \ar[ru] & & {} \ar@{}[rr]|-{\cdots} \ar[ru] & & {} \ar[ru] 
& & {} \ar[ru] & & {\vctDn{0}{1}{1}{1}{1}{1}{0}} \ar[ru] & & {\vctDn{1}{0}{1}{1}{1}{1}{1}} \ar[ru]
}$
\caption{Auslande-Reiten quiver of the Dynkin diagram $\mathbf{D_n}$ corresponding to the given orientation of arrows.
The positive roots shown correspond to the left side $\mathcal{X}_0$ and right side $\mathcal{Y}_0$ 
of the wings of vertex $[W_0]$.} \label{(DE)F:Dn}
\end{center}
\end{sidewaysfigure}
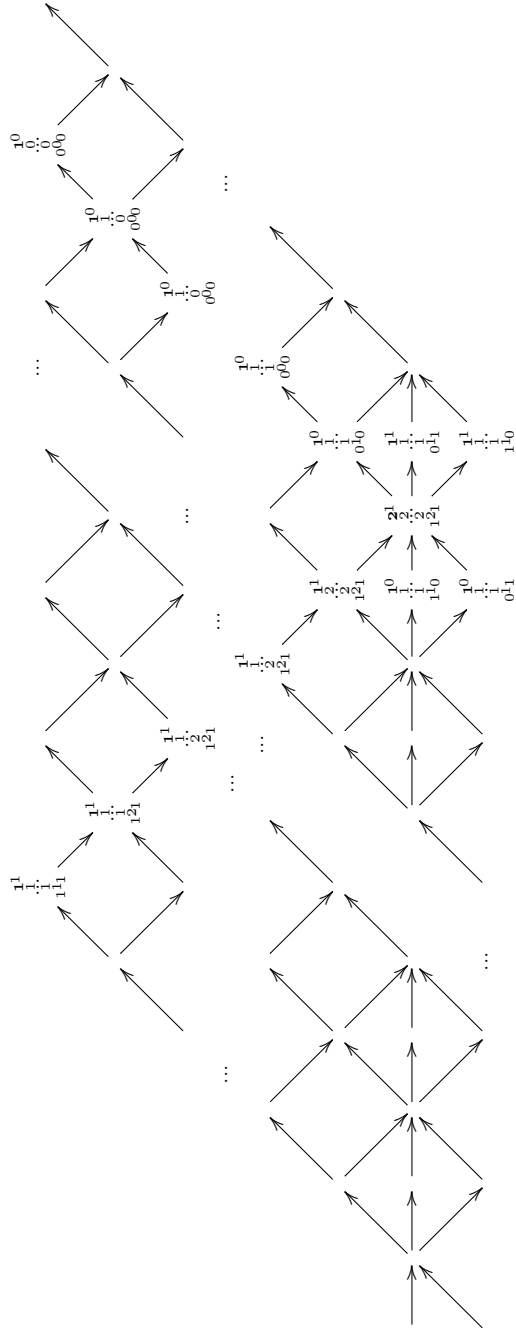
In figure~\ref{(DE)F:Dn} we show the Auslander-Reiten quiver $\Gamma(k\mathbf{D_n})$ of the Dynkin quiver corresponding to
$\mathbf{D_n}$. We order the wings of vertex $[W_0]$ starting with the wing that contains the sincere projective module 
(the wing of bigger order in cases $n>4$).  We fix indecomposable representations for the positive roots that appear
in figure~\ref{(DE)F:Dn} and show coefficient quivers corresponding to the maximal module and those representations in 
$\mathcal{X}_0$.
\begin{displaymath}
\xy
( -19,-32)="D1" *{};
( -16,-29.5)="D2" *{};
( -11,-5)="D3" *{};
( -14,-7.5)="D4" *{};
(-10,-1)="Col" *{\xymatrix@!0@C=15pt@R=9.5pt{{} \\ {} & & {}_{n}  
\\ & {}_{n-1} \ar@{-}[dd] \\ \\ & {}_{n-2} \\ & {\cdots} \\ & {}_{4} \ar@{-}[dd] 
\\ \\ & {}_{3} \ar@{-}[dd] \\ {}_{2} \\ & {}_{1} }};
( 0, 5)="UF" *{\CoefAlaU{1}}; 
( 0, 0)="UFEt" *{\scriptstyle \mathcal{C}(X^1_1)};
( 0, -20)="U"  *{\CoefXUDn}; 
(10, 5)="VF" *{\CoefAlaV{1}}; 
(10, 0)="VFEt" *{\scriptstyle \mathcal{C}(X^2_1)};
(10, -20)="V"  *{\CoefXDDn};
(20, 5)="WF" *{\CoefAlaW{1}};
(20, 0)="WFEt" *{\scriptstyle \mathcal{C}(X^3_1)}; 
(20, -20)="W"  *{\CoefXTuDn};
(30, 5)="W2F" *{\CoefAlaW{2}}; 
(30, 0)="W2FEt" *{\scriptstyle \mathcal{C}(X^3_2)};
(30, -20)="W2"  *{\CoefXTdDn};
(40, 5)="dots"  *{\cdots};
(50, 5)="WnF" *{\CoefAlaW{{}_{n_1\text{-}1}}}; 
(50, 0)="WnFEt" *{\scriptstyle \mathcal{C}(X^3_{n_1\text{-}1})};
(50, -20)="Wn"  *{\CoefXTnDn};
(60, 5)="MF" *{\CoefMax}; 
(60, 0)="MFEt" *{\scriptstyle \mathcal{C}(W_0)};
(60, -20)="M"  *{\CoefMaxDnO};
%-----------------------
"D3";"D4" **@{-};
"D1";"D2" **@{-};
\endxy
\end{displaymath}
For the computation of the reduced ditalgebra $\mathcal{A}^X$ consider the representations of $\mathbf{D_n}$ above
and observe that $n_2=n_3=2$ and that we can choose vectors $a^2_1$ and $a^3_1$ such that
$a^2_{2} =\left[ \begin{smallmatrix}1\\0\end{smallmatrix} \right]$ and $a^3_{2} =\left[\begin{smallmatrix}0\\1\end{smallmatrix}\right]$.
Notice that in this case we can take $a^1_1$ such that $a^1_{n_1}=\left[ \begin{smallmatrix}1\\1\end{smallmatrix} \right]$.
Then, defining $a=a^2_{2}$ and $a'=a^3_{2}$ we have $a^1_{n_1}=a+a'$ (cf. lemma~\ref{(DE)L:AXcero}).
%------------------------------------------------
%------------------------------------------------
\begin{displaymath}
 \xymatrix@C=4pc@R=2.5pc{
& & & & *+++[]{\bullet_{\omega}} \\
\\
& & & \bullet \ar@{.>}[rd]^-{\gamma^2_{1,2}} \ar@{<-}@/^.5pc/[ruu]^(.2){\xi^2_1} \\
\bullet \ar@{.>}[r]^(.7){\gamma^1_{1,2}} \ar@{.>}@/^1.5pc/[rrr]^(.3){\gamma^1_{1,n_1-1}} 
\ar@{.>}@/_2pc/[rrrr]_-{\gamma^1_{1,n_1}} \ar@{<-}@/^2pc/[rrrruuu]^(.27){\xi^1_{1}} 
 & \bullet \ar@{.>}@/_.5pc/[rr]_-{\gamma^1_{2,n_1-1}} 
\ar@{.>}@/^1.5pc/[rrr]^-{\gamma^1_{2,n_1}} \ar@{<-}@/^1pc/[rrruuu]^(.33){\xi^1_{2}} 
& \ldots & \bullet \ar@{.>}[r]^-{\gamma^1_{n_1-1,n_1}} \ar@{<-}@/^.8pc/[ruuu]^(.2){\xi^1_{n_1-1}}
& *++[]{\bullet_{m}} \ar@{<-}[uuu]^(.4){\xi'} \ar@{<-}@<-1ex>[uuu]_(.4){\xi} \\
& & & \bullet \ar@{.>}[ru]_-{\gamma^3_{1,2}} \ar@{<-}@/^.7pc/[ruuuu]^(.1){\xi^3_1} &
}
\end{displaymath}
The differential $\delta^x$ is described in the following table ($s \in \{1,2,3\}$).
%---------------------------------------------
\begin{center}
 \begin{tabular}{c | l }
Arrow & Its differential $\delta^x$ \\
\hline
$\xi$   & $\gamma^2_{1,2} \otimes \xi^2_1+\sum_{1\leq i<n_1}\gamma^1_{i,n_1} \otimes \xi^1_i$,    \\
$\xi' $ & $\gamma^3_{1,2} \otimes \xi^3_1+\sum_{1\leq i<n_1}\gamma^1_{i,n_1} \otimes \xi^1_i$,  \\
$\xi^s_j $ & $\sum_{1 \leq i<j}\gamma^s_{i,j} \otimes \xi^s_i$.  \\
 \end{tabular}
\end{center}
In particular observe that $\delta^x(\xi^1_1)=\delta^x(\xi^2_1)=\delta^x(\xi^3_1)=0$.

We describe now the reduced ditalgebra $\mathcal{A}^Y$. By lemma~\ref{(DE)L:AYcero} we must choose vectors
$b^1_1$, $b^2_1$ and $b^3_1$ such that
\[
 b^1_1 \neq a^1_{n_1}=\left[ \begin{smallmatrix}1\\1\end{smallmatrix} \right], \quad
 b^2_1 \neq a^2_{2}=\left[ \begin{smallmatrix}1\\0\end{smallmatrix} \right], \quad 
 b^3_1 \neq a^3_{2}=\left[ \begin{smallmatrix}0\\1\end{smallmatrix} \right], \quad
\]
(considered as elements in $\Hom_{A_0}(R,W_0)$) we take then $b^1_1 =\left[ \begin{smallmatrix}1\\0\end{smallmatrix} \right]$,
 $b^2_1 =\left[ \begin{smallmatrix}0\\1\end{smallmatrix} \right]$,
$b^3_1 =\left[ \begin{smallmatrix}1\\1\end{smallmatrix} \right]$
and make $b=b^1_1$ and $b'=b^2_1$. In this way we have $b^3_1=b+b'$.
%------------------------------------------------
%------------------------------------------------
\begin{displaymath}
 \xymatrix@C=4pc@R=2.5pc{
*+++[]{\bullet_{\omega}} \ar@<-.5ex>[ddd]_(.7){\theta} \ar@<.5ex>[ddd]^(.7){\theta'} 
\ar@/^.5pc/[ddr]^(.8){\theta^2_2} \ar@/^.8pc/[dddr]^(.8){\theta^1_2} \ar@/^.7pc/[ddddr]^(.9){\theta^3_2} 
\ar@/^1pc/[dddrrr]^(.7){\theta^1_{n_1-1}} \ar@/^2pc/[dddrrrr]^(.76){\theta^1_{n_1}} \\
\\
& \bullet \\
*++[]{\bullet_m} \ar@{.>}[r]^(.7){\varepsilon^1_{1,2}} \ar@{.>}[ru]^(.5){\varepsilon^2_{1,2}} \ar@{.>}[rd]_(.5){\varepsilon^3_{1,2}} 
\ar@{.>}@/^1.5pc/[rrr]^(.5){\varepsilon^1_{1,n_1-1}} \ar@{.>}@/_2pc/[rrrr]_(.5){\varepsilon^1_{1,n_1}}
& \bullet \ar@{.>}@/^1.5pc/[rrr]^(.7){\varepsilon^1_{2,n_1}} \ar@{.>}@/_.5pc/[rr]_(.5){\varepsilon^1_{2,n_1-1}} 
& \cdots & \bullet \ar@{.>}[r]^(.35){\varepsilon^1_{n_1-1,n_1}} & \bullet \\
& \bullet
}
\end{displaymath}
Then the reduced differential $\delta^y$ has the following form,
%---------------------------------------------
\begin{center}
 \begin{tabular}{c | l }
Arrow & Its differential $\delta^y$ \\
\hline
$\theta^1_j $ & $\varepsilon^1_{1,j} \otimes \theta+\sum_{1 < i<j}\varepsilon^1_{i,j} \otimes \theta^1_i$,  \\
$\theta^2_2 $ & $\varepsilon^2_{1,2} \otimes \theta'$,  \\
$\theta^3_2 $ & $\varepsilon^3_{1,2} \otimes \theta+\varepsilon^3_{1,2} \otimes \theta'$,  \\
$\theta$ & $0$,    \\
$\theta' $ & $0$.  \\
 \end{tabular}
\end{center}
To give preinjective representations of $\widetilde{\mathbf{D_n}}$ we choose first modules of the Dynkin quiver 
$\mathbf{D_n}$ corresponding to the elements of $\mathcal{Y}_0$ (see figure~\ref{(DE)F:Dn}).
Their coefficient quivers are shown next.
\begin{displaymath}
\xy
( -19,-32)="D1" *{};
( -16,-29.5)="D2" *{};
( -11,-5)="D3" *{};
( -14,-7.5)="D4" *{};
(-10,-1)="Col" *{\xymatrix@!0@C=15pt@R=9.5pt{{} \\ {} & & {}_{n}  
\\ & {}_{n-1} \ar@{-}[dd] \\ \\ & {}_{n-2} \\ & {\cdots} \\ & {}_{4} \ar@{-}[dd] 
\\ \\ & {}_{3} \ar@{-}[dd] \\ {}_{2} \\ & {}_{1} }};
( 0, 5)="UF" *{\CoefAlaW{2}}; 
( 0, 0)="UFEt" *{\scriptstyle \mathcal{C}(Y^3_2)};
( 0, -20)="U"  *{\CoefYDDn}; 
(10, 5)="VF" *{\CoefAlaV{2}};
(10, 0)="VFEt" *{\scriptstyle \mathcal{C}(Y^2_2)}; 
(10, -20)="V"  *{\CoefYUDn};
(20, 5)="WF" *{\CoefAlaU{2}};
(20, 0)="WFEt" *{\scriptstyle \mathcal{C}(Y^1_2)}; 
(20, -20)="W"  *{\CoefYTuDn};
(30, 5)="W2F" *{\CoefAlaU{3}};
(30, 0)="W2FEt" *{\scriptstyle \mathcal{C}(Y^1_3)}; 
(30, -20)="W2"  *{\CoefYTdDn};
(40, 5)="dots"  *{\cdots};
(50, 5)="WnF" *{\CoefAlaU{n_1}}; 
(50, 0)="WnFEt" *{\scriptstyle \mathcal{C}(Y^1_{n_1})};
(50, -20)="Wn"  *{\CoefYTnDn};
(60, 5)="MF" *{\CoefMax}; 
(60, 0)="MFEt" *{\scriptstyle \mathcal{C}(W_0)};
(60, -20)="M"  *{\CoefMaxDnO};
%-----------------------
"D3";"D4" **@{-};
"D1";"D2" **@{-};
\endxy
\end{displaymath}

For the construction of indecomposable $\widetilde{\mathbf{D_n}}$-representations we use these coefficient quivers
and the series of exceptional representations of the reduced ditalgebra $\mathcal{A}^Y$ having rank one 
(lemma~\ref{(DE)L:excRedRnkUno}) and rank two (lemma~\ref{(DE)L:excRedRnkDos}). In the following figures we show
the representation in matricial form (left) and the coefficient quiver for a particular case in the parameter $\ell$ (right).

%------------------
Rank one representations $M^{(1)}_{1,1,1}[\ell]=F^Y(\widetilde{M}^{(1)}_{1,1,1}[\ell])$ ($\ell\geq 0$).
\begin{equation} \label{(C)F:Dn1}
 \xy
(-60,0)  *{}; %Se pone para centrar la figura
(0,12)="Et"   *{\scriptstyle \mathcal{C}(M^{(1)}_{1,1,1}[4])};
(-23,4)="Col" *{\xymatrix@!0@L=6pt@R=10pt@C=35pt{
& k^{\ell+1} \ar[dd]_*+<1em>{^{\left[ \begin{smallmatrix} I_{\ell}^{\rightarrow}\\I_{\ell}^{\leftarrow} \end{smallmatrix} \right]}} \\
& & k^{\ell} \ar[ld]^(.3){\left[ \begin{smallmatrix} I_{\ell}\\I_{\ell} \end{smallmatrix} \right]} \\
& k^{2\ell} \ar[dd]_-{=} \\ \\
& k^{2\ell} \ar@{}[dd]|-{\cdots} \\ \\
& k^{2\ell} \ar[dd]_-{=} \\ \\
& k^{2\ell} \ar[dd]^-{\left[ \begin{smallmatrix} 0&I_{\ell} \end{smallmatrix} \right]} 
\ar[dl]_(.7){\left[ \begin{smallmatrix} I_{\ell}&0 \end{smallmatrix} \right]} \\
k^{\ell} \\
& k^{\ell}}};
(-16,8)="A1"  *{};
( -8,8)="A2"  *{};
(  0,8)="A3"  *{};
( 8, 8)="A4"  *{};
(16, 8)="A5"  *{};
(-13,0)="B1" *{};
( -5,0)="B2" *{};
( 3, 0)="B3" *{};
(11, 0)="B4" *{};
(-11,  0)="C1" *{};
(-3,  0)="C2" *{};
( 5,  0)="C3" *{};
(13,  0)="C4" *{};
(-12,-12)="D1" *{\CoefMaxDnO};
(-4, -12)="D2" *{\CoefMaxDnO};
( 4, -12)="D3" *{\CoefMaxDnO};
(12, -12)="D3" *{\CoefMaxDnO};
%-----------------------
"A1";"B1" **@{-};
"A2";"B2" **@{-};
"A3";"B3" **@{-};
"A4";"B4" **@{-};
"A2";"C1" **@{-};
"A3";"C2" **@{-};
"A4";"C3" **@{-};
"A5";"C4" **@{-};
\endxy
\end{equation}

%------------------
Rank one representations $M^{(1)}_{2,1,1}[\ell]_2=F^Y(\widetilde{M}^{(1)}_{2,1,1}[\ell]_2)$ ($\ell\geq 0$).
\begin{equation} \label{(C)F:Dn2}
 \xy
(-60,0)  *{}; %Se pone para centrar la figura
(3,12)="Et"  *{\scriptstyle \mathcal{C}(M^{(1)}_{2,1,1}[4]_2)};
(-23,4)="Col"  *{\xymatrix@!0@L=6pt@R=10pt@C=35pt{
& k^{\ell+1} \ar[dd]_*+<1em>{^{\left[ \begin{smallmatrix} I_{\ell+1}\\I_{\ell}^{\leftarrow} \end{smallmatrix} \right]}} \\
& & k^{\ell+1} \ar[ld]^(.3){\left[ \begin{smallmatrix} I_{\ell+1}\\I_{\ell}^{\rightarrow} \end{smallmatrix} \right]} \\
& k^{2\ell+1} \ar[dd]_-{=} \\ \\
& k^{2\ell+1} \ar@{}[dd]|-{\cdots} \\ \\
& k^{2\ell+1} \ar[dd]_-{=} \\ \\
& k^{2\ell+1} \ar[dd]^-{\left[ \begin{smallmatrix} 0&I_{\ell} \end{smallmatrix} \right]} 
\ar[dl]_(.7){\left[ \begin{smallmatrix} I_{\ell+1}&0 \end{smallmatrix} \right]} \\
k^{\ell+1} \\
& k^{\ell}}};
(-16,8)="A1"  *{};
( -8,8)="A2"  *{};
(  0,8)="A3"  *{};
( 8, 8)="A4"  *{};
(16, 8)="A5"  *{};
(-13,0)="B1" *{};
( -5,0)="B2" *{};
( 3, 0)="B3" *{};
(11, 0)="B4" *{};
(-11,  0)="C1" *{};
(-3,  0)="C2" *{};
( 5,  0)="C3" *{};
(13,  0)="C4" *{};
(19,  0)="D" *{};
(-12,-12)="D1" *{\CoefMaxDnO};
(-4, -12)="D2" *{\CoefMaxDnO};
( 4, -12)="D3" *{\CoefMaxDnO};
(12, -12)="D3" *{\CoefMaxDnO};
(19, -12)="D4" *{\CoefYUDn};
%-----------------------
"A1";"B1" **@{-};
"A2";"B2" **@{-};
"A3";"B3" **@{-};
"A4";"B4" **@{-};
"A2";"C1" **@{-};
"A3";"C2" **@{-};
"A4";"C3" **@{-};
"A5";"C4" **@{-};
"A5";"D" **@{-};
\endxy
\end{equation}

%------------------
Rank one representations $M^{(1)}_{2,1,1}[\ell]^2=F^Y(\widetilde{M}^{(1)}_{2,1,1}[\ell]^2)$ ($\ell\geq 0$).
\begin{equation} \label{(C)F:Dn3}
 \xy
(-60,0)  *{}; %Se pone para centrar la figura
(-3,12)="Et"  *{\scriptstyle \mathcal{C}(M^{(1)}_{2,1,1}[4]^2)};
(-27,4)="Col"  *{\xymatrix@!0@L=6pt@R=10pt@C=35pt{
& k^{\ell+1} \ar[dd]_*+<1em>{^{\left[ \begin{smallmatrix} I_{\ell}^{\rightarrow}\\I_{\ell+1} \end{smallmatrix} \right]}} \\
& & k^{\ell+1} \ar[ld]^(.3){\left[ \begin{smallmatrix} I_{\ell+1}\\I_{\ell}^{\rightarrow} \end{smallmatrix} \right]} \\
& k^{2\ell+1} \ar[dd]_-{=} \\ \\
& k^{2\ell+1} \ar@{}[dd]|-{\cdots} \\ \\
& k^{2\ell+1} \ar[dd]_-{=} \\ \\
& k^{2\ell+1} \ar[dd]^-{\left[ \begin{smallmatrix} 0&I_{\ell+1} \end{smallmatrix} \right]} 
\ar[dl]_(.7){\left[ \begin{smallmatrix} I_{\ell}&0 \end{smallmatrix} \right]} \\
k^{\ell} \\
& k^{\ell+1}}};
(-16,8)="A1"  *{};
( -8,8)="A2"  *{};
(  0,8)="A3"  *{};
( 8, 8)="A4"  *{};
(16, 8)="A5"  *{};
(-13,0)="B1" *{};
( -5,0)="B2" *{};
( 3, 0)="B3" *{};
(11, 0)="B4" *{};
(-11,  0)="C1" *{};
(-3,  0)="C2" *{};
( 5,  0)="C3" *{};
(13,  0)="C4" *{};
(-19,  0)="D" *{};
(-12,-12)="D1" *{\CoefMaxDnO};
(-4, -12)="D2" *{\CoefMaxDnO};
( 4, -12)="D3" *{\CoefMaxDnO};
(12, -12)="D3" *{\CoefMaxDnO};
(-19, -12)="D4" *{\CoefYDDn};
%-----------------------
"A1";"B1" **@{-};
"A2";"B2" **@{-};
"A3";"B3" **@{-};
"A4";"B4" **@{-};
"A2";"C1" **@{-};
"A3";"C2" **@{-};
"A4";"C3" **@{-};
"A5";"C4" **@{-};
"A1";"D" **@{-};
\endxy
\end{equation}

%------------------
Rank one representations ${}_iM^{(1)}_{2,1,1}[\ell]=F^Y({}_i\widetilde{M}^{(1)}_{2,1,1}[\ell])$ ($\ell\geq 0$).
\begin{equation} \label{(C)F:Dn4}
 \xy
(-56,-9.5)="Cuad" *{ \xy
( 8, 7.5)="A"  *{};
(-8, 7.5)="B"  *{};
( 8,-7.5)="C"  *{};
(-8,-7.5)="D"  *{};
%-----------------------
"A";"B" **@{-};
"C";"D" **@{-};
"A";"C" **@{-};
"B";"D" **@{-};
 \endxy};
%(-60,0)  *{}; %Se pone para centrar la figura
(-3,12)="Et"  *{\scriptstyle \mathcal{C}({}_iM^{(1)}_{2,1,1}[4])};
(-27,4)="Col"  *{\xymatrix@!0@L=6pt@R=10pt@C=35pt{
& k^{\ell+1} \ar[dd]_*+<1em>{^{\left[ \begin{smallmatrix} I_{\ell+1}\\I_{\ell}^{\rightarrow} \end{smallmatrix} \right]}} \\
& & k^{\ell} \ar[ld]^(.3){\left[ \begin{smallmatrix} I_{\ell}^{\uparrow}\\I_{\ell} \end{smallmatrix} \right]} \\
& k^{2\ell+1} \ar[dd]_-{=} \\ \\
& k^{2\ell+1} \ar@{}[dd]|-{\cdots} \\ \\
& k^{2\ell} \ar[dd]_-{=} \\ \\
& k^{2\ell} \ar[dd]^-{\left[ \begin{smallmatrix} I_{\ell}&0 \end{smallmatrix} \right]} 
\ar[dl]_(.7){\left[ \begin{smallmatrix} 0&I_{\ell} \end{smallmatrix} \right]} \\
k^{\ell} \\
& k^{\ell}}};
(-27,-3)="Col2"  *{\xymatrix@!0@L=6pt@R=10pt@C=35pt{
k^{2\ell+1} \ar[dd]_*+<1em>{^{I_{2\ell}^{\leftarrow}}} \\ \\
k^{2\ell}}};
(-54,-3)="dots1"  *{\scriptstyle \cdots};
(-54,-16)="dots2"  *{\scriptstyle \cdots};
(-16,8)="A1"  *{};
( -8,8)="A2"  *{};
(  0,8)="A3"  *{};
( 8, 8)="A4"  *{};
(16, 8)="A5"  *{};
(-13,0)="B1" *{};
( -5,0)="B2" *{};
( 3, 0)="B3" *{};
(11, 0)="B4" *{};
(-11,  0)="C1" *{};
(-3,  0)="C2" *{};
( 5,  0)="C3" *{};
(13,  0)="C4" *{};
(-19,  0)="D" *{};
(-12,-12)="D1" *{\CoefMaxDnO};
(-4, -12)="D2" *{\CoefMaxDnO};
( 4, -12)="D3" *{\CoefMaxDnO};
(12, -12)="D3" *{\CoefMaxDnO};
(-19, 0)="D4" *{\CoefAlaU{i}};
%-----------------------
"A1";"B1" **@{-};
"A2";"B2" **@{-};
"A3";"B3" **@{-};
"A4";"B4" **@{-};
"A2";"C1" **@{-};
"A3";"C2" **@{-};
"A4";"C3" **@{-};
"A5";"C4" **@{-};
"A1";"D" **@{-};
\endxy
\end{equation}

%------------------
Rank one representations $M^{(1)}_{2,2,1}[\ell]^2_2=F^Y(\widetilde{M}^{(1)}_{2,2,1}[\ell]^2_2)$ ($\ell\geq 0$).
\begin{equation} \label{(C)F:Dn5}
 \xy
(-60,0)  *{}; %Se pone para centrar la figura
(0,12)="Et"  *{\scriptstyle \mathcal{C}(M^{(1)}_{2,2,1}[4]^2_2)};
(-27,4)="Col"  *{\xymatrix@!0@L=6pt@R=10pt@C=35pt{
&   k^{\ell+1} \ar[dd]_*+<1em>{^{\left[ \begin{smallmatrix} I_{\ell+1}\\I_{\ell+1} \end{smallmatrix} \right]}} \\
& & k^{\ell+2} \ar[ld]^(.3){\left[ \begin{smallmatrix} I_{\ell+1}^{\rightarrow}\\I_{\ell+1}^{\leftarrow} \end{smallmatrix} \right]} \\
& k^{2\ell+2} \ar[dd]_-{=} \\ \\
& k^{2\ell+2} \ar@{}[dd]|-{\cdots} \\ \\
& k^{2\ell+2} \ar[dd]_-{=} \\ \\
& k^{2\ell+2} \ar[dd]^-{\left[ \begin{smallmatrix} I_{\ell+1}&0 \end{smallmatrix} \right]} 
\ar[dl]_(.7){\left[ \begin{smallmatrix} 0&I_{\ell+1} \end{smallmatrix} \right]} \\
k^{\ell+1} \\
& k^{\ell+1}}};
(-16,8)="A1"  *{};
( -8,8)="A2"  *{};
(  0,8)="A3"  *{};
( 8, 8)="A4"  *{};
(16, 8)="A5"  *{};
(-13,0)="B1" *{};
( -5,0)="B2" *{};
( 3, 0)="B3" *{};
(11, 0)="B4" *{};
(-11,  0)="C1" *{};
(-3,  0)="C2" *{};
( 5,  0)="C3" *{};
(13,  0)="C4" *{};
(-19,  0)="D" *{};
( 19,  0)="E" *{};
(-12,-12)="D1" *{\CoefMaxDnO};
(-4, -12)="D2" *{\CoefMaxDnO};
( 4, -12)="D3" *{\CoefMaxDnO};
(12, -12)="D3" *{\CoefMaxDnO};
(-19, -12)="D4" *{\CoefYDDn};
(19, -12)="D5" *{\CoefYUDn};
%-----------------------
"A1";"B1" **@{-};
"A2";"B2" **@{-};
"A3";"B3" **@{-};
"A4";"B4" **@{-};
"A2";"C1" **@{-};
"A3";"C2" **@{-};
"A4";"C3" **@{-};
"A5";"C4" **@{-};
"A1";"D" **@{-};
"A5";"E" **@{-};
\endxy
\end{equation}

%------------------
Rank one representations ${}_iM^{(1)}_{2,2,1}[\ell]_2=F^Y({}_i\widetilde{M}^{(1)}_{2,2,1}[\ell]_2)$ ($\ell\geq 0$).
\begin{equation} \label{(C)F:Dn6}
 \xy
(-57,-9.5)="Cuad" *{ \xy
( 8, 7.5)="A"  *{};
(-8, 7.5)="B"  *{};
( 8,-7.5)="C"  *{};
(-8,-7.5)="D"  *{};
%-----------------------
"A";"B" **@{-};
"C";"D" **@{-};
"A";"C" **@{-};
"B";"D" **@{-};
 \endxy};
%(-60,0)  *{}; %Se pone para centrar la figura
(0,12)="Et"  *{\scriptstyle \mathcal{C}({}_iM^{(1)}_{2,2,1}[4]_2)};
(-27,4)="Col"  *{\xymatrix@!0@L=6pt@R=10pt@C=35pt{
&   k^{\ell+1} \ar[dd]_*+<1em>{^{\left[ \begin{smallmatrix} I_{\ell+1}\\I_{\ell+1} \end{smallmatrix} \right]}} \\
& & k^{\ell+1} \ar[ld]^(.3){\left[ \begin{smallmatrix} I_{\ell}^{\rightarrow \uparrow}\\I_{\ell+1} \end{smallmatrix} \right]} \\
& k^{2\ell+2} \ar[dd]_-{=} \\ \\
& k^{2\ell+2} \ar@{}[dd]|-{\cdots} \\ \\
& k^{2\ell+1} \ar[dd]_-{=} \\ \\
& k^{2\ell+1} \ar[dd]^-{\left[ \begin{smallmatrix} I_{\ell}&0 \end{smallmatrix} \right]} 
\ar[dl]_(.7){\left[ \begin{smallmatrix} 0&I_{\ell+1} \end{smallmatrix} \right]} \\
k^{\ell+1} \\
& k^{\ell}}};
(-27,-3)="Col2"  *{\xymatrix@!0@L=6pt@R=10pt@C=35pt{
k^{2\ell+2} \ar[dd]_*+<1em>{^{I_{2\ell+1}^{\leftarrow}}} \\ \\
k^{2\ell+1}}};
(-54,-3)="dots1"  *{\scriptstyle \cdots};
(-54,-16)="dots2"  *{\scriptstyle \cdots};
(-16,8)="A1"  *{};
( -8,8)="A2"  *{};
(  0,8)="A3"  *{};
( 8, 8)="A4"  *{};
(16, 8)="A5"  *{};
(-13,0)="B1" *{};
( -5,0)="B2" *{};
( 3, 0)="B3" *{};
(11, 0)="B4" *{};
(-11,  0)="C1" *{};
(-3,  0)="C2" *{};
( 5,  0)="C3" *{};
(13,  0)="C4" *{};
(-19,  0)="D" *{};
( 19,  0)="E" *{};
(-12,-12)="D1" *{\CoefMaxDnO};
(-4, -12)="D2" *{\CoefMaxDnO};
( 4, -12)="D3" *{\CoefMaxDnO};
(12, -12)="D3" *{\CoefMaxDnO};
(-19, 0)="D4" *{\CoefAlaU{i}};
(19, -12)="D5" *{\CoefYUDn};
%-----------------------
"A1";"B1" **@{-};
"A2";"B2" **@{-};
"A3";"B3" **@{-};
"A4";"B4" **@{-};
"A2";"C1" **@{-};
"A3";"C2" **@{-};
"A4";"C3" **@{-};
"A5";"C4" **@{-};
"A1";"D" **@{-};
"A5";"E" **@{-};
\endxy
\end{equation}

%------------------
Rank one representations ${}_iM^{(1)}_{2,2,1}[\ell]^2=F^Y({}_i\widetilde{M}^{(1)}_{2,2,1}[\ell]^2)$ ($\ell\geq 0$).
\begin{equation} \label{(C)F:Dn7}
 \xy
(-57,-9.5)="Cuad" *{ \xy
( 8, 7.5)="A"  *{};
(-8, 7.5)="B"  *{};
( 8,-7.5)="C"  *{};
(-8,-7.5)="D"  *{};
%-----------------------
"A";"B" **@{-};
"C";"D" **@{-};
"A";"C" **@{-};
"B";"D" **@{-};
 \endxy};
%(-60,0)  *{}; %Se pone para centrar la figura
(0,12)="Et"  *{\scriptstyle \mathcal{C}({}_iM^{(1)}_{2,2,1}[4]^2)};
(-27,4)="Col"  *{\xymatrix@!0@L=6pt@R=10pt@C=35pt{
&   k^{\ell+1} \ar[dd]_*+<1em>{^{\left[ \begin{smallmatrix} I_{\ell+1}\\I_{\ell+1} \end{smallmatrix} \right]}} \\
& & k^{\ell+1} \ar[ld]^(.3){\left[ \begin{smallmatrix} I_{\ell+1}^{\uparrow}\\I_{\ell}^{\rightarrow} \end{smallmatrix} \right]} \\
& k^{2\ell+2} \ar[dd]_-{=} \\ \\
& k^{2\ell+2} \ar@{}[dd]|-{\cdots} \\ \\
& k^{2\ell+1} \ar[dd]_-{=} \\ \\
& k^{2\ell+1} \ar[dd]^-{\left[ \begin{smallmatrix} I_{\ell+1}&0 \end{smallmatrix} \right]} 
\ar[dl]_(.7){\left[ \begin{smallmatrix} 0&I_{\ell} \end{smallmatrix} \right]} \\
k^{\ell} \\
& k^{\ell+1}}};
(-27,-3)="Col2"  *{\xymatrix@!0@L=6pt@R=10pt@C=35pt{
k^{2\ell+2} \ar[dd]_*+<1em>{^{I_{2\ell+1}^{\leftarrow}}} \\ \\
k^{2\ell+1}}};
(-54,-3)="dots1"  *{\scriptstyle \cdots};
(-54,-16)="dots2"  *{\scriptstyle \cdots};
(-16,8)="A1"  *{};
( -8,8)="A2"  *{};
(  0,8)="A3"  *{};
( 8, 8)="A4"  *{};
(16, 8)="A5"  *{};
(-13,0)="B1" *{};
( -5,0)="B2" *{};
( 3, 0)="B3" *{};
(11, 0)="B4" *{};
(-11,  0)="C1" *{};
(-3,  0)="C2" *{};
( 5,  0)="C3" *{};
(13,  0)="C4" *{};
(-19,  0)="D" *{};
( 19,  0)="E" *{};
(-12,-12)="D1" *{\CoefMaxDnO};
(-4, -12)="D2" *{\CoefMaxDnO};
( 4, -12)="D3" *{\CoefMaxDnO};
(12, -12)="D3" *{\CoefMaxDnO};
(-19, 0)="D4" *{\CoefAlaU{i}};
(19, -12)="D5" *{\CoefYDDn};
%-----------------------
"A1";"B1" **@{-};
"A2";"B2" **@{-};
"A3";"B3" **@{-};
"A4";"B4" **@{-};
"A2";"C1" **@{-};
"A3";"C2" **@{-};
"A4";"C3" **@{-};
"A5";"C4" **@{-};
"A1";"D" **@{-};
"A5";"E" **@{-};
\endxy
\end{equation}

%------------------
Rank one representations ${}_iM^{(1)}_{2,2,2}[\ell]^2_2=F^Y({}_i\widetilde{M}^{(1)}_{2,2,2}[\ell]^2_2)$ ($\ell\geq 0$).
\begin{equation} \label{(C)F:Dn8}
 \xy
(-63,-9.5)="Cuad" *{ \xy
( 8, 7.5)="A"  *{};
(-8, 7.5)="B"  *{};
( 8,-7.5)="C"  *{};
(-8,-7.5)="D"  *{};
%-----------------------
"A";"B" **@{-};
"C";"D" **@{-};
"A";"C" **@{-};
"B";"D" **@{-};
 \endxy};
%(-60,0)  *{}; %Se pone para centrar la figura
(0,12)="Et"  *{\scriptstyle \mathcal{C}({}_iM^{(1)}_{2,2,2}[4]_2^2)};
(-30,4)="Col"  *{\xymatrix@!0@L=6pt@R=10pt@C=35pt{
&   k^{\ell+1} \ar[dd]_*+<1em>{^{\left[ \begin{smallmatrix} I_{\ell+1}\\I_{\ell+1}\\ \Renglon{1}{0}{\cdots}{0}{0} \end{smallmatrix} \right]}} \\
& & k^{\ell+2} \ar[ld]^(.3){\left[ \begin{smallmatrix} I_{\ell+2}\\I_{\ell}^{\rightrightarrows} \end{smallmatrix} \right]} \\
& k^{2\ell+3} \ar[dd]_-{=} \\ \\
& k^{2\ell+3} \ar@{}[dd]|-{\cdots} \\ \\
& k^{2\ell+2} \ar[dd]_-{=} \\ \\
& k^{2\ell+2} \ar[dd]^-{\left[ \begin{smallmatrix} 0&I_{\ell+1} \end{smallmatrix} \right]} 
\ar[dl]_(.7){\left[ \begin{smallmatrix} I_{\ell+1}&0 \end{smallmatrix} \right]} \\
  k^{\ell+1} \\
& k^{\ell+1}}};
(-30,-3)="Col2"  *{\xymatrix@!0@L=6pt@R=10pt@C=35pt{
k^{2\ell+3} \ar[dd]_*+<1em>{^{I_{2\ell+2}^{\rightarrow}}} \\ \\
k^{2\ell+2}}};
(-60,-3)="dots1"  *{\scriptstyle \cdots};
(-60,-16)="dots2"  *{\scriptstyle \cdots};
(-16,8)="A1"  *{};
( -8,8)="A2"  *{};
(  0,8)="A3"  *{};
( 8, 8)="A4"  *{};
(16, 8)="A5"  *{};
(-13,0)="B1" *{};
( -5,0)="B2" *{};
( 3, 0)="B3" *{};
(11, 0)="B4" *{};
(-11,  0)="C1" *{};
(-3,  0)="C2" *{};
( 5,  0)="C3" *{};
(13,  0)="C4" *{};
(-19,  0)="D" *{};
( 19,  0)="E" *{};
(-12,-12)="D1" *{\CoefMaxDnO};
(-4, -12)="D2" *{\CoefMaxDnO};
( 4, -12)="D3" *{\CoefMaxDnO};
(12, -12)="D3" *{\CoefMaxDnO};
(-19,-12)="D4" *{\CoefYDDn};
(19, -12)="D5" *{\CoefYUDn};
(-26,  0)="D6" *{\CoefAlaU{i}};
%-----------------------
"A1";"B1" **@{-};
"A2";"B2" **@{-};
"A3";"B3" **@{-};
"A4";"B4" **@{-};
"A2";"C1" **@{-};
"A3";"C2" **@{-};
"A4";"C3" **@{-};
"A5";"C4" **@{-};
"A1";"D"  **@{-};
"A5";"E"  **@{-};
"A1";"D6" **\crv{(-24,7)};
\endxy
\end{equation}

%------------------
Rank two representations ${}_iN^{(2)}_{2,2,2}[\ell]^2_2=F^Y({}_i\widetilde{N}^{(2)}_{2,2,2}[\ell]^2_2)$ ($\ell=2k$).
\begin{equation} \label{(C)F:Dn9}
 \xy
(-69,-9.5)="Cuad" *{ \xy
( 8, 7.5)="A"  *{};
(-8, 7.5)="B"  *{};
( 8,-7.5)="C"  *{};
(-8,-7.5)="D"  *{};
%-----------------------
"A";"B" **@{-};
"C";"D" **@{-};
"A";"C" **@{-};
"B";"D" **@{-};
 \endxy};
(0,12)="Et"  *{\scriptstyle \mathcal{C}({}_{i}N^{(2)}_{2,2,2}[6]_2^2)};
(-33,4)="Col"  *{\xymatrix@!0@L=6pt@R=10pt@C=35pt{
&   k^{\ell+2} \ar[dd]_*+<1em>{^{A}} \\
& & k^{\ell+2} \ar[ld]^(.3){B} \\
& k^{2\ell+3} \ar[dd]_-{=} \\ \\
& k^{2\ell+3} \ar@{}[dd]|-{\cdots} \\ \\
& k^{2\ell+2} \ar[dd]_-{=} \\ \\
& k^{2\ell+2} \ar[dd]^-{\left[ \begin{smallmatrix} I_{\ell+1}&0 \end{smallmatrix} \right]} 
\ar[dl]_(.7){\left[ \begin{smallmatrix} 0&I_{\ell+1} \end{smallmatrix} \right]} \\
  k^{\ell+1} \\
& k^{\ell+1}}};
(-33,-3)="Col2"  *{\xymatrix@!0@L=6pt@R=10pt@C=35pt{
k^{2\ell+3} \ar[dd]_*+<1em>{^{I_{2\ell+2}^{\leftarrow}}} \\ \\
k^{2\ell+2}}};
(-66,-3)="dots1"  *{\scriptstyle \cdots};
(-66,-16)="dots2"  *{\scriptstyle \cdots};
( -5,8)="A1"  *{};
(-13,8)="A2"  *{};
(-21,8)="A3"  *{};
(-29,8)="A4"  *{};
( -9,-12)="D1" *{\CoefMaxDnO};
(-17,-12)="D2" *{\CoefMaxDnO};
(-25,-12)="D3" *{\CoefMaxDnO};
(-10,0)="C1" *{};
(-18,0)="C2" *{};
(-26,0)="C3" *{};
( -8,0)="B1" *{};
(-16,0)="B2" *{};
(-24,0)="B3" *{};
( 5,8)="A'1"  *{};
(13,8)="A'2"  *{};
(21,8)="A'3"  *{};
(29,8)="A'4"  *{};
( 9,-12)="D'1" *{\CoefMaxDnO};
(17,-12)="D'2" *{\CoefMaxDnO};
(25,-12)="D'3" *{\CoefMaxDnO};
(10,0)="C'1" *{};
(18,0)="C'2" *{};
(26,0)="C'3" *{};
( 8,0)="B'1" *{};
(16,0)="B'2" *{};
(24,0)="B'3" *{};
(0, 0)="D" *{};
(32, 0)="E" *{};
(-32,0)="F" *{};
%(12, -12)="D3" *{\CoefMaxDnO};
(0,-12)="D4" *{\CoefYDDnPrime};
(32, -12)="D5" *{\CoefYUDn};
(-32,  0)="D6" *{\CoefAlaU{i}};
%-----------------------
"A1";"B1" **@{-};
"A2";"B2" **@{-};
"A3";"B3" **@{-};
"A2";"C1" **@{-};
"A3";"C2" **@{-};
"A4";"C3" **@{-};
"A'1";"B'1" **@{-};
"A'2";"B'2" **@{-};
"A'3";"B'3" **@{-};
"A'2";"C'1" **@{-};
"A'3";"C'2" **@{-};
"A'4";"C'3" **@{-};
"A1";"D"   **@{-};
"A'1";"D"  **@{-};
"A'4";"E"  **@{-};
"A4";"F"  **@{-};
%-----------------------
(-20,-45)="MA" *{ \xy
(-11, 0)="Et"  *{A=};
( 6, 12)="A"  *{};
(-6, 12)="B"  *{};
( 6,-12)="C"  *{};
(-6,-12)="D"  *{};
( 0,  0)="O"  *{};
( 0, 12)="AB"  *{};
( 0,-12)="CD"  *{};
( 6,  0)="AC"  *{};
(-6,  0)="BD"  *{};
( 6,  6)="AAC"  *{};
(-6,  6)="BBD"  *{};
( 6, -6)="ACC"  *{};
(-6, -6)="BDD"  *{};
(-3,  9)="Et0"  *{\scriptstyle I_{k+1}};
( 3,  3)="Et1"  *{\scriptstyle I_{k+1}};
(-3, -3)="Et2"  *{\scriptstyle I_{k}^{\rightarrow}};
( 3, -9)="Et3"  *{\scriptstyle I_{k+1}};
(-3,  5)="reng"  *{\renglon{0}{\cdots}{1}};
%-----------------------
"A";"B" **@{-};
"C";"D" **@{-};
"A";"C" **@{-};
"B";"D" **@{-};
"AC";"BD" **@{-};
"AB";"CD" **@{-};
"ACC";"BDD" **@{-};
"AAC";"BBD" **@{-};
 \endxy};
%-----------------------
(10,-45)="MB" *{ \xy
(-11, 0)="Et"  *{B=};
( 6, 12)="A"  *{};
(-6, 12)="B"  *{};
( 6,-12)="C"  *{};
(-6,-12)="D"  *{};
( 0,  0)="O"  *{};
( 0, 12)="AB"  *{};
( 0,-12)="CD"  *{};
( 6,  0)="AC"  *{};
(-6,  0)="BD"  *{};
( 6, -6)="ACC"  *{};
(-6, -6)="BDD"  *{};
(0 ,  6)="Et1"  *{I_{\ell+1}^{\rightarrow \uparrow}};
(-3, -3)="Et2"  *{\scriptstyle I_{k}^{\rightarrow}};
( 3, -9)="Et3"  *{\scriptstyle I_{k+1}};
%-----------------------
"A";"B" **@{-};
"C";"D" **@{-};
"A";"C" **@{-};
"B";"D" **@{-};
"AC";"BD" **@{-};
"O";"CD" **@{-};
"ACC";"BDD" **@{-};
 \endxy};
\endxy
\end{equation}

%------------------
Rank two representations ${}_{j,i}N^{(2)}_{3,2,2}[\ell]^2_2=F^Y({}_{j,i}\widetilde{N}^{(2)}_{3,2,2}[\ell]^2_2)$ ($\ell=2k$).
\begin{equation} \label{(C)F:Dn10}
 \xy
(0,12)="Et"  *{\scriptstyle \mathcal{C}({}_{j,i}N^{(2)}_{3,2,2}[6]_2^2)};
(-33,4)="Col"  *{\xymatrix@!0@L=6pt@R=10pt@C=35pt{
&   k^{\ell+2} \ar[dd]_*+<1em>{^{A}} \\
& & k^{\ell+2} \ar[ld]^(.3){B} \\
& k^{2\ell+4} \ar[dd]_-{=} \\ \\
& k^{2\ell+4} \ar@{}[dd]|-{\cdots} \\ \\
& k^{2\ell+2} \ar[dd]_-{=} \\ \\
& k^{2\ell+2} \ar[dd]^-{\left[ \begin{smallmatrix} I_{\ell+1}&0 \end{smallmatrix} \right]} 
\ar[dl]_(.7){\left[ \begin{smallmatrix} 0&I_{\ell+1} \end{smallmatrix} \right]} \\
  k^{\ell+1} \\
& k^{\ell+1}}};
%---------------------------
(-69,-9.5)="Cuad" *{ \xy
( 8, 7.5)="A"  *{};
(-8, 7.5)="B"  *{};
( 8,-7.5)="C"  *{};
(-8,-7.5)="D"  *{};
%-----------------------
"A";"B" **@{-};
"C";"D" **@{-};
"A";"C" **@{-};
"B";"D" **@{-};
 \endxy};
(-33,-3)="Col2"  *{\xymatrix@!0@L=6pt@R=10pt@C=35pt{
k^{2\ell+4} \ar[dd]_*+<1em>{^{I_{2\ell+3}^{\leftarrow}}} \\ \\
k^{2\ell+3}}};
(-66,-3)="dots1"  *{\scriptstyle \cdots};
(-66,-16)="dots2"  *{\scriptstyle \cdots};
%---------------------------
(-41,-13.5)="Cuad'" *{ \xy
( 8, 7.5)="A"  *{};
(-8, 7.5)="B"  *{};
( 8,-7.5)="C"  *{};
(-8,-7.5)="D"  *{};
%-----------------------
"A";"B" **@{-};
"C";"D" **@{-};
"A";"C" **@{-};
"B";"D" **@{-};
 \endxy};
(-22,-5)="Col'2"  *{\xymatrix@!0@L=6pt@R=10pt@C=35pt{
k^{2\ell+3} \ar[dd]^*+<1em>{^{I_{2\ell+2}^{\rightarrow}}} \\ \\
k^{2\ell+2}}};
(-44,-7)="dots'1"  *{\scriptstyle \cdots};
(-44,-20)="dots'2"  *{\scriptstyle \cdots};
%---------------------------
( -5,8)="A1"  *{};
(-13,8)="A2"  *{};
(-21,8)="A3"  *{};
(-29,8)="A4"  *{};
( -9,-12)="D1" *{\CoefMaxDnO};
(-17,-12)="D2" *{\CoefMaxDnO};
(-25,-12)="D3" *{\CoefMaxDnO};
(-10,0)="C1" *{};
(-18,0)="C2" *{};
(-26,0)="C3" *{};
( -8,0)="B1" *{};
(-16,0)="B2" *{};
(-24,0)="B3" *{};
( 5,8)="A'1"  *{};
(13,8)="A'2"  *{};
(21,8)="A'3"  *{};
(29,8)="A'4"  *{};
( 9,-12)="D'1" *{\CoefMaxDnO};
(17,-12)="D'2" *{\CoefMaxDnO};
(25,-12)="D'3" *{\CoefMaxDnO};
(10,0)="C'1" *{};
(18,0)="C'2" *{};
(26,0)="C'3" *{};
( 8,0)="B'1" *{};
(16,0)="B'2" *{};
(24,0)="B'3" *{};
(-2,0)="D" *{};
(32, 0)="E" *{};
(-32,0)="F" *{};
( 3, 0)="G" *{};
%(12, -12)="D3" *{\CoefMaxDnO};
(-2,-12)="D4" *{\CoefYDDnPrime};
(32, -12)="D5" *{\CoefYUDn};
(  3,  0)="D6" *{\CoefAlaU{j}};
(-32,  0)="D7" *{\CoefAlaU{i}};
%-----------------------
"A1";"B1" **@{-};
"A2";"B2" **@{-};
"A3";"B3" **@{-};
"A2";"C1" **@{-};
"A3";"C2" **@{-};
"A4";"C3" **@{-};
"A'1";"B'1" **@{-};
"A'2";"B'2" **@{-};
"A'3";"B'3" **@{-};
"A'2";"C'1" **@{-};
"A'3";"C'2" **@{-};
"A'4";"C'3" **@{-};
%"A5";"C4" **@{-};
"A1";"D"   **@{-};
"A'1";"D"  **\crv{(-1,7)};
"A'4";"E"  **@{-};
"A4";"F"  **@{-};
"A'1";"G"  **@{-};
%-----------------------
(-20,-45)="MA" *{  \xy
(-11, 0)="Et"  *{A=};
( 6, 12)="A"  *{};
(-6, 12)="B"  *{};
( 6,-14)="C"  *{};
(-6,-14)="D"  *{};
( 6,-12)="C'"  *{};
(-6,-12)="D'"  *{};
( 0,  0)="O"  *{};
( 0, 12)="AB"  *{};
( 0,-14)="CD"  *{};
( 6,  0)="AC"  *{};
(-6,  0)="BD"  *{};
( 6,  6)="AAC"  *{};
(-6,  6)="BBD"  *{};
( 6, -6)="ACC"  *{};
(-6, -6)="BDD"  *{};
(-3,  9)="Et0"  *{\scriptstyle I_{k+1}};
( 3,  3)="Et1"  *{\scriptstyle I_{k+1}};
(-3, -3)="Et2"  *{\scriptstyle I_{k}^{\rightarrow}};
( 3, -9)="Et3"  *{\scriptstyle I_{k+1}};
(-3,  5)="reng1"  *{\renglon{0}{\cdots}{1}};
( 3,-13.1)="reng2"  *{\renglon{1}{\cdots}{0}};
%-----------------------
"A";"B" **@{-};
"C";"D" **@{-};
"C'";"D'" **@{-};
"A";"C" **@{-};
"B";"D" **@{-};
"AC";"BD" **@{-};
"AB";"CD" **@{-};
"ACC";"BDD" **@{-};
"AAC";"BBD" **@{-};
 \endxy};
%-----------------------
(10,-45)="MB" *{ \xy
(-11, 0)="A"  *{B=};
( 6, 12)="A"  *{};
(-6, 12)="B"  *{};
( 6,-14)="C"  *{};
(-6,-14)="D"  *{};
( 0,  0)="O"  *{};
( 0, 12)="AB"  *{};
( 0,-14)="CD"  *{};
( 6,  0)="AC"  *{};
(-6,  0)="BD"  *{};
( 6, -6)="ACC"  *{};
(-6, -6)="BDD"  *{};
(0 ,  6)="Et1"  *{I_{\ell+1}^{\rightarrow \uparrow}};
(-3, -3)="Et2"  *{\scriptstyle I_{k}^{\rightarrow}};
( 3, -9)="Et3"  *{\scriptstyle I_{k+1}^{\downarrow}};
%-----------------------
"A";"B" **@{-};
"C";"D" **@{-};
"A";"C" **@{-};
"B";"D" **@{-};
"AC";"BD" **@{-};
"O";"CD" **@{-};
"ACC";"BDD" **@{-};
 \endxy};
\endxy
\end{equation}

%------------------
Rank two representations ${}_iM^{(2)}_{2,2,2}[\ell]^2_2=F^Y({}_i\widetilde{M}^{(2)}_{2,2,2}[\ell]^2_2)$ ($\ell=2k+1$).
\begin{equation} \label{(C)F:Dn11}
 \xy
(0,12)="Et"  *{\scriptstyle \mathcal{C}({}_{i}M^{(2)}_{2,2,2}[7]_2^2)};
(-36,4)="Col"  *{\xymatrix@!0@L=6pt@R=10pt@C=35pt{
&   k^{\ell+2} \ar[dd]_*+<1em>{^{A}} \\
& & k^{\ell+2} \ar[ld]^(.3){B} \\
& k^{2\ell+3} \ar[dd]_-{=} \\ \\
& k^{2\ell+3} \ar@{}[dd]|-{\cdots} \\ \\
& k^{2\ell+2} \ar[dd]_-{=} \\ \\
& k^{2\ell+2} \ar[dd]^-{\left[ \begin{smallmatrix} I_{\ell+1}&0 \end{smallmatrix} \right]} 
\ar[dl]_(.7){\left[ \begin{smallmatrix} 0&I_{\ell+1} \end{smallmatrix} \right]} \\
  k^{\ell+1} \\
& k^{\ell+1}}};
%---------------------------
(-75,-9.5)="Cuad" *{ \xy
( 8, 7.5)="A"  *{};
(-8, 7.5)="B"  *{};
( 8,-7.5)="C"  *{};
(-8,-7.5)="D"  *{};
%-----------------------
"A";"B" **@{-};
"C";"D" **@{-};
"A";"C" **@{-};
"B";"D" **@{-};
 \endxy};
(-36,-3)="Col2"  *{\xymatrix@!0@L=6pt@R=10pt@C=35pt{
k^{2\ell+3} \ar[dd]_*+<1em>{^{I_{2\ell+2}^{\leftarrow}}} \\ \\
k^{2\ell+2}}};
(-72,-3)="dots1"  *{\scriptstyle \cdots};
(-72,-16)="dots2"  *{\scriptstyle \cdots};
%---------------------------
( -5,8)="A1"  *{};
(-13,8)="A2"  *{};
(-21,8)="A3"  *{};
(-29,8)="A4"  *{};
(-37,8)="A5"  *{};
( -9,-12)="D1" *{\CoefMaxDnO};
(-17,-12)="D2" *{\CoefMaxDnO};
(-25,-12)="D3" *{\CoefMaxDnO};
(-33,-12)="D4" *{\CoefMaxDnO};
(-10,0)="C1" *{};
(-18,0)="C2" *{};
(-26,0)="C3" *{};
(-34,0)="C4" *{};
( -8,0)="B1" *{};
(-16,0)="B2" *{};
(-24,0)="B3" *{};
(-32,0)="B4" *{};
( 5,8)="A'1"  *{};
(13,8)="A'2"  *{};
(21,8)="A'3"  *{};
(29,8)="A'4"  *{};
( 9,-12)="D'1" *{\CoefMaxDnO};
(17,-12)="D'2" *{\CoefMaxDnO};
(25,-12)="D'3" *{\CoefMaxDnO};
(10,0)="C'1" *{};
(18,0)="C'2" *{};
(26,0)="C'3" *{};
( 8,0)="B'1" *{};
(16,0)="B'2" *{};
(24,0)="B'3" *{};
(-3,0)="D" *{};
(32, 0)="E" *{};
(-40,0)="F" *{};
( 2, 0)="G" *{};
(32, -12)="D5" *{\CoefYUDn};
( 2, -12)="D6" *{\CoefYDDnPrime};
(-40,  0)="D7" *{\CoefAlaU{i}};
%-----------------------
"A1";"B1" **@{-};
"A2";"B2" **@{-};
"A3";"B3" **@{-};
"A4";"B4" **@{-};
"A2";"C1" **@{-};
"A3";"C2" **@{-};
"A4";"C3" **@{-};
"A5";"C4" **@{-};
"A'1";"B'1" **@{-};
"A'2";"B'2" **@{-};
"A'3";"B'3" **@{-};
"A'2";"C'1" **@{-};
"A'3";"C'2" **@{-};
"A'4";"C'3" **@{-};
"A'4";"E"  **@{-};
"A5";"F"  **@{-};
"A'1";"F"  **\crv{(-45,18)};
"A'1";"G"  **@{-};
%-----------------------
(-20,-45)="MA" *{  \xy
(-11, 0)="Et"  *{A=};
( 6, 12)="A"  *{};
(-6, 12)="B"  *{};
( 6,-12)="C"  *{};
(-6,-12)="D"  *{};
( 6,-12)="C'"  *{};
(-6,-12)="D'"  *{};
( 0,  0)="O"  *{};
( 0, 12)="AB"  *{};
( 0,-12)="CD"  *{};
( 6,  0)="AC"  *{};
(-6,  0)="BD"  *{};
( 6,  6)="AAC"  *{};
(-6,  6)="BBD"  *{};
( 6, -6)="ACC"  *{};
(-6, -6)="BDD"  *{};
(-3,  9)="Et0"  *{\scriptstyle I_{k+2}};
( 3,  3)="Et1"  *{\scriptstyle I_{k+1}};
(-3, -3)="Et2"  *{\scriptstyle I_{k+1}^{\rightarrow}};
( 3, -9)="Et3"  *{\scriptstyle I_{k+1}};
( 3, 11)="reng1"  *{\renglon{1}{\cdots}{0}};
%-----------------------
"A";"B" **@{-};
"C";"D" **@{-};
"A";"C" **@{-};
"B";"D" **@{-};
"AC";"BD" **@{-};
"AB";"CD" **@{-};
"ACC";"BDD" **@{-};
"AAC";"BBD" **@{-};
 \endxy};
%-----------------------
(10,-45)="MB" *{ \xy
(-11, 0)="Et"  *{B=};
( 6, 12)="A"  *{};
(-6, 12)="B"  *{};
( 6,-12)="C"  *{};
(-6,-12)="D"  *{};
( 0,  0)="O"  *{};
( 0, 12)="AB"  *{};
( 0,-12)="CD"  *{};
( 6,  0)="AC"  *{};
(-6,  0)="BD"  *{};
( 6, -6)="ACC"  *{};
(-6, -6)="BDD"  *{};
(0 ,  6)="Et1"  *{I_{\ell+1}^{\rightarrow \uparrow}};
(-3, -3)="Et2"  *{\scriptstyle I_{k+1}^{\rightarrow}};
( 3, -9)="Et3"  *{\scriptstyle I_{k+1}};
%-----------------------
"A";"B" **@{-};
"C";"D" **@{-};
"A";"C" **@{-};
"B";"D" **@{-};
"AC";"BD" **@{-};
"O";"CD" **@{-};
"ACC";"BDD" **@{-};
 \endxy};
\endxy
\end{equation}

%------------------
Rank two representations ${}_{j,i}M^{(2)}_{3,2,2}[\ell]^2_2=F^Y({}_{j,i}\widetilde{M}^{(2)}_{3,2,2}[\ell]^2_2)$ ($\ell=2k+1$).
\begin{equation} \label{(C)F:Dn12}
 \xy
(0,12)="Et"  *{\scriptstyle \mathcal{C}({}_{j,i}M^{(2)}_{3,2,2}[7]_2^2)};
(-36,4)="Col"  *{\xymatrix@!0@L=6pt@R=10pt@C=35pt{
&   k^{\ell+2} \ar[dd]_*+<1em>{^{A}} \\
& & k^{\ell+2} \ar[ld]^(.3){B} \\
& k^{2\ell+4} \ar[dd]_-{=} \\ \\
& k^{2\ell+4} \ar@{}[dd]|-{\cdots} \\ \\
& k^{2\ell+2} \ar[dd]_-{=} \\ \\
& k^{2\ell+2} \ar[dd]^-{\left[ \begin{smallmatrix} I_{\ell+1}&0 \end{smallmatrix} \right]} 
\ar[dl]_(.7){\left[ \begin{smallmatrix} 0&I_{\ell+1} \end{smallmatrix} \right]} \\
  k^{\ell+1} \\
& k^{\ell+1}}};
%---------------------------
(-75,-9.5)="Cuad" *{ \xy
( 8, 7.5)="A"  *{};
(-8, 7.5)="B"  *{};
( 8,-7.5)="C"  *{};
(-8,-7.5)="D"  *{};
%-----------------------
"A";"B" **@{-};
"C";"D" **@{-};
"A";"C" **@{-};
"B";"D" **@{-};
 \endxy};
(-36,-3)="Col2"  *{\xymatrix@!0@L=6pt@R=10pt@C=35pt{
k^{2\ell+4} \ar[dd]_*+<1em>{^{I_{2\ell+3}^{\leftarrow}}} \\ \\
k^{2\ell+3}}};
(-72,-3)="dots1"  *{\scriptstyle \cdots};
(-72,-16)="dots2"  *{\scriptstyle \cdots};
%---------------------------
(-47,-13.5)="Cuad'" *{ \xy
( 8, 7.5)="A"  *{};
(-8, 7.5)="B"  *{};
( 8,-7.5)="C"  *{};
(-8,-7.5)="D"  *{};
%-----------------------
"A";"B" **@{-};
"C";"D" **@{-};
"A";"C" **@{-};
"B";"D" **@{-};
 \endxy};
(-25,-5)="Col'2"  *{\xymatrix@!0@L=6pt@R=10pt@C=35pt{
k^{2\ell+3} \ar[dd]^*+<1em>{^{I_{2\ell+2}^{\rightarrow}}} \\ \\
k^{2\ell+2}}};
(-50,-7)="dots'1"  *{\scriptstyle \cdots};
(-50,-20)="dots'2"  *{\scriptstyle \cdots};
%---------------------------
( -5,8)="A1"  *{};
(-13,8)="A2"  *{};
(-21,8)="A3"  *{};
(-29,8)="A4"  *{};
(-37,8)="A5"  *{};
( -9,-12)="D1" *{\CoefMaxDnO};
(-17,-12)="D2" *{\CoefMaxDnO};
(-25,-12)="D3" *{\CoefMaxDnO};
(-33,-12)="D4" *{\CoefMaxDnO};
(-10,0)="C1" *{};
(-18,0)="C2" *{};
(-26,0)="C3" *{};
(-34,0)="C4" *{};
( -8,0)="B1" *{};
(-16,0)="B2" *{};
(-24,0)="B3" *{};
(-32,0)="B4" *{};
( 5,8)="A'1"  *{};
(13,8)="A'2"  *{};
(21,8)="A'3"  *{};
(29,8)="A'4"  *{};
( 9,-12)="D'1" *{\CoefMaxDnO};
(17,-12)="D'2" *{\CoefMaxDnO};
(25,-12)="D'3" *{\CoefMaxDnO};
(10,0)="C'1" *{};
(18,0)="C'2" *{};
(26,0)="C'3" *{};
( 8,0)="B'1" *{};
(16,0)="B'2" *{};
(24,0)="B'3" *{};
(-3,0)="D" *{};
(32, 0)="E" *{};
(-40,0)="F" *{};
( 2, 0)="G" *{};
%(12, -12)="D3" *{\CoefMaxDnO};
(-3,0)="D4" *{\CoefAlaU{j}};
(32, -12)="D5" *{\CoefYUDn};
( 2, -12)="D6" *{\CoefYDDnPrime};
(-40,  0)="D7" *{\CoefAlaU{i}};
%-----------------------
"A1";"B1" **@{-};
"A2";"B2" **@{-};
"A3";"B3" **@{-};
"A4";"B4" **@{-};
"A2";"C1" **@{-};
"A3";"C2" **@{-};
"A4";"C3" **@{-};
"A5";"C4" **@{-};
"A'1";"B'1" **@{-};
"A'2";"B'2" **@{-};
"A'3";"B'3" **@{-};
"A'2";"C'1" **@{-};
"A'3";"C'2" **@{-};
"A'4";"C'3" **@{-};
"A'1";"D"  **\crv{(-1,7)};
"A'4";"E"  **@{-};
"A5";"F"  **@{-};
"A'1";"F"  **\crv{(-45,18)};
"A'1";"G"  **@{-};
%-----------------------
(-20,-45)="MA" *{  \xy
(-11, 0)="Et"  *{A=};
( 6, 12)="A"  *{};
(-6, 12)="B"  *{};
( 6,-14)="C"  *{};
(-6,-14)="D"  *{};
( 6,-12)="C'"  *{};
(-6,-12)="D'"  *{};
( 0,  0)="O"  *{};
( 0, 12)="AB"  *{};
( 0,-14)="CD"  *{};
( 6,  0)="AC"  *{};
(-6,  0)="BD"  *{};
( 6,  6)="AAC"  *{};
(-6,  6)="BBD"  *{};
( 6, -6)="ACC"  *{};
(-6, -6)="BDD"  *{};
(-3,  9)="Et0"  *{\scriptstyle I_{k+2}};
( 3,  3)="Et1"  *{\scriptstyle I_{k+1}};
(-3, -3)="Et2"  *{\scriptstyle I_{k+1}^{\rightarrow}};
( 3, -9)="Et3"  *{\scriptstyle I_{k+1}};
( 3, 11)="reng1"  *{\renglon{1}{\cdots}{0}};
( 3,-13.1)="reng2"  *{\renglon{1}{\cdots}{0}};
%-----------------------
"A";"B" **@{-};
"C";"D" **@{-};
"C'";"D'" **@{-};
"A";"C" **@{-};
"B";"D" **@{-};
"AC";"BD" **@{-};
"AB";"CD" **@{-};
"ACC";"BDD" **@{-};
"AAC";"BBD" **@{-};
 \endxy};
%-----------------------
(10,-45)="MB" *{ \xy
(-11, 0)="A"  *{B=};
( 6, 12)="A"  *{};
(-6, 12)="B"  *{};
( 6,-14)="C"  *{};
(-6,-14)="D"  *{};
( 0,  0)="O"  *{};
( 0, 12)="AB"  *{};
( 0,-14)="CD"  *{};
( 6,  0)="AC"  *{};
(-6,  0)="BD"  *{};
( 6, -6)="ACC"  *{};
(-6, -6)="BDD"  *{};
(0 ,  6)="Et1"  *{I_{\ell+1}^{\rightarrow \uparrow}};
(-3, -3)="Et2"  *{\scriptstyle I_{k+1}^{\rightarrow}};
( 3, -9)="Et3"  *{\scriptstyle I_{k+1}^{\downarrow}};
%-----------------------
"A";"B" **@{-};
"C";"D" **@{-};
"A";"C" **@{-};
"B";"D" **@{-};
"AC";"BD" **@{-};
"O";"CD" **@{-};
"ACC";"BDD" **@{-};
 \endxy};
\endxy
\end{equation}
%%%%%%%%%%%%%%%%%%%%%%%%%%%%%%%%

\begin{proposicion} \label{(C)P:Dn}
With the orientation given to the quiver $\widetilde{\mathbf{D_n}}$, almost all indecomposable preinjective 
$\widetilde{\mathbf{D_n}}$-representations are isomorphic to one of the modules given in figures~(\ref{(C)F:Dn1})-(\ref{(C)F:Dn12}).
\end{proposicion}
\bproof
By theorem~\ref{(DE)T:finito}($b$) and the choice of the reduction module $Y$, almost all indecomposable preinjective
module is isomorphic to an object in the image of the functor $F^Y$. 

By proposition~\ref{(DE)P:rango6}($a$), if $\widetilde{M}$ is an indecomposable $\mathcal{A}^Y$-module such that
$F^Y(\widetilde{M})$ is preinjective, then the dimension vector $\vdim \widetilde{M}$ is a rank one or rank two positive root.
By lemma~\ref{(DE)L:tablaComp}, $\vdim \widetilde{M}$ is one of the roots in table~\ref{(DE)T:raicesMenor}, excepting those
roots in the family $\mathcal{R}^{(2)}_{3,3,2}$, for all quivers
$\mathbf{D}_n$ are of type $n_1,2,2$. Then almost all indecomposable preinjective $\widetilde{\mathbf{D_n}}$-module 
is isomorphic to a representation of the form $F^Y(\widetilde{M})$ where $\widetilde{M}$ is some of the modules 
given in lemmas~\ref{(DE)L:excRedRnkUno} and~\ref{(DE)L:excRedRnkDos}, which are precisely those used in  
figures~(\ref{(C)F:Dn1})-(\ref{(C)F:Dn12}).
\eproof

\section{Reduced ditalgebras for the cases $\widetilde{\mathbf{E_m}}$.} \label{(B)S:mayor}
%--------------------------------------------------------------------------------------
%--------------------------------------------------------------------------------------
In this last section we give the reduced tensor algebras corresponding to the cases $\widetilde{\mathbf{E_6}}$, 
$\widetilde{\mathbf{E_7}}$ and $\widetilde{\mathbf{E_8}}$, for an arbitrary orientation of the arrows
(assuming that the extension vertex $\omega$ is a source). The explicit forms of the differentials $\delta^x$ and 
$\delta^y$ depend on the orientations of the arrows, and can be computed in each particular case from the formulae
given in lemmas~\ref{(DE)L:AX} and~\ref{(DE)L:AY}, and the choice of representations of the Dynkin quivers $\mathbf{E_m}$
corresponding to the elements in $\mathcal{X}_0$ and $\mathcal{Y}_0$ (cf. figure~\ref{(C)T:AREm} for the subspace 
orientation case).

Reduced tensor algebras of the path algebra of $\widetilde{\mathbf{E_6}}$ with respect to the admissible modules $X$ (left) 
and $Y$ (right).
%------------------------------------------------E6
\begin{displaymath}
 \xymatrix@C=4pc@R=2.5pc{
 & & *+++[]{\bullet_{\omega}} \ar[ddd]^(.6){\xi'} \ar@<-1ex>[ddd]_(.6){\xi} 
\ar@/_.7pc/[ddll]|(.5){\xi^1_1} \ar@/_.5pc/[ddl]|(.5){\xi^1_2} 
\ar@/_.7pc/[dddl]|(.5){\xi^2_2} \ar@/_.7pc/[dddll]|(.5){\xi^2_1} \ar@/_.5pc/[ddddl]|(.5){\xi^3_1} \\
\\
\bullet \ar@{.>}[r]^(.4){\gamma^1_{1,2}} \ar@{.>}[rrd]^(.4){\gamma^1_{1,3}} & \bullet \ar@{.>}[rd]^(.5){\gamma^1_{2,3}} \\
\bullet \ar@{.>}[r]^(.5){\gamma^2_{1,2}} \ar@{.>}@/_.7pc/[rr]_(.4){\gamma^2_{1,3}} & \bullet \ar@{.>}[r]^(.4){\gamma^2_{2,3}} 
& \bullet_m \\ & \bullet \ar@{.>}[ru]_-{\gamma^3_{1,2}}
} \quad
 \xymatrix@C=4pc@R=2.5pc{
*+++[]{\bullet_{\omega}} \ar[ddd]^(.6){\theta} \ar@<-1ex>[ddd]_(.6){\theta'} 
\ar@/^.7pc/[ddrr]|(.5){\theta^1_3} \ar@/^.5pc/[ddr]|(.5){\theta^1_2} 
\ar@/^.7pc/[dddr]|(.5){\theta^2_2} \ar@/^.7pc/[dddrr]|(.5){\theta^2_3} \ar@/^.5pc/[ddddr]|(.5){\theta^3_2} \\
\\
& \bullet \ar@{.>}[r]^(.6){\varepsilon^1_{2,3}} & \bullet  \\
\bullet_m \ar@{.>}@/_.7pc/[rr]_(.6){\varepsilon^2_{1,3}} \ar@{.>}[rd]_-{\varepsilon^3_{1,2}} \ar@{.>}[ru]^(.5){\varepsilon^1_{1,2}} 
\ar@{.>}[rru]^(.6){\varepsilon^1_{1,3}} \ar@{.>}[r]^(.5){\varepsilon^2_{1,2}}
& \bullet \ar@{.>}[r]^(.6){\varepsilon^2_{2,3}} & \bullet   \\ 
& \bullet 
}
\end{displaymath}

Reduced tensor algebra of the path algebra of $\widetilde{\mathbf{E_7}}$ with respect to $X$.
%------------------------------------------------E7X
\begin{displaymath}
 \xymatrix@C=4pc@R=2.5pc{
& & & *+++[]{\bullet_{\omega}} \ar[ddd]^(.6){\xi'} \ar@<-1ex>[ddd]_(.6){\xi} 
\ar@/_.7pc/[ddlll]|(.5){\xi^1_1} \ar@/_.7pc/[ddll]|(.5){\xi^1_2} \ar@/_.5pc/[ddl]|(.5){\xi^1_3} 
\ar@/_.7pc/[dddl]|(.5){\xi^2_2} \ar@/_.7pc/[dddll]|(.5){\xi^2_1} \ar@/_.5pc/[ddddl]|(.5){\xi^3_1} \\
\\
\bullet \ar@{.>}[r]^(.5){\gamma^1_{1,2}} \ar@{.>}@/^1.5pc/[rr]^(.5){\gamma^1_{1,3}} \ar@{.>}[rrrd]_(.3){\gamma^1_{1,4}} 
& \bullet \ar@{.>}[r]^(.4){\gamma^1_{2,3}} \ar@{.>}[rrd]^(.4){\gamma^1_{2,4}} & \bullet \ar@{.>}[rd]^(.5){\gamma^1_{3,4}} \\
& \bullet \ar@{.>}[r]_(.5){\gamma^2_{1,2}} \ar@{.>}@/_.7pc/[rr]_(.4){\gamma^2_{1,3}} & \bullet \ar@{.>}[r]_(.4){\gamma^2_{2,3}} 
& *+++[]{\bullet_m} \\ 
& & \bullet \ar@{.>}[ru]_-{\gamma^3_{1,2}}
}
\end{displaymath}

Reduced tensor algebra of the path algebra of $\widetilde{\mathbf{E_7}}$ with respect to $Y$.
%------------------------------------------------E7Y
\begin{displaymath}
 \xymatrix@C=4pc@R=2.5pc{
*+++[]{\bullet_{\omega}} \ar[ddd]^(.6){\theta} \ar@<-1ex>[ddd]_(.6){\theta'} 
\ar@/^.7pc/[ddrrr]|(.5){\theta^1_4} \ar@/^.7pc/[ddrr]|(.5){\theta^1_3} \ar@/^.5pc/[ddr]|(.5){\theta^1_2} 
\ar@/^.7pc/[dddr]|(.5){\theta^2_2} \ar@/^.7pc/[dddrr]|(.5){\theta^2_3} \ar@/^.5pc/[ddddr]|(.5){\theta^3_2} \\
\\
& \bullet \ar@{.>}@/^1.5pc/[rr]^(.5){\varepsilon^1_{2,4}} \ar@{.>}[r]^(.6){\varepsilon^1_{2,3}} & \bullet \ar@{.>}[r]^(.5){\varepsilon^1_{3,4}} & \bullet  \\
\bullet_m \ar@{.>}[rrru]_(.7){\varepsilon^2_{1,4}} \ar@{.>}@/_.7pc/[rr]_(.6){\varepsilon^2_{1,3}} \ar@{.>}[rd]_-{\varepsilon^3_{1,2}} \ar@{.>}[ru]^(.5){\varepsilon^1_{1,2}} 
\ar@{.>}[rru]^(.6){\varepsilon^1_{1,3}} \ar@{.>}[r]_(.5){\varepsilon^2_{1,2}}
& \bullet \ar@{.>}[r]_(.6){\varepsilon^2_{2,3}} & \bullet   \\ 
& \bullet 
}
\end{displaymath}

Reduced tensor algebra of the path algebra of $\widetilde{\mathbf{E_8}}$ with respect to $X$.
%------------------------------------------------E8X
\begin{displaymath}
 \xymatrix@C=4pc@R=2.5pc{
& & & & *+++[]{\bullet_{\omega}} \ar[ddd]^(.6){\xi'} \ar@<-1ex>[ddd]_(.6){\xi} 
\ar@/_.7pc/[ddllll]|(.5){\xi^1_1} \ar@/_.7pc/[ddlll]|(.5){\xi^1_2} \ar@/_.7pc/[ddll]|(.5){\xi^1_3} \ar@/_.5pc/[ddl]|(.5){\xi^1_4} 
\ar@/_.7pc/[dddl]|(.5){\xi^2_2} \ar@/_.7pc/[dddll]|(.5){\xi^2_1} \ar@/_.5pc/[ddddl]|(.5){\xi^3_1} \\
\\
\bullet \ar@{.>}[rrrrd]_(.3){\gamma^1_{1,5}} \ar@{.>}@/^3pc/[rrr]^(.3){\gamma^1_{1,4}} \ar@{.>}@/^1.5pc/[rr]^(.6){\gamma^1_{1,3}} \ar@{.>}[r]^(.5){\gamma^1_{1,2}} 
& \bullet \ar@{.>}[r]^(.5){\gamma^1_{2,3}} \ar@{.>}@/^1.5pc/[rr]^(.5){\gamma^1_{2,4}} \ar@{.>}[rrrd]_(.3){\gamma^1_{2,5}} 
& \bullet \ar@{.>}[r]^(.4){\gamma^1_{3,4}} \ar@{.>}[rrd]^(.4){\gamma^1_{3,5}} & \bullet \ar@{.>}[rd]^(.5){\gamma^1_{4,5}} \\
& & \bullet \ar@{.>}[r]_(.5){\gamma^2_{1,2}} \ar@{.>}@/_.7pc/[rr]_(.4){\gamma^2_{1,3}} & \bullet \ar@{.>}[r]_(.4){\gamma^2_{2,3}} 
& *+++[]{\bullet_m} \\ 
& & & \bullet \ar@{.>}[ru]_-{\gamma^3_{1,2}}
}
\end{displaymath}

Reduced tensor algebra of the path algebra of $\widetilde{\mathbf{E_8}}$ with respect to $Y$.
%------------------------------------------------E8Y
\begin{displaymath}
 \xymatrix@C=4pc@R=2.5pc{
*+++[]{\bullet_{\omega}} \ar[ddd]^(.6){\theta} \ar@<-1ex>[ddd]_(.6){\theta'} 
\ar@/^.7pc/[ddrrrr]|(.5){\theta^1_5} \ar@/^.7pc/[ddrrr]|(.5){\theta^1_4} 
\ar@/^.7pc/[ddrr]|(.5){\theta^1_3} \ar@/^.5pc/[ddr]|(.5){\theta^1_2} 
\ar@/^.7pc/[dddr]|(.5){\theta^2_2} \ar@/^.7pc/[dddrr]|(.5){\theta^2_3} \ar@/^.5pc/[ddddr]|(.5){\theta^3_2} \\
\\
& \bullet \ar@{.>}@/^3pc/[rrr]^(.7){\varepsilon^1_{2,5}} \ar@{.>}@/^1.5pc/[rr]^(.5){\varepsilon^1_{2,4}} \ar@{.>}[r]^(.6){\varepsilon^1_{2,3}} 
& \bullet \ar@{.>}@/^1.5pc/[rr]^(.4){\varepsilon^1_{3,5}} \ar@{.>}[r]^(.5){\varepsilon^1_{3,4}} & \bullet \ar@{.>}[r]^(.5){\varepsilon^1_{4,5}} & \bullet \\
\bullet_m \ar@{.>}[rrrru]_(.7){\varepsilon^1_{1,5}} \ar@{.>}[rrru]_(.7){\varepsilon^1_{1,4}} 
\ar@{.>}@/_.7pc/[rr]_(.6){\varepsilon^2_{1,3}} \ar@{.>}[rd]_-{\varepsilon^3_{1,2}} \ar@{.>}[ru]^(.5){\varepsilon^1_{1,2}} 
\ar@{.>}[rru]^(.6){\varepsilon^1_{1,3}} \ar@{.>}[r]_(.5){\varepsilon^2_{1,2}}
& \bullet \ar@{.>}[r]_(.6){\varepsilon^2_{2,3}} & \bullet   \\ 
& \bullet 
}
\end{displaymath}

%------------------
\begin{center}
\begin{sidewaysfigure}
$\xymatrix@!0@C=18pt@R=24pt{
& & {} \ar[rd] & & {\vctEse{1}{1}{0}{0}{0}{1}} \ar[rd] & & {} \ar[rd] & & {\vctEse{1}{0}{1}{0}{1}{1}} \ar[rd] & & {} \ar[rd] & & {} \\
& {} \ar[ru] \ar[rd] & & {} \ar[ru] \ar[rd] & & {\vctEse{2}{2}{1}{1}{0}{1}} \ar[ru] \ar[rd] 
& & {\vctEse{2}{1}{2}{1}{1}{1}} \ar[ru] \ar[rd] & & {} \ar[ru] \ar[rd] & & {} \ar[ru] \\
{} \ar[ru] \ar[rd] \ar[r] & {} \ar[r] & {} \ar[ru] \ar[rd] \ar[r] & {} \ar[r] & {} \ar[ru] \ar[rd] \ar[r] 
& {\vctEse{2}{1}{1}{1}{1}{1}} \ar[r] & {\vctEse{3}{2}{2}{1}{1}{2}} \ar[ru] \ar[rd] \ar[r] 
& {\vctEse{1}{1}{1}{0}{0}{1}} \ar[r] & {} \ar[ru] \ar[rd] \ar[r] & {} \ar[r] & {} \ar[ru] \ar[rd] \ar[r] & {} \\
& {} \ar[ru] \ar[rd] & & {} \ar[ru] \ar[rd] & & {\vctEse{2}{1}{2}{0}{1}{1}} \ar[ru] \ar[rd] 
& & {\vctEse{2}{2}{1}{1}{1}{1}} \ar[ru] \ar[rd] & & {} \ar[ru] \ar[rd] & & {} \ar[rd] \\
& & {} \ar[ru] & & {\vctEse{1}{0}{1}{0}{0}{1}} \ar[ru] & & {} \ar[ru] & & {\vctEse{1}{1}{0}{1}{0}{1}} \ar[ru] & & {} \ar[ru] & & {}
}\xymatrix@!0@C=19.5pt@R=24pt{
& & {} \ar[rd] & & {} \ar[rd] & & {} \ar[rd] & & {\vctEsi{2}{1}{1}{1}{1}{1}{1}} \ar[rd] & & {} \ar[rd] 
& & {\vctEsi{1}{1}{1}{1}{1}{0}{0}} \ar[rd] & & {} \ar[rd] & & {} \ar[rd] & & {} \\
& {} \ar[ru] \ar[rd] & & {} \ar[ru] \ar[rd] & & {} \ar[ru] \ar[rd] 
& & {} \ar[ru] \ar[rd] & & {\vctEsi{3}{2}{2}{1}{1}{2}{1}} \ar[ru] \ar[rd] & & {\vctEsi{2}{2}{2}{1}{1}{1}{0}} \ar[ru] \ar[rd] 
& & {} \ar[ru] \ar[rd] & & {} \ar[rd] \ar[ru] & & {} \ar[ru] \\
{} \ar[ru] \ar[rd] \ar[r] & {} \ar[r] & {} \ar[ru] \ar[rd] \ar[r] 
& {} \ar[r] & {} \ar[ru] \ar[rd] \ar[r] & {} \ar[r] 
& {} \ar[ru] \ar[rd] \ar[r] & {} \ar[r] & {} \ar[ru] \ar[rd] \ar[r] 
& {\vctEsi{2}{2}{1}{1}{1}{1}{0}} \ar[r] & {\vctEsi{4}{3}{3}{2}{2}{2}{1}} 
\ar[ru] \ar[rd] \ar[r] & {\vctEsi{2}{1}{2}{1}{1}{1}{1}} \ar[r] 
& {} \ar[ru] \ar[rd] \ar[r] & {} \ar[r] & {} \ar[ru] \ar[rd] \ar[r] 
& {} \ar[r] & {} \ar[ru] \ar[rd] \ar[r] & {} \\
& {} \ar[ru] \ar[rd] & & {} \ar[ru] \ar[rd] & & {} \ar[ru] \ar[rd]
& & {} \ar[ru] \ar[rd] & & {\vctEsi{3}{2}{3}{1}{2}{1}{1}} \ar[ru] \ar[rd] & & {\vctEsi{3}{2}{2}{1}{2}{2}{1}} \ar[ru] \ar[rd] 
& & {} \ar[ru] \ar[rd] & & {} \ar[ru] \ar[rd] & & {} \ar[rd] \\ 
& & {} \ar[rd] \ar[ru] & & {} \ar[rd] \ar[ru] & & {} \ar[rd] \ar[ru] 
& & {\vctEsi{2}{1}{2}{1}{1}{1}{0}} \ar[rd] \ar[ru] & & {} \ar[rd] \ar[ru] & & {\vctEsi{2}{2}{1}{1}{1}{1}{1}} \ar[ru] \ar[rd] 
& & {} \ar[ru] \ar[rd] & & {} \ar[ru] \ar[rd] & & {} \ar[rd]  \\ 
& & & {} \ar[ru] & & {} \ar[ru] & & {\vctEsi{1}{1}{1}{1}{0}{0}{0}} \ar[ru] & & {} \ar[ru] 
& & {} \ar[ru] & & {\vctEsi{1}{1}{0}{1}{0}{1}{0}} \ar[ru] & & {} \ar[ru] & & {} \ar[ru] & & {}
}$
\vspace{1cm}

$\xymatrix@!0@C=18.5pt@R=24pt{
& & {} \ar[rd] & & {} \ar[rd] & & {} \ar[rd] & & {} \ar[rd] & & {} \ar[rd] 
& & {} \ar[rd] & & {} \ar[rd] & & {\vctEo{2}{1}{2}{0}{2}{1}{1}{1}} \ar[rd] 
& & {} \ar[rd] & & {\vctEo{2}{1}{2}{1}{1}{1}{1}{1}} \ar[rd] 
& & {} \ar[rd] & & {} \ar[rd] & & {} \ar[rd] & & {} \ar[rd] & & {} \\
& {} \ar[ru] \ar[rd] & & {} \ar[ru] \ar[rd] & & {} \ar[ru] \ar[rd] 
& & {} \ar[ru] \ar[rd] & & {} \ar[ru] \ar[rd] & & {} \ar[ru] \ar[rd] 
& & {} \ar[ru] \ar[rd] & & {} \ar[rd] \ar[ru] & & {\vctEo{4}{3}{3}{1}{3}{2}{2}{1}} \ar[ru] \ar[rd] 
& & {\vctEo{4}{3}{3}{2}{2}{2}{2}{1}} \ar[ru] \ar[rd] & & {} \ar[ru] \ar[rd] & & {} \ar[ru] \ar[rd] 
& & {} \ar[ru] \ar[rd] & & {} \ar[ru] \ar[rd] & & {} \ar[ru] \\
{} \ar[ru] \ar[rd] \ar[r] & {} \ar[r] & {} \ar[ru] \ar[rd] \ar[r] 
& {} \ar[r] & {} \ar[ru] \ar[rd] \ar[r] & {} \ar[r] 
& {} \ar[ru] \ar[rd] \ar[r] & {} \ar[r] & {} \ar[ru] \ar[rd] \ar[r] 
& {} \ar[r] & {} \ar[ru] \ar[rd] \ar[r] & {} \ar[r] 
& {} \ar[ru] \ar[rd] \ar[r] & {} \ar[r] & {} \ar[ru] \ar[rd] \ar[r] 
& {} \ar[r] & {} \ar[ru] \ar[rd] \ar[r] & {\vctEo{3}{2}{3}{1}{2}{1}{2}{1}} \ar[r] 
& {\vctEo{6}{4}{5}{2}{4}{3}{3}{2}} \ar[ru] \ar[rd] \ar[r] 
& {\vctEo{3}{2}{2}{1}{2}{2}{1}{1}} \ar[r] & {} \ar[ru] \ar[rd] \ar[r] & {} \ar[r] & {} \ar[ru] \ar[rd] \ar[r] 
& {} \ar[r] & {} \ar[ru] \ar[rd] \ar[r] & {} \ar[r] & {} \ar[ru] \ar[rd] \ar[r] 
& {} \ar[r] & {} \ar[ru] \ar[rd] \ar[r] & {} \\
& {} \ar[ru] \ar[rd] & & {} \ar[ru] \ar[rd] & & {} \ar[ru] \ar[rd]
& & {} \ar[ru] \ar[rd] & & {} \ar[ru] \ar[rd] & & {} \ar[ru] \ar[rd] 
& & {} \ar[ru] \ar[rd] & & {} \ar[ru] \ar[rd] & & {\vctEo{5}{3}{4}{2}{3}{3}{2}{1}} \ar[ru] \ar[rd] 
& & {\vctEo{4}{3}{4}{1}{3}{2}{2}{1}} \ar[ru] \ar[rd] & & {} \ar[ru] \ar[rd] & & {} \ar[ru] \ar[rd] 
& & {} \ar[ru] \ar[rd] & & {} \ar[ru] \ar[rd] & & {} \ar[rd] \\ 
& & {} \ar[ru] \ar[rd] & & {} \ar[ru] \ar[rd] & & {} \ar[ru] \ar[rd] 
& & {} \ar[ru] \ar[rd] & & {} \ar[ru] \ar[rd] & & {} \ar[ru] \ar[rd] 
& & {} \ar[ru] \ar[rd] & & {\vctEo{4}{3}{3}{2}{2}{2}{1}{1}} \ar[ru] \ar[rd] & & {} \ar[ru] \ar[rd] 
& & {\vctEo{3}{2}{3}{1}{3}{1}{2}{1}} \ar[ru] \ar[rd] & & {} \ar[ru] \ar[rd] & & {} \ar[ru] \ar[rd] 
& & {} \ar[ru] \ar[rd] & & {} \ar[ru] \ar[rd] & & {} \ar[rd] \\ 
& & & {} \ar[ru] \ar[rd] & & {} \ar[ru] \ar[rd] & & {} \ar[ru] \ar[rd] & & {} \ar[ru] \ar[rd]
& & {} \ar[ru] \ar[rd] & & {} \ar[ru] \ar[rd] & & {\vctEo{3}{2}{2}{1}{1}{2}{1}{1}} \ar[ru] \ar[rd] & & {} \ar[ru] \ar[rd] 
& & {} \ar[ru] \ar[rd] & & {\vctEo{2}{1}{2}{0}{2}{1}{2}{1}} \ar[ru] \ar[rd] & & {} \ar[ru] \ar[rd] & & {} \ar[ru] \ar[rd] 
& & {} \ar[ru] \ar[rd] & & {} \ar[ru] \ar[rd] & & {} \ar[rd] \\
& & & & {} \ar[ru] & & {} \ar[ru] & & {} \ar[ru] & & {} \ar[ru] 
& & {} \ar[ru] & & {\vctEo{2}{1}{1}{1}{1}{1}{1}{1}} \ar[ru] & & {} \ar[ru] 
& & {} \ar[ru] & & {} \ar[ru] & & {\vctEo{1}{1}{1}{0}{1}{0}{1}{1}} \ar[ru] 
& & {} \ar[ru] & & {} \ar[ru] & & {} \ar[ru] & & {} \ar[ru] & & {}
}$ 
\caption{Auslander-Reiten quivers of the Dynkin diagrams $\mathbf{E_6}$ (top left),
$\mathbf{E_7}$ (top right) and $\mathbf{E_8}$ (bottom) with subspace orientation. 
In each case we show the positive roots corresponding to the left side $\mathcal{X}_0$ and right side $\mathcal{Y}_0$ of the wings
of vertex $[W_0]$.} \label{(C)T:AREm}
\end{sidewaysfigure}
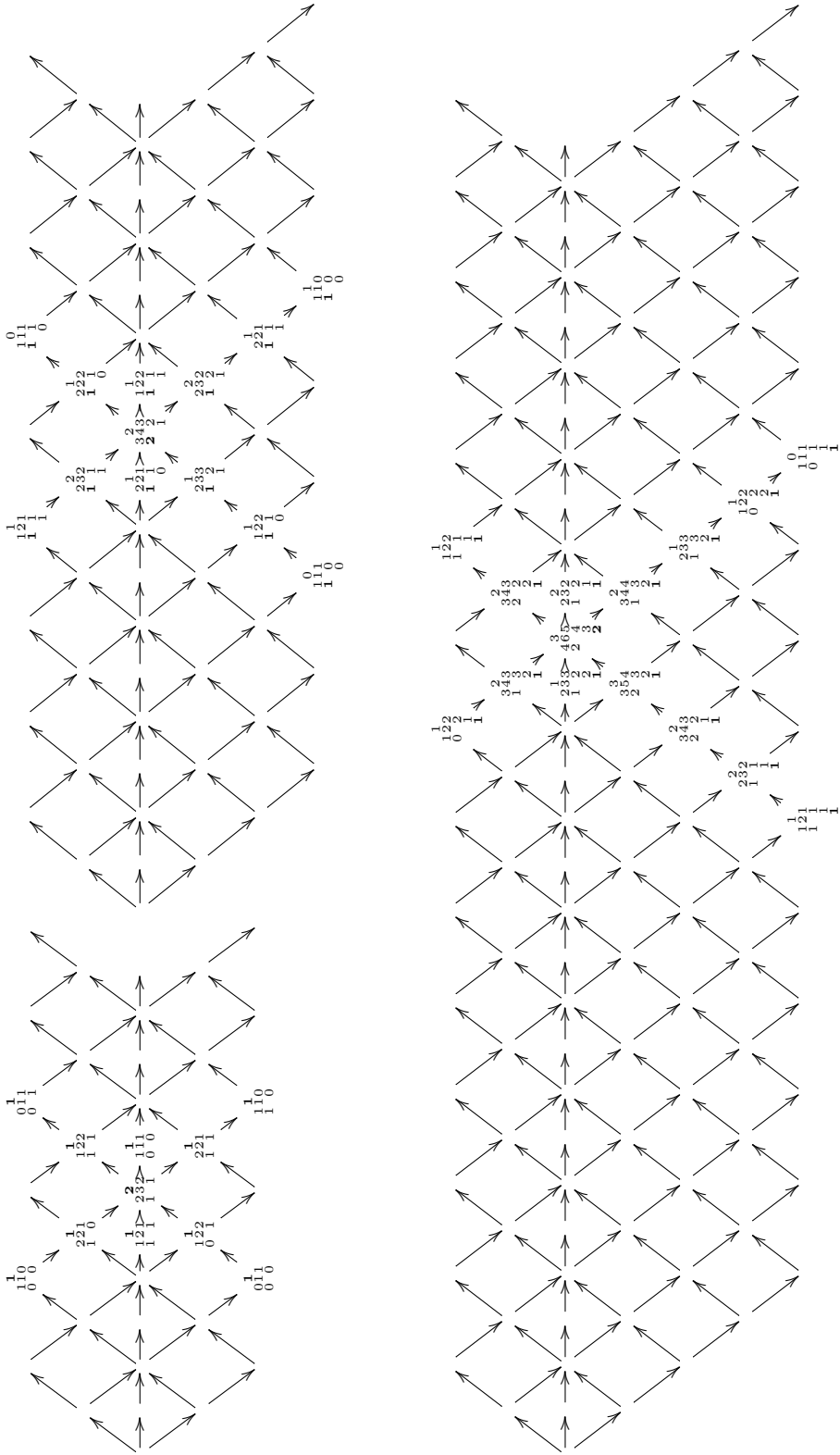
\end{center}

%\section{Planteamiento de problema.} \label{(B)S:Pro}
%--------------------------------------------------------------------------------------
%--------------------------------------------------------------------------------------

\appendix
%----------------------------------------------------------------------
%----------------------------------------------------------------------
%----------------------------------------------------------------------
\chapter{Differential tensor algebras.}
\label{Cap(A)}
%----------------------------------------------------------------------
%----------------------------------------------------------------------
\section*{Ditalgebras and their modules.} \label{(A)S:ditalg}
%--------------------------------------------------------------------------------------
%--------------------------------------------------------------------------------------
Let $k$ be an arbitrary field. In what follows all $k$-algebras will have identity elements and morphisms between $k$-algebras
will send indentities onto identities. All (bi)modules will be finite dimensional over $k$ and all modules will be left modules,
unless otherwise stated.

\begin{definicion} \label{(A)D:tens}
A $k$-algebra $T$ is called \textbf{positively graded} if there exist a decomposition of $k$-vector spaces
\[
 T=\bigoplus_{i=0}^{\infty}T_i,
\]
such that if $a \in T_i$ and $b \in T_j$, then $ab \in T_{i+j}$. Elements in $T_i$ are called \textbf{homogeneous of degree $i$}.
Let $T$ be a $k$-algebra, $A$ a subalgebra of $T$ and $V$ an $A$-$A$-subbimodule of $T$. We say that $T$ is 
\textbf{freely generated by $(A,V)$} if for any $k$-algebra $B$, any algebra morphism
$\varphi_0:A \to B$ and any morphism of $A$-$A$-bimodules $\varphi_1:V \to B$ (where $B$ is considered as an $A$-$A$-bimodule 
through $\varphi_0$), there is a unique algebra morphism $\varphi:T \to B$ which extends $\varphi_0$ and $\varphi_1$.
\[
 \xymatrix{
[A, \; V] \ar@<-.5ex>[d]_-{\varphi_0} \ar@<.5ex>[d]^-{\varphi_1} \ar@{^{(}->}[r] & T \ar@{.>}[dl]^-{\varphi}\\
B 
}
\]
\end{definicion}

\begin{lema} \label{(A)L:ejem}
Let $A$ be a $k$-algebra and $V$ an $A$-$A$-bimodule. Then the tensor algebra $T_A(V)$ is positively graded and freely 
generated by $(A,V)$.
\end{lema}

\bproof
Recall that 
\[
 T_A(V)=\bigoplus_{i=0}^{\infty}V^{\otimes i},
\]
where $V^{\otimes 0}:=A$, $V^{\otimes 1}:=V$ and $V^{\otimes i+1}:=V^{\otimes i} \otimes_A V$ for $i \geq 1$.
The multiplication is given by the tensor product, through the canonical isomorphisms $A \otimes V \cong V$ and
\[
\xymatrix{
V^{\otimes i} \otimes V^{\otimes j} \ar[r] & V^{\otimes i+j}.
}
\]
In this way it is clear that $T_A(V)$ is positively graded. Let $\varphi_0:A \to B$ and $\varphi_1:V \to B$ be as in 
the definition of freely generated algebras. Define $\varphi_i: V^{\otimes i} \to B$ for $i > 1$ as follows
\[
 \xymatrix@R=1pt{
{\sf X}_{i=1}^n V \ar[r] & B \\
(v_1,\ldots,v_n) \ar@{|->}[r] & \varphi_{1}(v_1)\cdots\varphi_1(v_n).
}
\]
Observe that this mapping is $A$-balanced through $\varphi_0$, hence it induces a morphism
$\varphi_i:V^{\otimes i} \to B$. The family $\{ \varphi_i \}$ induces a morphism of $A$-$A$-bimodules
\[
 \varphi: T_A(V)= \bigoplus_{i=0}^{\infty} V^{\otimes i} \longrightarrow B.
\]
It is easy to see that $\varphi$ is an algebra morphism. Unicity is clear.
\eproof

\begin{lema} \label{(A)L:isoTens}
 \begin{enumerate}
  \item[a)] If $T$ is freely generated by $(A,V)$ then the morphism $T \to T_A(V)$, determined by the inclusions
of $A$ and $V$ in $T$, is an isomorphism.
 \item[b)] $T_R(W_0 \oplus W_1) \cong T_{T_R(W_0)}(T_R(W_0) \otimes W_1 \otimes  T_R(W_0))$.
 \end{enumerate}
\end{lema}

\bproof
Since both $T$ and $T_A(V)$ are freely generated by $(A,V)$, there exist unique morphisms $\varphi$, $\psi$ which extend 
the inclusions of $A$ and $V$ in $T$, $T_A(V)$ respectively. 
\[
 \xymatrix{
[A, \; V] \ar@{^{(}->}[d] \ar@{^{(}->}[r] & T \ar@{.>}[dl]^-{\varphi} &
[A, \; V] \ar@{^{(}->}[d] \ar@{^{(}->}[r] & T_A(V) \ar@{.>}[dl]^-{\psi}\\
T_A(V) & & T
}
\]
Then $\psi \circ \varphi$ extends the inclusions of $A$ and $V$ in $T$, and by unicity $\psi \circ \varphi=Id_{T}$.
\[
 \xymatrix{
[A, \; V] \ar@{^{(}->}[d] \ar@{^{(}->}[r] & T \ar@{.>}[dl]^-{\psi \circ \varphi} &
[A, \; V] \ar@{^{(}->}[d] \ar@{^{(}->}[r] & T_A(V) \ar@{.>}[dl]^-{\varphi \circ \psi}\\
T & & T_A(V)
}
\]
Similarly $\varphi \circ \psi=Id_{T_A(V)}$, and hence we have $(a)$. 

To show $(b)$ define $A=T_R(W_0)$ and consider the inclusions $\sigma_0:R \to T_R(W_0 \oplus W_1)$ and
$\sigma_1:W_0 \to T_R(W_0 \oplus W_1)$, which determine an algebra morphism
\[
 \sigma:A \to T_R(W_0 \oplus W_1).
\]
Define also $\varphi_0=\sigma$ and
\[
\varphi_1: \xymatrix@C=3pc@R=.5pc{
A \otimes_R W_1 \otimes_R A \ar[r] & T_R(W_0 \oplus W_1) \\
a \otimes w \otimes b \ar[r] & \varphi_0(a)w\varphi_0(b),
}
\]
which is well defined (for $\varphi_0$ is morphism of $R$-$R$-bimodules) and by definition is a morphism of 
$A$-$A$-bimodules, considering the action of $A$ in $T_R(W_0\oplus W_1)$ through $\varphi_0$. 
Then $\varphi_0$ and $\varphi_1$ extend in a unique way to a morphism of $k$-algebras
\[
 \varphi:T_A(A \otimes_R W_1 \otimes_R A) \longrightarrow T_R(W_0 \oplus W_1).
\]

Take now the inclusion morphisms $\psi_0:R \to A$ and $\psi^0_1:W_0 \to A$ and consider the morphism of
$R$-$R$-bimodules
\[
\psi^1_1: \xymatrix@C=3pc@R=.5pc{
W_1 \ar[r] & A \otimes_R W_1 \otimes_R A  \\
w \ar[r] & 1 \otimes w \otimes 1.
}
\]
Make $\psi_1=\left[ \begin{smallmatrix} \psi^0_1\\ \psi^1_1 \end{smallmatrix} \right]:
W_0 \oplus W_1 \to T_A(A \otimes_R W_1 \otimes_R A)$ and notice that $\psi_1$ is also a morphism of
$R$-$R$-bimodules. Then $\psi_0$ and $\psi_1$ extend to a morphism of $k$-algebras 
$\psi:T_R(W_0 \oplus W_1) \to T_A(A \otimes_R W_1 \otimes_R A)$, which happens to be inverse of $\varphi$.
\eproof

\begin{definicion} \label{(A)D:ditalg}
\begin{enumerate}
 \item[a)] Let $T$ be a positively graded $k$-algebra (with degree $|a|$ for an homogeneous element $a$). 
A \textbf{differential} $\delta$ in $T$ is a $k$-linear transformation 
$\delta: T \longrightarrow T$ such that $\delta(T_i) \subset T_{i+1}$ and that satisfies the Leibniz rule:
if $a$ and $b$ are homogeneous elements, then
\begin{equation} \label{(A)EQ:leibniz}
 \delta(ab)=\delta(a)b + (-1)^{|a|}a\delta(b).
\end{equation} 

\item[b)] A \textbf{differential tensor algebra} (\textbf{ditalgebra}) $\mathcal{A}=(T,\delta)$ consists in a positively graded
$k$-algebra $T$ which is freely generated by $(T_0,T_1)$, together with a differential $\delta$ in $T$
such that $\delta^2=0$. A ditalgebra morphism $\varphi:(T,\delta) \to (T',\delta')$ is an algebra morphism 
$\varphi:T \to T'$ such that $\varphi(T_i)\subseteq T'_i$ for all $i$, and $\delta'\varphi=\varphi \delta$.
Denote by $A$ the subalgebra $T_0$ of $T$, and by $V$ the $A$-$A$-subbimodule $T_1$.
\end{enumerate}
\end{definicion}

Observe that if $\delta$ is a differential then $\delta(c)=0$ for any $c \in k \subset T$.
If $A$ is a $k$-algebra, then $\mathcal{A}=(T_A(0),0)$ is a ditalgebra (called \textbf{trivial} over $A$).
Recall that if $M$ and $N$ are left $A$-modules then $\Hom_k(M,N)$ is an $A$-$A$-bimodule with actions
\[
 [a\cdot f \cdot b](m)=af(bm).
\]

\begin{definicion} \label{(A)D:hom}
For $A$-modules $M$ and $N$ define the subspace
\[
 \Hom_{\mathcal{A}}(M,N) \subset \Hom_k(M,N) \oplus \Hom_{A \text{-}A}(V,\Hom_k(M,N))
\]
given by pairs $(f^0,f^1)$ which satisfy 
\begin{equation} \label{(A)EQ:hom}
 af^0(m)=f^0(am)+f^1(\delta(a))(m),
\end{equation} 
for any $a \in A$ and $m \in M$. Moreover, if $L$ is an $A$-modulo, define the composition function
\[
 \circ: \Hom_{\mathcal{A}}(M,N) \times \Hom_{\mathcal{A}}(L,M) \longrightarrow \Hom_{\mathcal{A}}(L,N)
\]
as $(f,g) \mapsto f \circ g=(f^0 \circ g^0, (f \circ g)^1)$, where $(f \circ g)^1$ is given for $v \in V$ with 
$\delta(v)=\sum_iv_i'v_i''$ by
\begin{equation} \label{(A)EQ:comp}
 (f \circ g)^1=f^0 \circ g^1(v) + f^1(v) \circ g^0 + \sum_i f^1(v_i')g^1(v_i'').
\end{equation} 
\end{definicion}

The last composition is well defined. In fact it determines, together with the $A$-modules as objects and
morphism groups given by $\Hom_{\mathcal{A}}(M,N)$, an additive $k$-category which we denote by $\mathcal{A}$-mod and is 
called category of (finite dimensional) modules over the ditalgebra $\mathcal{A}$ \cite[proposition 2.3]{BSZ09}.

\begin{lema} \label{(A)L:induc}
A ditalgebra morphism $\varphi : \mathcal{A} \longrightarrow \mathcal{A}'$ induces a functor
\[
 F_{\varphi} : \mathcal{A}'\text{-mod} \longrightarrow \mathcal{A} \text{-mod}.
\]
Moreover, if $\varphi$ is surjective, then $F_{\varphi}$ is faithful.
\end{lema}

\bproof
The restriction $\varphi_0:A \to A'$ induces a restriction functor $F_{\varphi_0}:A'\text{-mod} \longrightarrow A\text{-mod} $. 
Define $F_{\varphi}=F_{\varphi_0}$ in objects. For the definition of $F_{\varphi}$ in a morphism $f=(f^0,f^1):M \to N$ 
consider the composition
\[
 \xymatrix{
V \ar[r]^-{\varphi_1} & V' \ar[r]^-{f^1} & \Hom_k(M,N),
}
\]
which is clearly of $A$-$A$-bimodules, for $\varphi_1$ is a morphism of $A$-$A$-bimodules through the map $\varphi_0$.
Define then $F_{\varphi}(f^0,f^1)=(f^0,f^1 \circ \varphi_1)$. We verify that this is a morphism: for
$a \in A$ and $m \in M$ we have
\begin{eqnarray}
 a\cdot f^0(m) & = & \varphi_0(a)f^0(m)=f^0(\varphi_0(a)m)+f^1(\delta(\varphi_0(a)))(m)= \nonumber \\
& = & f^0(a \cdot m)+(f^1\circ \varphi_1)(\delta(a))(m). \nonumber
\end{eqnarray}
Clearly the mapping $F_{\varphi}$ sends identities onto identities.
To see that $F_{\varphi}$ respects compositions it is enough to verify that 
\[
 (g \circ f)^1 \circ \varphi_1= (F_{\varphi}(g) \circ F_{\varphi}(g))^1.
\]
Let $v \in V$ and make $\delta(v)=\sum_iv_i'v_i''$. Then
\begin{eqnarray}
 (F_{\varphi}(g) \circ F_{\varphi}(g))^1(v) & = & g^0f^1(\varphi_1(v))+g^1(\varphi_1(v))f^0+ \nonumber \\
& & + \sum_ig^1(\varphi_1(v_i'))f^1(\varphi_1(v_i''))=  \nonumber \\
& = & (g \circ f)^1(\varphi_1(v)). \nonumber
\end{eqnarray}
Assume now that $f=(f^0,f^1):M \to N$ is a morphism such that
$F(f)=(f^0,f^1 \circ \varphi_1)=0$. Then $f^0=0$ and if $\varphi $ is surjective then $f^1=0$, 
that is, $F_{\varphi}$ is a faithful functor.
\eproof

\section*{Layers, triangularity and extensions.} \label{(A)S:triang}
%--------------------------------------------------------------------------------------
%--------------------------------------------------------------------------------------
\begin{definicion} \label{(A)D:estrato}
Let $R$ be a $k$-algebra and $W$ a $R$-$R$-bimodule with decomposition $W=W_0\oplus W_1$.
Then $(R,W)$ is called \textbf{layer} of the ditalgebra $\mathcal{A}=(T,\delta)$ if $R,W_0 \subseteq T_0$, $W_1 \subseteq T_1$,
the algebra $T$ is freely generated by $(R,W)$ and $\delta(R)=0$.
Then $\delta:T \to T$ is a morphism of $R$-$R$-bimodules and to determine the diffe\-ren\-tial it is enough to give its values on 
$W_0$ and $W_1$
\cite[lema 4.4]{BSZ09}. 
\end{definicion}

\begin{lema} \label{(A)L:extend}
Let $\mathcal{A}=(T,\delta)$ be a ditalgebra with layer $(R,W_0\oplus W_1)$ and $A=T_0$. 
Assume that $A'$ is a subalgebra of $A$ and that $M$ and $N$ are $A'$-modules. Then there exists a natural isomorphism
in $M$ and $N$
\[
\xymatrix@C=1pc{
 \Hom_{R \text{-} R}(W_1,\Hom_k(M,N)) \ar[r]^-{\widehat{}} & \Hom_{A' \text{-} A'}(A' \otimes_R W_1 \otimes_R A',\Hom_k(M,N)).
}\]
\end{lema}
\bproof
We directly observe that the mappings 
\[
 \xymatrix@R=1pc@C=1pc{
\Hom_{R\text{-}R}(W_1,\Hom_k(M,N)) \ar[r]^-{\widehat{}} & \Hom_{A'\text{-}A'}(A' \otimes_R W_1 \otimes_R A',\Hom_k(M,N)) \\
h \ar@{|->}[r] & [a \otimes u \otimes a' \mapsto ah(u)a']
}
\]
and
\[
 \xymatrix@R=1pc@C=1pc{
\Hom_{A'\text{-}A'}(A' \otimes_R W_1 \otimes_R A',\Hom_k(M,N)) \ar[r] & \Hom_{R\text{-}R}(W_1,\Hom_k(M,N))  \\
f \ar@{|->}[r] & [u \mapsto f(1  \otimes u \otimes 1)].
}
\]
are inverse from each other. The naturality can be also proved in a direct way.
\eproof

By the lemma above, since $A \otimes_R W_1 \otimes_R A \cong T_1=V$, there is an isomorphism
\begin{equation*}
\xymatrix@R=.3pc{
\txt{$\Hom_R(M,N)$ \\ $\oplus$ \\ $\Hom_{R\text{-}R}(W_1,\Hom_k(M,N))$} \ar[r] & 
\txt{$\Hom_R(M,N)$ \\ $\oplus$ \\ $\Hom_{A\text{-}A}(V,\Hom_k(M,N))$}
}
\end{equation*} 
given by $(f^0,f^1) \mapsto (f^0,\widehat{f^1})$.

\begin{definicion} \label{(A)D:triang} 
We say that a layer $(R,W)$ of a ditalgebra $\mathcal{A}=(T,\delta)$ is \textbf{triangular} if
\begin{itemize}
 \item[1)] there is a sequence of $R$-$R$-subbimodules $0=W_0^0\subseteq W_0^1 \subseteq \cdots \subseteq W_0^r=W_0$
such that $\delta(W_0^{i+1}) \subset A_i W_1A_i$, for $0 \leq i < r$, where $A_i$ is the subalgebra of
$A$ generated by $R$ and $W_0^i$,
 \item[2)] there is a sequence of $R$-$R$-subbimodules $0=W_1^0\subseteq W_1^1 \subseteq \cdots \subseteq W_1^s=W_1$
such that $\delta(W_1^{j+1}) \subset AW_1^jAW_1^jA$, for $0 \leq j < s$.
\end{itemize}
If moreover each bimodule $W_0^i$ is a direct summand of $W_0^{i+1}$ for $0 \leq i < r$, we say that the layer $(R,W)$
is additive triangular.
\end{definicion}

Any ditalgebra with triangular layer which satisfies the conclusion of the following theorem is called
\textbf{Roiter ditalgebra}.

\begin{teorema}[Roiter property] \label{(A)T:roiter}
Let $\mathcal{A}=(T,\delta)$ be a ditalgebra with additiv triangular layer $(R,W)$ such that $R$ is semi-simple.
Consider a pair of $R$-modules $M$ and $N$ and morphisms 
\[
f^0 \in \Hom_R(M,N)  \qquad \text{and} \qquad f^1 \in \Hom_{R \text{-}R}(W_1,\Hom_k(M,N)).
\]
\begin{enumerate}
 \item[a)] If $f^0$ is a retraction (surjectivity is enough, for $R$ is semi-simple) and $N$ has structure of $A$-module, 
then there exists an structure of $A$-module for $M$ (which extends the structure of $R$-module) such that
\[
 (f^0,\widehat{f^1}):M \longrightarrow N,
\]
is a morphism in $\mathcal{A}$-mod.
 \item[b)] If $f^0$ is a section (injectivity is enough, for $R$ is semi-simple) and $M$ has structure of $A$-module, 
then there exists an structure of $A$-module for $N$ (which extends the structure of $R$-module) such that
\[
 (f^0,\widehat{f^1}):M \longrightarrow N,
\]
is a morphism in $\mathcal{A}$-mod.
\end{enumerate}
\end{teorema}

\bproof
$a)$ Since $f^0:M \to N$ is a retraction, there is a morphism of $R$-modules $g^0:N \to M$ such that $f^0g^0=Id_N$.
To give an action of $A=T_R(W_0)$ in $M$ it is enough to give a morphism $\varphi :W_0 \otimes M \to M$. 
We use the filtration of $W_0$
\[
 0 = W_0^0 \subset W_0^1 \subset \ldots \subset W_0^{r-1} \subset W_0^r = W_0,
\]
and inductively define $\varphi_j:W_0^j \otimes M \to M$ using the following formula 
(where the factor $f_{j-1}$ will be defined later)
\[
 \varphi_j(w \otimes m)=g^0(wf^0(m))-g^0[f^1_{j-1}(\delta(w))(m)], \quad \text{for } w \in W_0^j.
\]
If we have defined $\varphi_{j-1}$, then $M$ is an $A_{j-1}$-module and $\Hom_k(M,N)$ is an $A_{j-1}$-$A_{j-1}$-bimodule.
In this way, using lemma~\ref{(A)L:extend} and since the bimodule $A_{j-1} \otimes_R W_1 \otimes_R A_{j-1}$ is isomorphic to
$A_{j-1}W_1A_{j-1}$ for $W_0^i$ is direct summand of $W_0^{i+1}$ (cf. 5.2 in~\cite{BSZ09}), 
the morphism $f^1 \in \Hom_{R \text{-}R}(W_1,\Hom_k(M,N))$ can be extended to a morphism
\[
 f^1_{j-1} \in \Hom_{A_{j-1} \text{-}A_{j-1}}(A_{j-1}W_1A_{j-1},\Hom_k(M,N)).
\]
As base case, if $w \in W_0^1$ we have that $\delta(w) \in A_0W_1A_0$. But $A_0=R$ thus $\delta(w) \in W_1$.
Define 
\[
 \varphi_1(w \otimes m)=g^0(wf^0(m))-g^0[f^1(\delta(w))(m)].
\] 
Observe that in the last step we get a morphism $\varphi=\varphi_r$ such that
\[
 \varphi(w \otimes m)=g^0(wf^0(m))-g^0[\widehat{f^1}(\delta(w))(m)], \quad \text{for } w \in W_0.
\]
Notice that for $w \in R$ or $w \in W_0$ we have
\begin{eqnarray}
 f^0(w\cdot m) & = & f^0(\varphi(w \otimes m))= f^0(g^0(wf^0(m)))-f^0(g^0[\widehat{f^1}(\delta(w))(m)]) = \nonumber \\
& = & wf^0(m)-f^1(\delta(w))(m), \nonumber
\end{eqnarray}
that is,
\[
 wf^0(m)=f^0(w \cdot m)+f^1(\delta(w))(m).
\]
The last formula extends to $A=T_R(W_0)$. Indeed, assume that the formula is valid for elements in
$W_0^{\otimes \ell-1}$ and let $w \in W_0^{\otimes \ell-1}$ and $w_{\ell} \in W_0$. Then
\begin{eqnarray}
 (w_{\ell} w)f^0(m) & = & w_{\ell}[f^0(w\cdot m)+\widehat{f^1}(\delta(w))(m)] = \nonumber \\
& = & w_{\ell}f^0(w\cdot m)+\widehat{f^1}(w_{\ell}\delta(w))(m) = \nonumber \\
& = & f^0(w_{\ell}\cdot(w\cdot m)) + \widehat{f^1}(\delta(w_{\ell}))(w\cdot m)+\widehat{f^1}(w_{\ell}\delta(w))(m) = \nonumber \\
& = & f^0((w_{\ell}w)\cdot m)+\widehat{f^1}(\delta(w_{\ell}w)), \nonumber 
\end{eqnarray}
hence (by linearity) the formula is also valid in $W_0^{\otimes \ell}$, and thus in all $A$.
The proof of $(b)$ is similar.
\eproof

With the same idea in mind we can prove that (when $R$ is semi-simple) a morphism
$(f^0,f^1)$ in $\mathcal{A}$-mod is an isomorphism if and only if $f^0$ is an isomorphism in $R$-mod \cite[lemma 5.8]{BSZ09}
Assuming that $W_1$ is finitly generated, it can be also proved that $\mathcal{A}$-mod is a Krull-Schmidt category
\cite[theorem~5.13]{BSZ09}. We will assume from now on that $\mathcal{A}$ is a ditalgebra as we have described in this
paragraph. The following lemma can be found in \cite[lemma 5.14]{BSZ09} 

\begin{lema} \label{(A)L:sucExact}
Assume that there exists an exact pair  $\xymatrix{M \ar[r]^-{f} & E \ar[r]^-{g} & N}$ in $\mathcal{A}$-mod such that $gf=0$ and
$\xymatrix{0 \ar[r] & M \ar[r]^-{f^0} & E \ar[r]^-{g^0} & N \ar[r] & 0}$ is an (split) exact sequence in $R$-mod.
Then there is an isomorphism $h :E' \longrightarrow E$ such that
\[
 (h^{-1}f)^1=0 \quad \text{ and } \quad (gh)^1=0.
\]
\end{lema}

\begin{definicion} \label{(A)D:ext}
For two $\mathcal{A}$-modules $M$ and $N$ define the collection $\mathcal{E}(M,N)$ of pairs of composable morphisms $(f,g)$
\[
 \xymatrix{
N \ar[r]^-{f} & E \ar[r]^-{g} & M,
}
\]
such that $gf=0$ and the sequence of $R$-modules
\[
 \xymatrix{
0 \ar[r] & N \ar[r]^-{f^0} & E \ar[r]^-{g^0} & M \ar[r] & 0,
}
\]
is an exact split sequence.
Define in $\mathcal{E}(N,M)$ the relation $(f,g) \sim (f',g')$ if there exists an isomorphism $h:E \to E'$ such
that the following diagram is commutative
\[
 \xymatrix{
N \ar[r]^-{f} & E \ar[r]^-{g} \ar[d]^-{h} & M \\
N \ar[r]_-{f'} \ar@{=}[u] & E' \ar[r]_-{g'} & M \ar@{=}[u].
}
\]
Let $\Ext_{\mathcal{A}}(M,N)=\mathcal{E}(M,N)/\sim $.
\end{definicion}

\begin{lema} \label{(A)L:exacta}
Let $\mathcal{A}$ be a Roiter ditalgebra with layer $(R,W)$. There exists an exact sequence
\[
 \xymatrix@C=1pc{
0 \ar[r] &  \Hom_{\mathcal{A}}(M,N) \ar[r] &  \Hom_R(M,N) \oplus \Hom_{R \text{-}R}(W_1,\Hom_k(M,N)) \ar[r]^-{\sigma} &
}
\]
\[
\xymatrix@C=1pc{
\ar[r]^-{\sigma} &  \Hom_R(W_0 \otimes_R M,N) \ar[r]^-{\eta} &  \Ext_{\mathcal{A}}(M,N) \ar[r]  & 0,
}
\]
where $\sigma$ is the morphism given by 
\[
\sigma(f^0,f^1)(w \otimes m)=wf^0(m)-f^0(wm)-f^1(\delta(w))(m).
\]
\end{lema}

\bproof
Since $f^0$, $f^1$ and the restriction of $\delta$ to $W_1$ are morphisms of $R$-(bi)modules, $\sigma $ is well defined.
Let $M$ and $N$ be a pair of $A$-modules. By the isomorphism given in lemma~\ref{(A)L:extend}, the vector space 
$\Hom_{\mathcal{A}}(M,N)$ is kernel of $\sigma$. We define $\eta $ as follows. 
For a morphism of $R$-modules $h:W_0 \otimes M \to N$ define an $A$-module $N \oplus_h M$, whose underlying $R$-module 
is $N \oplus M$, in the following way. Since $A=T_R(W_0)$ is freely generated by $(R,W_0)$, it is enough to give an action 
$W_0 \otimes_R(N \oplus M) \to N \oplus M$. Take
\[
 \xymatrix@R=1pc{
W_0 \otimes_R (N \oplus M) \ar[r] & N \oplus M \\
w \otimes (n,m) \ar@{|->}[r] & (wn+h(w \otimes m), wm),
}
\]
which is a well defined morphism of $R$-modules. The split exact sequence of $R$-modules
$\xymatrix{0 \ar[r] & N \ar[r]^-{s^0=\left[ \substack{1 \\ 0} \right]} & N \oplus M \ar[r]^-{p^0=[0 \; 1]} & M \ar[r] & 0}$
can be considered as an exact sequence of $A$-modules,
\[
 \xymatrix{
0 \ar[r] & N \ar[r]^-{s^0} & N \oplus_h M \ar[r]^-{p^0} & M \ar[r] & 0,
}
\]
which in turn determines an exact sequence in $\mathcal{A}$-mod given by
\[
\eta_h= \xymatrix{
N \ar[r]^-{(s^0,0)} & N \oplus_h M \ar[r]^-{(p^0,0)} & M.
}
\]
Define $\eta(h)$ as the equivalence class of $\eta_h$ in $\Ext_{\mathcal{A}}(M,N)$. 

Conversely, if $M$, $N$ and $E$ are $A$-modules such that the underlying $R$-module of $E$ is $N\oplus M$ and
$\xymatrix{0 \ar[r] & N \ar[r]^-{s^0} & N \oplus_h M \ar[r]^-{p^0} & M \ar[r] & 0}$ is an exact sequence in $A$-mod,
then $E=N\oplus_h M$ for some morphism $h$ in $\Hom_R(W_0 \otimes_R M,N)$. \\
\underline{Step 1.} \textit{$\eta$ is surjective.} 
By lemma~\ref{(A)L:sucExact} for each element in $\Ext_{\mathcal{A}}(M,N)$ we can take a class representative
of the form $e=\xymatrix{N \ar[r]^-{(f^0,0)} & E \ar[r]^-{(g^0,0)} & M}$. Since we assume that $R$ is semi-simple, the
exact sequence of $R$-modules 
\[
\xymatrix{0 \ar[r] & N \ar[r]^-{f^0} & E \ar[r]^-{g^0} & M \ar[r] & 0}
\]
 splits, 
and hence there is an isomorphism $h^0$ of $R$-modules such that the following is a commutative diagram in $R$-mod,
\[
 \xymatrix{
0 \ar[r] & N \ar@{=}[d] \ar[r]^-{f^0} & E \ar[d]^-{h^0} \ar[r]^-{g^0} & M \ar@{=}[d] \ar[r] & 0 \\
0 \ar[r] & N \ar[r]^-{s^0} & N \oplus M \ar[r]^-{p^0} & M \ar[r] & 0.
}
\]
Since $\mathcal{A}$ is a Roiter ditalgebra, the $R$-module $N \oplus M$ admits an structure of $A$-module $\widehat{N \oplus M}$,
in such a way that $(h^0,0):E \to \widehat{N \oplus M}$ is a morphism in $\mathcal{A}$-mod.
By commutativity of the last diagram, $s^0$ and $p^0$ are morphisms of $A$-modules, and by the observation above
there is a morphism of $R$-modules $h:W_0 \otimes_R M \to N$ such that $\widehat{N \oplus M}=N \oplus_h M$.
Then there is a commutative diagram in $\mathcal{A}$-mod
\[
 \xymatrix{
N \ar@{=}[d] \ar[r]^-{(f^0,0)} & E \ar[d]^-{(h^0,0)} \ar[r]^-{(g^0,0)} & M \ar@{=}[d] \\
N \ar[r]_-{(s^0,0)} & N \oplus M \ar[r]_-{(p^0,0)} & M,
}
\]
and since $(h^0,0)$ is an isomorphism, $\eta(h)$ is the equivalence class of $e$. \\
\underline{Step 2.} \textit{The kernel of $\eta$ is the image of $\sigma$.}
Let $e$ be an element in $\Ext_{\mathcal{A}}(M,N)$ given as $e=\eta(h)$ for some morphism $h:W_0 \otimes M \to N$.
Then $e$ splits if and only if there exists an isomorphism
\[
 \xymatrix{
N \ar@{=}[d] \ar[r]^-{(s^0,0)} & N \oplus_0 M \ar[d]^-{\lambda} \ar[r]^-{(p^0,0)} & M \ar@{=}[d] \\
N \ar[r]_-{(s^0,0)} & N \oplus_h M \ar[r]_-{(p^0,0)} & M.
}
\]
By commutativity $\lambda$ has the form
\begin{displaymath}
 \left[ \begin{array}{cc}
	  I_N & (f^0,f^1) \\
	  0 & I_M 
        \end{array}
\right],
\end{displaymath}
with $f^0 \in \Hom_R(M,N)$ and $f^1 \in \Hom_{R\text{-}R}(W_0,\Hom_k(M,N))$.
Hence $e$ splits if and only if $\lambda$ is a morphism in $\mathcal{A}$-mod, that is
\[
 w\lambda^0(n+m)=\lambda^0(w(n+m))+\lambda^1(\delta(w))(n+m),
\]
or in matricial form, 
\[
 w\left[ \begin{array}{cc}
	  I_N & f^0 \\
	  0 & I_M 
        \end{array}
\right](n,m)=\left[ \begin{array}{cc}
	  I_N & f^0 \\
	  0 & I_M 
        \end{array}
\right](w(n,m))+\left[ \begin{array}{cc}
	  0 & f^1(\delta(w)) \\
	  0 & 0
        \end{array}
\right](n,m),
\]
if and only if
\[
 \left[ \begin{array}{c}
	  wn+wf^0(m) \\
	  wm 
        \end{array}
\right]=\left[ \begin{array}{cc}
	  I_N & f^0 \\
	  0 & I_M 
        \end{array}
\right](wn+h(w \otimes m), wm)+\left[ \begin{array}{c}
	   f^1(\delta(w))(m) \\
	   0
        \end{array}
\right],
\]
if and only if
\[
 h(w \otimes m) = wf^0(m)-f^0(wm)-f^1(\delta(w))(m),
\]
that is, if and only if $h=\sigma(f^0,f^1)$. The linearity of $\eta$ can be shown using the additive structure
in $\Ext_{\mathcal{A}}(M,N)$ (cf. proposition~6.13 in \cite{BSZ09}).
\eproof

%----------------------------------------------------------------------
%----------------------------------------------------------------------
\section*{Some properties of bimodules and its duals.} \label{(A)S:bimod}

Recall that a \textbf{right dual basis} $\{(u_i,\lambda_i)\}_{i\in I}$ of a right $S$-module $U$ is a collection
of elements $u_i \in U$ and $\lambda_i \in U^*=\Hom_S(U_S,S_S)$ such that for any $u \in U$ the set of indexes $i$
for which $\lambda_i(u) \neq 0$ is finite and
\[
 u=\sum_{i \in I} u_i\lambda_i(u).
\]
The proof of point $(a)$ in the following lemma can be found in \cite[lemma~11.3]{BSZ09}.
Point $(b)$ can be proved in a similar way.

\begin{lema} \label{(A)L:phiBimod}
Let $W$ be a $T$-$S$-bimodule and assume that $U$ is a $R$-$S$-bimodule which admits a finite dual basis 
$\{(u_i,\lambda_i)\}_{i \in I}$ as right $S$-module. Then there exists a natural isomorphism of $T$-$R$-bimodules
\[
\varphi: {}_TW_S \otimes_S U^*_R \longrightarrow \Hom_S(U,W), 
\]
given by $\varphi(w \otimes \lambda)(u)=w\lambda(u)$. The inverse $\varphi^{-1}$ has the form
\[
 \varphi^{-1}(h)=\sum_{i \in I} h(u_i) \otimes \lambda_i.
\]
\end{lema}

For the following result confer lemma~11.4 in \cite{BSZ09}.
\begin{lema} \label{(A)L:psiBimod}
Consider bimodules ${}_TA_S$ and ${}_SB_R$ such that $A_S$ has a finite right dual basis $\{(a_i,\lambda_i)\}_{i \in I}$.
Then there is a natural isomorphism of $R$-$T$-bimodules $\psi:B^* \otimes_S A^* \to (A \otimes_S B)^*$ given by the rule
$\psi(g \otimes f)[a \otimes b]=g(f(a)b)$. Its inverse satisfies, for an element $h \in (A \otimes_S B)^*$,
\[
 \psi^{-1}(h)=\sum_{i \in I} \rho_i(h) \otimes \lambda_i,
\]
where $\rho_i(h)[b]=h(a_i \otimes b)$.
\end{lema}

\begin{definicion} \label{(A)D:escision}
Let $\Gamma$ be a $k$-algebra, $S$ a subalgebra of $\Gamma$ and $P$ an ideal of $\Gamma$.
The pair $(S,P)$ is called splitting of $\Gamma$ if the module $P_S$ admits a finite dual basis
$\{(p_j,\gamma_j)\}_{j \in J}$ and there exists a decomposition of $S$-$S$-bimodules $\Gamma=S\oplus P$. 
In that case the product in $\Gamma$ induces a morphism of $S$-$S$-bimodules $m:P\otimes_S P \to P$.
The is a morphism of $S$-$S$-bimodules $\mu:=\psi^{-1}m^*: P^* \to P^* \otimes P^*$, 
where $P^*=\Hom_S(P_S,S_S)$. The mapping $\mu$ is called comultiplication of $P$ and satisfies the coassociativity rule,
that is, the following diagram is commutative,
\[
 \xymatrix{
P^* \ar[r]^-{\mu} \ar[d]_-{\mu} & P^* \otimes P^* \ar[d]^-{1 \otimes \mu} \\
P^* \otimes P^* \ar[r]_-{\mu \otimes 1} & P^* \otimes P^* \otimes P^*. 
}
\]
We have the following explicit expression for the coproduct in terms of the dual basis of $P$,
\[
 \mu(\gamma)= \sum_{i,j \in J}\gamma(p_ip_j)\gamma_j \otimes \gamma_i, \qquad \text{for $\gamma \in P^*$}.
\]
\end{definicion}

\begin{definicion} \label{(A)D:coacciones}
Assume that $B$ and $\Gamma$ are $k$-algebras and that $(S,P)$ is an splitting of $\Gamma$. 
Assume moreover that $X$ is a $B$-$\Gamma$-bimodule such that $X_S$ admits a finite dual basis
$\{(u_i,\nu_i)\}_{i \in I}$. There are morphisms of $B$-$\Gamma$-bimodules 
\[
m_{\ell}: \xymatrix@R.4pc{
X \otimes_S P \ar[r] & X \\
x \otimes p \ar@{|->}[r] & xp,
} \quad \text{and} \quad
m_{r}: \xymatrix@R.4pc{
X \ar[r] & \Hom_S(P,X) \\
x \ar@{|->}[r] & [p \mapsto xp],
}
\]
which determine the coactions of $P$ in $X$ given by
\[
\lambda:=\psi^{-1}m_{\ell}^*:X^* \to P^* \otimes_S X^* \quad \text{and} \quad \rho:=\varphi^{-1}m_r:X \to X \otimes P^*,
\]
where $\varphi$ and $\psi$ are the morphisms given in lemmas~\ref{(A)L:phiBimod} and~\ref{(A)L:psiBimod}.
Considering the bases of $X_S$ and $P_S$ we have the following expressions of $\lambda$ and $\rho$ in elements
$\nu \in X^*$ and $x \in X$,
\[
 \lambda(\nu)=\sum_{i \in I,j \in I}\nu(x_ip_j)\gamma_j \otimes \nu_i \quad \text{and} \quad
 \rho(x)=\sum_{j \in J}xp_j \otimes \gamma_j.
\]
Moreover, the morphisms $\lambda$ and $\rho$ are compatible with the coproduct $\mu$, that is,
the following diagrams commute,
\[
\xymatrix{
X^* \ar[r]^-{\lambda} \ar[d]_-{\lambda} & P^* \otimes X^* \ar[d]^-{1 \otimes \lambda} \\
P^* \otimes X^* \ar[r]_-{\mu \otimes 1} & P^* \otimes P^* \otimes X^*
} \quad \text{and} \quad  \xymatrix{
X \ar[r]^-{\rho} \ar[d]_-{\rho} & X \otimes P^* \ar[d]^-{1 \otimes \mu} \\
X \otimes P^* \ar[r]_-{\rho \otimes 1} & X \otimes P^* \otimes P^*.
}
\]
\end{definicion}

\section*{Reduction by admissible modules.} \label{(A)S:algorit}
%--------------------------------------------------------------------------------------
%--------------------------------------------------------------------------------------

\begin{definicion} \label{(A)D:subditalg}
Let $\mathcal{A}=(T,\delta)$ be a ditalgebra with layer $(R,W)$ and assume that $W_0=W_0'\oplus W_0''$ is
a decomposition of $R$-$R$-bimodules such that $\delta(W'_0)=0$. Consider the subalgebra $B$ of $T$
generated by $R$ and $W_0'$ and notice that $B$ is a subalgebra of $T_0$ freely generated by $(R,W_0')$.
Observe also that $\mathcal{A}'=(B,0)$ is a subditalgebra of $\mathcal{A}$ and that $T_0$ is freely generated by
$(B,BW_0''B)$. Then the algebra $T$ has layer (besides $(R,W)$) the layer $(B,\underline{W})$ 
(see lemma~\ref{(A)L:isoTens}), where $\underline{W}=\underline{W}_0 \oplus \underline{W}_1$ with 
$\underline{W}_0=BW_0''B$ and $\underline{W}_1=BW_1B$.
\end{definicion}

\begin{definicion} \label{(A)D:admisible}
Let $\mathcal{A}=(T,\delta)$ be a ditalgebra with layer $(R,W)$. Given an $\mathcal{A}$-module $X$, 
the algebra $\Gamma=\End_{\mathcal{A}}(X)^{op}$ acts on the right of $X$ through the rule $xp=p^0(x)$
for $p=(p^0,p^1)$ and $x \in X$. We say that $X$ is an admissible module if it satisfies the following conditions
\begin{itemize}
 \item[1)] $\Gamma$ admits an splitting $(S,P)$,
 \item[2)] the right $S$-module $X$ admits a finite dual basis, and
 \item[3)] the morphisms $f$ in $S$ have the form $f=(f^0,0)$.
\end{itemize}
\end{definicion}

Since in this work we are interested only in admissible $\mathcal{A}'$-modules for subditalgebras
$\mathcal{A}'$ as in~\ref{(A)D:subditalg}, the point $(3)$ in the definition above is automatically satisfied. 
In all cases we are interested in, $S$ is a semi-simple algebra, hence point $(2)$ holds as long as $X$ is a finite
dimensional module.

\begin{definicion} \label{(A):tensRed}
Assume that $X$ is an admissible $\mathcal{A}'$-module, where $\mathcal{A}'$ is a subditalgebra of
$\mathcal{A}$ as in definition~\ref{(A)D:subditalg}. Having in consideration the notation in the definition above,
define the $S$-$S$-bimodule $W^X=W^X_0 \oplus W^X_1$ by the formula
\[
 W^X_0=X^* \otimes_B \underline{W}_0 \otimes_B X \quad \text{and} \quad W^X_1=[X^* \otimes_B \underline{W}_1\otimes_BX]\oplus P^*.
\]
Take the \textbf{reduced tensor algebra} $T^X=T_S(W^X)$ with layer $(S,W^X)$.
\end{definicion}
We want to define a differential $\delta^X$ in $T^X$. For that purpose define the following auxiliar morphisms.

\begin{lema} \label{(A)L:sigma}
Assume that $X$ is an admissible $\mathcal{A}'$-module and that $T^X$ is the reduced tensor algebra as in definition 
above. Let $D(X)=(X^* \times X) \sqcup (P^* \otimes_S X^* \times X) \sqcup (X^* \times X \otimes_S P^*)$.
Recall that $T$ is freely generated by $(B,\underline{W})$. For any element $(u,v) \in D(X)$ there exists a 
linear morphism 
\[
 \sigma_{u,v}: T_B(\underline{W}) \longrightarrow T^X,
\]
such that for $w_1, \ldots,w_n \in \underline{W}$ we have $\sigma_{u,v}(w_1 \otimes \cdots \otimes w_n)$ given by
\begin{equation*}
\begin{array}{r}
\sum_{i_1,\ldots,i_{n-1}}u \otimes w_1 \otimes x_{i_1} \otimes \nu_{i_1} \otimes w_2 \otimes x_{i_2} \otimes \nu_{i_2} \otimes \\
\cdots \otimes x_{i_{n-1}} \otimes \nu_{i_{n-1}} \otimes w_n \otimes v,\\
\end{array}
\end{equation*}
and for $b \in B$, 
\begin{equation*}
\left\{
\begin{array}{ll}
\sigma_{\nu,x}(b)=\nu(bx), & \text{if } \nu \in X^*, x \in X,\\
\sigma_{u,x}(b)=\sum_r \varphi_r \psi_r (bx), & \text{if } x \in X,u=\sum_r \varphi_r \otimes \psi_r, \; \varphi_r \in P^*,\psi_r \in X^*,\\
\sigma_{\nu,v}(b)=\sum_s \nu(by_s)\varphi'_s, & \text{if } \nu \in X^*, v=\sum_s y_s \otimes \varphi'_s,\; y_s \in X,\varphi'_s \in P^*.
\end{array} \right.
\end{equation*}
\end{lema}

The proof, together with some properties of the morphisms $\sigma_{u,v}$, can be found in \cite[lemma 12.8]{BSZ09}.
Of particular interest, the morphism $\sigma_{\nu,x}$ preserves degree for any $x \in X$ and $\nu \in X^*$.

\begin{lema} \label{(A)L:diferen}
There exists a morphism of $S$-$S$-bimodules 
\[
\delta^X:[X^* \otimes_B \underline{W} \otimes_B X]\oplus P^* \to T^X,
\]
whose restriction to $P^*$ is the coproduct $\mu$ of $P$, 
and for $w \in \underline{W}_0 \cup \underline{W}_1$, $x \in X$ and $\nu \in X^*$, 
\[
 \delta^X(\sigma_{\nu,x}(w))=\sigma_{\lambda(\nu),x}(w) + \sigma_{\nu,x}(\delta(w)) + (-1)^{|w|+1}\sigma_{\nu,\rho(x)}(w),
\]
where $\lambda$ and $\rho$ are the coactions given in~\ref{(A)D:coacciones}. This morphism can be extended
to a differential in the tensor algebra $T^X$, which we also denote by $\delta^X$,
in such a way that $\mathcal{A}^X=(T^X,\delta^X)$ is a ditalgebra with layer $(S,W^X)$.
\end{lema}
The ditalgebra $\mathcal{A}^X$ is called \textbf{reduction of $\mathcal{A}$ with respect to $X$}. 
One of the main tools in this work is the \textbf{reduction functor} $F^X : \mathcal{A}^X\text{-mod} \to  
\mathcal{A}\text{-mod}$, as constructed by Bautista, Salmer\'on and Zuazua in proposition~12.10 in~\cite{BSZ09}. 
\begin{proposicion} \label{(A)P:funRed}
Let $\mathcal{A}=(T,\delta)$ be a ditalgebra with layer $(R,W)$ and $A=T_0$.
Let $X$ be an admissible $\mathcal{A}'$-module, where $\mathcal{A}'$ is a subditalgebra of $\mathcal{A}$ as in
definition~\ref{(A)D:subditalg}. Then there is a functor
\[
 F^X : \mathcal{A}^X\text{-mod} \to  \mathcal{A}\text{-mod},
\]
defined in the following way.

Given an $\mathcal{A}^X$-module $M$ denote by $\star$ the left action of $T^X=T_S(W^X)$ in $M$.
Make $F^X(M)=X \otimes_S M$ (considered as left $B$-module) to which we want to give structure of 
$A$-module. For that recall that $A$ is freely generated by $B$ and $\underline{W}_0$. Hence,
in order to give to $F^X(M)$ structure of $A$-module, it is enough to give a morphism of left $B$-modules 
$\underline{W}_0 \otimes_B X \otimes _S M \to X \otimes_S M$. Consider the composition
\[
\xymatrix@R=1pc{ 
\underline{W}_0 \otimes_B (X \otimes_S M) \ar[r] & \End_S(X) \otimes [\underline{W}_0 \otimes_B (X \otimes_S M)] 
\ar[d]^-{\varphi^{-1} \otimes 1} \\
& (X \otimes_S X^*) \otimes_B \underline{W}_0 \otimes_B (X \otimes_S M) \ar[d]^{\cong} \\
& X \otimes_S [(X^* \otimes_B \underline{W}_0 \otimes_B X) \otimes_S M] \ar@{=}[d] \\
& X \otimes_S (\underline{W}_0^X \otimes_S M) \ar[d]^-{1 \otimes \star} \\
& X \otimes_S M,
}
\]
where $\varphi$ is the morphism given in lemma~\ref{(A)L:phiBimod}. Using the expression of $\varphi^{-1}$ given in that
lemma, one shows the following formula for $a \in A$,
\[
a \cdot (x \otimes m)=\sum_{i}x_i \otimes \sigma_{\nu_i,x}(a) \star m.
\]
If $f=(f^0,f^1) \in \Hom_{\mathcal{A}^X}(M,N)$ then $F^X(f)$ is given by
\begin{eqnarray}
 (F^X(f))^0[x \otimes m] & = & x \otimes f^0(m)+\sum_j xp_j \otimes f^1(\gamma_j)[m], \nonumber \\
(F^X(f))^1(v)[x \otimes m] & = & \sum_i x_i \otimes f^1(\sigma_{\nu_i,x}(v))[m], \nonumber
\end{eqnarray}
where $x \in X$, $m \in M$ and $v \in V=T_1$.
\end{proposicion}

We give now conditions which guarantee that the reduction functor $F^X$ is full and faithful. 
Observe that for an admissible $\mathcal{A'}$-module $X$ we can perform two reductions
$\mathcal{A}^X$ and $\mathcal{A}'^X$.
The reduction functors $F^X:\mathcal{A}^X \text{-mod} \to \mathcal{A}\text{-mod}$ and
$F'^X:\mathcal{A}'^X \text{-mod} \to \mathcal{A}'\text{-mod}$ can be related in the following way.
Recall that there is a decomposition $W_0=W_0' \oplus W_0''$ with $\delta(W_0')=0$, and that
$B$ is the subalgebra of $T$ generated by $R$ and $W_0'$. Then the inclusion $r=B \to T$ determines
a ditalgebra morphism $\mathcal{A}' \to \mathcal{A}$.

On the other hand, recall that if $(S,P)$ is a splitting of $\End_{\mathcal{A}'}(X)^{op}$ then
the reduced ditalgebra $\mathcal{A}^X$ has layer $(S,W^X)$ where 
\[
 W_0^X = X^* \otimes_B \underline{W}_0 \otimes_B X \quad \text{and} \quad 
W_1^*=X^* \otimes_B \underline{W}_1 \otimes_B X \oplus P^*,
\]
$\underline{W}_0=BW''_0B$ and $\underline{W}_1=BW_1B$. Now, for the reduction $\mathcal{A}'^X$ consider
the trivial decomposition of bimodules $W_0'=W_0' \oplus 0$ in the layer $(R,W'_0)$ of $\mathcal{A}'$.
Then the reduced ditalgebra $\mathcal{A}'^X$ has layer $(S,P^*)$. The bimodule inclusion
$P^* \to W^X$ induces a morphism of ditalgebras $\widehat{r}:\mathcal{A}'^X \to \mathcal{A}^X$.
Moreover, the following diagram is commutative,
\[
 \xymatrix{
\mathcal{A}^X \text{-mod} \ar[r]^-{F^X} \ar[d]_-{F_{\widehat{r}}} & \mathcal{A} \text{-mod} \ar[d]^-{F_r} \\
\mathcal{A}'^X \text{-mod} \ar[r]_-{F'^X} & \mathcal{A}' \text{-mod}.
}
\]

\begin{definicion} \label{(A)D:completo}
Let $\mathcal{A}'=(B,0)$ be a ditalgebra with layer $(R,W'_0)$ as in de\-fi\-nition~\ref{(A)D:subditalg}.
An admissible $\mathcal{A}'$-module $X$ is called complete if the reduction functor 
$F'^X:\mathcal{A}'^X \text{-mod} \to \mathcal{A}'\text{-mod}$ is full and faithful.
\end{definicion}
The proof of the following results can be found in \cite[13.3 and 13.5]{BSZ09}.

\begin{lema} \label{(A)L:finita}
Let $\mathcal{A}'=(B,0)$ be a ditalgebra with layer $(R,W_0')$ as in de\-fi\-nition~\ref{(A)D:subditalg}. 
Assume that $X$ is an admissible $\mathcal{A}'$-module, where the algebra $\End_{\mathcal{A}'}(X)^{op}$ admits as
splitting the pair $(S,P)$. Assume moreover that $X$ is finite dimensional over the base field $k$ and that
$S$ is semi-simple. Then $X$ is a complete admissible $\mathcal{A}'$-module.
\end{lema}

\begin{proposicion} \label{(A)P:completo}
Assume that $X$ is a complete admissible $\mathcal{A}'$-module, where 
$\mathcal{A}'=(B,0)$ is a subditalgebra of the layered ditalgebra $\mathcal{A}$ 
as in definition~\ref{(A)D:subditalg}. Then the associated functor 
$F^X:\mathcal{A}^X \text{-mod} \to \mathcal{A}\text{-mod}$ is full and faithful.
\end{proposicion}

The following lemma describes the $\mathcal{A}$-modules which are in the image of the functor $F^X$ \cite[lemma 16.1]{BSZ09}.
\begin{lema} \label{(A)L:imagen}
Assume that $X$ is an admissible $\mathcal{A}'$-module, where $\mathcal{A}'=(B,0)$ is a subditalgebra of the layered ditalgebra 
$\mathcal{A}$ as in definition~\ref{(A)D:subditalg}. Assume that $M$ is a $S$-module such that
there exists an object $\overline{X \otimes_SM}$ in $\mathcal{A}$-mod whose underlying $B$-module structure
is the canonical structure in $X \otimes_SM$. Then there is a unique $\mathcal{A}'$-module $\overline{M}$
with underlying $S$-module structure $M$ such that $F^X(\overline{M})=\overline{X \otimes_SM}$.
\end{lema}

We end this section with the description of the effect of the reduction functor $F^X$ in extension groups.
Confer~\cite[16.6]{BSZ09}.

\begin{lema} \label{(A)L:extRed}
Assume that $X$ is an admissible $\mathcal{A}'$-module, where $\mathcal{A}'=(B,0)$ is a subditalgebra of the ditalgebra 
with triangular layer $\mathcal{A}$ as in definition~\ref{(A)D:subditalg}. Assume also that the decomposition 
$W_0=W_0' \oplus W_0''$ of the bimodule $W_0'$ appears in the filtration of $W_0$ in the definition of triangular layer. 
Then $\mathcal{A}'$ is a ditalgebra with triangular layer. Assume further that in the splitting $(S,P)$ of 
$\Gamma=\End_{\mathcal{A}'}(X)^{op}$, the algebra $S$ is semi-simple and $P=\rad \Gamma$
and that $X$ is finite dimensional over the field $k$. Then $F_X:\mathcal{A}^X \text{-mod} \to \mathcal{A}\text{-mod}$ 
is an exact functor and for any pair of $\mathcal{A}^X$-modules $M$ and $N$ it induces an exact sequence of vector spaces
\[
 \xymatrix@C=1pc@R=.5pc{
0 \ar[r] & \Ext^1_{\mathcal{A}^X}(M,N) \ar[r]^-{F_X^*} & \Ext^1_{\mathcal{A}}(F_X(M),F_X(N)) \\ 
& {} \ar[r]^-{F_r^*} & \Ext^1_{\mathcal{A}'}(F_rF_X(M),F_rF_X(N)) \ar[r] & 0.
}
\]
\end{lema}
%--------------------------------------------------------------------------------------

%-------------------------------------------------------------------
%\nocite{*}
\bibliographystyle{plain}
\bibliography{myrefs}
\addcontentsline{toc}{chapter}{Bibliography}

\end{document}